\RequirePackage{fix-cm} 
\RequirePackage[reqno]{amsmath}
\RequirePackage[l2tabu, orthodox]{nag}
\documentclass[final, pdftex, a4paper, twoside, reqno, 12pt]{amsart}

\pdfoutput=1

\usepackage[T1]{fontenc}
\usepackage[latin1]{inputenc}
\usepackage{lmodern}

\usepackage[activate, final, kerning, spacing, factor = 1100, stretch = 10, shrink = 10, spacing=nonfrench]{microtype}
\SetTracking{encoding = {*}, shape = sc}{60}

\usepackage[nonewpage, xindy]{imakeidx}
\makeindex[title=Alphabetical index, intoc]

\usepackage{bm}
\usepackage{esint}
\usepackage{amsbsy}
\usepackage{amsthm}
\usepackage{cancel}
\usepackage{dsfont}
\usepackage{empheq}
\usepackage{manfnt}
\usepackage{pifont}
\usepackage{xspace}
\usepackage{amssymb}
\usepackage{upgreek}
\usepackage{stmaryrd}
\usepackage{etoolbox}
\usepackage{extarrows}
\usepackage{textgreek}

\usepackage[nice]{nicefrac}
\usepackage{siunitx}
\usepackage{fourier-orns}

\usepackage{relsize}
\usepackage{mathrsfs}
\usepackage{textcomp}

\usepackage{expdlist}
\usepackage{etaremune}

\usepackage{mathtools}
\mathtoolsset{centercolon, mathic}

\usepackage[a4paper, margin=1.0in, bottom=1.33in]{geometry}
\usepackage[dvipsnames, table]{xcolor}
\usepackage[most]{tcolorbox}

\definecolor{bckg}{RGB}{20.8, 20.8, 20.8}
\definecolor{oneblue}{rgb}{0.0, 0.0, 0.85}
\definecolor{Lightblue}{RGB}{214, 214, 214}
\definecolor{bluepigment}{rgb}{0.2, 0.2, 0.6}
\definecolor{charcoal}{rgb}{0.21, 0.27, 0.31}
\definecolor{denimblue}{rgb}{0.08, 0.38, 0.74}
\definecolor{darkelectricblue}{rgb}{0.33, 0.41, 0.47}
\definecolor{katyblue}{rgb}{0.129412, 0.137255, 0.63}

\usepackage{booktabs}
\usepackage{array}
\usepackage{ellipsis}
\usepackage{multicol}

\usepackage{paralist}
\usepackage{csquotes}

\usepackage{listings}

\lstdefinelanguage{Maple}%
{morekeywords={and, assuming, break, by, catch, description, do, done,%
elif, else, end, error, export, fi, finally, for, from, global, if,%
implies, in, intersect, local, minus, mod, module, next, not, od,%
option, options, or, proc, quit, read, return, save, stop, subset, then,%
to, try, union, use, uses, while, xor},%
sensitive=true,%
morecomment=[l]\#,%
morestring=[b]",%
morestring=[d]"%
}[keywords,comments,strings]%

\definecolor{mymauve}{rgb}{0.58, 0.0, 0.82}
\usepackage{float}
\usepackage{graphicx}
\graphicspath{{figs/}}
\DeclareGraphicsExtensions{.pdf, .eps}

\usepackage{tikz}
\usepackage{tikz-cd}
\usepackage{subfigure}
\usepackage{morefloats}
\usepackage{indentfirst}

\usepackage[labelsep=period,
            font=normalsize,
            labelformat=simple,
            labelfont={bf,sf,color=bluepigment},
            justification=raggedright]{caption}

\usepackage[bottom, multiple, perpage]{footmisc}
\setfnsymbol{wiley}

\newcommand*{\Title}{\textcolor{bluepigment}{On complex algebraic singularities of nonlinear PDEs}}
\newcommand*{\Longtitle}{On complex algebraic singularities of some genuinely nonlinear PDEs}
\newcommand*{\Authors}{\textcolor{bluepigment}{D.~Dutykh and E.~Leichtnam}}
\newcommand*{\plogo}{\textcolor{gray}{{\texttt{arXiv.org} / \textsc{hal}}}}
\newcommand*{\Keywords}{nonlinear PDEs; initial value problem; complex singularities; ramification}

\usepackage[sort&compress, comma, square, numbers]{natbib}
\makeatletter
\renewcommand{\@biblabel}[1]{\textbf{[#1]}}
\makeatother

\usepackage{url}  
\usepackage[final=true,
           linktoc=all,
           pdftex,
           unicode=true,
           bookmarks=true,
           colorlinks=true,
           breaklinks=true,
           urlcolor=katyblue,
           linkcolor=MidnightBlue,
           citecolor=NavyBlue,
           pdffitwindow=false,
           bookmarksopen=false,
           pdfpagemode=UseNone,
           plainpages=false,
           hyperfigures=true,
           pdfstartview={FitV},
           pdftitle={\Longtitle},
           pdfauthor={Denys DUTYKH},
           pdfsubject={Applied Mathematics},
           pdfcreator={Dr. D},
           pdfproducer={pdfLaTeX},
           pdfkeywords={\Keywords},
           pdfnewwindow=true,
           backref,
           pagebackref=true, pdfa]{hyperref}

\usepackage[sort&compress, capitalise, noabbrev]{cleveref}

\usepackage[all]{hypcap}
\usepackage{doi}

\usepackage[explicit]{titlesec}

\newcommand\invisiblesection[1]{%
  \addcontentsline{toc}{section}{#1}%
  \sectionmark{#1}}

\titleformat{\section}
  {\color{NavyBlue}\Large\sffamily\bfseries}
  {}
  {0em}
  {\hspace*{-0.5em}\colorbox{bckg!5}{\parbox{\dimexpr\linewidth+0\fboxsep\relax}{\centering\Large\thesection. #1}}}
  [\vspace*{0.33em}]

\titleformat{name=\section, numberless}
  {\color{NavyBlue}\Large\sffamily\bfseries}
  {}
  {0.0em}
  {\hspace*{-0.5em}\colorbox{bckg!5}{\parbox{\dimexpr\linewidth+0\fboxsep\relax}{\centering#1}}}
  [\vspace*{0.33em}]

\titleformat{\subsection}
  {\color{NavyBlue}\large\sffamily\bfseries}
  {}
  {0.0em}
  {\hspace*{-0.5em}\colorbox{bckg!5}{\parbox{\dimexpr\linewidth+0\fboxsep\relax}{\centering\thesubsection. #1}}}
  [\vspace*{0.33em}]

\titleformat{name=\subsection, numberless}
  {\color{NavyBlue}\Large\sffamily\bfseries}
  {}
  {0em}
  {\hspace*{-0.5em}\colorbox{bckg!5}{\parbox{\dimexpr\linewidth+0\fboxsep\relax}{\centering#1}}}
  [\vspace*{0.33em}]

\titleformat{\subsubsection}
  {\color{bluepigment}\sffamily\normalsize\bfseries}
  {}
  {0.0em}
  {\hspace*{-0.5em}\colorbox{bckg!5}{\parbox{\dimexpr\linewidth+0\fboxsep\relax}{\centering\thesubsubsection. #1}}}
  [\vspace*{0.33em}]

\titleformat{name=\subsubsection, numberless}
  {\color{bluepigment}\sffamily\normalsize\bfseries}
  {}
  {0.0em}
  {\hspace*{-0.5em}\colorbox{bckg!5}{\parbox{\dimexpr\linewidth+0\fboxsep\relax}{\centering#1}}}
  [\vspace*{0.33em}]

\titleformat{\paragraph}[runin]
  {\color{bluepigment}\sffamily\small\bfseries}
  {}
  {0em}
  {#1}

\titlespacing{\section}{1.0em}{1.5em plus 2pt minus 2pt}%
{1.0em plus 2pt minus 2pt}[0em]
\titlespacing{\subsection}{1.0em}{1.5em plus 2pt minus 2pt}%
{1.0em plus 2pt minus 2pt}[0em]
\titlespacing{\subsubsection}{1.0em}{1.5em plus 2pt minus 2pt}%
{1.0em plus 2pt minus 2pt}[0em]

\usepackage{titletoc}

\contentsmargin{0.5em}
\newlength{\tocsep} 
\titlecontents{section}[\tocsep]
  {\addvspace{10pt}\bfseries\sffamily}
  {\contentslabel[\thecontentslabel]{\tocsep}}
  {}
  {\ \titlerule*[0.75pc]{.}\ \thecontentspage}
  []

\titlecontents{subsection}[\tocsep]
  {\addvspace{8pt}\sffamily}
  {\contentslabel[\thecontentslabel]{\tocsep}}
  {}
  {\ \titlerule*[0.5pc]{.}\ \thecontentspage}
  []

\titlecontents{subsubsection}[\tocsep]
  {\addvspace{4pt}\footnotesize\sffamily}
  {\contentslabel[\thecontentslabel]{\tocsep}\hspace*{0.5em}}
  {}
  {\ \titlerule*[0.35pc]{.}\ \thecontentspage}
  []

\makeatletter
\def\@setauthors{%
  \begingroup
  \def\thanks{\protect\thanks@warning}%
  \trivlist
  \centering\footnotesize \@topsep30\p@\relax
  \advance\@topsep by -\baselineskip
  \item\relax
  \author@andify\authors
  \def\\{\protect\linebreak}%
  \textsc{\normalsize\textcolor{charcoal}{\authors}}%
  \ifx\@empty\contribs
  \else
    ,\penalty-3 \space \@setcontribs
    \@closetoccontribs
  \fi
  \endtrivlist
  \endgroup
}
\def\@settitle{\begin{center}%
  \baselineskip14\p@\relax
    \bfseries
    \textsc{\Large\textcolor{charcoal}{\@title}}
  \end{center}%
}
\makeatother

\usepackage{enumitem}
\setlist[description]{%
  topsep = 9pt,               
  itemsep = 7pt,               
  labelsep = 10pt,
  font={\bfseries\color{NavyBlue}}, 
}
\setlist[itemize]{%
  itemsep = 5pt,
  font={\color{NavyBlue}}
}

\usepackage{fancyhdr}
\usepackage{lastpage}

\usepackage{xpatch}
\xpretocmd\headrule{\color{lightgray}}{}{\PatchFailed}

\numberwithin{equation}{section}

\theoremstyle{plain}
\newtheorem{theorem}{Theorem}[section]
\newtheorem{lemma}{Lemma}[section]
\newtheorem{prop}{Proposition}[section]
\newtheorem{corollary}{Corollary}[section]
\newtheorem{conj}{Conjecture}[section]

\theoremstyle{definition}
\newtheorem{definition}{Definition}[section]
\newtheorem{example}{Example}[section]

\theoremstyle{remark}
\newtheorem{remark}{Remark}[section]

\usepackage{setspace}
\vfuzz \hfuzz
\newcommand{\up}[1]{$^{\mathrm{\small\textsf{#1}}}$} 

\newtcbox{\mymath}[1][]{%
    nobeforeafter, math upper, tcbox raise base,
    enhanced, colframe = black!35,
    colback = black!5, boxrule = 1pt, arc = 0mm,
    #1}

\newcommand{\vO}{\bm{0}}
\newcommand{\ah}{\hat{a}}
\newcommand{\bh}{\hat{b}}
\newcommand{\ph}{\hat{p}}
\newcommand{\qh}{\hat{q}}
\newcommand{\Y}{\Upsilon}
\newcommand{\cz}{\breve{z}}
\newcommand{\C}{\mathds{C}}
\newcommand{\E}{\mathds{E}}
\newcommand{\N}{\mathds{N}}
\newcommand{\R}{\mathds{R}}
\newcommand{\V}{\mathds{V}}
\newcommand{\Z}{\mathds{Z}}
\newcommand{\kk}{\mathds{k}}
\newcommand{\Id}{\mathds{1}}
\newcommand{\Aa}{\mathds{A}}
\newcommand{\Bb}{\mathds{B}}
\newcommand{\Ma}{\mathds{M}}
\newcommand{\A}{\mathscr{A}}
\newcommand{\B}{\mathscr{B}}
\newcommand{\U}{\mathscr{U}}
\newcommand{\ud}{\mathrm{d}}
\newcommand{\ui}{\mathrm{i}}
\newcommand{\ue}{\mathrm{e}}
\newcommand{\D}{\mathcal{D}}
\newcommand{\F}{\mathcal{F}}
\newcommand{\I}{\mathcal{I}}
\newcommand{\W}{\mathscr{W}}
\newcommand{\X}{\mathscr{X}}
\newcommand{\xib}{\bar{\xi}}
\newcommand{\jj}{\mathcal{j}}
\newcommand{\Cl}{\mathcal{C}}
\newcommand{\Cr}{\mathscr{C}}
\newcommand{\Hh}{\mathcal{H}}
\newcommand{\Mm}{\mathcal{M}}
\newcommand{\Nn}{\mathcal{N}}
\newcommand{\Qq}{\mathcal{Q}}
\newcommand{\Rr}{\mathcal{R}}
\newcommand{\Ss}{\mathcal{S}}
\newcommand{\Tt}{\mathcal{T}}
\newcommand{\Yy}{\mathscr{Y}}
\newcommand{\Wu}{\mathscr{W}}
\newcommand{\Cc}{\mathfrak{C}}
\renewcommand{\H}{\mathscr{H}}
\renewcommand{\L}{\mathscr{L}}
\renewcommand{\P}{\mathcal{P}}
\renewcommand{\O}{\mathcal{O}}
\newcommand{\WF}{\mathsf{WF}\,}
\newcommand{\x}{\boldsymbol{x}}
\newcommand{\y}{\boldsymbol{y}}
\newcommand{\z}{\boldsymbol{z}}
\newcommand{\fO}{\underline{0}}
\newcommand{\fW}{\underline{1}}
\renewcommand{\S}{\mathfrak{S}}
\newcommand{\FF}{\boldsymbol{F}}

\renewcommand{\mapsto}{\longmapsto}

\renewcommand{\emptyset}{\varnothing}
\newcommand{\M}{\boldsymbol{\mathsf{M}}}
\newcommand{\T}{\boldsymbol{\mathsf{T}}}
\newcommand{\FFt}{\tilde{\boldsymbol{F}}}
\newcommand{\Nc}{\boldsymbol{\mathsf{N}}}
\newcommand{\Vv}{\boldsymbol{\mathsf{V}}}
\newcommand{\Xx}{\mathsf{\boldsymbol{X}}}
\newcommand{\gr}{\boldsymbol{\mathsf{gr}}\,}
\newcommand{\Br}{\boldsymbol{\mathsf{Br}}\,}
\newcommand{\GL}{\boldsymbol{\mathsf{GL}}\,}
\newcommand{\SL}{\boldsymbol{\mathsf{SL}}\,}
\newcommand{\Ann}{\boldsymbol{\mathsf{Ann}}\,}
\newcommand{\Gal}{\boldsymbol{\mathsf{Gal}}\,}
\newcommand{\PSL}{\boldsymbol{\mathsf{PSL}}\,}
\renewcommand{\eta}{\hspace{0.05em}\mbox{\texteta}}
\newcommand{\powerset}{\raisebox{.15\baselineskip}{\Large\ensuremath{\wp}}}

\DeclareMathOperator{\cod}{cod}
\DeclareMathOperator{\dom}{dom}

\DeclareMathOperator{\comp}{\circ}
\DeclareMathOperator{\coker}{coker}
\DeclareMathOperator{\singsupp}{sing\:supp}

\newcommand{\eps}{\varepsilon}
\renewcommand{\leq}{\leqslant}
\renewcommand{\geq}{\geqslant}

\newcommand{\nm}[1]{\textsc{#1}}
\newcommand{\cf}{\emph{cf.}\xspace}
\newcommand{\ie}{\emph{i.e.}\xspace}
\newcommand{\eg}{\emph{e.g.}\xspace}

\renewcommand{\sim}{\thicksim}

\newcommand{\bydef}{:\Rightarrow}

\newcommand{\defeq}{\mathop{\stackrel{\,\mathrm{def}}{\eqcolon}\,}}
\newcommand{\eqdef}{\mathop{\stackrel{\,\mathrm{def}}{\coloneq}\,}}
\newcommand{\pd}[2]{\frac{\partial\hspace{0.0556em} #1}{\partial\/ #2}}

\newcommand{\od}[2]{\frac{\mathrm{d}\hspace{0.0556em} #1}{\mathrm{d}\/#2}}

\newcommand{\scal}{\;\raisebox{0.25ex}{\tikz\filldraw[black, x=1pt, y=1pt](0,0) circle (1);}\;}

\DeclareFontFamily{U}{MnSymbolC}{}
\DeclareSymbolFont{MnSyC}{U}{MnSymbolC}{m}{n}
\DeclareFontShape{U}{MnSymbolC}{m}{n}{
    <-6>  MnSymbolC5
   <6-7>  MnSymbolC6
   <7-8>  MnSymbolC7
   <8-9>  MnSymbolC8
   <9-10> MnSymbolC9
  <10-12> MnSymbolC10
  <12->   MnSymbolC12}{}
\DeclareMathSymbol{\intprod}{\mathbin}{MnSyC}{'270}

\renewcommand{\emptyset}{\varnothing}

\DeclarePairedDelimiterX\abs[1]\lvert\rvert{
  \ifblank{#1}{\:\cdot\:}{\,#1\,}
}
\DeclarePairedDelimiterX\norm[1]\lVert\rVert{
  \ifblank{#1}{\:\cdot\:}{\,#1\,}
}
\DeclarePairedDelimiterX\set[1]{\lbrace}{\rbrace}{\,#1\,}
\DeclarePairedDelimiterX\Span[1]{\langle}{\rangle}{\,#1\,}
\DeclarePairedDelimiterX\floor[1]{\lfloor}{\rfloor}{\,#1\,}
\DeclarePairedDelimiterX\Inner[2]{\langle}{\rangle}{\,#1,\,#2\,}
\DeclarePairedDelimiterX\Set[2]{\lbrace}{\rbrace}{\,#1\ \delimsize\vert\ #2\:}
\DeclarePairedDelimiterX\lb[2]{[}{]}{\,\ifblank{#1}{-}{#1}\,,\,\ifblank{#2}{-}{#2}\,}
\DeclarePairedDelimiterX\pb[2]{\lbrace}{\rbrace}{\,\ifblank{#1}{-}{#1}\,,\,\ifblank{#2}{-}{#2}\,}

\makeatletter
\newcommand{\etabar}{\text{\eta@bar}}
\newcommand{\eta@bar}{%
  \vphantom{$\m@th \eta$}%
  \ooalign{%
    $\m@th \eta$\cr
    \hidewidth\kern.25em\smash{\raisebox{-0.7ex}{$\m@th\mathchar'55$}}\hidewidth\cr}%
}
\makeatother

\newcommand{\half}{{\textstyle{1\over2}}}
\newcommand{\fifth}{{\textstyle{1\over5}}}
\newcommand{\third}{{\textstyle{1\over3}}}

\newcommand{\twothirds}{{\textstyle{2\over3}}}
\newcommand{\fourthirds}{{\textstyle{4\over3}}}
\newcommand{\fivethirds}{{\textstyle{5\over3}}}

\newcommand{\eightthirds}{{\textstyle{8\over3}}}

\makeatletter
\renewcommand*\env@matrix[1][\arraystretch]{%
  \edef\arraystretch{#1}%
  \hskip -\arraycolsep
  \let\@ifnextchar\new@ifnextchar
  \array{*\c@MaxMatrixCols c}}
\makeatother

\makeatletter
\renewenvironment{abstract}{%
    \small\thispagestyle{empty}
    \null\vfil
    {\textcolor{RoyalBlue}{\scshape\abstractname.}}
    \quotation
    }
{\endquotation\vfil\null\clearpage}
\makeatother

\newcommand{\TheEnd}{
\bigskip\bigskip
\begin{center}
  \Large
  \decofourleft\hspace*{0.5em}\floweroneleft\hspace*{0.5em}\decoone\hspace*{0.5em}\floweroneright\hspace*{0.5em}\decofourright
\end{center}
\bigskip\bigskip}

\usepackage[acronym, toc]{glossaries}

\makeglossaries

\newacronym{cp}{CP}{\nm{Cauchy} problem}
\newacronym{cpu}{CPU}{Central Processing Unit}
\newacronym{pde}{PDE}{Partial Differential Equation}
\newacronym{ivp}{IVP}{Initial Value Problem}
\newacronym{ibe}{iBE}{inviscid \nm{Burgers} Equation}
\newacronym{2d}{2D}{Two dimensions}
\newacronym{3d}{3D}{Three dimensions}

\begin{document}

\title[\Title]{\Longtitle}

\author[D.~Dutykh]{Denys Dutykh\textcolor{denimblue}{$^*$}}
\address{\textcolor{denimblue}{\bf D.~Dutykh:} Mathematics Department, Khalifa University of Science and Technology, PO Box 127788, Abu Dhabi, United Arab Emirates}
\email{\href{mailto:Denys.Dutykh@ku.ac.ae}{Denys.Dutykh@ku.ac.ae}}
\urladdr{\url{http://www.denys-dutykh.com/}}
\thanks{\textcolor{denimblue}{$^*$}\it Corresponding author}

\author[E.~Leichtnam]{\'Eric Leichtnam}
\address{\textcolor{denimblue}{\bf E.~Leichtnam:} CNRS, Institut de Math\'ematiques de Jussieu--Paris Rive Gauche UMR 7586, UP7D -- B\^atiment Sophie Germain, 56-58 avenue de France --- BC 7012, 75205 Paris Cedex 13}
\email{\href{mailto:eric.leichtnam@imj-prg.fr}{Eric.Leichtnam@imj-prg.fr}}
\urladdr{\url{https://webusers.imj-prg.fr/~eric.leichtnam/}}



\begin{titlepage}
\clearpage
\pagenumbering{arabic}
\thispagestyle{empty} 
\noindent
{\Large Denys \textsc{Dutykh}}
\\[0.001\textheight]
{\textit{\textcolor{gray}{Khalifa University of Science and Technology, Abu Dhabi, UAE}}}
\\[0.02\textheight]
{\Large \'Eric \textsc{Leichtnam}}
\\[0.001\textheight]
{\textit{\textcolor{gray}{CNRS, Institut de Math\'ematiques de Jussieu--Paris Rive Gauche, France}}}
\\[0.16\textheight]

\vspace*{2.99cm}

\colorbox{Lightblue}{
  \parbox[t]{1.0\textwidth}{
    \centering\huge
    \vspace*{0.75cm}
    
    \textsc{\textcolor{katyblue}{\Longtitle}}
    
    \vspace*{0.75cm}
  }
}

\vfill 

\raggedleft     
{\large \plogo} 
\end{titlepage}


\clearpage
\thispagestyle{empty} 
\par\vspace*{\fill}   
\begin{flushright} 
{\textcolor{RoyalBlue}{\textsc{Last modified:}} \today}
\vspace*{1.0em}
\end{flushright}


\clearpage
\maketitle
\thispagestyle{empty}


\begin{abstract}

In this manuscript, we highlight a new phenomenon of complex algebraic singularities formation for solutions of a large class of genuinely nonlinear Partial Differential Equations (\acrshort{pde}s). We start from a unique \nm{Cauchy} datum which is holomorphic ramified like $x_{\,1}^{\,\frac{1}{k\,+\,1}}$ (and its powers) around the \textbf{smooth} locus $x_{\,1}\ =\ 0$ and is sufficiently singular. Then, we expect the existence of a solution which should be holomorphic ramified around the \textbf{singular} locus $\S$ defined by the vanishing of the discriminant of an algebraic equation of degree $k\ +\ 1\,$. Notice, moreover, that the monodromy of the \nm{Cauchy} datum is \textbf{Abelian}, whereas one of the solutions is \textbf{non-Abelian} and that $\S$ depends on the \nm{Cauchy} datum in contrast to the \nm{Leray} principle (stated for linear problems only). This phenomenon is due to the fact that the \acrshort{pde} is genuinely nonlinear and that the \nm{Cauchy} datum is sufficiently singular. First, we investigate the case of the \acrfull{ibe}. Later, we state a general \cref{conj:main}, which describes the expected phenomenon. We view this \cref{conj:main} as a working programme allowing us to develop interesting new Mathematics. We also state \cref{conj:2}, which is a particular case of the general \cref{conj:main} but keeps all the flavour and difficulty of the subject. Then, we propose a new algorithm with a map $\FF$ such that a fixed point of $\FF$ would give a solution to the problem associated with \cref{conj:2}. Then, we perform convincing, elaborate numerical tests which suggest that a \nm{Banach} norm should exist for which the mapping $\FF$ should be a contraction so that the solution (with the above specific algebraic structure) should be unique. This work is a continuation of \cite{Leichtnam1993}.

\bigskip\bigskip
\noindent \textcolor{RoyalBlue}{\textbf{\keywordsname:}} \Keywords \\

\smallskip
\noindent \textcolor{RoyalBlue}{\textbf{MSC:}} 
\smallskip \\
\noindent \textcolor{RoyalBlue}{\textbf{PACS:}} 

\end{abstract}


\newpage
\pagestyle{empty}
\tableofcontents
\clearpage


\newpage
\listoffigures
\listoftables
\clearpage
\pagestyle{fancy}





\clearpage
\bigskip\bigskip
\section*{Preface}
\addcontentsline{toc}{section}{Preface}
\thispagestyle{empty}

In this manuscript, we investigate the process of singularity formation in two examples of fully (or genuinely) nonlinear \acrfull{pde}s starting from a single algebraic singular \nm{Cauchy} datum. Namely, we consider specific ramified initial value (or \textsc{Cauchy}-type) problems. In order to better understand the real singularities, our approach consists of looking at what happens in the complex domain. Thus, we formally assume that space and time variables take complex values. The analysis of the \acrshort{pde}s in the complex domain sheds new light on the origin of familiar real singularities. In the case of the \acrfull{ibe} (with an algebraic singular \nm{Cauchy} datum), we are able to carry out all computations analytically using two completely different methods: \nm{Cauchy}--\nm{Kovalevskaya} theory and contact geometry on the space of $1-$jets. We propose a generalisation to the case of genuinely nonlinear \acrshort{pde}s of orders two or higher, and we state several conjectures. Moreover, for a second order \acrshort{pde}, we propose an iterative scheme, which allows us to construct efficiently approximate complex solutions. The convergence of this scheme is empirically demonstrated on a sufficient number of numerical examples, which is a good indication in favour of our conjectures. We stress that the proposed algorithm is entirely new. This scheme might also be used to show the existence and uniqueness of solutions to the ramified \acrfull{ivp} corresponding to our conjectures.

In more precise mathematical terms, in this manuscript, we highlight a new phenomenon of complex algebraic singularities formation for solutions of genuinely nonlinear \acrfull{pde}
\begin{equation*}
  P\,(\,u\,)\,(\,t,\,\x\,)\ =\ 0\,, \qquad (\,t,\,\x\,)\ \in\ \C\times\C^{\,n}
\end{equation*}
starting from a unique \nm{Cauchy} datum, which is sufficiently singular. More precisely, we start from a unique \nm{Cauchy} datum $u\,(\,0,\,\x\,)$ which is ramified near the origin like $x_{\,1}^{\,\frac{1}{k\,+\,1}}$ around the smooth hyper surface $x_{\,1}\ =\ 0$ and assume that $P$ is simply characteristic with respect to the co-normal of $x_{\,1}\ =\ 0\,$. Then, we expect the existence of a solution to $P\,(\,u\,)\ =\ \fO$ of a specific form (ansatz) $u\,(\,t,\,\x\,)\ =\ \A\,\bigl(\,t,\,\x,\,z\,(\,t,\,\x\,)\,\bigr)\,$, where $\A$ is a holomorphic function near the origin and $z\,(\,t,\,\x\,)$ is a solution of an algebraic equation
\begin{equation}\label{eq:1.1}
  z^{\,k\,+\,1}\,(\,t,\,\x\,)\ -\ a_{\,k\,-\,1}\,(\,t,\,\x\,)\,z^{\,k\,-\,1}\,(\,t,\,\x\,)\ -\ \ldots\ -\ a_{\,0}\,(\,t,\,\x\,)\ =\ 0\,.\tag{$\ast$}
\end{equation}
Thus, the solution is holomorphic ramified around the singular locus defined by the discriminant (swallow-tail) of \cref{eq:1.1}. Since we search for a solution under a special form which is ramified along a single (singular) hyper-surface, we need only a single \nm{Cauchy} datum. Actually, ``morally'' the solution is unique once we have fixed the choice of a root to the polynomial equation defined by the simple characteristic hypothesis. The fact that the solution is ramified around a singular locus and has a non-\nm{Abelian} monodromy, whereas the \nm{Cauchy} datum is ramified around a smooth locus and has \nm{Abelian} monodromy is a (new) phenomenon which is due to the facts that the \nm{Cauchy} datum is sufficiently singular and that the \acrshort{pde} $P\,(\,u\,)\ =\ \fO$ is genuinely nonlinear. Let us provide two examples to illustrate the concept of genuine nonlinearity, while the more general definition will be given in the formulation of \cref{conj:main}. Let $u$ be a holomorphic function
\begin{align*}
  u\,:\ \C^{\,n\,+\,1}\ &\longrightarrow\ \C \\
  (\,t,\,x_{\,1},\,x_{\,2},\,\ldots,\,x_{\,n}\,)\ &\mapsto\ u\,(\,t,\,x_{\,1},\,x_{\,2},\,\ldots,\,x_{\,n}\,)\,.
\end{align*}
Then, the following `initial' value problem is genuinely nonlinear
\begin{equation*}
  u_{\,t\,t}\ -\ u_{\,x_{\,1}}\,(\,u_{\,x_{\,1}\,x_{\,1}}\ +\ u_{\,x_{\,2}\,x_{\,2}}\ +\ \ldots\ u_{\,x_{\,n}\,x_{\,n}})\ =\ \fO\,, \qquad u\,(\,0,\,\x\,)\ =\ c_{\,1}\,x_{\,1}\ +\ c_{\,2}\,x_{\,1}^{1\,+\,\frac{1}{3}}\,,
\end{equation*}
and the next one is not
\begin{equation*}
  u_{\,t\,t}\ -\ u_{\,x_{\,2}}\,(\,u_{\,x_{\,1}\,x_{\,1}}\ +\ u_{\,x_{\,2}\,x_{\,2}}\ +\ \ldots\ u_{\,x_{\,n}\,x_{\,n}})\ =\ \fO\,, \qquad u\,(\,0,\,\x\,)\ =\ c_{\,1}\,x_{\,1}\ +\ c_{\,2}\,x_{\,1}^{1\,+\,\frac{1}{3}}\,,
\end{equation*}
where $c_{\,1,\,2}\ \in\ \C^{\,\times}\,$.

The full generalization of this phenomenon is stated in \cref{conj:main} for genuinely nonlinear \acrshort{pde} of order greater or equal to two. Simpler versions of this Conjecture are stated in \cref{conj:2} and \cref{conj:m}. But \cref{conj:2} keeps all the flavour and the conceptual difficulty of \cref{conj:main}. \cref{conj:2} considers the following \acrshort{pde}  
\begin{equation*}
  u_{\,t\,t}\ -\ u_{\,x}\,u_{\,x\,x}\ =\ \fO
\end{equation*}
along with the \nm{Cauchy} datum
\begin{equation*}
  u\,(\,0,\,x\,)\ =\ \sum_{j\,=\,1}^{N_{\,0}}\,c_{\,j}\,x^{\,1\,+\,\frac{j\,-\,1}{3}}\,,
\end{equation*}
where $N_{\,0}\ \in\ \N_{\,>\,2}\,$, $c_{\,1,\,2}\ \in\ \C^{\,\times}$ and $\set*{c_{\,j}}_{\,j\,=\,3}^{\,N_{\,0}}\ \subseteq\ \C$ are constants. Moreover, a root of the polynomial equation associated with the simple characteristic hypothesis is chosen. Then, \cref{conj:2} states the existence of an algebraic solution $u\,(\,t,\,x\,)$ ramified along a cusp
\begin{equation*}
  4\,p^{\,3}\,(\,t,\,x\,)\ -\ 27\,q^{\,2}\,(\,t,\,x\,)\ =\ 0
\end{equation*}
depending on the \nm{Cauchy} datum and the previous choice. This conjecture is a working programme allowing the development of interesting new Mathematics. The algebraic and geometric studies carried out in Sections~\ref{sec:ring} and \ref{sec:dmod} allow us to construct, in \cref{sec:sec}, a new algorithm with a map $\FF$ such that a fixed point of $\FF$ provides a solution to \cref{conj:2}. In \cref{sec:num}, we carry out various elaborate numerical experiments which show quite convincingly that the iterations $\FF^{\,\comp\,m}\,(\,w_{\,0}\,)$ seem to converge (as $m\ \to\ +\,\infty$) to a unique fixed point. This suggests that a family of semi-norms should exist, allowing the construction of a \nm{Banach} norm for which $\FF$ would be contracting. However, this is a very difficult problem, much harder than the one treated successfully in \cite{Leichtnam1987} corresponding to the case of algebraic equations of degree two (see also \cref{sec:hope} for a short reminder).

In the case of the \acrfull{ibe}, we are able to carry out all the computations analytically using two completely different methods: the \nm{Cauchy}--\nm{Kovalevskaya} theorem and the contact geometry approach on the space of $1-$jets. This new phenomenon is exhibited very clearly in the case of the \acrshort{ibe}. For the most general case of genuinely nonlinear \acrshort{pde}s of the order higher or equal to two, we state the general \cref{conj:main}, which describes the new phenomenon previously mentioned very precisely. \cref{conj:main} is a very difficult problem. We view it as a motivation to develop new Mathematics. Indeed, in this paper, we try to lay the foundations of proof. We also state Conjectures~\ref{conj:2} and \ref{conj:m}, which are particular cases of the general \cref{conj:main}. The \cref{conj:2} deals with the case $k\ =\ 2$ and
\begin{align}\label{eq:1.2a}
  P\,(\,u\,)\ \coloneq\ u_{\,t\,t}\ -\ u_{\,x}\,u_{\,x\,x}\ &=\ \fO\,, \tag{$\ast\ast_{1}$} \\
  u\,(\,0,\,x\,)\ &=\ \sum_{j\,=\,1}^{N_{\,0}}\,c_{\,j}\,x^{\,1\,+\,\frac{j\,-\,1}{3}}\,, \tag{$\ast\ast_{2}$}\label{eq:1.2b}
\end{align}
where $N_{\,0}\ \geq\ 3\,$, $c_{\,1,\,2}\ \in\ \C^{\,\times}$ and $c_{\,j}\ \in\ \C\,$, $(\,3\ \leq\ j\ \leq\ N_{\,0}\,)$ are constants, and $(\,t,\,x\,)\ \in\ \C^{\,2}\,$. Morally, \cref{conj:2} is not much easier than \cref{conj:main}, but its simpler formulation allows for direct numerical verifications.

We revisit and improve the geometric study of the algebra $\O\,[\,z\,]\,$, where $z^{\,3}\ =\ p\,z\ +\ q$ and $\O$ denotes the algebra of germs of holomorphic functions $\set*{a\,(\,p,\,q\,)}$ near the origin $\vO\ \in\ \C^{\,2}\,$. This study allows us to construct a new algorithm (and the corresponding ansatz) along with a mapping $\FF$ such that a fixed point of $\FF$ gives a solution to the Problem~\eqref{eq:1.2a}, \eqref{eq:1.2b}. We provide various numerical experiments show that that the iterations $\FF^{\comp\,m}\,(\,w_{\,0}\,)$ seem to converge numerically (when $m\ \to\ +\,\infty$) to a fixed point which does not depend on $w_{\,0}\,$. This provides strong empirical evidence that this algorithm will provide a unique solution to the Problem~\eqref{eq:1.2a}, \eqref{eq:1.2b} of the form $\A\,\bigl(\,t,\,\x,\,z\,(\,t,\,\x\,)\,\bigr)$ once a choice of a root to the simply characteristic equation (of degree two here) has been chosen. It remains a very difficult open problem to find appropriate a \nm{Banach} (normed) algebra for which
\begin{equation*}
  \norm{\FF\,(\,w\,)\ -\ \FF\,(\,v\,)}\ \leq\ c\,\norm{w\ -\ v}\,,
\end{equation*}
with $0\ <\ c\ <\ 1\,$, which would give a proof of \cref{conj:2}.

\vspace*{4em}
\begin{flushright}\noindent
Abu Dhabi,\hfill {\it Denys Dutykh} \\
Saint-Jean-de-Chevelu,\hfill {\it Eric Leichtnam} \\
December 2023\hfill { }
\end{flushright}


\clearpage
\bigskip\bigskip
\section{Introduction}
\label{sec:intro}

The process of formation of singularities to solutions of nonlinear \acrlong{pde}s is a huge and important topic which has stimulated a lot of interesting works (see \eg \cite{Popivanov2003, Weideman2003}). We begin by recalling some known results in this topic in order to describe briefly the landscape in which our results and goals will be realized. 

Even infinitely smooth initial data under the dynamics of a nonlinear \acrshort{pde} will not remain in general smooth for all times. Henceforth, the topic of the singularity formation in solutions to nonlinear \acrshort{pde}s has been central in the study of these equations. Perhaps the wave-breaking phenomenon is the most familiar and far from being completely studied, \emph{natural} singularity formation processes \cite{Banner1993}. Moreover, the wave-breaking process can be easily seen by everyone.

An attempt to classify various \emph{finite time} singularities in \acrshort{pde}s has been made in \cite{Eggers2007}. Physical (and several other) applications have clearly motivated this classification. The question of eventual blow-up in \acrshort{3d} incompressible \nm{Euler} and \nm{Navier}--\nm{Stokes} equations is open and central to many current theoretical and numerical investigations \cite{Hou2009}. We also refer to \cite{Caflisch2015} as an excellent review of available numerical approaches to detect complex singularities of \acrshort{pde}s. Complex singularities of the \nm{Lorenz} dynamical system (with complex time) have been studied in \cite{Viswanath2010}.

The present manuscript is devoted in the first place to the study of complex algebraic singularities in the \acrfull{ibe}. Namely, we study the algebraic singularity blow-up under the dynamics of some first and second-order genuinely nonlinear \acrshort{pde}s. Among the second-order \acrshort{pde}s, we focus on a particular second-order \acrshort{pde} belonging to the celebrated family of $p-$systems described in some detail below. 

The \acrshort{ibe} can be written as \cite{Burgers1948}:
\begin{equation}\label{eq:burg}
  u_{\,t}\ -\ u\,u_{\,x}\ =\ \fO\,,
\end{equation}
where in the (classical) real case 
\begin{align*}
  u\,:\ \R_{\,\geq\,t_{\,0}}\,\times\,\R\ &\longrightarrow\ \R\,, \\
  (\,t,\,x\,)\ &\mapsto\ u\,(\,t,\,x\,)\,.
\end{align*}
The subscripts $(\,-\,)_{\,t}$ and $(\,-\,)_{\,x}$ denote the usual partial derivatives $\pd{}{t}$ and $\pd{}{x}$ respectively. In order to obtain an \acrfull{ivp}, Equation~\eqref{eq:burg} has to be completed by an appropriate initial condition:
\begin{equation*}
  u\,(\,t_{\,0},\,x\,)\ =\ u_{\,0}\,(\,x\,)\,, \qquad x\ \in\ \R\,.
\end{equation*}
It is well-known that the \acrshort{ibe} will develop the gradient catastrophe in finite time from a generic initial condition. More precisely, the existence time $t_{\,0}\ \leq\ t\ <\ t_{\,s}$ depends on the initial condition:
\begin{equation*}
  t_{\,s}\ =\ \frac{1}{\sup\limits_{x\,\in\,\R} u_{\,0,\,x}\,(\,x\,)}\,.
\end{equation*}
Please note that the real singularity never occurs (in forward time) if $u_{\,0,\,x}$ takes strictly negative values for all finite $x\ \in\ \R\,$. Below, for the sake of simplicity, we shall take $t_{\,0}\ =\ 0\,$.

The real singularity is point-like in the sense that it occurs locally in a single point at a fixed time $t_{\,s}\,$. In contrast to real singularities, we shall demonstrate below that complex algebraic singularities take place on sets of positive dimensions (\ie analytical hyper-surfaces in $\C^{\,2}$ in our case). For the sake of completeness, we would like to mention that some singularities may even happen on sets of fractional \nm{Hausdorff} dimensions \cite{Fontelos2000}. This mechanism is conjectured for the \nm{Navier}--\nm{Stokes} equations in \acrshort{3d}. We can also remember an old idea of J.~\nm{Leray} that singularities in the \nm{Navier}--\nm{Stokes} equations could be related to the phenomenon of turbulence \cite{Boldrighini2012}. However, this idea was not followed, up to now and to the best of our knowledge, by any significant results.

In a similar line of thinking, we shall consider also in the present manuscript the following second-order nonlinear hyperbolic equation:
\begin{equation}\label{eq:sec}
  u_{\,t\,t}\ -\ u_{\,x}\,u_{\,x\,x}\ =\ \fO\,,
\end{equation}
whose \acrshort{ivp} requires traditionally two data to be specified:
\begin{align}\label{eq:ivp0}
  u\,(\,0,\,x\,)\ &=\ u_{\,0}\,(\,x\,)\,, \qquad x\ \in\ \R\,, \\
  u_{\,t}\,(\,0,\,x\,)\ &=\ u_{\,1}\,(\,x\,)\,, \qquad x\ \in\ \R\,.
\end{align}
However, we shall consider a ramified \nm{Cauchy} problem with only one (algebraic) \nm{Cauchy} datum \eqref{eq:ivp0} and we shall seek a solution $(\,t,\,x\,)\ \mapsto\ u\,(\,t,\,x\,)$ of a very special (algebraic) form so that only the first \nm{Cauchy} datum is needed. Both considered Equations~\eqref{eq:burg} and \eqref{eq:sec} come from the same family of \acrshort{pde}s discussed in Section~\ref{sec:scope}. Mathematically speaking, we discuss in the present study the hyperbolic sub-family of genuinely nonlinear \acrshort{pde}s. The general notion of genuinely nonlinear \acrshort{pde} that we use is stated as Assumption (2) in \cref{conj:main}. For example, the following \acrshort{pde}
\begin{equation*}
  u_{\,t\,t}\ -\ u_{\,x}\,u_{\,x\,x}\ =\ \fO
\end{equation*}
is genuinely nonlinear, whereas
\begin{equation*}
  u_{\,t\,t}\ -\ u\,u_{\,x\,x}\ =\ \fO
\end{equation*}
is not genuinely nonlinear but only quasilinear. Other examples will be given in \cref{sec:scope}.

\begin{example}\label{eq:toy}
Let us consider the following \nm{Cauchy} `toy' problem in $\C^{\,2}\,$:
\begin{align*}
  \pd{u}{t}\ &=\ u^{\,2}\,, \\
  u\,(\,0,\,x\,)\ &=\ \frac{1}{x}\,.
\end{align*}
It is not difficult to see that the exact solution to this problem is 
\begin{equation*}
  u\,(\,t,\,x\,)\ =\ \frac{1}{x\ -\ t}\,.
\end{equation*}
This simple example shows that it is absolutely crucial to assume that \nm{Cauchy} data \eqref{eq:ivp0} is bounded if we want to establish the existence of a holomorphic solution ramified along a characteristic hyper-surface. However, we mention that this toy example is not representative of our study for the following reasons:
\begin{enumerate}
  \item The singular locus $\S\ \coloneq\ \Set*{t\ \in\ \C}{(\,t,\,t\,)}$ is a regular hyper-surface.
  \item The hyper-surface $\S$ is not characteristic for the differential operator $\pd{}{t}\,$.
\end{enumerate}
Henceforth, our example has some serious limitations.
\end{example}

It is well-known that the behaviour of the real-valued and complex-valued solutions is quite different. We shall mention the example of the \acrshort{2d} viscid \nm{Burgers} system. For real-valued solutions O.~\nm{Ladyzhenskaya} proved\footnote{This result is a kind of folklore theorem because the original reference \cite{Ladyzhenskaya2014} does not contain a very detailed proof.} in $1963$ global existence and uniqueness result in \nm{Sobolev} spaces \cite{Ladyzhenskaya2014}. Later, it turned out that complex-valued solutions differ drastically from real ones. In particular, \nm{Li} and \nm{Sinai} proved in \cite{Li2006a} using the renormalization group method that complex-valued solutions (in $n$D) can develop finite time singularities. This was proven earlier in \cite{Polacik2008} for gradient-like solutions. The \acrshort{2d} complex \nm{Burgers} equation was studied in \cite{Li2010}, where an open set of a six-parameter family of initial conditions is constructed such that the corresponding solutions exhibit blowups in finite time. Finally, complex-valued initial conditions have been showing \emph{numerically} to develop a singularity in finite time in complete agreement with theory. These results have also been extended to \acrshort{3d} complex-valued \nm{Navier}--\nm{Stokes} equations \cite{Li2008c}. The complex singularities in \acrshort{2d} \nm{Euler} equations have been studied in \cite{Gilbert2011}.

Our study differs at least in two important respects from the previously described line of thinking:
\begin{enumerate}
  \item We consider hyperbolic and inviscid equations starting with the \acrshort{ibe}.
  \item We complexify not only the dependent variables but also the independent ones.
\end{enumerate}
The viscous \nm{Burgers} equation with real-time and complex spatial variables has been considered in \cite{Gal2013} in the context of semi-groups. Similarly, the linear wave (or \nm{D'Alembert}) and telegraph equations were treated in \cite{Gal2012}. Perhaps the complex view on \acrshort{pde}s can be traced back to the work of S.~\nm{Kovalevskaya} who considered the \acrfull{cp} for the usual (linear) heat equation with complex temporal and spatial variables \cite{Kovalevskaya1875}. In particular, she showed that there are examples of holomorphic non-entire initial conditions such that the power series solution does not converge in any neighbourhood of $t_{\,0}\ \in\ \C\,$. The \nm{Borel} summability of these divergent solutions has been studied much later \cite{Lutz1999}. However, our study's goal is to shed new light from the complex geometry angle on the formation of real singularities in fully nonlinear \acrshort{pde}s. The complex \acrshort{ibe} appeared surprisingly also in random matrix theory \cite{Matytsin1994} and random surface models \cite{Kenyon2007}.

In the present work, we focus first on the classical \acrshort{ibe} and the following ramified \acrshort{ivp}:
\begin{align*}
  u_{\,t}\ -\ u\,u_{\,x}\ &=\ \fO\,, \\
  u\,(\,0,\,x\,)\ &=\ x^{\,\third}\,.
\end{align*}
We show that the last \acrshort{ivp} possesses solutions of the form $u\ =\ a\ +\ z\,$, where $a$ is a regular holomorphic function and $z$ verifies the cubic equation $z^{\,3}\ =\ p\,z\ +\ q\,$. The obtained solution can be analytically continued along any path in a small neighbourhood of the origin $\vO\ \in\ \C^{\,2}$ and originating from $(\,0,\,x_{\,0}\,)\ \in\ \C^{\,2}$ and avoiding the cusp singularity $4\,p^{\,3}\ -\ 27\,q^{\,2}\ =\ \fO\,$. This is explained in detail in \cref{sec:burg}. Thus, the solution $u$ is ramified along the \emph{singular} locus $4\,p^{\,3}\ -\ 27\,q^{\,2}\ =\ \fO$ and its monodromy group\footnote{Informally speaking, the monodromy group is a group of transformations that encodes what happens with data when we turn around the singularity. More formally, a monodromy is a representation of the fundamental group. Let $\Vv$ be a holomorphic connected variety. Let us take $x_{\,0}\ \in\ \Vv$ and $u$ be a holomorphic germ defined at point $x_{\,0}\,$, which can be prolongated along any closed loop issued from $x_{\,0}\,$. If $\gamma$ is any such closed loop at $x_{\,0}\,$, by $\gamma\cdot u$ we denote the result of this prolongation. Finally, by $\E$, we denote the vector space generated by all such elements $\gamma\cdot u\,$. If $\dim\,\E\ <\ +\infty\,$, we say that $u$ is of \emph{finite determination}. The homotopy equivalence classes $[\,\gamma\,]$ form the fundamental group $\pi_{\,1}\,(\,\Vv\,)\,$. Finally, the monodromy is defined as the following map:
\begin{align*}
  \pi_{\,1}\,(\,\Vv\,)\ &\longrightarrow\ \GL\,(\,\E\,)\,, \\
  \gamma\ &\mapsto\ (\,u\ \mapsto\ \gamma\cdot u\,)\,.
\end{align*}
Of course, it has to be shown that all these maps are well-defined, and the result does not depend on the representative $\gamma$ of the homotopy class $[\,\gamma\,]\,$. Fortunately, it can be done without any difficulties.} (\cf Remarks~\ref{rem:4.1} and \ref{rem:4.2}) is the (non-commutative) permutation group $S_{\,3}$ (\cf \cref{rem:4.1}) whereas the \nm{Cauchy} datum $x^{\,\frac{1}{3}}$ is ramified around a \emph{smooth} locus with a commutative monodromy group, namely $\Z\,/\,3\,\Z$\,. This phenomenon (including the changing of the monodromy group) is due both to the genuine nonlinearity of the \acrshort{ibe} and the fact the algebraic \nm{Cauchy} datum is sufficiently singular (see Section~\ref{sec:lin} and \cref{thm:2.2}). \textbf{However, the main goal of our study is to formulate a convincing numerical convergence evidence along with the generalisation of this phenomenon for genuinely nonlinear \acrshort{pde}s of the order two and higher (see Conjectures~\ref{conj:2}, \ref{conj:m} and \ref{conj:main}).} A first generalization is stated in \cref{conj:2} which considers the \acrfull{pde} \eqref{eq:sec} with the initial condition
\begin{equation*}
  u\,(\,0,\,x\,)\ =\ \sum_{j\,=\,1}^{N_{\,0}}\,c_{\,j}\,x^{\,1\,+\,\frac{j\,-\,1}{3}}\,,
\end{equation*}
where $N_{\,0}\ \in\ \N_{\,>\,2}\,$, $c_{\,1,\,2}\ \in\ \C^{\,\times}$ and $\set*{c_{\,j}}_{\,j\,=\,3}^{\,N_{\,0}}\ \subseteq\ \C$ are constants. It states the existence of an algebraic solution $u\,(\,t,\,x\,)$ ramified along a cusp $4\,p^{\,3}\,(\,t,\,x\,)\ -\ 27\,q^{\,2}\,(\,t,\,x\,)\ =\ 0$ \emph{depending} on the \nm{Cauchy} datum. We strongly believe that Conjectures~\ref{conj:2}, \ref{conj:m} and the most general \cref{conj:main} are important problems requiring the development of interesting new Mathematics. 

This manuscript is the first step in which we shall provide solid mathematical foundations and numerical evidence to believe that they are true. Especially regarding the \cref{conj:2}, we shall construct in \cref{sec:sec} a new (and highly non-trivial) algorithm, allowing us to obtain quite convincing numerical convergence results in \cref{sec:num}. Moreover, we sketch a tentative strategy to theoretically address the general \cref{conj:main}. In Sections \ref{sec:ring} and \ref{sec:dmod}, we revisit the general results of \cite{Leichtnam1993} in the specific case of the cubic equation $z^{\,3}\ =\ p\,z\ +\ q$ and make them more precise. The big advantage here is that most computations can be done almost explicitly, which allows us to highlight the underlying structures, which will allow us to construct a new practical algorithm from \cref{sec:sec}. More precisely, in Section~\ref{sec:ring} we revisit the ring\footnote{Throughout this article, the word ``ring'' means ``ring with identity''.} $\O\,\llbracket\,z\,\rrbracket$ (with $z^{\,3}\ =\ p\,z\ +\ q$) and we introduce the primitive operator $\partial_{\,q}^{\,-\,1}\,$. The reason why the derivation $\partial_{\,q}$ operator governs the structure of the ring $\O\,\llbracket\,z\,\rrbracket$ is elucidated in Section~\ref{sec:dmod}. In the same Section we establish the fact that $(\,z,\,z^{\,2}\,)$ is a solution of a holonomic $\D-$module with the characteristic variety $\Vv$ which is included in the union of the zero section and the co-normal\footnote{Let us remind that the co-normal is constituted of all co-tangent vectors which annihilate on the tangent space to a given manifold (here, an analytic curve or surface).} to the cusp $4\,p^{\,3}\ -\ 27\,q^{\,2}\ =\ \fO\,$:
\begin{equation*}
  \Vv\ \bigcap\ \T_{\,\vO}^{\,\ast}\,\C^{\,2}\ =\ \bigl(\,0,\,0\,;\,0,\,\C\,\bigr)\,.
\end{equation*}
In fact, using a deep result of \nm{Kashiwara}, we are able to prove an even stronger result stating that the characteristic variety $\Vv$ \emph{is} the union of the zero section and the co-normal to the cusp. The characteristic variety is an important geometric invariant of a $\D-$module. The integrability of the characteristic variety is a central result in the theory of $\D-$modules. Loosely speaking, this result says that the singular support of a $\D-$module is an involutive sub-variety in the co-tangent bundle. The involutivity property may be seen under the sheaf angle: the ideal sheaf defining the singular support is closed with respect to the natural \nm{Poisson} bracket on the co-tangent bundle \cite{Singh2014}.

Let us remind you of some additional relevant background material which will help you understand the sequel of this Section. We denote by $\D$ the sheaf of holomorphic differential operators on $\C^{\,2}\ \cong\ \set*{(\,p,\,q\,)}\,$. The vector space of sections of $\D$ over an open subset $\U$ is denoted by $\D\,(\,\U\,)\,$.
\begin{definition}
Let $M$ be a sheaf of $\D-$modules ($\D-$module for short) over an open subset $\U$ of $\C^{\,2}\,$. The $\D-$module $\Mm$ is said to be coherent if for any sheaf homomorphism of $\D-$modules
\begin{equation*}
  \phi\,:\ \D^{\,m}\ \longrightarrow\ \Mm
\end{equation*}
and for any point $z$ of $\U\,$, one can find an open neighbourhood $\U_{\,z}$ of $z$ and a finite number of sections $s_{\,1}\,$, $s_{\,2}\,$, \ldots, $s_{\,r}$ (over $\U_{\,z}$) of the sheaf $(\,\ker\,\phi\,)$ such that for any open subset $\W\ \subseteq\ \U_{\,z}\,$, the restrictions $s_{\,1}\,\vert_{\,\W}\,$, $s_{\,2}\,\vert_{\,\W}\,$, \ldots, $s_{\,r}\,\vert_{\,\W}$ generate $(\,\ker\,\phi\,)\,(\,\W\,)$ as a module over $\D\,(\,\W\,)\,$.
\end{definition}

The fact that above the origin, there is only one line (the co-normal to $\partial_{\,q}$) in $\Vv$ is absolutely crucial for our constructions. In particular, it explains why the primitive $\partial_{\,q}^{\,-\,1}$ exists and why the ``micro-local singularities''\footnote{We take this expression in the quotes because the distributions $z$ and $z^{\,2}$ are not yet well-defined in the real case \emph{stricto sensu}. It will be a topic of our future works.} of the product $z^{\,2}$ remain under the control above the origin $\vO\ \in\ \C^{\,2}\,$. This observation will constitute one of the key ingredients in the future proof of \cref{conj:2}. However, before engaging in this huge endeavour, we thought obtaining convincing numerical tests to converge our new iterative algorithm was indispensable. Indeed, its convergence was not obvious \emph{a priori}, and its proof will constitute a big challenge. We end with a short heuristic explanation of why we believe that if we impose a second \nm{Cauchy} datum to \cref{conj:2}, then there will be no solution in general. Indeed, from the linear result \cite{Leichtnam1993}, two cusps will show up. The products of $u$ and its derivatives will provide a spreading of singularities (associated with the co-normal of each of the cusps) over the origin, which will not be possible to control. We shall return to this question in the forthcoming works.

The present manuscript is organized as follows. In Section~\ref{sec:hope}, we briefly review some results of \cite{Leichtnam1987}, which proves \cref{conj:main} in the case of the algebraic equation of degree two and which provides some foundational material for an approach of Conjectures~\ref{conj:2} and \ref{conj:m}. In Section~\ref{sec:scope}, we review an interesting family of genuinely nonlinear \acrshort{pde} to which our Conjectures should apply. First, we understand the formation of complex singularities in the \acrshort{ibe} using the \nm{Cauchy}--\nm{Kovalevskaya} theory in Section~\ref{sec:ck}, then, for the sake of completeness, we obtain the same results using completely different methods of contact geometry in Section~\ref{sec:cg}. The reason to tackle the same problem from two completely different angles allows us to explore different ideas underlying the general conjecture formulated below. Some brief reminder on the ring of formal series with holomorphic germs\footnote{The notion of a germ of a mathematical object captures the local properties of that object. The germ is precisely the equivalence class of objects that share the same local property. In order to implement the idea of germs, the space has to be at least topological to give sense to the word `local'. The name `germ' was introduced into Mathematics in the continuation of the \emph{sheaf} metaphor.} coefficients are made in Section~\ref{sec:ring}, and some remarks on holonomic $\D-$modules are presented in Section~\ref{sec:dmod}. The second order nonlinear \acrshort{pde} \eqref{eq:sec} is analyzed in Section~\ref{sec:sec}, and a completely new algorithm to approximate ramified solutions is presented in Section~\ref{sec:it}. The convergence of this algorithm would prove a particular case of the general \cref{conj:main}. We demonstrate numerically in Section~\ref{sec:num} that the proposed algorithm converges in practice, which constitutes one of the main achievements of the present work. Some possible further generalizations are discussed in Section~\ref{sec:gen}. This study's main conclusions and perspectives are outlined in Section~\ref{sec:concl}. Finally, in Appendix~\ref{app:a}, we show how our theory relates to the classical theory of shock waves in an elementary example.


\section{Review of existing results}
\label{sec:hope}

In this Section, we remind the general theory of singular solutions to genuinely nonlinear \acrshort{pde}s, which was initiated by the second author some thirty years ago, and the construction of holomorphic ramified solutions have been considered even earlier, see \eg \cite{Wagschal1974, Hamada1976, Ishii1980}. Certain linear ramified \nm{Cauchy} problems were considered in \cite{Leichtnam1990} and second-order semi-linear problems in \cite{Leichtnam1991}. The linear ramified differential operators were considered in \cite{DAgnolo1991} in the setting of sheaves theory. They demonstrated independently a weaker version of the results presented in \cite{Leichtnam1993}. \nm{Hamada} considered the homogeneous linear case even earlier in \cite{Hamada1969}. The micro-local existence result was established in \cite{Bony1976} for certain classes of linear differential operators. The same homogeneous case was investigated later using the contact transformation approach in \cite{Kashiwara1978}. The linear homogeneous case was also considered using the tools of the micro-local theory of sheaves in \cite{DAgnolo1991a}.

The reminders will allow us to understand better the context and set up of Conjectures~\ref{conj:2}, \ref{conj:m} and \ref{conj:main}. Actually, \cref{thm:2} (to be recalled below) is the particular case of Conjectures~\ref{conj:2} and \ref{conj:main} obtained by restricting oneself to algebraic equations of the second order. This is an additional good reason in favour of \cref{conj:2}.

Let us consider an operator $\P$ of the order $m\,$:
\begin{equation*}
  u\ \mapsto\ \P\,(\,u\,)\ \eqdef\ \sum_{\abs{\sigma}\,=\,m}\,P_{\,\sigma}\,(\,\x,\,u,\,\partial_{\,\x}^{\,\beta}\,u\,)\,\partial_{\,\x}^{\,\sigma}\,u\ +\ \Rr\,(\,\x,\,\partial_{\,\x}^{\,\beta}\,u\,)\,
\end{equation*}
where we use the multi-indices\footnote{If $\beta\ =\ (\,\beta_{\,0},\,\beta_{\,1}\,\ldots,\,\beta_{\,n}\,)\ \in\ \Z^{\,n\,+\,1}$ be a multi-index. By $\abs{\cdot}$ we denote the \emph{height} of the multi-index defined as
\begin{equation*} 
  \abs{\beta}\ \eqdef\ \sum_{j\,=\,0}^{\,n}\,\abs{\beta_{\,j}}\,.
\end{equation*}
The derivation operator with respect to the variable $x_{\,j}\,$, $j\ \in\ (\,n\,+\,1)^{\,\sqsubset}$ is denoted by $\partial_{\,x_{\,j}}\,$. If $\beta\ \in\ \N^{\,n\,+\,1}$ is a multi-index of derivation, then we shall write:
\begin{equation*}
  \partial^{\,\beta}\ \eqdef\ \partial^{\,\beta_{\,0}}_{\,x_{\,0}}\comp\partial^{\,\beta_{\,1}}_{\,x_{\,1}}\comp\ldots\partial^{\,\beta_{\,n}}_{\,x_{\,n}}\,.
\end{equation*}} $\sigma$ and $\beta$ with $\abs{\beta}\ \leq\ m\ -\ 1\,$. The use of multi-index notation is often attributed in the literature to Laurent~\nm{Schwartz}. All functions $P_{\,\sigma}$ and $\Rr$ are supposed to be holomorphic in their arguments $\forall\,\x$ in some vicinity $\U_{\,\vO}$ of $\vO\ \in\ \C^{\,n+1}\,$, $\forall\,\partial_{\,\x}^{\,\beta}\,u$ in some vicinity of $u_{\,0}^{\,\beta}\ \eqdef\ \partial_{\,\x}^{\,\beta}\,u\,(\,\vO\,)\,$, $\beta\ \in\ \N^{\,n+1}\,$. The question of local solutions to equation $\P\,u\ =\ \fO$ can be naturally asked. The solution $u\,(\,\x\,)$ is holomorphic and ramified around a hyper-surface $\S$ defined by equation $s\,(\,\x\,)\ =\ 0$ and passing by $\vO\ \in\ \C^{\,n+1}\,$. This solution has the following analytical form \cite{Leichtnam1987}:
\begin{equation}\label{eq:ser0}
  u\,(\,\x\,)\ =\ a\,(\,\x\,)\ +\ \sum_{k\,=\,0}^{+\infty}\,b_{\,k}\,(\,\x\,)\,s^{\,\gamma_{\,k}}\,(\,\x\,)\,,
\end{equation}
where the sequence $\set*{\gamma_{\,k}}_{\,k\,=\,0}^{\,+\infty}$ takes its values in $\R_{\,>\,0}$ and is strictly increasing and tending to $+\,\infty\,$. The functions $a$ and $b_{\,k}\,$, $\forall\,k\ \in\ \N$ are holomorphic in $\U_{\,\vO}\ \ni\ \vO\,$. The series \eqref{eq:ser0} is convergent on the covering space of $\U_{\,\vO}\setminus\S\,$. The following assumption is also adopted: the hyper-surface $\S$ is simply characteristic for the linearized equation $\P_{\,\mathrm{lin}}$ of the genuinely nonlinear operator $\P$ in every point $a\,(\,\x\,)\,$. In other words, the principal symbol $p_{\,m}$ of $\P$
\begin{equation*}
  p_{\,m}\,(\,\x,\,\xib\,)\ \eqdef\ \sum_{\abs{\sigma}\,=\,m}\,P_{\,\sigma}\,\bigl(\,\x,\,a\,(\,\x\,),\,\partial_{\,\x}^{\,\beta}\,a\,(\,\x\,)\,\bigr)\,\xib^{\,\sigma}\,, \qquad \abs{\beta}\ \leq\ m\ -\ 1\,,
\end{equation*}
verifies two conditions:
\begin{itemize}
  \item $p_{\,m}\,\bigl(\,\x,\,\ud\,s\,(\,\x\,)\,\bigr)\,\bigr\vert_{\,\S}\ \equiv\ 0\,$,
  \item $\partial_{\,\xib}\,p_{\,m}\,\bigl(\,\vO,\,\ud\,s\,(\,\vO\,)\,\bigr)\,\bigr\vert_{\,\S}\ \neq\ 0\,$.
\end{itemize}
The hypothesis of simple characteristics is a holomorphic analogue of the strictly hyperbolic operator relative to $\S$ in the real case.

We also suppose that we know the `initial' data $\partial_{\,t}^{\,\beta}\,a\,\vert_{\,t\,=\,0}$ ($\forall\,\abs{\beta}\ \leq\ m\ -\ 1$), $b_{\,k}\,\vert_{\,t\,=\,0}$ ($\forall\,k\ \geq\ 0$) and $s\,\vert_{\,t\,=\,0}\,$, where $t\ =\ 0$ is the equation a complex smooth hyper-surface $\Tt$ transversal to the vector field $\partial_{\,\xib}\,p_{\,m}\,\bigl(\,\x,\,\ud\,s\,(\,\x\,)\,\bigr)\,$. Below, we consider two separate cases.


\subsection{Weakly singular solutions}
\label{sec:wss}

In this case, we take a holomorphic function $a$ and a hyper-surface $\S$ (given by equation $s\,(\,\x\,)\ =\ 0$) verifying preceding conditions, we assume $\P\,(\,a\,)\ \equiv\ \fO$ and we seek for weakly singular solutions of type \eqref{eq:ser0} with
\begin{equation*}
  \gamma_{\,k}\ \coloneq\ m\ +\ (\,k\ +\ 1\,)\,\mu\,, \qquad \mu\ \in\ ]\,0,\,1\,[
\end{equation*}
so that the constructed solution $u$ admits $m$ bounded derivatives. A change of variables allows us to assume that $\S$ is defined by the equation
\begin{equation*}
  s\,(\,\x\,)\ =\ y_{\,1}\ =\ 0\,,
\end{equation*}
where we decompose $\x\ \equiv\ (\,t,\,\y\,)$ and we introduce the following notation:
\begin{equation*}
  1\ -\ \frac{\y}{R}\ \eqdef\ \prod_{j\,=\,1}^{n}\,\Bigl(\,1\ -\ \frac{y_{\,j}}{R}\,\Bigr)\,, \qquad R\ >\ 0\,.
\end{equation*}
By analogy,
\begin{equation*}
  \frac{1}{1\ -\ \frac{\y}{R}}\ \eqdef\ \prod_{j\,=\,1}^{n}\,\frac{1}{1\ -\ \frac{y_{\,j}}{R}}\,, \qquad R\ >\ 0\,.
\end{equation*}

\begin{definition}
Let $u\ =\ \sum_{j\,=\,1}^{+\,\infty}\,u_{\,j}\,y^{\,j}$ and $v\ =\ \sum_{j\,=\,1}^{+\,\infty}\,v_{\,j}\,y^{\,j}$ be two formal power series with complex coefficients. We shall say that $u\ \preccurlyeq\ v$ if and only if
\begin{equation*}
  \abs{u_{\,j}}\ \leq\ \abs{v_{\,j}}\,, \qquad \forall\,j\ \in\ \N\,.
\end{equation*}
\end{definition}
The last definition can be easily extended to the multi-variable case:
\begin{definition}
Let $p\ \in\ \N^{\,\times}$ and $u\ =\ \sum_{\beta\,\in\,\N^{\,p}}\,u_{\,\beta}\,y^{\,\beta}$ and $v\ =\ \sum_{\beta\,\in\,\N^{\,p}}\,v_{\,\beta}\,y^{\,\beta}$ be two formal power series with complex coefficients. We shall say that $u\ \preccurlyeq\ v$ if and only if
\begin{equation*}
  \abs{u_{\,\beta}}\ \leq\ \abs{v_{\,\beta}}\,, \qquad \forall\,\beta\ \in\ \N^{\,p}\,.
\end{equation*}
\end{definition}

\begin{definition}
Let $\Y$ be a formal variable. Given $R\ >\ 0\,$, $\eps\ >\ 0\,$, $p\ \in\ \N$ and $d\ \in\ \N\,$, we denote by $\Aa^{\,d}\,(\,R,\,\eps,\,p\,)$ the set of formal power series of the type:
\begin{equation*}
  u\ =\ u\,(\,\x,\,\Y\,)\ =\ \sum_{j\,=\,0}^{+\,\infty}\,\sum_{k\,=\,1}^{+\,\infty}\,u_{\,k,\,j}\,(\,t,\,\y\,)\,\Y^{\,k\,\mu\,+\,j\,+\,d}\,, \qquad 0\ <\ \mu\ <\ 1\,,
\end{equation*}
where functions $u_{\,k,\,j}$ are holomorphic in some neighbourhood of $\vO\ \in\ \C^{\,n\,+\,1}$ and verify
\begin{equation}\label{eq:bound}
  u_{\,k,\,j}\ \preccurlyeq\ \sum_{l\,=\,0}^{+\,\infty}\,\frac{t^{\,l}}{\bigl(\,1\ -\ \frac{\y}{R}\,\bigr)^{\,k\,+\,j\,+\,l}}\cdot\frac{(\,k\,+\,j\,+\,l\,)!}{l!\,(\,k\,+\,j\,+\,p\,)!}\cdot C_{\,k\,j\,l}^{\,p}
\end{equation}
such that
\begin{equation}\label{eq:norm}
  \norm{u}_{\,p,\,d}\ \eqdef\ \sum_{k,\,j,\,l\,\in\,\N}\,C_{\,k\,j\,l}^{\,p}\,\eps^{\,k\,+\,2\,j\,+\,l}\ <\ +\,\infty\,,
\end{equation}
where $C_{\,k\,j\,l}^{\,p}\,(\,u\,)$ are the smallest non-negative real constants allowing to bound functions $u_{\,k,\,j}$ in \eqref{eq:bound}. We also introduce the following related definitions:
\begin{equation*}
  \Aa^{\,d}\,(\,R\,)\ \eqdef\ \bigcup_{\eps\,>\,0}\,\Aa^{\,d}\,(\,R,\,\eps,\,p\,)\,, \qquad
  \Aa^{\,d}\ \eqdef\ \bigcup_{R\,>\,0}\,\Aa^{\,d}\,(\,R\,)\,.
\end{equation*}
\end{definition}
To each element $u\ \in\ \Aa^{\,d}\,(\,R,\,\eps,\,p\,)$ we put non-injectively in correspondence a function
\begin{equation*}
  u\,(\,\x\,)\ =\ \sum_{j\,=\,0}^{+\,\infty}\,\sum_{k\,=\,1}^{+\,\infty}\,u_{\,k,\,j}\,y^{\,k\,\mu\,+\,j\,+\,d}\,.
\end{equation*}
It is not difficult to see that if $u_{\,1}\ \in\ \Aa^{\,d_{\,1}}$ and $u_{\,2}\ \in\ \Aa^{\,d_{\,2}}\,$, then $u_{\,1}\cdot u_{\,2}\ \in\ \Aa^{\,d_{\,1}\,+\,d_{\,2}}\,$. If $\beta$ is a multi-index in $\N^{\,n\,+\,1}$ of the height $\leq\ d$ then the operator $\partial^{\,\beta}$ takes the elements from $\Aa^{\,d}$ and sends them to $\Aa^{\,d\,-\,\abs{\beta}}\,$. Since $\S$ is characteristic, then $\P_{\,\mathrm{lin}}$ sends $\Aa^{\,d}$ to $\Aa^{\,d\,-\,m\,+\,1}\,$. The following property specifies the geometric domain of the convergence associated with the norm $\norm{\cdot}_{\,p,\,d}$ defined in \eqref{eq:norm}:
\begin{lemma}[\cite{Leichtnam1987}]
Let $u\ \in\ \Aa^{\,d}\,(\,R,\,\eps,\,p\,)\,$. Then, every $u_{\,k,\,j}$ is a holomorphic function in the domain
\begin{equation*}
  \Omega\ \eqdef\ \Set[\bigg]{(\,t,\,\y\,)\ \in\ \C^{\,n\,+\,1}}{\abs{y_{\,j}}\ <\ R\,, \quad \abs{t}\ <\ \eps\,\prod_{j\,=\,1}^{n}\,\Bigl(\,1\,-\,\frac{\abs{y_{\,j}}}{R}\,\Bigr)}\,,
\end{equation*}
and, for every compact $K\ \subset\ \Omega$ there exists a constant $C_{\,K}\ >\ 0$ such that
\begin{equation*}
  \sup_{K}\,\abs{u_{\,k,\,j}}\ \leq\ C_{K}^{\,k\,+\,j\,+\,1}\,.
\end{equation*}
Conversely, let $\set{u_{\,k,\,j}}$ be a sequence of holomorphic functions defined on an open neighbourhood $\U_{\,\vO}$ of the point $\vO\ \in\ \C^{\,n\,+\,1}$ such that
\begin{equation*}
  \sup_{\U_{\,\vO}}\,\abs{u_{\,k,\,j}}\ \leq\ C^{\,k\,+\,j\,+\,1}\,.
\end{equation*}
Then,
\begin{equation*}
  \sum_{j\,=\,0}^{+\,\infty}\, \sum_{k\,=\,1}^{+\,\infty}\, u_{\,k,\,j}\, \Y^{\,k\,\mu\,+\,j\,+\,d}\ \in\ \Aa^{\,d}\,.
\end{equation*}
\end{lemma}

Then, finally, we may state the first result:
\begin{theorem}[\cite{Leichtnam1987}]\label{thm:2.1}
Take $u_{\,0}\ \in\ \Aa^{\,m}\,(\,\C^{\,n}\,)$ and let $\Tt$ be an analytic hyper-surface passing by $\vO\ \in\ \C^{\,n\,+\,1}$ being defined by the analytical equation $t\ =\ 0\,$. This hyper-surface is transversal to the field $\partial_{\,\xib}\,p_{\,m}\,(\,\x,\,\ud\,y_{\,1}\,)$ (and, thus, to $\S$). Then, there exists $r\ \in\ \Aa^{\,m}\,(\,\C^{\,n\,+\,1}\,)$ such that
\begin{equation*}
  \P\,\bigl(\,a\,(\,\x\,)\ +\ r\,(\,\x\,)\,\bigr)\ \equiv\ \fO \qquad \text{and} \qquad r\,\vert_{\,t\,=\,0}\ =\ u_{\,0}\,.
\end{equation*}
\end{theorem}
The proof consists, roughly speaking, of applying the fixed point iterations.

\begin{remark}
The obtained solution $a\,(\,\x\,)\ +\ r\,(\,\x\,)$ in the previous Theorem is ramified around a smooth locus (independent of the value of $\mu$) because the \nm{Cauchy} datum is not sufficiently singular.
\end{remark}


\subsection{Strongly singular solutions}

In this Section we are looking for solutions of type \eqref{eq:ser0} with
\begin{equation*}
  \gamma_{\,k}\ \eqdef\ m\ -\ \frac{1}{2}\ +\ \frac{k}{2}\,.
\end{equation*}
We stress out that the function $a$ is not necessarily a solution of equation $\P\,(\,a\,)\ =\ \fO$ anymore if $b_{\,0}\ \neq\ \fO\,$. In this Section, $a$ and $\S$ will be considered as problem unknowns on the same footing with $\set*{b_{\,k}}_{\,k\,=\,0}^{\,+\,\infty}\,$.

\begin{definition}
Let $S$ be a formal variable. By $\Bb^{\,m\,-\,1}\,(\,\C^{\,n\,+\,1}\,)$ we designate the algebra of formal series of the form:
\begin{equation*}
  r\,(\,\x,\,S\,)\ =\ \sum_{j\,=\,0}^{+\,\infty}\,b_{\,k}\,(\,\x\,)\,S^{\,\frac{k\,+\,2\,m\,-\,1}{2}}\,.
\end{equation*}
We assume also that all $\set{b_{\,j}}_{\,j\,\in\,\N}$ are holomorphic functions defined on the same neighbourhood $\U_{\,\vO}$ of the point $\vO\ \in\ \C^{\,n\,+\,1}$ and there exists a constant $C\ >\ 0$ such that
\begin{equation*}
  \abs{b_{\,j}\,(\,\x\,)}\ \leq\ C^{\,k\,+\,1}\,, \qquad \forall\,j\ \in\ \N\,, \qquad \forall\,\x\ \in\ \U_{\,\vO}\,.
\end{equation*}
\end{definition}
Below, we shall substitute the formal variable $S$ by the holomorphic function $s$ vanishing at the origin $\vO\ \in\ \C^{\,n\,+\,1}$ and defining the characteristic hyper-surface for the operator $\P_{\,\mathrm{lin}}\,$.

For any $\abs{\beta}\ =\ 1\,$, we define the action of $\partial^{\,\beta}$ on $\Bb^{\,m\,-\,1}$ by conveying that
\begin{equation*}
  \partial^{\,\beta}\,\bigl(\,S^{\,\frac{k\,+\,2\,m\,-\,1}{2}}\,\bigr)\ =\ \frac{k\,+\,2\,m\,-\,1}{2}\;S^{\,\frac{k\,+\,2\,m\,-\,3}{2}}\,\partial^{\,\beta}\,s\,(\,\x\,)\,.
\end{equation*}
The action of $\partial^{\,\beta}$ on $\Bb^{\,m\,-\,1}$ is defined recursively for any $\abs{\beta}\ \leq\ m\,$. To each element $r\,(\,\x,\,S\,)\ \in\ \Bb^{\,m\,-\,1}\,(\,\C^{\,n\,+\,1}\,)$ we associate non-injectively a function $\x\ \mapsto\ r\,\bigl(\,\x,\,s\,(\,\x\,)\,\bigr)$ defined near the origin.

Let $\x\ \equiv\ (\,t,\,\y\,)\ \in\ \C\times\C^{\,n}$ and we assume that the operator $\P$ has the form $\P\ =\ \partial_{\,t}^{\,m}\ +\ \Qq\,(\,t,\,\y,\,\partial^{\,\beta}\,)$ with $\abs{\beta}\ \leq\ m$ and $\beta\ \neq\ (\,m,\,0,\,\ldots,\,0\,)\,$. The principal symbol of the linearized operator $\P_{\,\mathrm{lin}}\,(\,a\,)$ of $\P$ at $a$ has the form $p_{\,m}\,(\,t,\,\y,\,\partial^{\,\beta}\,\tau,\,\xib\,)$ with $\abs{\beta}\ \leq\ m\ -\ 1\,$. Finally, we consider a point $(\,\tau_{\,0},\,\xib^{\,0}\,)\ \in\ (\,\C\times\C^{\,n}\,)\setminus\vO$ together with \nm{Cauchy} data $\set{a_{\,j}\,(\,\y\,)}_{\,j\,=\,0}^{\,m\,-\,1}$ holomorphic in some neighbourhood of $\vO\ \in\ \C^{\,n}\,$. We assume that
\begin{equation*}
  p_{\,m}\,\bigl(\,0,\,\partial^{\,\sigma} a_{\,j}\,(\,\vO\,),\,\tau^{\,0},\,\xib^{\,0}\,\bigr)\ =\ \fO\,, \qquad
  \partial_{\,\tau}\,p_{\,m}\,\bigl(\,0,\,\partial^{\,\sigma} a_{\,j}\,(\,\vO\,),\,\tau^{\,0},\,\xib^{\,0}\,\bigr)\ \neq\ \fO\,,
\end{equation*}
where $j\ +\ \abs{\sigma}\ \leq\ m\ -\ 1\,$. We choose also a root $\tau_{\,1}\,(\,\x,\,u_{\,\beta},\,\xib\,)$ of the equation $\tau\ \mapsto\ p_{\,m}\,(\,\x,\,u,\,\tau,\,\xib\,)\ =\ 0$ holomorphic near the point $(\,0,\,\partial^{\,\sigma} a_{\,j}\,(\,\vO\,),\,\xib^{\,0}\,)\,$. Finally, let $s_{\,1}$ be a holomorphic function near $\vO\ \in\ \C^{\,n}$ such that
\begin{equation*}
  s_{\,1}\,(\,\vO\,)\ =\ 0\,, \qquad \ud\,s_{\,1}\,(\,\vO\,)\ =\ \xib^{\,0}\,.
\end{equation*}
We can state the main result:
\begin{theorem}[\cite{Leichtnam1987}]\label{thm:2}
Let $\sum_{\,k\,=\,0}^{\,+\,\infty}\,b_{\,k}^{\,0}\,(\,\y\,)\,S^{\,\frac{k\,+\,2\,m\,-\,1}{2}}\ \in\ \Bb^{\,m\,-\,1}\,(\,\C^{\,n}\,)\,$. Then, there exists a neighbourhood of the point $\vO\ \in\ \C^{\,n\,+\,1}$ along with holomorphic functions $s\,$, $a\,$, $\set{b_{\,k}}_{\,k\,=\,0}^{\,+\,\infty}$ defined on it such that $s$ is characteristic for $\P_{\,\mathrm{lin}}$ and
\begin{align*} 
  b_{\,k}\,(\,0,\,\y\,)\ &\equiv\ b_{\,k}^{\,0}\,, \quad \forall\,k\ \in\ \N\,, \\
  \partial_{\,t}^{\,j}\,a\,(\,0,\,\y\,)\ &\equiv\ a_{\,j}\,, \quad j\ \in\ (\,m\,+\,1\,)^{\,\sqsubset}\,, \\
  s\,(\,0,\,\y\,)\ &\equiv\ s_{\,1}\,, \\
  \partial_{\,t}\,s\,(\,0\,)\ &\equiv\ \tau^{\,0}\,.
\end{align*}
The series
\begin{equation*}
  r\,(\,\x,\,S\,)\ \eqdef\ \sum_{k\,=\,0}^{\,+\,\infty}\,b_{\,k}\,(\,\x\,)\,S^{\,\frac{k\,+\,2\,m\,-\,1}{2}}\ \in\ \Bb^{\,m\,-\,1}
\end{equation*}
and $\x\ \mapsto\ u\,(\,\x\,)\ \equiv\ a\,(\,\x\,)\ +\ r\,\bigl(\,\x,\,s\,(\,\x\,)\,\bigr)$ is a solution of the equation $\P\,(\,u\,)\ =\ \fO\,$.
\end{theorem}
The space $\Bb^{\,m\,-\,1}\,(\,\C^{\,n\,+\,1}\,)$ is endowed with appropriate norm to define a \nm{Banach} algebra. This Theorem is proved by applying the fixed point \nm{Picard} iteration, which also gives a practical algorithm converging to the required solution \cite{Leichtnam1991}.

\begin{remark}
Since the \nm{Cauchy} datum $u_{\,0}\ \equiv\ u\,(\,0,\,\y\,)$ is sufficiently singular, the hyper-surface $s\ =\ \fO$ around which the solution is ramified depends on $u_{\,0}$ and cannot be prescribed in advance. Moreover, since $u_{\,0}$ is ramified like $s_{\,1}^{\frac{1}{2}}\,$, the algebraic equation $z^{\,2}\ +\ b\,z\ +\ c\ =\ \fO$ will inevitably show up in the construction of the solution. The discriminant of this equation is equal to $\Delta\ \coloneq\ b^{\,2}\ -\ 4\,c\,$. Observe that $\Delta\ =\ \fO$ defines a smooth locus with respect to the variables $(\,b,\,c\,)\,$. That is why \cref{thm:2} is easier to prove than its generalizations to the algebraic equations of order $k\ \geq\ 3\,$, \cf Conjectures \ref{conj:2}, \ref{conj:m} and \ref{conj:main}.
\end{remark}


\section{An interesting family of genuinely nonlinear \acrshort{pde}s}
\label{sec:scope}

We consider a family of \acrshort{pde}s with the unknown function $u\,:\ \R\,\times\,\R_{\,\geq\,0}\ \longrightarrow\ \R\,$:
\begin{equation}\label{eq:fam}
  \L_m\,(\,u\,)\ \eqdef\ \partial_{t}^{\,m}u\ -\ \partial_{x}^{\,m-1}\,u\,\partial_{x}^{\,m}\,u\ =\ \fO\,,
\end{equation}
where $m\ \geq\ 1$ is an integer parameter and $\partial_{(-)}^{\,m}$ denotes the $m$\up{th} order partial derivatives with respect to independent variables $t$ or $x\,$. Let us introduce the auxiliary variables:
\begin{equation*}
  v_{\,1}\,(\,t,\,x\,)\ \eqdef\ \partial_{\,x}^{\,m-1}\,u\,(\,t,\,x\,)\,, \qquad
  v_{\,m}\,(\,t,\,x\,)\ \eqdef\ \partial_{\,t}^{\,m-1}\,u\,(\,t,\,x\,)\,.
\end{equation*}
Then, we have
\begin{align*}
  \partial_{\,t}^{\,m-1}\,v_{\,1}\ -\ \partial_{\,x}^{\,m-1}\,v_{\,m}\ &=\ \fO\,, \\
  \partial_{\,t}\,v_{\,m}\ -\ \partial_{\,x}\,\bigl(\,\half\;v_{\,1}^{\,2}\,\bigr)\ &=\ \fO\,.
\end{align*}
We shall also introduce the functions $\set{v_{\,j}}_{\,j\,=\,2}^{\,m-1}$ satisfying the following system:
\begin{align*}
  \partial_{\,t}\,v_{\,1}\ -\ \partial_{\,x}\,v_{\,2}\ &=\ \fO\,, \\
  \partial_{\,t}\,v_{\,2}\ -\ \partial_{\,x}\,v_{\,3}\ &=\ \fO\,, \\
  \vdots\ & \ \vdots \\
  \partial_{\,t}\,v_{\,m}\ -\ \partial_{\,x}\,\bigl(\,\half\;v_{\,1}^{\,2}\,\bigr)\ &=\ \fO\,.
\end{align*}
For regular solutions (\ie of the class $C^{1}\,$), the last system can be rewritten under the following matrix form:
\begin{equation}\label{eq:sys}
  \partial_{\,t}\,\begin{pmatrix}
    v_{\,1} \\
    v_{\,2} \\
    \vdots \\
    v_{\,m}
  \end{pmatrix}\ -\
  \begin{pmatrix}
    0 & 1 & 0 & \cdots & 0 \\
    0 & 0 & 1 & \cdots & 0 \\
    \vdots & \vdots & \vdots & \ddots & \vdots \\
    v_{\,1} & 0 & 0 & \cdots & 0
  \end{pmatrix}\cdot\partial_{\,x}\,\begin{pmatrix}
    v_{\,1} \\
    v_{\,2} \\
    \vdots \\
    v_{\,m}
  \end{pmatrix}\ =\ \vO\,,
\end{equation}
where $\vO$ is a function taking a constant (zero) vector value, \ie
\begin{align*}
  \vO\,:\ \R_{\,\geq\,0}\times\R\ &\longrightarrow\ \R^{\,m}\,, \\
  (\,t,\,x\,)\ &\mapsto\ (\,0,\,0,\,\ldots,\,0\,)^{\,\top}\,.
\end{align*}
Formulations \eqref{eq:fam} and \eqref{eq:sys} are equivalent to each other on the set of sufficiently regular solutions.

It is not difficult to see that this system matrix is nothing else but the transposed companion matrix of the polynomial
\begin{equation*}
  \lambda^{\,m}\ -\ v_{\,1}\ \in\ \R\,[\,\lambda\,]\,.
\end{equation*}
If we assume that $v_{\,1}\ >\ 0$ (just for definiteness), the eigenvalues of this matrix are all simple and given by:
\begin{equation*}
  \lambda_{\,k}\ \eqdef\ v_{1}^{\,\frac{1}{m}}\,\ue^{\,\ui\,\frac{2\,\pi\,k}{m}}\,, \qquad k\ \in\ m^{\,\sqsupset}\,.
\end{equation*}
The eigenvector associated with the eigenvalue $\lambda_{\,k}$ is a direction in the kernel of the following linear operator:
\begin{equation*}
  \begin{pmatrix}
    -\lambda_{\,k} & 1 & 0 & \cdots & 0 \\
    0 & -\lambda_{\,k} & 1 & \cdots & 0 \\
    \vdots & \vdots & \vdots & \ddots & \vdots \\
    v_{\,1} & 0 & 0 & \cdots & -\lambda_{\,k}
  \end{pmatrix}
\end{equation*}
or simply the following vector:
\begin{equation*}
  \omega_{\,k}\ \eqdef\ (\,1,\,\lambda_{\,k},\,\ldots,\,\lambda_{\,k}^{\,m-1}\,)^{\,\top}
\end{equation*}
We can easily verify that
\begin{equation*}
  \Inner{\ud\,\lambda_{\,k}}{\omega_{\,k}}\ =\ \frac{\lambda_{\,k}}{m\,v_{\,1}}\ \neq\ \fO\,,
\end{equation*}
which means that the $k$\up{th} characteristic field is \emph{genuinely nonlinear} according to the definition from \cite[Section~\textsection5]{Lax1973}. It is also obvious that System~\eqref{eq:sys} is not hyperbolic since it admits complex characteristic speeds if $m\ >\ 2\,$. Below, we shall consider in some detail the two hyperbolic cases $m\ \in\ 1$ and $m\ =\ 2\,$. The former corresponds to the familiar \acrshort{ibe}. The latter deserves more attention. It can be written as the following system of two equations:
\begin{align}\label{eq:psys}
  \partial_{\,t}\,v_{\,1}\ -\ \partial_{\,x}\,v_{\,2}\ &=\ \fO\,, \\
  \partial_{\,t}\,v_{\,2}\ -\ \partial_{\,x}\,\bigl(\,\sigma\,(\,v_{\,1}\,)\,\bigr)\ &=\ \fO\,.
\end{align}
In the last system (also known as the $p-$system), we can recognize a particular case of nonlinear elasticity written in \nm{Lagrangian} variables, where $v_{\,1}$ is the deformations field and $v_{\,2}$ corresponds to the velocity field. Finally, the mechanical constraint $\sigma\,(w)\ \eqdef\ \half\;w^{\,2}$ is induced by the deformation $v_{\,1}\,$. See \cite[Chapter~7]{Dafermos2010} and, in particular, Equation~(7.1.11) therein for more details. It is known that System~\eqref{eq:psys} is hyperbolic when $\od{\sigma}{w}\ \equiv\ w\ >\ 0\,$, which becomes the condition $v_{\,1}\ >\ 0$ in our case. The two characteristic fields corresponding to speeds
\begin{equation*}
  \lambda_{\,\pm}\ =\ \pm\sqrt{v_{\,1}}
\end{equation*}
and to eigenvectors
\begin{equation*}
  \omega_{\,\pm}\ =\ (\,1,\,\pm\sqrt{v_{\,1}}\,)^{\,\top}
\end{equation*}
are genuinely nonlinear as in the general case (\ie $m\ >\ 2$).


\section{Inviscid Burgers equation}
\label{sec:burg}

In this Section we consider a particular \acrshort{ivp} for the \acrshort{ibe} \eqref{eq:burg}:
\begin{subequations}\label{eq:cauchy}
\begin{align}\label{eq:cauchy1}
  u_{\,t}\ -\ u\,u_{\,x}\ &=\ \fO\,, \\
  u\,(\,0,\,x\,)\ &=\ x^{\,\third}\,. \label{eq:cauchy02}
\end{align}
\end{subequations}
The main particularity of our study is that we consider the space and time variables to be complex, including the unknown function $u\,$:
\begin{equation*}
  u\,:\ \C\,\times\,\C\ \longrightarrow\ \C\,.
\end{equation*}
The derivatives with respect to $x$ and $t$ are understood from now on in the sense of complex analysis. The initial condition containing an algebraic branching point in the complex domain was taken on purpose. Consequently, the \acrshort{ivp} \eqref{eq:cauchy} is called the \emph{ramified} \nm{Cauchy} \emph{problem}. Solving the ramified \nm{Cauchy} problem \eqref{eq:cauchy} means by definition to find an open neighbourhood $\U_{\,\vO}$ of $\vO\ \in\ \C^{\,2}$ along with a germ of the holomorphic function $u\,(\,t,\,x\,)$ at a point of $\U_{\,\vO}\,\cap\,\set{t\,=\,0}$ (different from the origin) satisfying \eqref{eq:cauchy1} and \eqref{eq:cauchy02} such that $u\,(\,t,\,x\,)$ can be continued holomorphically along any path in $\U_{\,\vO}$ which does not meet a certain characteristic hyper-surface which restricts to $x\ =\ 0$ when $t\ =\ 0\,$. Below, we describe two different approaches to study the \acrshort{ivp} \eqref{eq:cauchy} because both of them contain the foundational ideas to approach the \cref{conj:2}.


\subsection{A linear digression}
\label{sec:lin}

It is often helpful to consider first the linear problem before tackling the fully nonlinear formulation. It often brings some useful insight. The linear counterpart to \acrshort{ivp} \eqref{eq:cauchy} reads
\begin{align*}
  u_{\,t}\ -\ u_{\,x}\ &=\ \fO\,, \\
  u\,(\,0,\,x\,)\ &=\ x^{\,\third}\,.
\end{align*}
It is clear to see that the last \acrshort{ivp} has the following \emph{exact} solution:
\begin{equation*}
  u\,(\,t,\,x\,)\ =\ (\,t\ +\ x\,)^{\,\third}\,.
\end{equation*}
Moreover, the singular locus of this solution is the following affine algebraic variety:
\begin{equation*}
  \S_{\,\mathrm{s}}\ \eqdef\ \Set*{(\,t,\,x\,)\ \in\ \C^{\,2}}{t\ +\ x\ =\ 0}\,.
\end{equation*}
It is clear to see that this locus is completely regular as a geometric object. This is a characteristic property of linear equations. Below in \cref{sec:ck}, we shall demonstrate that in the genuinely \textbf{nonlinear} case \eqref{eq:cauchy}, the corresponding singular locus is \textbf{singular} itself, and we shall describe its algebraic singularity explicitly. We emphasize that \cref{eq:toy} was (simply) nonlinear and also had a completely regular locus. Also, the \nm{Galois} group of the field extension of the \nm{Cauchy} data and of the solution is the same \nm{Abelian} group:
\begin{equation*}
  \Gal\,\bigl(\,\C\,(\,x^{\,\frac{1}{3}}\,)\,:\,\C\,(\,x\,)\,\bigr)\ =\ U_{\,3}\,.
\end{equation*}
This observation is to be compared with the genuinely nonlinear case described below. Moreover, the solution $u\,(\,t,\,x\,)$ and the \nm{Cauchy} datum $u\,(\,0,\,x\,)$ share the same monodromy group.

Let us underline one more aspect in which linear and genuinely nonlinear problems differ. If we remove from $\C^{\,2}$ the singular locus and we compute the fundamental group of the obtained domain, we shall obtain\footnote{Indeed, this result is not difficult to obtain:
\begin{equation*}
  \pi_{\,1}\,(\,\C^{\,2}\,\setminus\,\S_{\,\mathrm{s}}\,)\ =\ \pi_{\,1}\,(\,\C^{\,2}\,\setminus\,\Set*{x\ =\ 0}{(\,t,\,x\,)\ \in\ \C^{\,2}}\,)\ =\ \pi_{\,1}\,(\,\C\,\times\,\C\,\setminus\,\set*{0}\,)\ =\ \pi_{\,1}\,(\,\C\,\setminus\,\set*{0}\,)\ =\ \Z\,.
\end{equation*}}:
\begin{equation*}
  \pi_{\,1}\,(\,\C^{\,2}\,\setminus\,\S_{\,\mathrm{s}}\,)\ =\ \Z\,,
\end{equation*}
which is an \nm{Abelian} group. In the genuinely nonlinear case, the fundamental group will be non-commutative.


\subsection{The Cauchy--Kovalevskaya approach}
\label{sec:ck}

In the present Section, we employ the \nm{Cauchy}--\nm{Kovalevskaya} theory in order to study the \acrfull{ivp} \eqref{eq:cauchy}. We shall even obtain a slightly more general result by considering the following \acrshort{ivp}:
\begin{subequations}\label{eq:cauchy2}
\begin{align}
  u_{\,t}\ -\ u\,u_{\,x}\ &=\ \fO\,, \label{eq:cauchy12} \\
  u\,(\,0,\,x\,)\ &=\ a_{\,0}\,(\,x\,)\ +\ x^{\,\third}\,. \label{eq:cauchy22}
\end{align}
\end{subequations}

\begin{theorem}\label{thm:ck}
Consider the \acrshort{ivp} \eqref{eq:cauchy} for the \acrshort{ibe}. Then, there exist three holomorphic functions $(\,t,\,x\,)\ \mapsto\ a\,(\,t,\,x\,)\,$, $(\,t,\,x\,)\ \mapsto\ p\,(\,t,\,x\,)\,$ and $(\,t,\,x\,)\ \mapsto\ q\,(\,t,\,x\,)$ defined in some neighbourhood of the point $(\,0,\,0\,)\ \in\ \C^{\,2}$ such that
\begin{equation}\label{eq:conds}
  a\,(\,0,\,x\,)\ =\ a_{\,0}\,(\,x\,)\,, \qquad
  p\,(\,0,\,x\,)\ =\ 0\,, \qquad
  q\,(\,0,\,x\,)\ =\ x
\end{equation}
and
\begin{equation}\label{eq:ans}
  u\,(\,t,\,x\,)\ \eqdef\ a\,(\,t,\,x\,)\ +\ z\,(\,t,\,x\,)
\end{equation}
being a unique solution to \cref{eq:cauchy} with the complex auxiliary function $(\,t,\,x\,)\ \mapsto\ z\,(\,t,\,x\,)$ verifying the following algebraic relation\footnote{This cubic relation can be beneficially seen as an algebraic equation $f\,(\,z\,)\ \coloneq\ z^{\,3}\ -\ p\,z\ -\ q\ =\ \fO$ over some field $\kk$ of characteristic different from $2$ and $3\,$. We would like to mention here that the \nm{Galois} group associated with this equation depends on the fact whether its discriminant $4\,p^{\,3}\ -\ 27\,q^{\,2}$ is a square or not in $\kk\,$. In the former case, $\Gal_{\,\kk}\,(\,f\,)\ =\ \Z\,/\,3\,\Z$ and in the latter case $\Gal_{\,\kk}\,(\,f\,)\ =\ S_{\,3}\,$, \ie the symmetric group over a finite set of $3$ symbols. Since we are working over the field of complex numbers $\kk\ \coloneq\ \C\,(\,p,\,q\,)\,$, the \nm{Galois} group is the latter.}
\begin{equation}\label{eq:cube}
  z^{\,3}\ =\ p\,z\ +\ q\,.
\end{equation}
\end{theorem}

\begin{proof}
First of all, we check that the ansatz \eqref{eq:ans} verifies the initial condition \eqref{eq:cauchy2}. Indeed, by taking the limit $t\ \to\ 0$ in \cref{eq:cube}, we obtain
\begin{equation*}
  z^{\,3}\,(\,0,\,x\,)\ =\ p\,(\,0,\,x\,)\cdot z\,(\,0,\,x\,)\ +\ q\,(\,0,\,x\,)\ \stackrel{\text{\eqref{eq:conds}}}{\equiv}\ x\,.
\end{equation*}
Thus, $z\,(\,0,\,x\,)\ =\ x^{\,\third}\,$. Substituting this result into the solution ansatz \eqref{eq:ans}, we obtain:
\begin{equation*}
  u\,(\,0,\,x\,)\ \bydef\ a\,(\,0,\,x\,)\ +\ z\,(\,0,\,x\,)\ \stackrel{\text{\eqref{eq:conds}}}{\equiv}\ a_{\,0}\,(\,x\,)\ +\ x^{\,\third}\,.
\end{equation*}
By differentiating the algebraic relation \eqref{eq:cube} with respect to $t$ and $x\,$, we may compute the partial derivatives of the function $z\,$:
\begin{equation}\label{eq:pd}
  z_{\,t}\ =\ \frac{p_{\,t}\,z\ +\ q_{\,t}}{3\,z^{\,2}\ -\ p}\,, \qquad
  z_{\,x}\ =\ \frac{p_{\,x}\,z\ +\ q_{\,x}}{3\,z^{\,2}\ -\ p}\,.
\end{equation}
Let us substitute the solution ansatz \eqref{eq:ans} into the \acrshort{ibe} \eqref{eq:cauchy1}:
\begin{equation*}
  a_{\,t}\ +\ z_{\,t}\ -\ (\,a\ +\ z\,)\,(\,a_{\,x}\ +\ z_{\,x}\,)\ =\ \fO\,.
\end{equation*}
By substituting partial derivatives \eqref{eq:pd} and multiplying both sides by the denominator $3\,z^{\,2}\ -\ p\,$, we obtain:
\begin{equation*}
  (\,3\,z^{\,2}\ -\ p\,)\,a_{\,t}\ +\ z\,p_{\,t}\ +\ q_{\,t}\ -\ (\,a\ +\ z\,)\,\bigl(\,a_{\,x}\,(\,3\,z^{\,2}\ -\ p\,)\ +\ z\,p_{\,x}\ +\ q_{\,x}\,\bigr)\ =\ \fO\,.
\end{equation*}
Now, by equating the coefficients in front of equal powers of $z\,$, we have the following sequence of equalities:
\begin{align*}
  z^{\,0}\,:&\qquad -\,p\,a_{\,t}\ +\ q_{\,t}\ +\ \cancel{a\,p\,a_{\,x}}\ -\ a\,q_{\,x}\ =\ \fO\,, \\
  z^{\,1}\,:&\qquad p_{\,t}\ +\ \cancel{p\,a_{\,x}}\ -\ a\,p_{\,x}\ -\ q_{\,x}\ =\ \fO\,, \\
  z^{\,2}\,:&\qquad a_{\,t}\ -\ \cancel{a\,a_{\,x}}\ -\ \third\;p_{\,x}\ =\ \fO\,, \\
  z^{\,3}\,:&\qquad 3\,a_{\,x}\ =\ \fO\,.
\end{align*}
The last equation can be substituted into three previous ones to obtain:
\begin{align}\label{eq:ck}
 \begin{split}
  q_{\,t}\ &=\ a\,q_{\,x}\ +\ p\,a_{\,t}\ \equiv\ a\,q_{\,x}\ +\ \third\,p\,p_{\,x}\,, \\
  p_{\,t}\ &=\ a\,p_{\,x}\ +\ q_{\,x}\,, \\
  a_{\,t}\ &=\ \third\;p_{\,x}\,.
 \end{split}
\end{align}
The first two equations are not incompatible because co-tangent vectors $\ud\,p$ and $\ud\,q$ are linearly independent by continuity in the vicinity of the point $(\,0,\,0\,)\ \in\ \C^{\,2}\,$. Indeed,
\begin{equation*}
  \bigl(\,p_{\,t},\,p_{\,x}\,\bigr)\,(\,0,\,0\,)\ \equiv\ (\,1,\,0\,)\,, \qquad
  \bigl(\,q_{\,t},\,q_{\,x}\,\bigr)\,(\,0,\,0\,)\ \equiv\ \bigl(\,a_{\,0}\,(\,0\,),\,1\,\bigr)\,.
\end{equation*}
It is not difficult to see that the last two co-vectors are linearly independent\footnote{Indeed, \begin{equation*}\begin{vmatrix} \ 1 & 0\;\;\; \\ \ a_{\,0}\,(\,0\,) & 1\;\;\; \end{vmatrix}\ \equiv\ 1\ \neq\ 0\,, \qquad \forall\,a_{\,0}\,(\,0\,)\ \in\ \C\,.\end{equation*}} for any value of $a_{\,0}\,(\,0\,)\ \in\ \C\,$. Thus, after applying the classical \nm{Cauchy}--\nm{Kovalevskaya} theorem to System \eqref{eq:ck}, we obtain the required local existence and uniqueness result.
\end{proof}

\begin{corollary}
The local existence and uniqueness for the \acrshort{ivp} \eqref{eq:cauchy} is obtained by simply choosing $a_{\,0}\ \coloneq\ \fO\,$.
\end{corollary}

\begin{remark}\label{rem:4.1}
It is instructive to have a look at the obtained solution from the point of view of \nm{Galois} theory. We already established above that the \nm{Galois} group associated with the initial data is the \nm{Abelian} group $U_{\,3}\,$. The dynamics of the \acrshort{ibe} transforms this initial data into the local analytic solution \eqref{eq:ans}. The \nm{Galois} group associated to the solution $u$ can be easily computed:
\begin{equation*}
  \Gal\,\bigl(\,\C\,(\,p,\,q\,)\,[\,z\,]\,:\,\C\,(\,p,\,q\,)\,\bigr)\ =\ S_{\,3}\,.
\end{equation*}
We stress that the last group is obviously different from $U_{\,3}$ and, additionally, it is non-commutative.
\end{remark}

Investigating the situations where the solution established in the last \cref{thm:ck} becomes singular is also very interesting. Let us fix some values of $(\,p_{\,0},\,q_{\,0}\,)$ of $(\,p,\,q\,)\ \in\ \C^{\,2}$ for which $z_{\,0}$ is a simple root of \cref{eq:cube}. Then, by the classical implicit function theorem, there exists a holomorphic function $(\,p,\,q\,)\ \mapsto\ z\,(\,p,\,q\,)$ defined in the vicinity of the point $(\,p_{\,0},\,q_{\,0}\,)$ of the complex plane. Henceforth, the singular locus for \cref{eq:cube} will consist of the following points:
\begin{equation*}
  \Delta\ \eqdef\ \Set*{(\,p,\,q\,)\ \in\ \C^{\,2}}{\text{\cref{eq:cube} has a double root}}\,.
\end{equation*}
It is not difficult to describe this set analytically. For this, we eliminate $z$ from the following system of equations:
\begin{equation*}
  \begin{dcases*}
    \ z^{\,3}\ -\ p\,z\ -\ q\ =\ 0\,, \\
    \ 3\,z^{\,2}\ -\ p\ =\ 0\,.
  \end{dcases*}
\end{equation*}
The last two equations can be seen as a system of linear equations with respect to parameters $(\,p,\,q\,)\,$. It admits the following parametric solution:
\begin{equation*}
  \begin{dcases*}
    \ q\ =\ -2\,z^{\,3}\,, \\
    \ p\ =\ 3\,z^{\,2}\,.
  \end{dcases*}
\end{equation*}
Henceforth, the singular locus is
\begin{equation*}
  \Delta\ \bydef\ \Set*{(\,p,\,q\,)\ \in\ \C^{\,2}}{4\,p^{\,3}\ -\ 27\,q^{\,2}\ =\ 0}\,.
\end{equation*}
It is not difficult to recognize here the swallow tail (cusp) singularity \cite{Woodcock1974}, which coincides with the cusp in this particular case ($n\ =\ 2\,$, see below).

\begin{remark}\label{rem:4.2}
The \cref{thm:ck} establishes the local existence of the solution in some polydisc $D$ around the origin $(\,0,\,0\,)\ \in\ \C^{\,2}\,$. Let us take the solution germ $u\,(\,t,\,x\,)$ at $(\,t_{\,0},\,x_{\,0}\,)\ \in\ D\,\setminus\,\Delta^{\,<}\,(\,\set{\vO}\,)\,$. If $\gamma$ denotes a closed loop around the point $(\,t_{\,0},\,x_{\,0}\,)$ which does not intersect the singular locus $\Delta^{\,<}\,(\,\set{\vO}\,)\,$, then one can obtain the analytical continuation of the germ $u\,(\,t,\,x\,)$ along the loop $\gamma$ \cite{Forster1981}. In other terms, this defines the action of the first homotopy (or fundamental) group $\pi_{\,1}\,\bigl(\,\U_{\,\vO}\,\setminus\,\Delta^{\,<}\,(\,\set{\vO}\,)\,\bigr)$ on the solution \eqref{eq:ans}:
\begin{equation*}
  \gamma\cdot u\ \coloneq\ a\ +\ \gamma\cdot z\,,
\end{equation*}
since the regular part $a$ is not ramified. We note that this group action turns out to be \emph{transitive} (see below \cref{thm:trans}). We say that $\gamma\cdot z$ is the monodromy action of $\gamma$ on $z\,$. Depending on the chosen loop $\gamma\ \in\ D\,\setminus\,\Delta^{\,<}\,(\,\set{\vO}\,)\,$, this action may bring us to a different root of the algebraic \cref{eq:cube}. This situation is schematically depicted in Figure~\ref{fig:cont}. We also mention that this action can be lifted to the universal covering.

\begin{figure}
  \centering
  \includegraphics[width=0.59\textwidth]{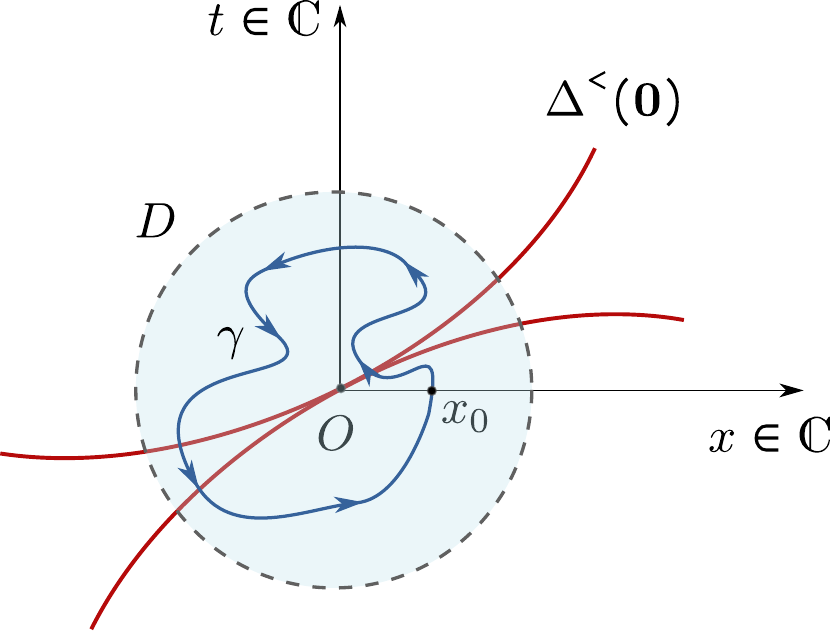}
  \caption{\small\em A schematic illustration for the analytical continuation of a solution germ defined at $x_{\,0}\ \in\ D$ along a closed path $\gamma\ \subseteq\ D\,$.}
  \label{fig:cont}
\end{figure}

It is well-known (see \eg \cite{Arnold1990}) that the fundamental group\footnote{Both Authors of this manuscript are deeply convinced that the mathematical reality is connected. We have some reasons to think that there are some deep hidden connections between the problem considered in our study and the number theory. Indeed, it is known that the braid group $\Br_{3}$ is the universal central extension of the modular group $\PSL_{2}\,(\,\Z\,)\ \eqdef\ \SL_{2}\,(\,\Z\,)\,/\,Z\,\bigl(\,\SL_{2}\,(\,\Z\,)\,\bigr)\ \equiv\ \SL_{2}\,(\,\Z\,)\,/\,\set*{\pm\Id}$ (here by $\Id$ we understand the unit element of $\SL_{2}\,(\,\Z\,)$). The last modular group is deeply rooted in the theory of modular forms and newforms in the number theory. However, at the current stage, we are unable to make this statement more precise.} of the complex plain without the cusp is
\begin{equation*}
  \pi_{\,1}\,\bigl(\,\C^{\,2}\,\setminus\,\Delta^{\,<}\,(\,\set{\vO}\,)\,\bigr)\ =\ \Br_{3}\,,
\end{equation*}
the braid group with three strands is an infinite non-commutative group (in contrast to the linear case discussed above).
\end{remark}

Let us clarify our terminology. Consider the following polynomial equation\footnote{We would like to mention that the \nm{Galois} group of this polynomial equation is the symmetric group $S_{\,n\,+\,1}$ over a finite set of $n\ +\ 1$ symbols, which is the most non-commutative finite group, in the sense that the centre of $Z\,(\,S_{\,n}\,)$ is trivial $\forall\,n\,\geq\,3\,$. The point is that if the underlying \nm{Galois} group is non-commutative, then the analytical prolongation of the solution will be non-commutative as well, \ie the solution value will depend on the order of loops starting and returning to a given point.}:
\begin{equation}\label{eq:pol}
  F\,(\,\x,\,z\,)\ \coloneq\ z^{\,n\,+\,1}\ -\ x_{\,n}\,z^{\,n\,-\,1}\ -\ x_{\,k\,-\,1}\,z^{\,n\,-\,2}\ -\ \ldots\ -\ x_{\,2}\,z\ -\ x_{\,1}\ =\ \fO\,,
\end{equation}
where $\x\ \eqdef\ (\,x_{\,1},\,x_{\,2},\,\ldots,\,x_{\,n}\,)\ \in\ \C^{\,n}\,$. By $\Delta\,(\,\x\,)$ we denote the discriminant of the polynomial \cref{eq:pol} \cite{Gelfand1994}. By swallow tail we shall designate the analytical hyper-surface $\Delta^{\,<}\,(\,\set{0}\,)\,$. We provide below the most general definition of the swallow tail which suits our purposes:

\begin{definition}[\cite{Leichtnam1993}]
An analytical hyper-surface $\S$ in $\C^{\,n\,+\,1}$ defined in some open neighbourhood $\U_{\,\vO}$ is called the swallow tail with the edge $\vO\ \in\ \C^{\,n\,+\,1}$ if there exists $k\ \in\ \N_{\,\geq\,1}$ and $k$ holomorphic functions
\begin{equation*}
  g_{\,j}\,:\ \U_{\,\vO}\ \longrightarrow\ \C\,, \qquad g_{\,j}\,(\,\vO\,)\ =\ 0\,, \qquad j\ \in\ k^{\,\sqsupset}\,,
\end{equation*}
such that the differentials $\set*{\ud\,g_{\,j}\,(\,\vO\,)}_{\,j\,=\,1}^{\,k}$ are linearly independent and $\S$ is defined as the locus of the equation
\begin{equation*}
  \Delta\,\bigl(\,g_{\,1}\,(\,\x\,),\,g_{\,2}\,(\,\x\,),\,\ldots,\,g_{\,k}\,(\,\x\,)\,\bigr)\ =\ \fO\,,
\end{equation*}
where, as above, $\Delta\,(\,g_{\,1},\,g_{\,2},\,\ldots,\,g_{\,k}\,)$ denotes the discriminant of the following polynomial equation in $z$ \cite{Gelfand1994}:
\begin{equation*}
  z^{\,k\,+\,1}\ -\ g_{\,k}\,z^{\,k\,-\,1}\ -\ \ldots\ -\ g_{\,2}\,z\ -\ g_{\,1}\ =\ \fO\,.
\end{equation*}
\end{definition}
We define the co-normal $\Nn\,(\,\S\,)$ of $\S$ to be the closure in $\T^{\,\ast}\,\U_{\,\vO}\ \setminus\ \vO$ of the co-normal $\Nn\,(\,\S_{\,\mathrm{reg}}\,)\,$, where $\S_{\,\mathrm{reg}}$ is the smooth (regular) part of $\S\,$, \ie
\begin{equation*}
  \S_{\,\mathrm{reg}}\ \eqdef\ \Set*{\y\ \in\ \C^{\,k}}{\Delta\,(\,\y\,)\ =\ 0, \quad \ud\,\Delta\,(\,\y\,)\ \neq\ 0}\,.
\end{equation*}
It can be also shown that for any polydisc $D$ centered around $\vO\ \in\ \C^{\,k}\,$, $D\,\bigcap\,\S_{\,\mathrm{reg}}$ is a connected domain. This property is important from the topological point of view.

\begin{remark}
Let $\x_{\,0}\ \in\ \U_{\,\vO}\ \setminus\ \S\,$. Consider a holomorphic germ\footnote{Consider a point $x$ of a complex topological manifold $\V\,$. Consider also two maps $f,\,g\,:\ \U_{\,x}^{\,(f),\,(g)}\ \longrightarrow\ \C$ defined on some neighbourhoods (not necessarily the same) of the point $x\ \in\ \V\,$. These functions define the same germ (or they belong to the same equivalence class) if there is an open neighbourhood $\U\ \subseteq\ \U_{\,x}^{\,(f)}\ \bigcap\ \U_{\,x}^{\,(g)}$ such that
\begin{equation*}
  f\,(\,u\,)\ =\ g\,(\,u\,)\,, \qquad \forall\,u\ \in\ \U\,.
\end{equation*}
This fact can be denoted as $f\ \sim_{\,x}\ g$ and the germ generated by the function $f$ is defined as
\begin{equation*}
  [\,f\,]_{\,x}\ \eqdef\ \Set*{g\,:\ \U_{\,x}\ \longrightarrow\ \C}{f\ \sim_{\,x}\ g}\,.
\end{equation*}
The set of germs at a point $x$ can be obviously endowed with the ring structure over the field $\C\,$. Here, we discussed the scalar case. The generalization to the vectorial case is straightforward. We would like to mention that the set of germs does not possess a non-trivial topology. Henceforth, it makes little sense to speak of the convergence of sequences of germs.} at the point $\bigl(\,g_{\,1}^{\,(\,0\,)},\,g_{\,2}^{\,(\,0\,)},\ldots,\,g_{\,n}^{\,(\,0\,)}\,\bigr)\ \coloneq\ \bigl(\,g_{\,1}\,(\,\x_{\,0}\,),\,g_{\,2}\,(\,\x_{\,0}\,),\ldots,\,g_{\,n}\,(\,\x_{\,0}\,)\,\bigr)\,$:
\begin{equation*}
  \bigl(\,g_{\,1},\,g_{\,2},\,\ldots,\,g_{\,n}\,\bigr)\ \mapsto\ z\,(\,g_{\,1},\,g_{\,2},\,\ldots,\,g_{\,n}\,)\,,
\end{equation*}
where $z$ is the solution to the polynomial \cref{eq:pol}. The classical implicit function theorem ensures the existence of this germ. The germ $z\,(\,g_{\,1},\,g_{\,2},\,\ldots,\,g_{\,n}\,)$ is a holomorphic function ramified around $\Delta^{\,<}\,(\,\set{0}\,)$ and it can be prolonged holomorphically along any path starting from the point $\bigl(\,g_{\,1}\,(\,\x_{\,0}\,),\,g_{\,2}\,(\,\x_{\,0}\,),\ldots,\,g_{\,n}\,(\,\x_{\,0}\,)\,\bigr)$ and avoiding the swallow tail $\Delta^{\,<}\,(\,\set{0}\,)\,$. Then, the germ in $\x_{\,0}\,$:
\begin{equation*}
  \x\ \mapsto\ z\,\bigl(\,g_{\,1}\,(\,\x\,),\,g_{\,2}\,(\,\x\,),\ldots,\,g_{\,n}\,(\,\x\,)\,\bigr)
\end{equation*}
can be prolongated holomorphically along any path starting from $\x_{\,0}$ and belonging to $\U_{\,\vO}\ \setminus\ \S\,$. These facts are rigorously established in \cite[Theorem~3.2]{Leichtnam1993}:
\begin{theorem}\label{thm:4.2}
Let $\set*{z_{\,j}}_{\,j\,=\,1}^{\,n\,+\,1}$ denote the set of solutions to polynomial \cref{eq:pol}. Then, for $\forall\,j\ \in\ (\,n\ +\ 1\,)^{\,\sqsupset}\,$, the germ $z_{\,j}$ defined at $\x_{\,0}$ can be holomorphically prolongated along any path starting from $\x_{\,0}$ and belonging to $\C^{\,n}\,\setminus\,\Delta^{\,<}\,(\,\set{0}\,)\,$.
\end{theorem}
\begin{proof}
See \cite[p.~40]{Leichtnam1993}.
\end{proof}

Basically, we came to the point where we may define a group epimorphism from the fundamental group on the permutation group of simple roots to the algebraic \cref{eq:pol}:
\begin{theorem}[\cite{Leichtnam1993}]\label{thm:trans}
Let $D$ be a polydisc of some positive radius centered around $\vO\ \in\ \C^{\,n}$ and $\gamma$ is a closed path drawn in $D\,\setminus\,\Delta^{\,<}\,(\,\set{\vO}\,)$ starting and ending at $\x_{\,0}\,$: 
\begin{align*}
  [\,0,\,1\,]\ &\longrightarrow\ D\,\setminus\,\Delta^{\,<}\,(\,\set{\vO}\,)\,, \\
  t\ &\mapsto\ \gamma\,(\,t\,)\,
\end{align*}
such that $\gamma\,(\,0\,)\ =\ \gamma\,(\,1\,)\ \equiv\ \x_{\,0}\,$. By $\gamma\cdot z\,(\,\x_{\,0}\,)$ we denote the holomorphic germ obtained from $z\,(\,\x_{\,0}\,)$ after the analytical continuation along the loop $\gamma\,$. Thanks to the implicit function theorem and the fact that the roots of \eqref{eq:pol} are simple in $D\,\setminus\,\Delta^{\,<}\,(\,\set{\vO}\,)$ allows us to define the mapping:
\begin{align*}
  S_{\,n\,+\,1}\ &\longrightarrow\ S_{\,n\,+\,1}\, \\
  z\,(\,\x_{\,0}\,)\ &\mapsto\ \gamma\cdot z\,(\,\x_{\,0}\,)\,,
\end{align*}
where $S_{\,n\,+\,1}$ is the permutation group on $n\ +\ 1$ symbols $\set*{z_{\,1},\,z_{\,2},\,\ldots,\,z_{\,n\,+\,1}}\,$. Then, the following mapping
\begin{align*}
  \pi_{\,1}\,\bigl(\,D\,\setminus\,\Delta^{\,<}\,(\,\set{\vO}\,)\,\bigr)\ &\longrightarrow\ S_{\,n\,+\,1}\,, \\
  \gamma\ &\mapsto\ \bigl(\,z_{\,j}\,(\,\x_{\,0}\,)\ \mapsto\ \gamma\,z_{\,j}\,(\,\x_{\,0}\,)\,\bigr)
\end{align*}
defines the desired group epimorphism.
\end{theorem}
\end{remark}


\subsubsection{Uniformization of the local solution}

It is known that the \nm{Cauchy}--\nm{Kovalevskaya} theory provides only the \emph{local} analytic solutions. However, we work in a complex analytic setting where a germ can be continued along any path, avoiding singularities. Henceforth, the local solution obtained in the previous Section can be continued analytically along any path starting in some neighbourhood in origin and avoiding the swallow tail singularity $\Delta^{\,<}\,(\,\set{0}\,)\,$. The analytical continuation theorem ensures that two solution germs at a point $x$ obtained from two \emph{homotopic} paths starting at $x_{\,0}$ are identical. This fundamental theorem of the analytical continuation guarantees that the \emph{uniform} solution is well-defined in each homotopy class of paths \cite{Forster1981}.

\begin{figure}[H]
  \centering
  \includegraphics[width=0.51\textwidth]{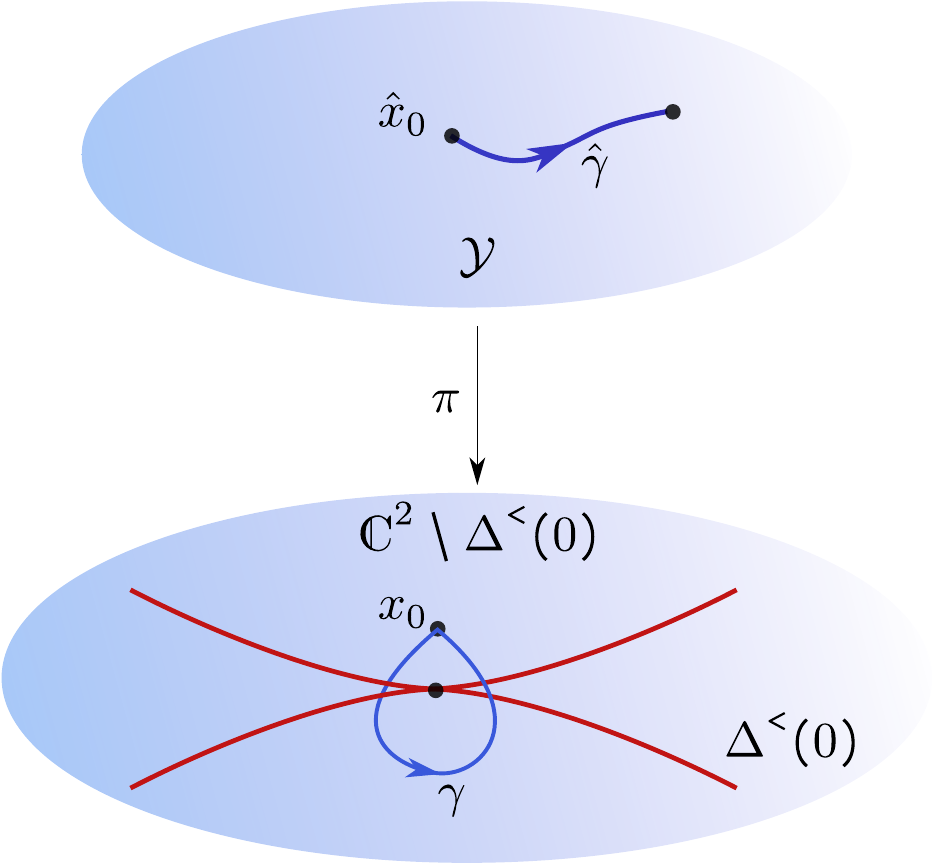}
  \caption{\small\em Illustration of the solution uniformization on the universal covering $\Yy\,$.}
  \label{fig:uni}
\end{figure}

We shall use this analytical continuation technique in order to lift the obtained local solution to the global one on the universal covering. Let us denote the universal covering\footnote{Let $\X$ be a connected topological manifold. We shall say that a topological manifold $\Yy$ is a covering of $\X$ if there exists a surjection $\pi\,:\ \X\ \twoheadrightarrow\ \Yy$ such that $\forall\,x\ \in\ \X\,$, $\exists\,\U_{\,x}$ an open neighbourhood with the property that $\pi^{\,<}\,\bigl(\,\set*{\U_{\,x}}\,\bigr)\ =\ \cup_{i\,\in\,I}\,V_{\,i}$ is open in $\Yy$ and $\forall\,i\ \in\ I$ the maps $\pi_{\,i}\,:\ V_{\,i}\ \longrightarrow\ \U_{\,x}$ are diffeomorphisms. A covering is said to be the \emph{universal covering} if $\Yy$ is also simply connected (we also remind that a space is called `simply connected' if any closed loop is homotopic to a point). The equivalence classes of homotopic paths give basically the universal covering. The term `\textbf{the} universal covering' is not an abuse of the language because all universal coverings of space are isomorphic to each other.} of $\X\ \eqdef\ \C^{\,2}\,\setminus\,\Delta^{\,<}\,(\,\set{0}\,)$ by $\Yy\,$. In the previous Section, we already mentioned that the fundamental group of $\X$ is $\Br_{\,3}$ and the \nm{Galois} group associated to \cref{eq:cube} is $S_{\,3}\,$. It is known from the theory of groups that there exists a group epimorphism\footnote{An epimorphism is simply a surjection subject to some additional properties.}
\begin{equation*}
  \chi\,:\ \Br_{\,3}\ \twoheadrightarrow\ S_{\,3}\,.
\end{equation*}
It is also a well-known fact that $\ker\chi\ \trianglelefteq\ \Br_{\,3}$ is a normal sub-group and $\Br_{\,3}\,/\,\ker\chi\ \cong\ S_{\,3}\,$. Moreover, for every $\gamma\ \in\ \ker\chi\,$, the monodromy action of $\chi\,(\,\gamma\,)$ fixes each solution $z$ of \eqref{eq:cube}. We introduce also an intermediate covering $\Rr\ \eqdef\ \Yy\,/\,\ker\chi\,$. It is composed of $\Yy$ elements orbits by the (sub-)group action $\ker\chi\,$. It is also the \nm{Riemann} surface which defines the \nm{Galois} covering of the structural automorphisms group (\ie $\ker\chi$). This situation is represented in the following commutative diagram:
\begin{equation*}
\begin{tikzcd}
  \pi_{\,1}\,(\,\X\,)\arrow[r] & \Yy \arrow[dd, "\pi", two heads]\arrow[rd] & \\
                               &                          & \Rr\arrow[ld, "\pi_{\,\chi}", two heads] \\
                               & \X                       &  
\end{tikzcd}
\end{equation*}
Hence, the solution germ $u$ on $\X$ induces a uniform holomorphic function $u\comp\pi_{\,\chi}$ on the intermediate covering $\Rr\,$. In a similar way, the solution germ $u$ can be lifted to the universal covering $\Yy$ by taking the composition $u\comp\pi\,$. We added a little illustration of the uniformization process in \cref{fig:uni}. This completes our short description of the solution uniformization.

\subsection{The contact geometry approach}
\label{sec:cg}

In this Section, we recover the same swallow tail singularity using the contact geometry approach. It can also be seen as a holomorphic method of characteristics.

Consider a complex space\footnote{We may think of it as a $1-$jet space of a complex function $(\,t,\,x)\ \mapsto\ u\,(\,t,\,x\,)\,$.}
\begin{equation*}
  \C^{\,5}\ \equiv\ \C^{\,2\,+\,2\,+\,1}\ \eqdef\ \set*{(\,t,\,x,\,u,\,\tau,\,\xi\,)\ \in\ \C^{\,5}}\,,
\end{equation*}
together with a holomorphic contact form
\begin{equation*}
  \omega\ \eqdef\ \ud\,u\ -\ \tau\,\ud\,t\ -\ \xi\,\ud\,x\,,
\end{equation*}
which defines the contact structure on $\C^{\,5}\,$. Let us also consider a first-order holomorphic \acrshort{pde}:
\begin{equation*}
  \F\,(\,t,\,x,\,u,\,u_{\,t},\,u_{\,x}\,)\ \equiv\ \F\,(\,t,\,x,\,u,\,\tau,\,\xi\,)\ =\ \fO\,.
\end{equation*}
Typically, the practical application we are interested in is given by the \acrfull{ibe}:
\begin{equation}\label{eq:cb}
  \F\,(\,t,\,x,\,u,\,\tau,\,\xi\,)\ \coloneq\ \tau\ -\ u\,\xi\ =\ \fO\,.
\end{equation}
It is assumed that the family $\set*{\ud\,\F,\,\omega}$ is free\footnote{In the contrary case, the contact field $\H_{\,\F}$ will be identically equal to zero.} and $\F$ is a real smooth function of its arguments.

The contact field $\H_{\,\F}$ is defined by the following condition \cite{Arnold1997}:
\begin{equation*}
  \ud\,\omega\,(\,\H_{\,\F},\,h\,)\ \equiv\ \ud\,\F\scal h\,, \qquad \forall\,h\ \in\ \ker\,\omega\,.
\end{equation*}
In our case, the contact field can be easily computed explicitly\footnote{The contact field can be easily computed explicitly for a general $\F$ as well \cite[Chapter~V, Section~\textsection1]{Leichtnam1985}:
\begin{equation*}
  \H_{\,\F}\ \coloneq\ \biggl(\,\pd{\F}{\tau},\,\pd{\F}{\xi},\,\tau\,\pd{\F}{\tau}\,+\,\xi\,\pd{\F}{\xi},\,-\,\Bigl(\,\pd{\F}{t}\,+\,\tau\,\pd{\F}{u}\,\Bigr),\,-\,\Bigl(\,\pd{\F}{x}\,+\,\xi\,\pd{\F}{u}\,\Bigr)\,\biggr)\,.
\end{equation*}}:
\begin{equation}\label{eq:cf}
  \H_{\,\F}\ \bydef\ \bigl(\,\F_{\,\tau},\,\F_{\,\xi},\,\tau\,\F_{\,\tau}\ +\ \xi\,\F_{\,\xi},\,-\,(\,\F_{\,\tau}\ +\ \tau\,\F_{\,u}\,),\,-\,(\,\F_{\,x}\ +\ \xi\,\F_{\,u}\,)\,\bigr)\,.
\end{equation}
By $\phi^{\,t}$ we denote the flow map along the vector field $\H_{\,\F}\,$:
\begin{equation*}
 \begin{dcases*}
  \ \od{}{t}\;\phi^{\,t}\,(\,z\,)\,\Bigr\vert_{\,t\,=\,0}\ =\ \H_{\,\F}\,(\,z\,)\,, \\
  \ \phi^{\,0}\,(\,z\,)\ =\ z\,.
 \end{dcases*}
\end{equation*}
The integral curves of this flow are called the \emph{characteristics}. By construction, the function $\F$ is constant along this flow. We can easily check it by direct computation:
\begin{multline*}
  \od{}{t}\;\Bigl(\,\F\,\bigl(\,\phi^{\,t}\,(\,z\,)\,\bigr)\,\Bigr)\ \stackrel{\text{\eqref{eq:cf}}}{=}\ \F_{\,t}\,\F_{\,\tau}\ +\ \F_{\,x}\,\F_{\,\xi}\ +\ \F_{\,u}\,(\,\tau\,\F_{\,\tau}\ +\ \xi\,\F_{\,\xi}\,)\ +\\ 
  \F_{\,\tau}\,(\,-\,\F_{\,t}\ -\ \tau\,\F_{\,u}\,)\ +\ \F_{\,\xi}\,(\,-\,\F_{\,x}\ -\ \xi\,\F_{\,u}\,)\ \equiv\ \fO\,.
\end{multline*}
Obviously, the last identity is the consequence of the fact that
\begin{equation*}
  \ud\,\F\cdot\H_{\,\F}\ \bydef\ \ud\,\omega\,(\,\H_{\,\F},\,\H_{\,\F}\,)\ =\ \fO\,.
\end{equation*}
It is not difficult to check that $\H_{\,\F}\ \in\ \ker\,\omega\,$.

We recall that the \nm{Legendrian} sub-variety $\V$ in our case is a holomorphic, possibly singular, integral sub-variety of $\ker\,\omega$ of dimension $\dim\,\V\ =\ 2$ such that
\begin{align}\label{eq:leg}
 \begin{split}
  \T\,\V\ &\subseteq\ \ker\,\omega\,, \\
  \V\ &\subseteq\ \F^{\,-1}\,(\,0\,)\,.
 \end{split}
\end{align}
We shall say that $\V$ is the generalized solution in the sense of \nm{Lie}. We know that a contact transformation allows obtaining a classical solution locally given by the function graph. If the canonical projection
\begin{align*}
  \pi\,:\ \V\ &\longrightarrow\ \C^{\,2}\,, \\
  (\,t,\,x,\,u,\,\tau,\,\xi\,)\ &\mapsto\ (\,t,\,x\,)
\end{align*}
is of (maximal) rank $2\,$, then $\V$ is a $1-$jet $\jj^{\,1}\,(\,u\,)$ of a complex holomorphic function $u\,$, \ie
\begin{equation*}
  \V\ \equiv\ \Set*{(\,t,\,x\,)\ \in\ \C^{\,2}}{\bigl(\,t,\,x,\,u\,(\,t,\,x\,),\,u_{\,t}\,(\,t,\,x\,),\,u_{\,x}\,(\,t,\,x\,)\,\bigr)}\,.
\end{equation*}
We remind that all considerations in this article are local in the vicinity of $(\,0,\,0\,)\ \in\ \C^{\,2}\,$. The \nm{Legendrian} $\V$ should be understood in the same sense.

The \nm{Legendrian} $\V$ can be also defined by a \emph{generating function} \cite{Arnold1997}, \ie $(\,\tau,\,x\,)\ \mapsto\ \Ss\,(\,\tau,\,x\,)\,$:
\begin{equation*}
  \V\ =\ \Set*{(\,\tau,\,x\,)\ \in\ \C^{\,2}}{\bigl(\,-\,\Ss_{\,\tau},\,\tau,\,\Ss\ -\ \tau \Ss_{\,\tau},\,x,\,\Ss_{\,x}\,\bigr)}\,.
\end{equation*}
It is not difficult to check that conditions \eqref{eq:leg} hold. For example,
\begin{equation*}
  \ud\,\bigl(\,\Ss\ -\ \tau\,\Ss_{\,\tau}\,\bigr)\ -\ \tau\,\ud\,\bigl(\,-\,\Ss_{\,\tau}\,\bigr)\ -\ \Ss_{\,x}\,\ud\,x\ =\ \ud\,\Ss\ -\ \Ss_{\,\tau}\,\ud\,\tau\ -\ \Ss_{\,x}\,\ud\,x\ \equiv\ \fO\,.
\end{equation*}


\subsubsection{Application to the \acrshort{ibe}}

In this case, the first order holomorphic \acrshort{pde} was defined above in \cref{eq:cb}. In this case, the contact field becomes:
\begin{equation}\label{eq:hfb}
  \H_{\,\F}\ \bydef\ \bigl(\,1,\,-\,u,\,\underbrace{\tau\ -\ u\,\xi}_{\equiv\ \fO\text{ by \eqref{eq:cb}}},\,\tau\,\xi,\,\xi^{\,2}\,\bigr)\,.
\end{equation}
Using the \acrshort{ibe} \eqref{eq:cb} and \nm{Cauchy} data \eqref{eq:cauchy2}, we may reconstruct the $1-$jet of the function $u$ at $t\ =\ 0\,$:
\begin{equation*}
  \jj^{\,1}\,(\,u\,)\ \bigcap\ \set*{t\ =\ 0}\ \equiv\ \biggl(\,0,\,x,\,x^{\,\third},\,\frac{1}{3}\;x^{\,-\,\third},\,\frac{1}{3}\;x^{\,-\,\twothirds}\,\biggr)\ \defeq\ z_{\,0}\,.
\end{equation*}
The last point $z_{\,0}$ will be used as an initial condition in finding the flow $\phi^{\,s}\,(\,z_{\,0}\,)$ (we change the flow parameter $t$ to $s\,$, because $t$ is already employed in the \acrshort{ibe}). To simplify computations, we notice also that the third component in $\H_{\,\F}$ vanishes (\cf \cref{eq:hfb}) since $\F$ conserves its value along the flow trajectories:
\begin{equation*}
 \begin{dcases*}
  \ \od{}{t}\;\phi^{\,s}\,(\,z_{\,0}\,)\,\Bigr\vert_{\,s\,=\,0}\ =\ \H_{\,\F}\,(\,z_{\,0}\,)\ \equiv\ \bigl(\,1,\,-\,u,\,0,\,\tau\,\xi,\,\xi^{\,2}\,\bigr)\,, \\
  \ \phi^{\,0}\,(\,z_{\,0}\,)\ =\ z_{\,0}\,.
 \end{dcases*}
\end{equation*}
The flow can be explicitly written in this case:
\begin{multline*}
  \phi^{\,s}\,(\,z_{\,0}\,)\ \equiv\ \bigl(\,t\,(\,s\,),\,x\,(\,s\,),\,u\,(\,s\,),\,\tau\,(\,s\,),\,\xi\,(\,s\,)\,\bigr)\ \coloneq\\ 
  \biggl(\,s,\;x\ -\ x^{\,\third}\,s,\;x^{\,\third},\;\frac{x^{\,\third}}{3\,x^{\,\twothirds}\ -\ s},\;\frac{1}{3\,x^{\,\twothirds}\ -\ s}\,\biggr)\,.
\end{multline*}
Hence, $u\ \equiv\ x^{\,\third}\,$. We introduce a new variable
\begin{equation*}
  y\ \eqdef\ x\ -\ x^{\,\third}\,t\ \equiv\ u^{\,3}\ -\ u\,t\,.
\end{equation*}
In these new variables, the flow becomes:
\begin{equation*}
  \phi^{\,s}\,(\,z_{\,0}\,)\ \equiv\ \biggl(\,t,\,y,\;u,\;\frac{u}{3\,u^{\,2}\ -\ t},\;\frac{1}{3\,u^{\,2}\ -\ t}\,\biggr)\,,
\end{equation*}
with an additional algebraic relation:
\begin{equation}\label{eq:alg}
  u^{\,3}\ -\ u\,t\ -\ y\ =\ 0\,.
\end{equation}
The singularity occurs when $t\ \to\ 3\,u^{\,2}\,$. It is not difficult to see that by definition of variably $y$, and we have in the limit $y\ \to\ -\,2\,u^{\,3}\,$. After eliminating $u\,$, we obtain an important algebraic relation between $t$ and $y$ on the singularity:
\begin{equation*}
  4\,t^{\,3}\ -\ 27\,y^{\,2}\ =\ 0\,.
\end{equation*}
It corresponds to the swallow tail (cusp) singularity seen above in Section~\ref{sec:ck}.

\begin{remark}
If in \nm{Cauchy} problem \eqref{eq:cauchy} we take another ramified initial data
\begin{equation*}
  u\,(\,0,\,x\,)\ =\ x^{\,\fifth}\,,
\end{equation*}
then direct computations show that the algebraic relation between $u\,$, $t$ and $y$ would be
\begin{equation*}
  u^{\,5}\ -\ u\,t\ -\ y\ =\ 0\,.
\end{equation*}
\end{remark}


\subsection{Generalizations}

A natural multidimensional generalization of the \acrshort{ivp} \eqref{eq:cauchy2} for an unknown function
\begin{align*}
  u\,:\ \C\times\C^{\,n}\ &\longrightarrow\ \C\,, \\
  (\,t,\,x_{1},\,x_{2},\,\ldots,\,x_{n}\,)\ &\mapsto\ u\,(\,t,\,x_{1},\,x_{2},\,\ldots,\,x_{n}\,)
\end{align*}
is the following one:
\begin{subequations}\label{eq:cauchy3}
\begin{align}
  u_{\,t}\ -\ u\,\bigl(\,u_{\,x_1}\ +\ u_{\,x_2}\ +\ \ldots\ +\ u_{\,x_n}\,\bigr)\ &=\ \fO\,, \\
  u\,(\,0,\,x_{1},\,x_{2},\,\ldots,\,x_{n}\,)\ &=\ u_{\,0}\,(\,x_{1},\,x_{2},\,\ldots,\,x_{n}\,)\ +\ x_{1}^{\,\frac{1}{p}}\,,
\end{align}
\end{subequations}
for some holomorphic function $u_{\,0}$ in the vicinity of $(\,0,\,0,\ldots,\,0\,)\ \in\ \C^{\,n}\,$. The right-hand side in the multi-dimensional \acrshort{pde} could be a polynomial in variables $t\,$, $x$ and $u\,$. In the homogeneous case, it is not difficult to show the following
\begin{theorem}
Consider the \acrshort{ivp} \eqref{eq:cauchy3}. Then, there exist $p$ germs of holomorphic functions $a_{\,0}\,$, $a_{\,1}\,$, \ldots, $a_{\,p-1}$ in the vicinity of $(\,t,\,x_{1},\,x_{2},\,\ldots,\,x_{n}\,)\ \equiv\ (\,t,\,\x\,)\ \coloneq\ (0,\,0,\ldots,\,0)\ \in\ \C^{\,n+1}$ along with a holomorphic function $(\,t,\,\x\,)\ \mapsto\ v\,(\,t,\,\x\,)$ defined in the vicinity of the origin such that $u\,(\,t,\,\x\,)\ \eqdef\ v\,(\,t,\,\x\,)\ +\ z\,(\,t,\,\x\,)$ is the unique solution to \eqref{eq:cauchy3}. The function $(\,t,\,\x\,)\ \mapsto\ z\,(\,t,\,\x\,)$ is defined as a solution to the following algebraic equation (with $t$ and $\x$ being seen as parameters):
\begin{equation*}
  z^{\,p}\ =\ a_{\,p-1}\,(\,t,\,\x\,)\,z^{\,p-1}\ +\ \ldots\ a_{\,1}\,(\,t,\,\x\,)\,z\ +\ a_{\,0}\,(\,t,\,\x\,)
\end{equation*}
together with `initial' conditions:
\begin{align*}
  a_{\,0}\,(\,0,\,x_{1},\,x_{2},\,\ldots,\,x_{n}\,)\ &=\ x_{\,1}\,, \\
  a_{\,j}\,(\,0,\,x_{1},\,x_{2},\,\ldots,\,x_{n}\,)\ &=\ \fO\,, \qquad j\ \in\ (\,p\,-\,1\,)^{\,\sqsupset}\,.
\end{align*}
\end{theorem}
\begin{proof}
This can be done by analogy with the proof of \cref{thm:ck}.
\end{proof}


\section{A ring of algebraic convergent power series}
\label{sec:ring}

In the present Section, we glimpse the very rich and complex geometry of the ring\footnote{If necessary, this ring can also be considered an algebra over the field of complex numbers $\C\,$.} of formal power series. The goal of this Section is to provide a differential study of the algebra $\O\,\llbracket\,z\,\rrbracket$ (with $z^{\,3}\ =\ p\,z\ +\ q$). This study will exhibit some of the tools allowing us to construct the Algorithm of Section~\ref{sec:it}, which will provide a test of our \cref{conj:2}.

To fix the ideas, let us consider the following ring \cite{VanderPut1977}:
\begin{definition}
We define
\begin{equation*}
  \O\,\llbracket\,z\,\rrbracket\ \eqdef\ \set*{\sum_{j\,=\,0}^{+\,\infty} a_{\,j}\,(\,p,\,q\,)\,z^{\,j}}\,,
\end{equation*}
where the coefficients $\set*{a_{\,j}\,(\,p,\,q\,)}_{\,j\,\in\,\N}$ are the germs of analytic functions at the origin\footnote{All germs in this study are implicitly supposed to be defined at the origin $(\,0,\,0\,)\ \in\ \C^{\,2}\,$.} $(\,0,\,0\,)\ \in\ \C^{\,2}$ and the complex variable $z$ satisfies the algebraic equation \eqref{eq:cube}. The algebra of these germs is denoted by $\O\,$. We require also that there exists $\exists\,\eps\ >\ 0$ (depending on the $\set*{a_{\,j}\,(\,p,\,q\,)}_{\,j\,\in\,\N}$) such that
\begin{equation*}
  \sum_{j\,=\,0}^{+\,\infty}\,\sup\limits_{\substack{\abs{p}\ <\ \eps \\ \abs{q}\ <\ \eps}} \abs{a_{\,j}\,(\,p,\,q\,)}\,\eps^{\,j}\ <\ +\,\infty\,.
\end{equation*}
\end{definition}

By applying the \nm{Weierstra}\ss{} division theorem \cite{Narasimhan1966}, we may show that this ring is isomorphic to the following reduced form:
\begin{equation*}
  \O\,[\,z\,]\ \eqdef\ \Set*{a\,(\,p,\,q\,)\ +\ b\,(\,p,\,q\,)\,z\ +\ c\,(\,p,\,q\,)\,z^{\,2}}{z^{\,3}\ =\ p\,z\ +\ q}\,,
\end{equation*}
for some germs of analytic functions $a\,$, $b$ and $c\,$. Our goal consists of building some form of the differential and integral calculus on $\O\,\llbracket\,z\,\rrbracket$ (or, in our case, equivalently on $\O\,[\,z\,]$).

\begin{remark}\label{rem:p}
Notice that we may remove the explicit dependence on variable $q$ in coefficients $a\,$, $b$ and $c$ since according to \cref{eq:cube}, $q\ =\ z^{\,3}\ -\ p\,z$ and the coefficients ($a\,$, $b\,$, $c$) may be re-expanded in the second variable $q\ \mapsto\ \set*{a,\,b,\,c}\,(\,p,\,q\,)$ after the substitution of the new expression for $q\,$. Finally, after collecting all the terms, we obtain the equivalent representation of the ring using a formal power series:
\begin{equation*}
  \O\,[\,z\,]\ \eqdef\ \Set*{\sum_{j\,=\,0}^{+\,\infty}\,b_{\,j}\,(\,p\,)\,z^{\,j}}{z^{\,3}\ =\ p\,z\ +\ q}\,.
\end{equation*}
It is not difficult to see that the same remark applies also to the ring $\O\,\llbracket\,z\,\rrbracket$ that we considered from the beginning:
\begin{equation*}
  \sum_{j\,=\,0}^{+\,\infty} a_{\,j}\,(\,p,\,q\,)\,z^{\,j}\ \equiv\ \sum_{j\,=\,0}^{+\,\infty} b_{\,j}\,(\,p\,)\,z^{\,j}\,.
\end{equation*}
Namely, we have demonstrated the following
\begin{prop}
Every element $\sum_{j\,=\,0}^{+\,\infty} a_{\,j}\,(\,p,\,q\,)\,z^{\,j}\ \in\ \O\,\llbracket\,z\,\rrbracket$ may be written under the form $\sum_{j\,=\,0}^{+\,\infty} b_{\,j}\,(\,p\,)\,z^{\,j}$ for some appropriately chosen holomorphic coefficients $\set*{b_{\,j}\,(\,p\,)}_{\,j\,=\,0}^{\,+\,\infty}\,$.
\end{prop}
\end{remark}

Now, we start to construct the differential calculus over the rings $\O\,\llbracket\,z\,\rrbracket$ and $\O\,[\,z\,]\,$:
\begin{lemma}\label{eq:lem1}
The partial derivatives of the solution $z$ with respect to coefficients $p$ and $q$ are given by:
\begin{equation*}
  \pd{z}{p}\ =\ \frac{z}{3\,z^{\,2}\ -\ p}\,, \qquad
  \pd{z}{\,q}\ =\ \frac{1}{3\,z^{\,2}\ -\ p}\,.
\end{equation*}
\end{lemma}

\begin{proof}
By taking the differential of Equation~\eqref{eq:cube}, we obtain:
\begin{equation*}
  (\,3\,z^{\,2}\ -\ p\,)\,\ud\,z\ =\ z\,\ud p\ +\ \ud q
\end{equation*}
or
\begin{equation*}
  \ud\,z\ =\ \frac{z}{3\,z^{\,2}\ -\ p}\,\ud p\ +\ \frac{1}{3\,z^{\,2}\ -\ p}\;\ud q\,.
\end{equation*}
After identifying the left and right-hand sides, we obtain the required result.
\end{proof}

\begin{remark}\label{rem:wf}
Lemma~\ref{eq:lem1} (and \cref{sec:ck}) shows that the derivatives of a branch of $z\,(\,p,\,q\,)$ blow up when the point $(\,p,\,q\,)$ tends to a point of the discriminant locus $\Delta^{\,<}\,(\,\set{0}\,)\,$. 
\end{remark}


\subsection{A lyrical digression}

\begin{figure}
  \centering
  \includegraphics[width=0.59\textwidth]{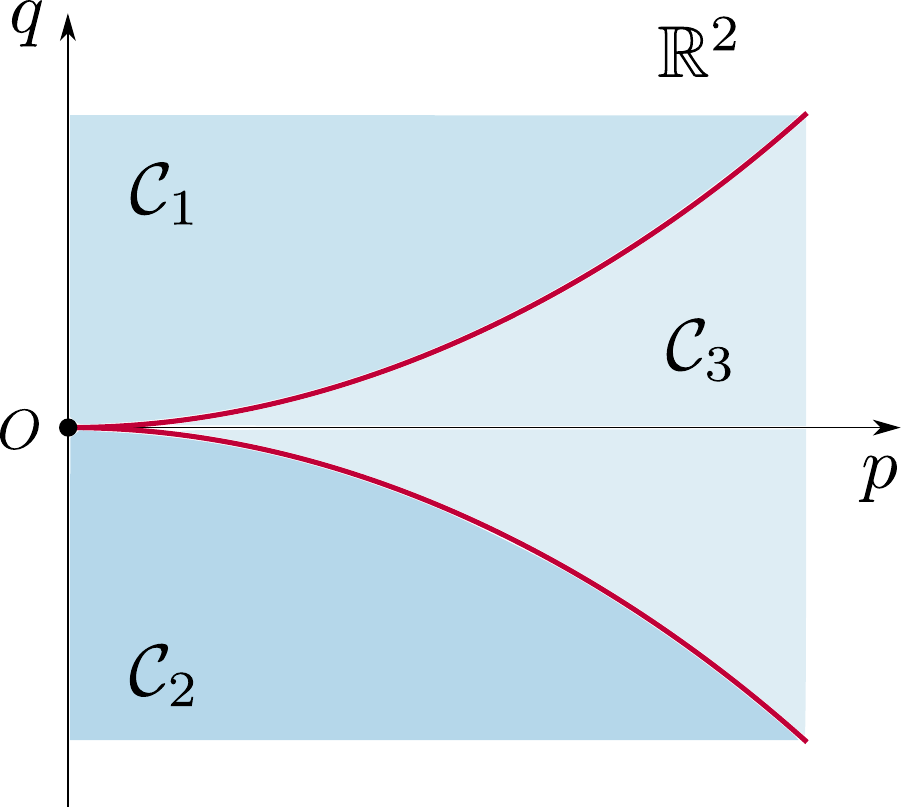}
  \caption{\small\em Division of the real domain $M\,\setminus\,\Delta^{\,<}\,(\,\set{0}\,)\ \subseteq\ \R^{\,2}$ in three regions $\Cl_{\,1,\,2,\,3}\,$.}
  \label{fig:real}
\end{figure}

Let us return to the real space $\R^{\,2}\ \subseteq\ \C^{\,2}$ for the moment. At the current stage, we conjecture the following intuitive, but precise statement about micro-local singularities of $z$ and $z^{\,2}$ (real) distributions:
\begin{prop}\label{prop:wf}
\begin{equation}\label{eq:wf}
  \WF\,(\,\cz\,)\ =\ \WF\,(\,\cz^{\,2}\,)\ =\ \overline{\Nn\,\bigl(\,\Delta^{\,<}\,(\,\set{0}\,)\,\setminus\,\set*{\vO\ \in\ \R^{\,2}}\,\bigr)}\,,
\end{equation}
where by the bar we denote the closure operation and $\Nn\,(\cdot)$ is the co-normal set and the distributions $\cz$ and $\cz^{\,2}$ are defined below.
\end{prop}
\begin{proof}
This proposition is left for our future works on the subject. Namely, in future work, we wish to apply the results of this manuscript to the construction, in the real domain, of weak solutions to some genuinely nonlinear \acrshort{pde}s. The first step would be the following investigation. Divide the domain $\set*{(\,p,\,q\,)\ \in\ \R_{\,\geq\,0}\,\times\,\R}\,\setminus\,\Delta^{\,<}\,(\,\set{0}\,)$ in three connected components $\Cl_{\,1,\,2,\,3}\,$ as depicted in \cref{fig:real}. Choose in each branch $\Cl_{\,j}$ a branch $z_{\,j}$ of solution to the \cref{eq:cube}. Let us denote by $\cz$ the distribution on $M\ \eqdef\ \set*{(\,p,\,q\,)\ \in\ \R_{\,\geq\,0}\,\times\,\R}$ obtained by gluing together all the local pieces $\set*{z_{\,j}}_{\,j\,=\,1}^{\,3}\,$. Then, in a separate work we plan to address the following question: would it be possible to choose the $z_{\,j}$ in each $\Cl_{\,j}$ so that the wave front sets $\WF\,(\,\cz\,)$ and $\WF\,(\,\cz^{\,2}\,)$ is included in the intersection of the co-tangent bundle $T^{\,\ast}\,M$ with the co-normal to $\Delta^{\,<}\,(\,\set{0}\,)\,$? It is in this sense that the proposition above has to be understood.
\end{proof}
In the last \cref{eq:wf}, there is an obvious part of the statement:
\begin{align*}
  \singsupp\,(\,z\,)\ &\subseteq\ \Delta^{\,<}\,(\,\set{0}\,)\,, \\
  \singsupp\,(\,z^{\,2}\,)\ &\subseteq\ \Delta^{\,<}\,(\,\set{0}\,)\,.
\end{align*}

The deeper sense of \cref{prop:wf} consists in making the connection between our theory and the more classical theory of shock waves with the first steps in this direction given in \cref{app:a}.


\subsection{Primitivization}

For the reasons which will become clearer below, the goal of this Section is to introduce and study the primitive of the ring $\O\,[\,z\,]$ with respect to the variable $q\,$:
\begin{equation*}
  \partial_{\,q}^{\,-\,1}\,\O\,[\,z\,]\ \eqdef\ \Set*{a\ +\ b\,z\ +\ c\,z^{\,2}\ \in\ \O\,[\,z\,]}{\partial_{\,q}\,(\,a\ +\ b\,z\ +\ c\,z^{\,2}\,)\ \in\ \O\,[\,z\,]}\,.
\end{equation*}
The following Lemma shows that the derivative $\partial_{\,q}$ of an element of $\O\,[\,z\,]$ belongs to $\O\,[\,z\,]$ only under certain conditions:
\begin{lemma}
  Let $z$ be a solution to \cref{eq:cube} and $\Delta\ \bydef\ 4\,p^{\,3}\ -\ 27\,q^{\,2}$ be its discriminant. Take an element $a\ +\ b\,z\ +\ c\,z^{\,2}\ \in\ \O\,[\,z\,]\,$. Then, $\partial_{\,q}\,(\,a\ +\ b\,z\ +\ c\,z^{\,2}\,)\ \in\ \O\,[\,z\,]$ if and only if the discriminant $\Delta$ divides $6\,c\,q\ -\ 2\,b p$ and $9\,b\,q\ -\ 4\,c p^{\,2}\,$.
\end{lemma}
\begin{proof}

Let us show the main steps of computations which yield the required result. First of all, we differentiate the element $a\ +\ b\,z\ +\ c\,z^{\,2}\ \in\ \O\,[\,z\,]$ with respect to $q$ and we use the result from \cref{eq:lem1}:
\begin{equation*}
  \partial_{\,q}\,(\,a\ +\ b\,z\ +\ c\,z^{\,2}\,)\ =\ a_{\,q}\ +\ b_{\,q}\,z\ +\ c_{\,q}\,z^{\,2}\ +\ \frac{b\ +\ 2\,c\,z}{3\,z^{\,2}\ -\ p}\,.
\end{equation*}
It is clear that the derivative $\partial_{\,q}\,(\,a\ +\ b\,z\ +\ c\,z^{\,2}\,)\ \in\ \O\,[\,z\,]$ if and only if the last fraction can be presented as
\begin{equation*}
  \frac{b\ +\ 2\,c\,z}{3\,z^{\,2}\ -\ p}\ \equiv\ \alpha\ +\ \beta\,z\ +\ \gamma\,z^{\,2}
\end{equation*}
for some holomorphic germs $\alpha\,$, $\beta$ and $\gamma\,$. The system of algebraic equations to determine these germs is obtained from the following obvious relation:
\begin{equation*}
  b\ +\ 2\,c\,z\ =\ (\,\alpha\ +\ \beta\,z\ +\ \gamma\,z^{\,2}\,)\,(\,3\,z^{\,2}\ -\ p\,)\,.
\end{equation*}
After expanding the powers of $z$ higher than two using \cref{eq:cube} and identifying the coefficients in front of $z^{\,0}\,$, $z^{\,1}$ and $z^{\,2}\,$, we obtain the following system of equations:
\begin{align*}
  3\,\beta\,q\ -\ \alpha\,p\ &=\ b\,, \\
  3\,\gamma\,q\ +\ 2\,\beta\,p\ &=\ 2\,c\,, \\
  3\,\alpha\ +\ 2\,\gamma\,p\ &=\ \fO\,.
\end{align*}
Simple algebraic manipulations yield:
\begin{align*}
  \alpha\ &\coloneq\ -\;\frac{2\,p\,(\,6\,c\,q\ -\ 2\,b\,p\,)}{\Delta}\,, \\
  \beta\ &\coloneq\ \frac{9\,b\,q\ -\ 4\,c\,p^{\,2}}{\Delta}\,, \\
  \gamma\ &\coloneq\ \frac{3\,(\,6\,c\,q\ -\ 2\,b\,p\,)}{\Delta}\,.
\end{align*}
This yields the conclusion of this Lemma.
\end{proof}
Similar conditions can be formulated for the derivation with respect to $p\,$. However, it is out of the scope of our study.

The next question we have to discuss is the integral calculus on $\O\,[\,z\,]$ in the sense of the operation inverse to the differentiation. Let us focus on the independent variable\footnote{The reason is that the integration only with respect to this variable $q$ can be shown to be well-defined.} $q$ and introduce the following
\begin{definition}
An element $w\,(\,z\,)\ \in\ \O\,[\,z\,]$ is called a primitive of $v\,(\,z\,)\ \in\ \O\,[\,z\,]$ if two conditions\footnote{The second condition is imposed to remove any ambiguity with respect to constant terms.} are satisfied:
\begin{itemize}
  \item $\partial_{\,q}\,w\ =\ v\,$,
  \item $w\,(\,0\,)\ =\ 0\,$.
\end{itemize}
This fact will be denoted as $w\ =\ \partial_{\,q}^{\,-\,1}\,v$ and, by definition, we have $\partial_{\,q}\comp\partial_{\,q}^{\,-\,1}\ =\ \Id_{\,\O\,[\,z\,]}\,$.
\end{definition}
To illustrate this definition, we provide one concrete example:
\begin{equation*}
  \partial_{\,q}^{\,-\,1}\,z\ =\ \frac{3\,q\,z\ +\ p\,z^{\,2}}{4}\,.
\end{equation*}
The last formula can be obtained by specifying \cref{eq:primq} to $z\,$. Otherwise, it can be checked by a direct verification\footnote{Indeed, the following sequence of equalities holds:
\begin{multline*}
   \partial_{\,q}\comp\partial_{\,q}^{\,-\,1}\,z\ =\ \partial_{\,q}\,\biggl(\,\frac{3\,q\,z\ +\ p\,z^{\,2}}{4}\,\biggr)\ =\ \frac{1}{4}\;\Bigl(\,\frac{2\,p\,z}{3\,z^{\,2}\, -\, p}\ +\ 3\,z\ +\ \frac{3\,q}{3\,z^{\,2}\, -\, p}\,\Bigr)\ = \\ 
   \frac{9\,z^{\,3}\ -\ p\,z\ +\ 3\,q}{4\,(3\,z^{\,2}\ -\ p)}\ =\ \frac{9\,z^{\,3}\ -\ p\,z\ +\ 3\,(\,z^{\,3}\ -\ p\,z\,)}{4\,(3\,z^{\,2}\ -\ p)}\ =\\ 
   \frac{12\,z^{\,3}\ -\ 4\,p\,z}{4\,(3\,z^{\,2}\ -\ p)}\ =\ \frac{4\,z\,(\,3\,z^{\,2}\ -\ p\,)}{4\,(3\,z^{\,2}\ -\ p)}\ \equiv\ z\,.
\end{multline*}}.

In the following two \cref{lem:unique,lem:prim}, we show that this definition is meaningful and effective. However, first of all, we make use of the \cref{rem:p} and switch back to the isomorphic ring object $\O\,\llbracket\,z\,\rrbracket$ since it is directly used in our construction presented in \cref{sec:it}:
\begin{lemma}\label{lem:unique}
The primitive $\partial_{\,q}^{\,-\,1}\,v\,$, $\forall\,v\ \in\ \O\,\llbracket\,z\,\rrbracket$ is unique.
\end{lemma}
\begin{proof}
In order to show this statement, it is enough to demonstrate that if $\partial_{\,q}\,w\ =\ 0$ and $w\,(\,0\,)\ =\ 0$ imply that $w\ =\ \fO\ \in\ \O\,\llbracket\,z\,\rrbracket\,$. Indeed, from the second condition, we conclude that
\begin{equation*}
  w\ =\ \sum_{j\,=\,1}^{\infty}\,c_{\,j}\,(\,p\,)\,z^{\,j}
\end{equation*}
for some germs $\set*{c_{\,j}\,(\,p\,)}_{\,j\,=\,1}^{\,\infty}\,$. Let us compute its derivative (which exists because $w$ is a primitive):
\begin{equation*}
  \partial_{\,q}\,w\ =\ \sum_{j\,=\,1}^{\infty}\,j\,c_{\,j}\,(\,p\,)\,\frac{z^{\,j\,-\,1}}{3\,z^{\,2}\,-\,p}\ =\ \fO\,.
\end{equation*}
We multiply both sides of the last equality by $3\,z^{\,2}\ -\ p\,$:
\begin{equation*}
  \sum_{j\,=\,1}^{\infty}\,j\,c_{\,j}\,(\,p\,)\,z^{\,j\,-\,1}\ =\ \fO\,.
\end{equation*}
By using the unicity theorem for power series, we conclude that
\begin{equation*}
  c_{\,j}\,(\,p\,)\ \equiv\ \fO\,, \qquad \forall\,j\ \geq\ 1\,.
\end{equation*}
\end{proof}

\begin{lemma}\label{lem:prim}
Let $\sum_{j\,=\,0}^{\infty}\,a_{\,j}\,(\,p\,)\,z^{\,j}\ \in\ \O\,\llbracket\,z\,\rrbracket\,$. Then, the primitive $\partial_{\,q}^{\,-\,1}$ of this element with respect to $q$ belongs also to $\O\,\llbracket\,z\,\rrbracket$ and can be explicitly computed (see \cref{eq:primq}).
\end{lemma}
\begin{proof}
The proof is constructive. Let us assume that there exists another element $\sum_{j\,=\,0}^{\infty}\,b_{\,j}\,(\,p\,)\,z^{\,j}\ \in\ \O\,\llbracket\,z\,\rrbracket$ such that
\begin{equation*}
  \partial_{\,q}\,\biggl(\,\sum_{j\,=\,0}^{\infty}\,b_{\,j}\,(\,p\,)\,z^{\,j}\,\biggr)\ =\ \sum_{j\,=\,0}^{\infty}\,a_{\,j}\,(\,p\,)\,z^{\,j}
\end{equation*}
and we shall construct it. By interchanging the summation and derivation operators, and using \cref{eq:lem1}, we obtain:
\begin{equation*}
  \sum_{j\,=\,1}^{\infty}\,b_{\,j}\,(\,p\,)\,\frac{j\,z^{\,j\,-\,1}}{3\,z^{\,2}\ -\ p}\ =\ \sum_{j\,=\,0}^{\infty}\,a_{\,j}\,(\,p\,)\,z^{\,j}\,.
\end{equation*}
By multiplying both sides by $3\,z^{\,2}\ -\ p$ and rearranging the terms, we obtain:
\begin{equation*}
  \sum_{j\,=\,0}^{\infty}\,(\,j\ +\ 1\,)\,b_{\,j\,+\,1}\,(\,p\,)\,z^{\,j}\ =\ \sum_{j\,=\,0}^{\infty}\,\Bigl(\,3\,a_{\,j\,-\,2}\ -\ p\,a_{\,j}\,(\,p\,)\,\Bigr)\,z^{\,j}\,,
\end{equation*}
where we assume the following convention $a_{\,-\,1}\ \equiv\ a_{\,-\,2}\ \eqdef\ \fO\,$. It can be readily deduced that
\begin{equation*}
  b_{\,j\,+\,1}\,(\,p\,)\ \coloneq\ \frac{1}{j\ +\ 1}\;\bigl(\,-\,p\,a_{\,j}\,(\,p\,)\ +\ 3\,a_{\,j\,-\,2}\,(\,p\,)\,\bigr)\,, \qquad \forall\,j\ \in\ \N\,.
\end{equation*}
Henceforth, we obtain the required result:
\begin{empheq}[box={\mymath}]{equation}\label{eq:primq}
  \partial_{\,q}^{\,-\,1}\,\biggl(\,\sum_{j\,=\,0}^{\infty}\,a_{\,j}\,(\,p\,)\,z^{\,j}\,\biggr)\ =\ \underbrace{\sum_{j\,=\,1}^{\infty}\,\frac{1}{j}\;\bigl(\,-\,p\,a_{\,j\,-\,1}\,(\,p\,)\ +\ 3\,a_{\,j\,-\,3}\,(\,p\,)\,\bigr)\,z^{\,j}}_{\displaystyle{\defeq\ w\ \in\ \O\,\llbracket\,z\,\rrbracket}}\,.
\end{empheq}
Once we constructed the primitive given in the right-hand side of the above equation, it is not difficult to check by the direct calculation that we have found the desired element in $\O\,\llbracket\,z\,\rrbracket\,$:
\begin{equation*}
  \partial_{\,q}\,w\ \bydef\ \partial_{\,q}\,\biggl(\,\sum_{j\,=\,1}^{\infty}\,\frac{1}{j}\;\bigl(\,-\,p\,a_{\,j\,-\,1}\,(\,p\,)\ +\ 3\,a_{\,j\,-\,3}\,(\,p\,)\,\bigr)\,z^{\,j}\,\biggr)\ =\ \sum_{j\,=\,0}^{\infty}\,a_{\,j}\,(\,p\,)\,z^{\,j}\,.
\end{equation*}
This is left to the reader as a simple exercise\footnote{For the sake of completeness of our manuscript, we provide here this simple computation:
\begin{multline*}
  \partial_{\,q}\,w\ =\ \sum_{j\,=\,1}^{\infty}\,\frac{1}{j}\;\bigl(\,3\,a_{\,j\,-\,3}\ -\ p\,a_{\,j\,-\,1}\,\bigr)\;\frac{j\,z^{\,j\,-\,1}}{3\,z^{\,2}\ -\ p}\ =\ \sum_{j\,=\,1}^{\infty}\,\bigl(\,3\,a_{\,j\,-\,3}\ -\ p\,a_{\,j\,-\,1}\,\bigr)\;\frac{z^{\,j\,-\,1}}{3\,z^{\,2}\ -\ p}\ =\\
  \frac{1}{3\,z^{\,2}\ -\ p}\;\biggl(\,3\,\sum_{j\,=\,1}^{\infty}\,a_{\,j\,-\,3}\,z^{\,j\,-\,1}\ -\ p\,\sum_{j\,=\,1}^{\infty}\,a_{\,j\,-\,1}\,z^{\,j\,-\,1}\,\biggr)\ =\\ 
  \frac{1}{3\,z^{\,2}\ -\ p}\;\biggl(\,3\,z^{\,2}\,\sum_{j\,=\,0}^{\infty}\,a_{\,j}\,z^{\,j}\ -\ p\,\sum_{j\,=\,0}^{\infty}\,a_{\,j}\,z^{\,j}\,\biggr)\ \equiv\ \sum_{j\,=\,0}^{\infty}\,a_{\,j}\,z^{\,j}\,.
\end{multline*}}. The found primitive $w$ in \cref{eq:primq} satisfies also obviously the second condition $w\,(\,0\,)\ =\ 0\,$.
\end{proof}

It is worthwhile to observe that one cannot define the primitive with respect to the other independent variable $p\,$. It is due, in particular, to the non-uniqueness as it follows from this simple example (employing again \cref{eq:lem1}):
\begin{equation*}
  \partial_{p}\,q\ \equiv\ \partial_{p}\,(\,z^{\,3}\ -\ p\,z\,)\ =\ 3\,z^{\,2}\;\frac{z}{3\,z^{\,2}\ -\ p}\ -\ z\ -\ p\;\frac{z}{3\,z^{\,2}\ -\ p}\ \equiv\ \fO\ =\ \partial_{p}\,\fO\,.
\end{equation*}
Hence, we just showed the following
\begin{prop}
Let $w\ \in\ \O\,\llbracket\,z\,\rrbracket\,$. Then, two conditions
\begin{itemize}
  \item $\partial_{p}\,w\ \equiv\ \fO\,$,
  \item $w\,(\,0\,)\ =\ 0$
\end{itemize}
do not necessarily imply that $w\ \equiv\ \fO\,$.
\end{prop}
\begin{proof}
See above.
\end{proof}
A conceptual explanation of the non-existence of the primitive $\partial_{\,p}^{\,-\,1}$ will be given at the very end of \cref{sec:dmod}.

The primitive with respect to $q$ has another good property that deserves to be mentioned here: $p-$differentiable in the ring $\O\,\llbracket\,z\,\rrbracket\,$. The following computation demonstrates this property and gives explicitly the formula for the $p-$derivative:
\begin{multline*}
  \partial_{p}\comp\partial_{\,q}^{\,-\,1}\,\biggl(\;\sum_{j\,=\,0}^{\infty}\,a_{\,j}\,z^{\,j}\;\biggr)\ =\ \partial_{p}\,\biggl(\;\sum_{j\,=\,1}^{\infty}\,\frac{1}{j}\;\bigl(\,3\,a_{\,j\,-\,3}\ -\ p\,a_{\,j\,-\,1}\,\bigr)\,z^{\,j}\;\biggr)\ =\\
  \underbrace{\sum_{j\,=\,1}^{\infty}\,\frac{1}{j}\;\bigl(\,3\,a_{\,j\,-\,3}^{\,\prime}\ +\ p\,a_{\,j\,-\,1}^{\,\prime}\ -\ a_{\,j\,-\,1}\,\bigr)\,z^{\,j}}_{\displaystyle{(\,\ast\,)}}\ +\ \underbrace{\sum_{j\,=\,1}^{\infty}\,\frac{z^{\,j}}{3\,z^{\,2}\ -\ p}\;\bigl(\,3\,a_{\,j\,-\,3}\ -\ p\,a_{\,j\,-\,1}\,\bigr)}_{\displaystyle{(\,\ast\,\ast\,)}}\,,
\end{multline*}
where by prime $(\cdot)^{\,\prime}$ we denote the derivative with respect to $p\,$. The first term $(\,\ast\,)$ is regular and belongs to $\O\,\llbracket\,z\,\rrbracket\,$. The second term $(\,\ast\,\ast\,)$ can be further transformed to the regular form:
\begin{multline*}
  (\,\ast\,\ast\,)\ \bydef\ \sum_{j\,=\,1}^{\infty}\,\frac{z^{\,j}}{3\,z^{\,2}\ -\ p}\;\bigl(\,3\,a_{\,j\,-\,3}\ -\ p\,a_{\,j\,-\,1}\,\bigr)\ =\\ 
  \frac{1}{3\,z^{\,2}\ -\ p}\;\biggl(\;3\,\sum_{j\,=\,0}^{\infty}\,a_{\,j}\,z^{\,j\,+\,3}\ -\ p\,\sum_{j\,=\,0}^{\infty}\,b_{\,j}\,z^{\,j\,+\,1}\;\biggr)\ =\\ 
  \frac{z}{3\,z^{\,2}\ -\ p}\;\biggl(\;3\,z^{\,2}\,\sum_{j\,=\,0}^{\infty}\,a_{\,j}\,z^{\,j}\ -\ p\,\sum_{j\,=\,0}^{\infty}\,b_{\,j}\,z^{\,j}\;\biggr)\ \equiv\ z\,\sum_{j\,=\,0}^{\infty}\,b_{\,j}\,z^{\,j}\,.
\end{multline*}
Hence, to summarize, we just demonstrated the following result:
\begin{empheq}[box={\mymath}]{equation*}
  \partial_{p}\comp\partial_{\,q}^{\,-\,1}\,\biggl(\;\sum_{j\,=\,0}^{\infty}\,a_{\,j}\,z^{\,j}\;\biggr)\ =\ \sum_{j\,=\,1}^{\infty}\,\frac{1}{j}\;\bigl(\,3\,a_{\,j\,-\,3}^{\,\prime}\ +\ p\,a_{\,j\,-\,1}^{\,\prime}\ -\ a_{\,j\,-\,1}\,\bigr)\,z^{\,j}\ +\ z\,\sum_{j\,=\,0}^{\infty}\,b_{\,j}\,z^{\,j}\,.
\end{empheq}
More generally, one can similarly show that
\begin{prop}
Let $u\ =\ \partial_{\,q}^{\,-\,k}\,v$ for some element $v\ \in\ \O\,\llbracket\,z\,\rrbracket$ and $k\ \in\ \N$. Then, $\partial_{p}^{\,\alpha}\comp\partial_{\,q}^{\,\beta}\,u\ \in\ \O\,\llbracket\,z\,\rrbracket$ provided that $\alpha\ +\ \beta\ \leq\ k\,$.
\end{prop}


\section{Holonomic coherent $\D-$modules}
\label{sec:dmod}

In this Section we perform the explicit construction of a coherent\footnote{We remind that a $\D-$module $\Mm$ is said to be \emph{coherent} if for any $m\ \in\ \N^{\,\times}$ and any $\D-$module homomorphism
\begin{equation*}
  \D^{\,m}\ \eqdef\ \underbrace{\D\,\times\,\D\,\times\ldots\times\,\D}_{\displaystyle{m\ \text{times}}}\ \longrightarrow\ \Mm
\end{equation*}
has a kernel which is locally a $\D-$module of finite type.} holonomic\footnote{As a reminder, a coherent $\D-$module is said to be \emph{holonomic} if its characteristic variety $\Vv$ is \nm{Lagrangian}. In other words, in every regular point $x\ \in\ \Vv\,$, we have
\begin{equation*}
  \T_{\,x}\,\Vv\ =\ (\,\T_{\,x}\,\Vv\,)^{\,\top}\,,
\end{equation*}
where the orthogonal complement $(\cdot)^{\,\top}$ is taken in the sense of the natural symplectic structure on the co-tangent bundle $\T^{\,\ast}\,\Vv\,$. We remind that in general one has only the \emph{involution} property of the characteristic variety \cite{Gabber1981}:
\begin{equation*}
  (\,\T_{\,x}\,\Vv\,)^{\,\top}\ \subseteq\ \T_{\,x}\,\Vv\,.
\end{equation*}
This is a symplectic geometry view on the involutivity. More algebraic points of view are also possible, see \eg \cite{Gabber1981, Singh2014}.} $\D-$module whose geometry is underlying our solution strategy of \cref{conj:2}. The goal is to study the singularities of the obtained solution. More precisely, we show that the vector\footnote{After some simplification, the solution ansatz we consider looks like $a\ +\ b\,z\ +\ c\,z^{\,2}\,$. So, one could be inclined to consider instead the vector $(\,\fW,\,z,\,z^{\,2}\,)^{\,\top}\,$. However, the first component is completely regular, and since we are interested in singularities, we retain only the two last components.} $(\,z,\,z^{\,2}\,)^{\,\top}$ is a solution of a holonomic $\D-$module whose characteristic variety $\Vv$ is included into the union of the zero section and of the co-normal to the swallow tail. In the case of a smooth hyper-surface defined by $x\ =\ \fO\,$, the multivalued function $x^{\,\frac{1}{3}}$ is a solution of the equation
\begin{equation*}
  x\,u_{\,x}\ -\ \frac{u}{3}\ =\ \fO\,.
\end{equation*}
So, morally, the ``micro-local'' singularities of $x^{\,\frac{1}{3}}$ live in the co-normal of $x\ =\ \fO$ and this property is one of the ingredients of the proof of \cref{thm:2}. Therefore, in order to understand the \cref{conj:2}, it is important to work out the analogue of this property in the case of $z$ satisfying $z^{\,3}\ =\ p\,z\ +\ q\,$. This is the goal of this Section.

So, we consider the algebraic \cref{eq:cube} whose discriminant $\Delta$ is given by
\begin{equation*}
  \Delta\ \bydef\ 4\,p^{\,3}\ -\ 27\,q^{\,2}\,.
\end{equation*}
The co-normal\footnote{Informally speaking, the co-normal $\Nn\,(\,\S\,)$ plays the r\^ole of the characteristic variety for the quasi-linear operator \eqref{eq:pde2} defined below.} $\Nn$ to $\S\ \eqdef\ \Delta^{\,<}\,(\,\set{0}\,)$ can be easily computed:
\begin{equation}\label{eq:conormal}
  \Nn\,(\,\S\,)\ \bydef\ \overline{\Nn\,(\,\S_{\,\mathrm{reg}}\,)}\ =\ \Set*{(\,3\,z^{\,2},\,-\,2\,z^{\,3},\,z\,\lambda,\,\lambda\,)}{(\,z,\,\lambda\,)\ \in\ \C^{\,2}}\ \subseteq\ \T^{\,\ast}\,\C^{\,2}\,,
\end{equation}
where we use the following co\"ordinates to parametrize the co-tangent bundle:
\begin{equation*}
  \T^{\,\ast}\,\C^{\,2}\ \equiv\ \set*{(\,p,\,q,\,\xi_{\,1},\,\xi_{\,2}\,)}\ \cong\ \C^{\,4}\,.
\end{equation*}
In the sequel, the co-normal $\Nn\,(\,\S\,)$ will be always understood as the closure $\overline{\Nn\,(\,\S_{\,\mathrm{reg}}\,)}$ according to its definition.

\begin{remark}
The fact that the co-normal $\Nn\,(\,\S\,)$ to the swallow tail contains only one complex direction\footnote{It would be more accurate to speak about a line, of course. In the real case the singular directions are $(\,0,\,\pm\,1\,)\,$, which lie along the same line.} above the origin is absolutely crucial in our study).
\end{remark}

\begin{lemma}\label{lem:p}
Let us define the following polynomials on the co-tangent bundle $\T^{\,\ast}\,\C^{\,2}\,$:
\begin{subequations}
\begin{align}
  P_{\,1}\,(\,p,\,q;\,\xi_{\,1},\,\xi_{\,2}\,)\ &\eqdef\ \frac{p}{3}\,\xi_{\,2}^{2}\ -\ \xi_{\,1}^{2}\,, \qquad m_{\,1}\ \eqdef\ \deg\,(\,P_{\,1}\,)\ =\ 2\,, \label{eq:pol1} \\
  P_{\,2}\,(\,p,\,q;\,\xi_{\,1},\,\xi_{\,2}\,)\ &\eqdef\ \frac{q}{2}\,\xi_{\,2}^{3}\ +\ \xi_{\,1}^{3}\,, \qquad m_{\,2}\ \eqdef\ \deg\,(\,P_{\,2}\,)\ =\ 3\,, \label{eq:pol2} \\
  P_{\,3}\,(\,p,\,q;\,\xi_{\,1},\,\xi_{\,2}\,)\ &\eqdef\ (\,4\,p^{\,3}\ -\ 27\,q^{\,2}\,)\,\xi_{\,1}^{2}\,, \qquad m_{\,3}\ \eqdef\ \deg\,(\,P_{\,3}\,)\ =\ 2\,, \label{eq:pol3}
\end{align}
\end{subequations}
whose respective degrees in variables $(\,\xi_{\,1},\,\xi_{\,2}\,)$ are denoted by $m_{\,k}\,$, $j\ \in\ 3^{\,\sqsupset}\,$. Then,
\begin{equation*}
  \bigcap_{j\,=\,1}^{3}\,P_{\,j}^{\,<}\,(\,\set{0}\,)\ =\ \Nn\,(\,\S\,)\ \bigcup\ \set*{(\,p,\,q,\,0,\,0\,)}\ \subseteq\ \T^{\,\ast}\,\C^{\,2}\,.
\end{equation*}
\end{lemma}

\begin{proof}
Due to the continuity of polynomial functions, it is enough to check the statement on the regular part of the co-normal $\Nn\,(\,\S\,)\,$. First, let us prove the inclusion $\subseteq\,$. So, consider a point
\begin{equation*}
  (\,p,\,q;\,\xi_{\,1},\,\xi_{\,2}\,)\ \in\ \bigcap_{j\,=\,1}^{3}\,P_{\,j}^{\,<}\,(\,\set{0}\,)\,.
\end{equation*}
Observe that if $\xi_{\,2}\ =\ 0\,$, then $P_{\,1}\,(\,p,\,q;\,\xi_{\,1},\,\xi_{\,2}\,)\ =\ 0$ implies that $\xi_{\,1}\ =\ 0\,$. Consequently, the point $(\,p,\,q;\,\xi_{\,1},\,\xi_{\,2}\,)$ belongs to the zero section. If we assume now that $4\,p^{\,3}\ -\ 27\,q^{\,2}\ \neq\ 0\,$, then $P_{\,3}\,(\,p,\,q;\,\xi_{\,1},\,\xi_{\,2}\,)\ =\ 0$ implies $\xi_{\,2}\ =\ 0\,$. We draw similarly the conclusion that the point $(\,p,\,q;\,\xi_{\,1},\,\xi_{\,2}\,)$ belongs to the zero section again. Now assume that $\xi_{\,2}\ \neq\ 0$ and $4\,p^{\,3}\ -\ 27\,q^{\,2}\ =\ 0\,$. Then, there exists $z\ \in\ \C$ such that $(\,p,\,q\,)\ \equiv\ (\,3\,z^{\,2},\,-2\,z^{\,3}\,)\,$. If $z\ =\ 0\,$, then from $P_{\,1}\,(\,p,\,q;\,\xi_{\,1},\,\xi_{\,2}\,)\ =\ 0$ we deduce that $\xi_{\,1}\ =\ 0$ and $(\,0,\,0;\,0,\,\xi_{\,2}\,)$ belongs to $\Nn\,(\,\S\,)$ as desired. If $z\ \neq\ 0\,$, then we can deduce from the equalities $P_{\,1,\,2}\,(\,p,\,q;\,\xi_{\,1},\,\xi_{\,2}\,)\ =\ 0\,$:
\begin{equation*}
  \xi_{\,1}\ =\ \frac{\xi_{\,1}^{3}}{\xi_{\,1}^{2}}\ =\ \frac{-\,\frac{q}{2}\;\xi_{\,2}^{3}}{\frac{p}{3}\;\xi_{\,2}^{2}}\ =\ z\,\xi_{\,2}\,.
\end{equation*}
This proves the inclusion. Now let us prove the reverse inclusion $\supseteq\,$. The zero section is clearly included in $\cap_{j\,=\,1}^{3}\,P_{\,j}^{\,<}\,(\,\set{0}\,)\,$. So, let us consider instead $(\,3\,z^{\,2},\,-\,2\,z^{\,3},\,z\,\lambda,\,\lambda\,)\ \in\ \Nn\,(\,\S\,)\,$. Then, an easy computation shows that
\begin{equation*}
  P_{\,j}\,\bigl(\,3\,z^{\,2},\,-\,2\,z^{\,3},\,z\,\lambda,\,\lambda\,\bigr)\ \equiv\ \fO\,, \qquad \forall\,j\ \in\ 3^{\,\sqsupset}\,.
\end{equation*}
This proves the Lemma.
\end{proof}

\begin{theorem}\label{thm:6.1}
There exists an open polydisc $D$ centered around $\vO\ \in\ \C^{\,2}$ such that $\forall\,k\ \in\ 3^{\,\sqsupset}$ and $\forall\,(\,i,\,j\,)\ \in\ 2^{\,\sqsupset}\times 2^{\,\sqsupset}$ we can find a differential operator $R_{\,i\,j}^{\,k}$ with holomorphic coefficients on $D$ such that:
\begin{enumerate}
  \item The operator $R_{\,i\,j}^{\,k}\,$, $\forall\,(\,i\ \neq\ j\,)$ is of the order $m_{\,k}\ -\ 1$ at most.
  \item The differential operator $R_{\,i\,i}^{\,k}$ is of the order $m_{\,k}$ and admits the polynomial $P_{\,k}\,(\,p,\,q;\,\xi_{\,1},\,\xi_{\,2}\,)$ as its principal symbol.
  \item The following identity holds:
  \begin{equation*}
    \begin{pmatrix}[1.2]
      R_{\,1\,1}^{\,k} & R_{\,1\,2}^{\,k} \\
      R_{\,2\,1}^{\,k} & R_{\,2\,2}^{\,k}
    \end{pmatrix}\cdot
    \begin{pmatrix}[1.2]
      z \\
      z^{\,2}
    \end{pmatrix}\ \equiv\ 
    \begin{pmatrix}[1.2]
      \fO \\
      \fO
    \end{pmatrix}\,.
  \end{equation*}
\end{enumerate}
\end{theorem}

\begin{proof}
This theorem can be proved along the same lines as \cite[Theorem~4.37]{Leichtnam1993}.
\end{proof}


\subsection{The holonomic $\D-$module construction}

The holonomic $\D-$module $\Mm$ will be constructed in the proof of the following
\begin{theorem}\label{thm:2.2}
The characteristic variety\footnote{It might be useful to say a few words about the characteristic variety of a $\D-$module in case a reader is not familiar with this notion. The annihilator of $\Mm$ is defined as
\begin{equation*}
  \Ann\,\Mm\ \eqdef\ \Set*{P\ \in\ \gr\,\D}{P\cdot m\ =\ \fO\,, \forall\,m\ \in\ \gr\,\Mm}\,.
\end{equation*}
It is an ideal in $\gr\,\D\,$. However, this ideal generally depends on the choice of a good filtration. So, one has to consider instead its \emph{radical}:
\begin{equation*}
  \sqrt{\Ann\,\Mm}\ \eqdef\ \Set*{P\ \in\ \gr\,\D}{P^{\,k}\ \in\ \Ann\,\Mm\ \text{for some}\ \exists\,k\ \in\ \N}\,.
\end{equation*}
It can be shown that the object $\sqrt{\Ann\,\Mm}$ does not depend on the filtration. The set of common zeros of $\sqrt{\Ann\,\Mm}$ is referred to as the characteristic variety $\Vv$ of $\Mm\,$:
\begin{equation*}
  \Vv\,\bigl(\,\sqrt{\Ann\,\Mm}\,\bigr)\ \eqdef\ \Set*{(\,\x_{\,0},\,\xib_{\,0}\,)}{P\,(\,\x_{\,0},\,\xib_{\,0}\,)\ \equiv\ 0\,,\ \forall\,P\ \in\ \sqrt{\Ann\,\Mm}}\,.
\end{equation*}} $\Vv$ of $\Mm$ is included in $\bigcap_{j\,=\,1}^{3}\,P_{\,j}^{\,<}\,(\,\set{0}\,)\ \equiv\ \Nn\,(\,\S\,)\ \bigcup\ \set*{(\,p,\,q,\,0,\,0\,)}$ and
\begin{enumerate}
  \item The $\D-$ module $\Mm$ is holonomic.
  \item The vector $(\,z,\,z^{\,2}\,)$ defines a holomorphic solution to $\Mm$ on every simply connected open subset of $\C^{\,2}\,\setminus\,\Delta^{\,<}\,(\,\set{0}\,)\,$.
\end{enumerate}
\end{theorem}

\begin{proof}
We denote by $\D$ (respectively by $\D^{\,(\,j\,)})$ the sheaf of holomorphic differential operators (of the order not greater than $j\,$, respectively) on the polydisc $D\,$. The grading of $\D$ is defined as
\begin{equation*}
  \gr\,\D\ \eqdef\ \bigoplus_{j\,=\,-1}^{+\,\infty}\,\frac{\D^{\,(\,j\,+\,1\,)}}{\D^{\,(\,j\,)}}\,, \qquad \D^{\,(\,-\,1\,)}\ \equiv\ \set*{0}\,.
\end{equation*}
Algebraically speaking, this object is isomorphic to the polynomial algebra $\O\,[\,\xi_{\,1},\,\xi_{\,2}\,]\,$.

We consider a $\D-$module morphism $\phi$ defined as:
\begin{align*}
  \phi\,:\ \D^{\,6}\ &\longrightarrow\ \D^{\,2} \\
  (\,Q_{\,\ell}^{\,k}\,)_{\substack{\,k\,\in\,3^{\,\sqsupset} \\ \ell\,\in\,2^{\,\sqsupset}}}\ &\mapsto\ \sum_{\substack{k\,\in\,3^{\,\sqsupset} \\ \ell\, \in\, 2^{\,\sqsupset}}} \bigl(\,Q_{\,\ell}^{\,k}\,R_{\,\ell,\,1}^{\,k},\,Q_{\,\ell}^{\,k}\,R_{\,\ell,\,2}^{\,k}\,\bigr)\,.
\end{align*}
We are particularly interested in the co-kernel of this morphism:
\begin{equation*}
  \Mm\ \eqdef\ \coker\,\phi\ \bydef\ \cod\,\phi\,/\,\phi_{\,>}\,(\,\set{\dom\,\phi}\,)\ \equiv\ \D^{\,2}\,/\,\phi_{\,>}\,(\,\set{\D^{\,6}}\,)
\end{equation*}
along with the canonical projection
\begin{equation*}
  \pi\,:\ \D^{\,2}\ \longrightarrow\ \Mm\,.
\end{equation*}
We stress out that $\Mm$ is a $\D-$module by construction. \cref{thm:6.1} allows us to define a $\D-$module sheaf homomorphism by the formula:
\begin{align*}
  \Mm\ &\longrightarrow\ \O\,[\,z\,] \\
  (\,\A,\,\B\,)\ &\mapsto\ \A\,z\ +\ \B\,z^{\,2}\,,
\end{align*}
where the right-hand side is a ramified solution existing on $D\,\setminus\,\Delta^{\,<}\,(\,\set{0}\,)$ defined by the pair $(\,z,\,z^{\,2}\,)\,$. Then part (2) of the Theorem follows from \cref{thm:4.2}. We introduce also a \emph{good}\footnote{We say that a filtration $\set*{\Mm_{\,j}}_{\,j\,\in\,\Z_{\,\geq\,-\,1}}$ on $\Mm$ is \emph{good} if $\Mm_{\,j}$ is coherent $\forall\,j\ \in\ \Z_{\,\geq\,-\,1}$ and there exists $\exists\,k_{\,0}\ \in\ \Z_{\,\geq\,-\,1}$ such that $\forall\,k\ \geq\ k_{\,0}\,$:
\begin{equation*}
  \D^{\,(\,k\,)}\cdot\Mm_{\,j}\ =\ \Mm_{\,k\,+\,j}\,.
\end{equation*}
We remind that, in general, one only has a weaker property:
\begin{equation*}
  \D^{\,(\,k\,)}\cdot\Mm_{\,j}\ \subseteq\ \Mm_{\,k\,+\,j}\,.
\end{equation*}} filtration\footnote{A filtration of $\Mm$ is an increasing sequence of sub-modules $\set*{\Mm_{\,j}}_{\,j\,\in\,\Z_{\,\geq\,-\,1}}$ of $\Mm$ verifying two properties:
\begin{itemize}
  \item $\Mm\ \equiv\ \bigcup_{j\,=\,-1}^{+\infty}\Mm_{\,j}\,$.
  \item $\D^{\,(\,k\,)}\cdot\Mm_{\,j}\ \subseteq\ \Mm_{\,j\,+\,k}\,$, $\forall\,(\,k,\,j\,)\ \in\ \N^{\,2}$\,.
\end{itemize}} $\set*{\Mm_{\,j}}_{\,j\,\in\,\Z_{\,\geq\,-\,1}}$ on $\Mm$ in the following way \cite{Malgrange1977}:
\begin{equation*}
  \Mm_{\,-\,1}\ \eqdef\ \set*{0}\,, \qquad \Mm_{\,j}\ \eqdef\ \pi\,(\,\D^{\,(\,j\,)},\,\D^{\,(\,j\,)}\,)\, \qquad \forall\,j\ \geq\ 0\,, \qquad \Mm\ \equiv\ \bigcup_{j\,=\,-1}^{+\infty}\Mm_{\,j}\,.
\end{equation*}
Consider $k\ \in\ 3^{\,\sqsupset}$ and $(\,\A,\,\B\,)\ \in\ \Mm_{\,j\,+\,1}$ with $j\ \geq\ -\,1\,$. The grading of $\Mm$ is
\begin{equation*}
  \gr\,\Mm\ \eqdef\ \bigoplus_{j\,=\,-1}^{+\,\infty}\,\frac{\Mm^{\,(\,j\,+\,1\,)}}{\Mm^{\,(\,j\,)}}\,.
\end{equation*}
The action of $\D$ on $\Mm$ induces the action of $\gr\,\D$ on $\gr\,\Mm$ in the obvious way. In $\Mm_{\,k\,+\,1\,+\,m_{\,k}}\,/\,\Mm_{\,k\,+\,m_{\,k}}$ we have:
\begin{equation*}
  P_{\,k}\cdot(\,\A,\,\B\,)\ =\ \A\cdot(\,R^{\,k}_{\,1,\,1},\,R_{\,1,\,2}^{\,k}\,)\ +\ \B\cdot(\,R^{\,k}_{\,2,\,1},\,R_{\,2,\,2}^{\,k}\,)\ =\ \fO\,.
\end{equation*}
Hence, the multiplication of $\Mm_{\,k\,+\,1}\,/\,\Mm_{\,k}$ by $P_{\,k}$ induces a zero mapping into $\Mm_{\,k\,+\,1\,+\,m_{\,k}}\,/\,\Mm_{\,k\,+\,m_{\,k}}\,$. 

The characteristic variety of $\Mm$ minus the zero section $\set*{(\,p,\,q,\,0,\,0\,)}$ is included into a \nm{Lagrangian} sub-variety of $\T^{\,\ast}\C^{\,2}\,\setminus\,\set*{(\,p,\,q,\,0,\,0\,)}\,$, so $\Mm$ is holonomic. This completes the proof of the Theorem.
\end{proof}
Of course, a couple of holomorphic ramified functions $(\,z,\,z^{\,2}\,)$ does not define a (uniform) distribution, but morally, \cref{thm:2.2} means that the ``micro-local singularities'' of $(\,z,\,z^{\,2}\,)$ live in the co-normal to the swallow tail $\Nn\,(\,\S\,)\,$.

We reiterate the important fact there exists only one direction\footnote{It would be more accurate to speak about a line, of course. In the real case the singular directions are $(\,0,\,\pm\,1\,)\,$, which lie along the same (complex) line.} in $\Nn\,(\,\S\,)$ above $(\,0,\,0\,)\ \in\ \C^{\,2}\,$, namely $(\,0,\,\lambda\,)\,$, $\lambda\ \in\ \C\,$. The principal symbol of $\partial_{\,q}$ (respectively, $\partial_{\,p}$) does not vanish (respectively, vanishes) on $\Set*{(\,0,\,0;\,0,\,\lambda\,)}{\lambda\ \in\ \C^{\,\times}}\,$. This explains why the primitive $\partial_{\,q}^{\,-\,1}$ exists on $\O\,\llbracket\,z\,\rrbracket$ whereas $\partial_{\,p}^{\,-\,1}$ does not exist as we saw above in \cref{sec:ring}. Lastly, \cref{thm:2.2} shows morally that the singularities of the product $z\cdot z$ do not spread over the origin and remain confined in $\Nn\,(\,\S\,)\,$.


\subsection{A computational digression}

Let $P_{\,i}$ and $P_{\,j}\,$ ($i,\,j\ \in\ 3^{\,\sqsupset}$) be any two polynomials from \cref{lem:p}. Theorem~\ref{thm:2.2} shows that the characteristic variety $\Vv$ is contained in
\begin{equation*}
  \Vv\ \subseteq\ P_{\,i}^{\,<}\,(\,\set{0}\,)\,, \qquad
  \Vv\ \subseteq\ P_{\,j}^{\,<}\,(\,\set{0}\,)\,.
\end{equation*}
Hence, $\T_{\,\x}\,\Vv\ \subseteq\ \ker\ud P_{\,i,\,j}\,(\,\x\,)\,$, where $\x\ \in\ \Vv$ is a smooth point of the characteristic variety $\Vv\,$. Then, if $\omega\ \eqdef\ \sum_{k\,=\,1}^{2}\,\ud x_{\,k}\,\wedge\,\ud \xi_{\,k}$ denotes the standard symplectic $2-$form on $\T_{\,\x}\,\Vv\,$, then $\forall\,h\ \in\ \T_{\,\x}\,\Vv\,$:
\begin{equation}\label{eq:orth}
  \ud P_{\,i,\,j}\,(\,\x\,)\cdot h\ \equiv\ \omega\,\bigl(\,\Hh_{\,P_{\,i,\,j}}\,(\,\x\,),\,h\,\bigr)\ =\ 0\,,
\end{equation}
where $\Hh_{\,P_{\,i,\,j}}$ is the following vector field:
\begin{equation*}
  \Hh_{\,P_{\,i,\,j}}\ \eqdef\ \sum_{k\,=\,1}^{2}\,\Bigl(\,\pd{P_{\,i,\,j}}{\xi_{\,k}}\,\pd{}{x_{\,k}}\ -\ \pd{P_{\,i,\,j}}{x_{\,k}}\,\pd{}{\xi_{\,k}}\,\Bigr)\,.
\end{equation*}
The second equality in \eqref{eq:orth} shows that $\Hh_{\,P_{\,i,\,j}}\,(\,\x\,)\ \in\ \bigl(\,\T_{\,\x}\,\Vv\,\bigr)^{\,\perp}\,$, where the orthogonal complement is taken in the sense of the natural symplectic structure on the (co-)tangent bundle $\T^{\,\ast}\,\Vv\,$. Due to the involutivity\footnote{We remind that if the characteristic variety $\Vv$ is \nm{Lagrangian}, then one has an even stronger property $\bigl(\,\T_{\,\x}\,\Vv\,\bigr)^{\,\perp}\ =\ \T_{\,\x}\,\Vv\,$.} property of the characteristic variety $\Vv$ \cite{Gabber1981}, we have that
\begin{equation*}
  \bigl(\,\T_{\,\x}\,\Vv\,\bigr)^{\,\perp}\ \subseteq\ \T_{\,\x}\,\Vv\,.
\end{equation*}
Henceforth, $\Hh_{\,P_{\,i,\,j}}\,(\,\x\,)\ \in\ \T_{\,\x}\,\Vv\,$. Consequently, we obtain the following result:
\begin{equation*}
  \omega\,\bigl(\,\Hh_{\,P_{\,i}}\,(\,\x\,),\,\Hh_{\,P_{\,j}}\,(\,\x\,)\,\bigr)\ \equiv\ \pb*{P_{\,i}}{P_{\,j}}\,(\,\x\,)\ =\ 0\,,
\end{equation*}
where the \nm{Poisson} bracket operator definition is implied by the symplectic structure:
\begin{equation*}
  \pb*{P}{Q}\ \eqdef\ \sum_{k\,=\,1}^{2}\,\Bigl(\,\pd{P}{\xi_{\,k}}\,\pd{Q}{x_{\,k}}\ -\ \pd{P}{x_{\,k}}\,\pd{Q}{\xi_{\,k}}\,\Bigr)\,.
\end{equation*}
Since $\x$ is a smooth point of $\Vv\,$, \cref{thm:2.2} allows to show that $\Vv$ and $\bigcap_{j\,=\,1}^{3}\,P_{\,j}^{\,<}\,(\,\set{0}\,)$ coincide in a small neighbourhood of $\x\,$. Then, by \nm{Hilbert} Nullstellensatz and \cref{thm:2.2}, one deduces that the \nm{Poisson} brackets $\pb*{P_{\,i}}{P_{\,j}}$ ($i,\,j\ \in\ 3^{\,\sqsupset}$) must belong at least to the radical ideal $\sqrt{\I}\,$, where $\I\ \eqdef\ \Span*{P_{\,1},\,P_{\,2},\,P_{\,3}}\,$. We are going to check this property using modern computer algebra methods which are deeply based on the \nm{Gr\"obner} bases methods \cite{Rosenkranz2007}.

First of all, let us initialize and load the library to work with polynomial ideals:
\begin{lstlisting}[language=Maple]
  restart:
  with(PolynomialIdeals):
\end{lstlisting}
Then, we define the polynomials from \cref{lem:p}:
\begin{lstlisting}[language=Maple]
  P[1] := (p, q, r, s) -> 1/3*p*s^2 - r^2;
  P[2] := (p, q, r, s) -> 1/2*q*s^3 + r^3;
  P[3] := (p, q, r, s) -> (4*p^3 - 27*q^2)*r^2;
\end{lstlisting}
along with the procedure which allows to compute the \nm{Poisson} bracket of two smooth functions:
\begin{lstlisting}[language=Maple]
  PBracket := proc(P,Q)
    description "Returns the Poisson bracket of two functions P and Q";
    return simplify(diff(P(p, q, r, s), r)*diff(Q(p, q, r, s), p) - diff(P(p, q, r, s), p)*diff(Q(p, q, r, s), r) + diff(P(p, q, r, s), s)*diff(Q(p, q, r, s), q) - diff(P(p, q, r, s), q)*diff(Q(p, q, r, s), s));
  end proc;
\end{lstlisting}
Then, we define the polynomial ideal $\I$:
\begin{lstlisting}[language=Maple]
  I := PolynomialIdeal(P[1](p, q, r, s), P[2](p, q, r, s), P[3](p, q, r, s));
\end{lstlisting}
and check whether the ideal $\I$ is radical:
\begin{lstlisting}[language=Maple]
  IsRadical(I);
\end{lstlisting}
The last function returns the logical value \texttt{false}. Henceforth, we construct its radical:
\begin{lstlisting}[language=Maple]
  RI := Radical(I);
\end{lstlisting}
We have to compute also the corresponding pairwise \nm{Poisson} brackets of polynomials $P_{\,1,\,2,\,3}\,$:
\begin{lstlisting}[language=Maple]
  P12 := (p, q, r, s) -> PBracket(P[1], P[2]);
  P13 := (p, q, r, s) -> PBracket(P[1], P[3]);
  P23 := (p, q, r, s) -> PBracket(P[2], P[3]);
\end{lstlisting}
Finally, we can check whether the \nm{Poisson} brackets belong to the radical ideal $\sqrt{\I}\,$:
\begin{lstlisting}[language=Maple]
  IdealMembership(P12(p, q, r, s), RI);
  IdealMembership(P13(p, q, r, s), RI);
  IdealMembership(P23(p, q, r, s), RI);
\end{lstlisting}
All these functions return the value \texttt{true}. This motivates us to establish an even stronger result comparing to \cref{thm:2.2}.

It is also interesting to note that $\pb*{P_{\,1}}{P_{\,2}}$ and $\pb*{P_{\,2}}{P_{\,3}}$ belong even to the ideal $\I\,$. However, it is not the case for $\pb*{P_{\,1}}{P_{\,3}}\,$. These facts can be simply checked as follows:
\begin{lstlisting}[language=Maple]
  IdealMembership(P12(p, q, r, s), I);
  IdealMembership(P23(p, q, r, s), I);
  IdealMembership(P13(p, q, r, s), I);
\end{lstlisting}
The first two calls return \texttt{true} while the last one returns the value \texttt{false}. However, it is not difficult to check that $\pb*{P_{\,1}}{P_{\,3}}^{\,2}$ already belongs to $\I\,$:
\begin{lstlisting}[language=Maple]
  IdealMembership(P13(p, q, r, s)^2, I);
\end{lstlisting}
which returns the value \texttt{true}.


\subsection{A stronger result}

In fact, we were able to prove a stronger result compared to \cref{thm:2.2}. One actually has the equality in \cref{thm:2.2}:
\begin{theorem}\label{thm:eq}
Under conditions of \cref{thm:2.2}, one has
\begin{equation*}
  \Vv\,\setminus\,\set*{(\,p,\,q,\,0,\,0\,)}\ =\ \Nn\,(\,\S\,)\,\setminus\,\set*{(\,p,\,q,\,0,\,0\,)}\,,
\end{equation*}
where $\set*{(\,p,\,q,\,0,\,0\,)}$ refers to the zero section.
\end{theorem}

\begin{proof}
In \cref{thm:2.2} we already established the inclusion $\Vv\,\setminus\,\set*{(\,p,\,q,\,0,\,0\,)}\ \subseteq\ \Nn\,(\,\S\,)\,\setminus\,\set*{(\,p,\,q,\,0,\,0\,)}$ and from it we deduce that $\Mm$ is holonomic. Let
\begin{align*}
  \pi\,:\ \T^{\,\ast}\C^{\,2}\ &\longrightarrow\ \C^{\,2} \\
  (\,\x,\,\xib\,)\ &\mapsto\ \x
\end{align*}
denote the standard projection. We shall prove our statement to be a contradiction. So, let us assume that $\Nn\,(\,\S\,)\,\setminus\,\set*{(\,p,\,q,\,0,\,0\,)}$ is not included in $\Vv\,\setminus\,\set*{(\,p,\,q,\,0,\,0\,)}$ and let us deduce a contradiction. Consider a point
\begin{equation*}
  \bigl(\,p_{\,0},\,q_{\,0};\,\xi_{\,1}^{\,0},\,\xi_{\,2}^{\,0}\,\bigr)\ \in\ \Nn\,(\,\S\,)\,\setminus\,\bigl(\,\Vv\,\cup\,\set*{(\,p,\,q,\,0,\,0\,)}\,\bigr)\,.
\end{equation*}
By homogeneity, for any $\mu\ \in\ \C^{\,\times}\,$, the point $\bigl(\,p_{\,0},\,q_{\,0};\,\mu\,\xi_{\,1}^{\,0},\,\mu\,\xi_{\,2}^{\,0}\,\bigr)$ does not belong to $\Vv\,$. Since $\pi^{\,-\,1}\,\set*{(\,p_{\,0},\,q_{\,0}\,)}\,\cap\,\Nn\,(\,\S\,)$ is a complex line and since $\Vv\,\setminus\,\set*{(\,p,\,q,\,0,\,0\,)}\ \subseteq\ \Nn\,(\,\S\,)\,\setminus\,\set*{(\,p,\,q,\,0,\,0\,)}\,$, we immediately deduce that $(\,p_{\,0},\,q_{\,0}\,)\ \not\in\ \pi\,\bigl(\,\Vv\,\setminus\,\set*{(\,p,\,q,\,0,\,0\,)}\,\bigr)\,$. Now recall that the characteristic variety minus the zero section is (always) a homogeneous closed subset of $\T^{\,\ast}\,\C^{\,2}\,\setminus\,\set*{(\,p,\,q,\,0,\,0\,)}\,$. Therefore, there exists an open neighbourhood $\U^{\,(\,0\,)}$ of $(\,p_{\,0},\,q_{\,0}\,)\ \in\ \C^{\,2}$ such that
\begin{equation*}
  \pi\,\bigl(\,\Vv\,\setminus\,\set*{(\,p,\,q,\,0,\,0\,)}\,\bigr)\,\bigcap\,\U^{\,(\,0\,)}\ =\ \emptyset\,.
\end{equation*}

Consider now a \nm{Whitney} stratification of $\C^{\,2}\,$:
\begin{equation*}
  \C^{\,2}\ =\ \bigcup_{\lambda\,\in\,\Lambda}\,\Xx_{\,\lambda}
\end{equation*}
such that $\Vv\ \subseteq\ \bigcup_{\lambda\,\in\,\Lambda}\,\Nc^{\,\ast}\Xx_{\,\lambda}\,$, where $\Xx_{\,\lambda}$ are smooth (open) sub-manifolds and $\Nc^{\,\ast}\Xx_{\,\lambda}$ denotes the co-normal bundle to $\Xx_{\,\lambda}\,$. Recall also that $\set*{\Xx_{\,\lambda}}_{\,\lambda\,\in\,\Lambda}$ form a partition of $\C^{\,2}\,$. Since $\pi\,\bigl(\,\Vv\,\setminus\,\set*{(\,p,\,q,\,0,\,0\,)}\,\bigr)\,\cap\,\U^{\,(\,0\,)}\ =\ \emptyset\,$, we can assume that $(\,p,\,q,\,0,\,0\,)\ \in\ \Xx_{\,\lambda_{\,0}}\,$, where $\Xx_{\,\lambda_{\,0}}$ is an open subset of $\C^{\,2}$ so that its co-normal is the zero section. Then, consider an open simply connected neighbourhood $\U^{\,(\,\lambda_{\,0}\,)}$ of $(\,p_{\,0},\,q_{\,0}\,)$ which is included in $\Xx_{\,\lambda_{\,0}}\,$. Consider also an open simply connected subset $\Wu^{\,(\,\lambda_{\,0}\,)}$ included in $\U^{\,(\,\lambda_{\,0}\,)}\,\setminus\,\Delta^{\,<}\,(\,\set{0}\,)\,$. Now we can apply to $\Mm$ the first theorem of \nm{Kashiwara} \cite[Theorem~3.1]{Kashiwara1975} (see also a brief reminder below). It states, in particular, that the restriction of
\begin{equation*}
  \mathsf{Ext\,}^{\,0}_{\,\D}\,(\,\Mm,\,\O\,)\ \equiv\ \mathsf{Hom}_{\,\D}\,(\,\Mm,\,\O\,)
\end{equation*}
to $\Xx_{\,\lambda_{\,0}}$ is a local system. Since by \cref{thm:2.2}, the vector $(\,z,\,z^{\,2}\,)$ is a holomorphic solution of $\Mm$ on $\Wu^{\,(\,\lambda_{\,0}\,)}\,$, \nm{Kashiwara}'s theorem implies that one can analytically continue $(\,z,\,z^{\,2}\,)$ to $\Wu^{\,(\,\lambda_{\,0}\,)}\,$. But this is impossible because if we make the analytic continuation of $\partial_{\,q}\,z\ =\ \dfrac{1}{3\,z^{\,2}\ -\ p}$ along a path which ends up at $(\,p_{\,0},\,q_{\,0}\,)\,$, it will blow up. This contradiction, therefore, shows that the reverse inclusion is satisfied, and this Theorem is proved.
\end{proof}

For the sake of the exposition completeness, we formulate the following fundamental result of M.~\nm{Kashiwara} on which our \cref{thm:eq} is based:
\begin{theorem}[\cite{Kashiwara1975}]
Let $\mathfrak{M}$ by a maximally overdetermined system on a complex manifold $\Xx$ and $\Xx\ =\ \bigcup\,\Xx_{\,\alpha}$ be a stratification of $\Xx$ satisfying the regularity conditions of H.~\nm{Whitney} such that the singular support of $\mathfrak{M}$ is contained in the union of co-normal projective bundles of the strata. Then, the restriction of $\mathsf{Ext\,}^{\,i}_{\,\D_{\,\Xx}}\,(\,\mathfrak{M},\,\O_{\,\Xx}\,)$ to each stratum is a locally constant sheaf of finite rank.
\end{theorem}
In the Theorem above the following notations were used:
\begin{description}
  \item[$\D_{\,\Xx}\ $] The sheaf of differential operators of \emph{finite order} on the complex manifold $\Xx\,$.
  \item[$\O_{\,\Xx}\ $] The sheaf of holomorphic functions on $\Xx\,$, which is a left coherent $\D_{\,\Xx}-$module.
  \item[$\mathfrak{M}\ $] The coherent $\D_{\,\Xx}-$module.
  \item[$\mathsf{Ext}\ $] The right derived functor of the hom-functor.
\end{description}


\section{Second order PDE}
\label{sec:sec}

In this Section, we study the next representative of the \acrshort{pde} family \eqref{eq:fam} with $m\ =\ 2\,$. Typically, we would like to develop the methods which would allow us to study the following ramified \acrshort{cp}\footnote{We underline the fact that the term \acrfull{cp} is understood here in the \emph{incomplete} (or loose) sense specified in the Introduction \cref{sec:intro} since our \nm{Cauchy} datum contains only one initial condition for a second order problem in time \eqref{eq:pde2}.} for a genuinely nonlinear second order \acrshort{pde}:
\begin{subequations}\label{eq:second}
 \begin{align}\label{eq:pde2}
   u_{\,t\,t}\ -\ u_{\,x}\,u_{\,x\,x}\ &=\ \fO\,, \\
   u\,(\,0,\,x\,)\ &=\ \sum_{j\,=\,1}^{N_{\,0}}\,c_{\,j}\,x^{\,1\,+\,\frac{j\,-\,1}{3}}\,,\label{eq:ic3}
 \end{align}
\end{subequations}
with $N_{\,0}\ \in\ \N_{\,>\,2}\,$, $c_{\,1,\,2}\ \in\ \C^{\,\times}$ and $\set*{c_{\,j}}_{\,j\,=\,3}^{\,N_{\,0}}\ \subseteq\ \C$ some constants. Namely, we would like to provide sufficient evidence that the following conjecture holds. This is the main goal of this Section:

\begin{conj}\label{conj:2}
Consider the ramified \acrshort{ivp} \eqref{eq:second}. Then, for each choice of a root to the algebraic equation\footnote{We have to assume that $u_{\,x}\,(\,0,\,0\,)\ \neq\ 0$ in order to have two distinct roots.} $\tau^{\,2}\ -\ u_{\,x}\,(\,0,\,0\,)\cdot 1^{2}\ =\ \fO$ (where we put $\xi\ \leftarrow\ 1\ \neq\ 0$ by homogeneity of the principal symbol), there exist holomorphic functions $p\,(\,t,\,x\,)\,$, $q\,(\,t,\,x\,)$ and $\set*{a_{\,j}\,(\,t,\,x\,)}_{\,j\,=\,0}^{\,2}$ defined in a neighbourhood of $(\,0,\,0\,)\ \in\ \C^{\,2}$ with
\begin{equation*}
  p\,(\,0,\,x\,)\ =\ \fO\,, \qquad
  q\,(\,0,\,x\,)\ =\ x
\end{equation*}
such that the following assertion holds: the problem \eqref{eq:second} admits a local solution of the form:
\begin{equation}\label{eq:7.N}
  u\,(\,t,\,x\,)\ =\ \partial_{\,q}^{\,-1}\,\biggl(\,\sum_{j\,=\,0}^{\,2}\,a_{\,j}\,(\,t,\,x\,)\,z^{\,j}\,\biggr)\,,
\end{equation}
where $z\,(\,t,\,x\,)$ satisfies the following algebraic equation:
\begin{equation*}
  z^{\,3}\ =\ p\,(\,t,\,x\,)\,z\ +\ q\,(\,t,\,x\,)\,.
\end{equation*}
\end{conj}

\begin{proof}
Open problem.
\end{proof}

\begin{remark}
Please, notice that the second order ramified \acrshort{ivp} \eqref{eq:second} contains only one \nm{Cauchy} data. To compensate this fact, we impose a special form (\ie ansatz) for the solution that we are seeking. The condition $c_{\,1}\ \neq\ 0$ is necessary in order to ensure that the equation $\tau^{\,2}\ -\ u_{\,x}\,(\,0,\,0\,)\cdot 1^{2}\ =\ \fO$ has two distinct roots. The condition $c_{\,2}\ \neq\ 0$ means that the \nm{Cauchy} datum is sufficiently singular, it is crucial to ensure that the matrix $\Ma$ defined in \cref{eq:matrix} in the eikonal \cref{eq:eiko} is invertible. Actually, the Theorem~\ref{thm:2.1} from \cref{sec:wss} shows that if the constant $c_{\,2}$ is equal to zero, then the solution will be ramified around a smooth hyper-surface.

The choice of a root in the algebraic equation will allow fixing the value of $q_{\,t}\,(\,0,\,0\,)\,$, then it is very likely that the solution under the form \eqref{eq:7.N} will be unique. Notice that the functions $p\,(\,t,\,x\,)$ and $q\,(\,t,\,x\,)$ will depend on the \nm{Cauchy} datum.
\end{remark}

Despite the fact that we do not have a proof for this conjecture, we shall present in this Section strong evidence that it is true. Our symbolic--numeric computations will be performed in the case of $c_{\,1}\ =\ 1$ and $c_{\,2}\ =\ \frac{3}{4}$ for the sake of convenience.

In the sequel we shall rather write the solution ansatz without applying the \nm{Weierstra}\ss{} division theorem \cite{Narasimhan1966}, \ie
\begin{equation*}
  \sum_{j\,=\,0}^{\,2}\,a_{\,j}\,(\,t,\,x\,)\,z^{\,j}\ \equiv\ \sum_{k\,=\,0}^{\,\infty}\,b_{\,k}\,(\,t,\,x\,)\,z^{\,k}\,.
\end{equation*}
Despite the apparent complication, in this way, it will be easier to express the iterative process.


\subsection{Construction of iterations}
\label{sec:it}

Motivated by \cref{lem:unique} and \cref{thm:eq}, we extract from the \acrshort{pde} \eqref{eq:pde2} the following characteristic equation:
\begin{equation*}
  q_{\,t}\,(\,0,\,0\,)^{\,2}\ -\ u_{\,x}\,(\,0,\,0\,)\,q_{\,x}\,(\,0,\,0\,)\ =\ 0\,.
\end{equation*}
Observe that $q_{\,x}\,(\,0,\,0\,)\ =\ 1$ and that the choice $c_{\,1}\ =\ 1$ necessarily implies $u_{\,x}\,(\,0,\,0\,)\ =\ 1\,$. Then, in the previous equation, we shall choose the root $q_{\,t}\,(\,0,\,0\,)\ =\ 1$ for the rest of this Section.

The approach we adopt in this Section is greatly inspired by \cref{sec:ring,sec:dmod}. We would like to construct a solution to \cref{eq:pde2} in the following form, motivated by \cref{lem:prim}:
\begin{equation}\label{eq:ser}
  u\,(\,t,\,x\,)\ =\ \sum_{k\,=\,1}^{\infty}\,\frac{1}{k}\;\Bigl(\,-\,p\,(\,t,\,x\,)\,b_{\,k-1}\,(\,t,\,x\,)\ +\ 3\,b_{\,k-3}\,(\,t,\,x\,)\,\Bigr)\,z^{\,k}\ \in\ \O\,\llbracket\,z\,\rrbracket\,,
\end{equation}
where $p\,(\,t,\,x\,)\,$, $q\,(\,t,\,x\,)$ and $\set*{b_{\,k}\,(\,t,\,x\,)}_{\,k\,=\,0}^{\,\infty}$ are holomorphic functions of its two variables $(\,t,\,x\,)\ \in\ \C^{\,2}$ defined in the vicinity $ \U_{\,(\,0,\,0\,)}$ of $(\,0,\,0\,)\ \in\ \C^{\,2}$ that verify the following `initial' conditions:
\begin{subequations}\label{eq:iconds}
\begin{equation}\label{eq:pconds}
  p\,(\,0,\,x\,)\ \equiv\ \fO\,, \qquad
  q\,(\,0,\,x\,)\ \equiv\ x\,, \qquad q_{\,t}\,(\,0,\,0\,)\ =\ 1\,,
\end{equation}
\begin{equation}\label{eq:bcond}
  b_{\,0}\,(\,0,\,x\,)\ \equiv\ \fW\,, \qquad
  b_{\,1}\,(\,0,\,x\,)\ \equiv\ \fW\,, \qquad
  b_{\,k}\,(\,0,\,x\,)\ \equiv\ \fO\,,\quad \forall\,k\ \geq\ 2\,,
\end{equation}
\end{subequations}
and
\begin{equation*}
  z^{\,3}\ =\ p\,z\ +\ q\,, \qquad \forall\,(\,t,\,x\,)\ \in\ \U_{\,(\,0,\,0\,)}\ \subseteq\ \C^{\,2}\,.
\end{equation*}
The `initial' conditions \eqref{eq:iconds} correspond to the \nm{Cauchy} datum \eqref{eq:ic3} with $c_{\,1}\ =\ 1$ and $c_{\,2}\ =\ \dfrac{3}{4}\,$:
\begin{equation}\label{eq:init1}
  u\,(\,0,\,x\,)\ =\ x\ +\ \frac{3}{4}\;x^{\;\fourthirds}\,.
\end{equation}
Several other initial conditions will be considered below in \cref{sec:num}.

It seems impossible to construct directly such a solution $u\,(\,t,\,x\,)$ having a property that the series \eqref{eq:ser} is convergent $\forall\,(\,t,\,x\,)\ \in\ \U_{\,(\,0,\,0\,)}\,$. Instead, we propose to construct this solution $u\,(\,t,\,x\,)$ using the fixed point iterative scheme.

Using the ansatz \eqref{eq:ser} and this scheme, we shall see that the solution is completely determined if we have the following collection of holomorphic functions:
\begin{equation}\label{eq:coll}
  \Cc\ \eqdef\ \Bigl(\,p\,(\,t,\,x\,),\,q\,(\,t,\,x\,),\,\set*{b_{\,k}\,(\,t,\,x\,)}_{\,k\,=\,0}^{\,\infty}\,\Bigr)\ \bigcap\ \text{\eqref{eq:iconds}}\,.
\end{equation}
By construction, the functions from this collection $\Cc$ verify the `initial conditions' \eqref{eq:iconds}. Thus, if we know these functions, we can easily recover the solution $u\,(\,t,\,x\,)$ thanks to \eqref{eq:ser}. Hence, the collection~\eqref{eq:coll} will be called the solution data. 

The first naive tentative to solve \cref{conj:2} consists in replacing $u\,(\,t,\,x\,)$ by \eqref{eq:ser} in the nonlinear differential operator $u\ \mapsto\ u_{\,t\,t}\ -\ u_{\,x}\,u_{\,x\,x}\,$. One obtains a certain expression $(\,\maltese\,)\,$. Unfortunately, it seems to be unrealistic to construct the solutions directly in this way. Henceforth, we introduce an iterative procedure where we replace in one specific place of the expression $(\,\maltese\,)$ the function $p\,(\,t,\,x\,)$ by $\ph\,(\,t,\,x\,)\,$, respectively, $q\,(\,t,\,x\,)$ by $\qh\,(\,t,\,x\,)$ and $\set*{b_{\,k}\,(\,t,\,x\,)}_{\,k\,=\,0}^{\,\infty}$ by $\set*{\bh_{\,k}\,(\,t,\,x\,)}_{\,k\,=\,0}^{\,\infty}\,$. This procedure is better explained in Appendix~\ref{app:0}. As a result, we obtain the new expression:
\begin{equation*}
  (\,\widehat{\maltese}\,)\ =\ \sum_{k\,=\,0}^{\infty}\,\alpha_{\,k}\,z^{\,k}\ +\ \frac{\sum_{k\,=\,0}^{\infty}\,\beta_{\,k}\,z^{\,k}}{3\,z^{\,2}\ -\ p}\,.
\end{equation*}
Then, by requiring that $3\,z^{\,2}\ -\ p$ divides $\sum_{k\,=\,0}^{\infty}\,\beta_{\,k}\,z^{\,k}$ we obtain the eikonal \cref{eq:eiko} which gives us the new components $(\,\ph,\,\qh\,)\,$. Roughly speaking, \eqref{eq:eiko} means that the cusp is characteristic for the \acrshort{pde}. By requiring that $(\,\widehat{\maltese}\,)$ vanishes, we obtain the following discrete `time' (but infinite dimensional) dynamical system:
\begin{equation*}
  \Bigl(\,\ph\,(\,t,\,x\,),\,\qh\,(\,t,\,x\,),\,\set{\bh_{\,k}\,(\,t,\,x\,)}_{\,k\,=\,0}^{\,\infty}\,\Bigr)\ \coloneq\ \FF\,\Bigl(\,p\,(\,t,\,x\,),\,q\,(\,t,\,x\,),\,\set*{b_{\,k}\,(\,t,\,x\,)}_{\,k\,=\,0}^{\,\infty}\,\Bigr)\,.
\end{equation*}
Henceforth, we may say that the ideas of geometric optics inspire the construction of iterations. A similar procedure was employed by \nm{Wagschal} in \cite{Wagschal1974}.

Then, by the construction of the mapping $\FF\,$, its fixed point gives us the desired solution to \cref{conj:2}. In this Section we shall construct the mapping $\FF$ and implement it in a computer algebra system \texttt{Maple}${}^{\text{\texttrademark}}$:
\begin{align}\label{eq:mapF}
  \FF\,:\ \Cc\ &\longrightarrow\ \Cc\,, \\
  \Bigl(\,p\,(\,t,\,x\,),\,q\,(\,t,\,x\,),\,\set*{b_{\,k}\,(\,t,\,x\,)}_{\,k\,=\,0}^{\,\infty}\,\Bigr)\ &\mapsto\ \Bigl(\,\ph\,(\,t,\,x\,),\,\qh\,(\,t,\,x\,),\,\set{\bh_{\,k}\,(\,t,\,x\,)}_{\,k\,=\,0}^{\,\infty}\,\Bigr)\,. \nonumber
\end{align}
As we already mentioned, the mapping $\FF$ is constructed such that if we have a fixed point
\begin{equation*}
  \FF\,\Bigl(\,p\,(\,t,\,x\,),\,q\,(\,t,\,x\,),\,\set*{b_{\,k}\,(\,t,\,x\,)}_{\,k\,=\,0}^{\,\infty}\,\Bigr)\ \equiv\ \Bigl(\,p\,(\,t,\,x\,),\,q\,(\,t,\,x\,),\,\set*{b_{\,k}\,(\,t,\,x\,)}_{\,k\,=\,0}^{\,\infty}\,\Bigr)\,,
\end{equation*}
then the solution \eqref{eq:ser} reconstructed from this data verifies \eqref{eq:pde2} (\cf \cref{eq:prop}). Additionally, since we impose Conditions~\eqref{eq:bcond}, this mapping will yield that $\bigl(\,q_{\,t}\,(\,0,\,0\,)\,\bigr)^{\,2}\ =\ 1\,$. The fact that we have chosen $q_{\,t}\,(\,0,\,0\,)\ =\ 1$ corresponds to the choice of a root to equation $\tau^{\,2}\ -\ \partial_{\,x}\,u\,(\,0,\,0\,)\cdot 1^{2}\ =\ \fO$ in \cref{conj:2} (where we replaced $\xi\ \leftarrow\ 1\ \neq\ 0$ by homogeneity of the principal symbol).

We define the following family of functions:
\begin{multline*}
  \B_{\,k}\ \eqdef\ \sum_{j\,=\,0}^{k}\,\biggl[\,\Bigl(\,1\, -\, \frac{1}{j}\,\Bigr)\,p_{\,x}\,b_{\,j\,-\,1}\ +\ q_{\,x}\,b_{\,j}\ +\ \frac{1}{j}\;\Bigl(\,-\,p\,b_{\,j\,-\,1,\,x}\ +\ 3\,b_{\,j\,-\,3,\,x}\,\Bigr)\,\biggr]\times\\
  \biggl[\,(\,k\, -\, j\, -\, 1\,)\,p_{x}^{\,2}\,b_{\,k\,-\,j\,-\,1}\ +\ 2\,(\,k\, -\, j\,)\,p_{\,x}\,q_{\,x}\,b_{\,k\,-\,j}\ +\ (\,k\, -\, j\, +\, 1\,)\,q_{x}^{\,2}\,b_{\,k\,-\,j\,+\,1}\,\biggr]\\ 
  -\ (\,k\, -\, 1\,)\cdot p_{t}^{\,2}\,b_{\,k\,-\,1}\,, \qquad \forall\,k\ \in\ \N\,,
\end{multline*}
where we additionally use the convention:
\begin{equation}\label{eq:conv}
  b_{\,-j}\ \equiv\ \fO\,, \forall\,j\ \geq\ 1\,.
\end{equation}
The \emph{eikonal equation}\footnote{Strictly speaking, the term `\emph{eikonal equation}' is not defined for genuinely nonlinear \acrshort{pde}s and singular hyper-surfaces. However, we have good reasons to think that what we are doing is the right generalization of this notion to our case.} is defined as
\begin{equation}\label{eq:eiko}
  \begin{pmatrix}
    \ph_{\,t} \\
    \addlinespace
    \qh_{\,t}
  \end{pmatrix}\ =\ \Ma^{\,-\,1}\cdot
  \begin{pmatrix}
    \sum_{k\,=\,0}^{+\,\infty}\; \Bigl(\,\dfrac{p}{3}\,\Bigr)^{\,k}\,\B_{\,2\,k} \\
    \addlinespace
    \sum_{k\,=\,0}^{+\,\infty}\; \Bigl(\,\dfrac{p}{3}\,\Bigr)^{\,k}\,\B_{\,2\,k\,+\,1}
  \end{pmatrix}\,.
\end{equation}
Since functions $\ph$ and $\qh$ belong to the class $\Cc\,$, they verify the following initial conditions:
\begin{equation*}
  \ph\,(\,0,\,x\,)\ =\ \fO\,, \qquad
  \qh\,(\,0,\,x\,)\ =\ x\,.
\end{equation*}
The functional $2\times2$ matrix $\Ma$ is defined as
\begin{equation}\label{eq:matrix}
  \Ma\ \eqdef\ \begin{pmatrix}
    4\,\sum_{k\,=\,0}^{+\,\infty}\; \Bigl(\,\dfrac{p}{3}\,\Bigr)^{\,k}\,k\,q_{\,t}\,b_{\,2\,k} & \sum_{k\,=\,0}^{+\,\infty}\; \Bigl(\,\dfrac{p}{3}\,\Bigr)^{\,k}\,(\,2\,k\,+\,1\,)\,q_{\,t}\,b_{\,2\,k\,+\,1} \\
    \addlinespace
    2\,\sum_{k\,=\,0}^{+\,\infty}\; \Bigl(\,\dfrac{p}{3}\,\Bigr)^{\,k}\,(\,2\,k\,+\,1)\,q_{\,t}\,b_{\,2\,k\,+\,1} & 2\,\sum_{k\,=\,0}^{+\,\infty}\; \Bigl(\,\dfrac{p}{3}\,\Bigr)^{\,k}\,(\,k\,+\,1\,)\,q_{\,t}\,b_{\,2\,k\,+\,2}
  \end{pmatrix}\,.
\end{equation}
We underline that the inversion of the last matrix is necessary to obtain an explicit form of the eikonal System~\eqref{eq:eiko} solved with respect to the time derivatives. It is not difficult to see (by substituting the necessary initial data and conventions) that the first Equation in \eqref{eq:eiko} yields $\ph_{\,t}\,(\,0,\,0\,)\ =\ \half\,$.

Let us define also the following family of functions:
\begin{multline*}
  \A_{\,k}\ \eqdef\ (\,k\,-\,1\,)\,b_{\,k\,-\,1}\,p_{t}^{\,2}\ +\ 2\,k\,b_{\,k}\,\ph_{\,t}\,q_{\,t}\ +\ (\,k\,+\,1\,)\;\qh_{\,t}\,q_{\,t}\,b_{\,k\,+\,1}\ -\\
  \sum_{j\,=\,0}^{k}\,\biggl[\,\Bigl(\,1\,-\,\frac{1}{j}\,\Bigr)\,b_{\,j\,-\,1}\,p_{\,x}\ +\ b_{\,j}\,q_{\,x}\ +\ \frac{1}{j}\;\bigl(\,-\,p\,b_{\,j\,-\,1,\,x}\ +\ 3\,b_{\,j\,-\,3,\,x}\,\bigr)\,\biggr]\times\\
  \biggl[\,(\,k\,-\,j\,-\,1\,)\,b_{\,k\,-\,j\,-\,1}\,p_{x}^{\,2}\ +\ 2\,(\,k\,-\,j\,)\,b_{\,k\,-\,j}\,p_{\,x}\,q_{\,x}\ +\ (\,k\,-\,j\,+\,1\,)\,b_{\,k\,-\,j\,+\,1}\,q_{x}^{\,2}\,\biggr]\,, \quad \forall\,k\ \in\ \N\,.
\end{multline*}
Also, for any $k\ \in\ \N$ we define:
\begin{equation*}
  \Cr_{\,k}\ \eqdef\ \sum_{j\,=\,0}^{\infty}\,\A_{\,k\,+\,2\,+\,2\,j}\,\Bigl(\,\frac{p}{3}\,\Bigr)^{\,j}\,.
\end{equation*}
Now, we are ready to write down equations which define the family of functions $\set*{\bh_{\,k}\,(\,t,\,x\,)}_{\,k\,=\,0}^{\,\infty}\,$:
\begin{multline}\label{eq:b}
  \pd{\bh_{\,k}}{t}\ =\ \frac{1}{2\,q_{\,t}}\;\Biggl\{\,-\,\Bigl(\,2\,-\,\frac{2}{k}\,\Bigr)\,b_{\,k\,-\,1,\,t}\,p_{\,t}\ -\ \Bigl(\,1\,-\,\frac{1}{k}\,\Bigr)\,b_{\,k\,-\,1}\,p_{\,t\,t}\ -\ b_{\,k}\,q_{\,t\,t}\\ 
  -\ \frac{1}{k}\;\bigl(\,-\,p\,b_{\,k\,-\,1,\,t\,t}\ +\ 3\,b_{\,k\,-\,3,\,t\,t}\,\bigr)\ -\ \Cr_{\,k}\ +\\ 
  \sum_{j\,=\,0}^{k}\,\biggl[\,\Bigl(\,1\,-\,\frac{1}{j}\,\Bigr)\,b_{\,j\,-\,1}\,p_{\,x}\ +\ b_{\,j}\,q_{\,x}\ +\ \frac{1}{j}\;\bigl(\,-\,p\,b_{\,j\,-\,1,\,x}\ +\ 3\,b_{\,j\,-\,3,\,x}\,\bigr)\,\biggr]\times\\
  \biggl[\,\Bigl(\,2\,-\,\frac{2}{k\,-\,j}\,\Bigr)\,b_{k\,-\,j\,-\,1,\,x}\,p_{\,x}\ +\ \Bigl(\,1\,-\,\frac{1}{k\,-\,j}\,\Bigr)\,b_{\,k\,-\,j\,-\,1}\,p_{\,t\,t}\ +\ b_{\,k\,-\,j}\,q_{\,x\,x}\ +\\ 
  2\,b_{\,k\,-\,j,\,x}\,q_{\,x}\ +\ \frac{1}{k\,-\,j}\;\bigl(\,-\,p\,b_{\,k\,-\,j\,-\,1,\,x\,x}\ +\ 3\,b_{\,k\,-\,j\,-\,3,\,x\,x}\,\bigr)\,\biggr]\,\Biggr\}\,, \qquad \forall\,k\ \in\ \N\,.
\end{multline}
We underline that functions $\set*{\bh_{\,k}\,(\,t,\,x\,)}_{\,k\,=\,0}^{\,\infty}$ satisfy the initial conditions \eqref{eq:bcond}. The convention \eqref{eq:conv} must be employed to interpret correctly the last equation.

The solution of differential Equations~\eqref{eq:eiko} and \eqref{eq:b} completes the construction of the fixed point mapping $\FF$ \eqref{eq:mapF}. Then, the connection between the mapping $\FF$ and the ramified \nm{Cauchy} problem \eqref{eq:second} is elucidated in the following
\begin{lemma}\label{eq:prop}
Let us assume that we have fixed point data:
\begin{equation*}
  \FF\,\Bigl(\,p\,(\,t,\,x\,),\,q\,(\,t,\,x\,),\,\set*{b_{\,k}\,(\,t,\,x\,)}_{\,k\,=\,0}^{\,\infty}\,\Bigr)\ =\ \Bigl(\,p\,(\,t,\,x\,),\,q\,(\,t,\,x\,),\,\set*{b_{\,k}\,(\,t,\,x\,)}_{\,k\,=\,0}^{\,\infty}\,\Bigr)
\end{equation*}
satisfying the `initial' conditions \eqref{eq:iconds}. Then, the function
\begin{equation}\label{eq:infsol}
  u\,(\,t,\,x\,)\ =\ \partial_{\,q}^{\,-1}\,\biggl(\,\sum_{k\,=\,0}^{+\,\infty}\,b_{\,k}\,(\,t,\,x\,)\,z^{\,k}\,\biggr)\,,
\end{equation}
with $z^{\,3}\ =\ p\,z\ +\ q$ satisfies \acrshort{pde} \eqref{eq:pde2} together with the initial condition \eqref{eq:ic3}.
\end{lemma}

\begin{proof}
By direct computations, we have the following expressions for ansatz \eqref{eq:infsol} derivatives:
\begin{equation*}
  u_{\,\alpha}\ =\ \sum_{k\,=\,0}^{+\,\infty}\,\biggl[\,\Bigl(\,1\,-\,\frac{1}{k}\,\Bigr)\,b_{\,k\,-\,1}\,p_{\,\alpha}\ +\ b_{\,k}\,q_{\,\alpha}\ +\ \frac{1}{k}\;\bigl(\,-\,p\,b_{\,k\,-\,1,\,\alpha}\ +\ 3\,b_{\,k\,-\,3,\,\alpha}\,\bigr)\,\biggr]\,z^{\,k}\,,
\end{equation*}
\begin{multline*}
  u_{\,\alpha\,\alpha}\ =\ \sum_{k\,=\,0}^{+\,\infty}\,\biggl[\,\Bigl(\,2\ -\ \frac{2}{k}\,\Bigr)\,b_{\,k\,-\,1,\,\alpha}\,p_{\,\alpha}\ +\ \Bigl(\,1\ -\ \frac{1}{k}\,\Bigr)\,b_{\,k\,-\,1}\,p_{\,\alpha\,\alpha}\ +\ b_{\,k}\,q_{\,\alpha\,\alpha}\ +\ 2\,b_{\,k,\,\alpha}\,q_{\,\alpha}\ +\\ 
  \frac{1}{k}\;\bigl(\-\,p\,b_{\,k\,-\,1,\,\alpha\,\alpha}\ +\ 3\,b_{\,k\,-\,3,\,\alpha\,\alpha}\,\bigr)\,\biggr]\,z^{\,k}\ +\\ 
  \frac{1}{3\,z^{\,2}\,-\,p}\;\sum_{k\,=\,0}^{+\,\infty}\,\biggl[\,(\,k\,-\,1\,)\,b_{\,k\,-\,1}\,p_{\alpha}^{\,2}\ +\ 2\,k\,b_{\,k}\,p_{\,\alpha}\,q_{\,\alpha}\ +\ (\,k\,+\,1\,)\,b_{\,k\,+\,1}\,q_{\alpha}^{\,2}\,\biggr]\,z^{\,k}\,,
\end{multline*}
where the symbol $\alpha\ \in\ \set*{t,\,x}\,$.

Then, using the definition of the mapping $\FF$ along with the fact that
\begin{equation*}
  \Bigl(\,\ph\,(\,t,\,x\,),\,\qh\,(\,t,\,x\,),\,\set{\bh_{\,k}\,(\,t,\,x\,)}_{\,k\,=\,0}^{\,\infty}\,\Bigr) \ =\ \Bigl(\,p\,(\,t,\,x\,),\,q\,(\,t,\,x\,),\,\set*{b_{\,k}\,(\,t,\,x\,)}_{\,k\,=\,0}^{\,\infty}\,\Bigr)\,,
\end{equation*}
one obtains the result stated in this Lemma. For more details, we refer also to Appendix~\ref{app:0}.
\end{proof}
To the best of our knowledge, the presented algorithm is completely new.


\subsection{Generalizations}

We may attempt to generalize the result stated in \cref{thm:ck} to the case of the \acrshort{pde}s family \eqref{eq:fam} with $m\ \geq\ 2\,$. However, we stress that at the current stage, it remains at the level of a conjecture:

\begin{conj}\label{conj:m}
Consider the following ramified \acrshort{ivp} for \cref{eq:fam}:
\begin{subequations}\label{eq:psysS}
 \begin{align}
   \L_{\,m}\,(\,u\,)\ \bydef\ \partial_{t}^{\,m}u\ -\ \partial_{x}^{\,m-1}\,u\,\partial_{x}^{\,m}\,u\ &=\ \fO\,, \\
   u\,(\,0,\,x\,)\ &=\ c_{\,1}\,x^{\,m\,-\,1}\ +\ c_{\,2}\,x^{\,m\,-\,1\,+\,\third} \nonumber \\ 
   &\quad +\ c_{\,3}\,x^{\,m\,-\,1\,+\,\frac{2}{3}}\ +\ x^{\,m}\cdot u_{\,0}\,(\,x\,)\,, \label{eq:psysb}
 \end{align}
\end{subequations}
where $c_{\,1,\,2}\ \in\ \C^{\,\times}\,$, $c_{\,3}\ \in\ \C$ and $u_{\,0}$ is a holomorphic function in some neighbourhood of $0\ \in\ \C\,$. Then, for each choice of a root to the algebraic equation\footnote{In order to have $m$ distinct roots, we have to choose the initial data so that $\partial_{\,x}^{\,m\,-\,1}\,u\,(\,0,\,0\,)\ \neq\ 0\,$. As a result, we may say that the solution procedure is parametrized by the elements of the multiplicative group of the roots of unity $U_{\,m}\ \eqdef\ \Set*{\xi^{\,m}\ =\ 1}{\xi\ \in\ \C}\,$.} $\tau^{\,m}\ -\ \partial_{\,x}^{\,m\,-\,1}\,u\,(\,0,\,0\,)\cdot 1^{m}\ =\ \fO$ (where we replaced $\xi\ \leftarrow\ 1\ \neq\ 0$ by homogeneity of the principal symbol), there exist holomorphic functions $p\,(\,t,\,x\,)\,$, $q\,(\,t,\,x\,)$ and $\bigl\{\,a_{\,j}\,(\,t,\,x\,)\,\bigr\}_{j\,=\,0}^{\,2}$ defined in the vicinity of $(\,0,\,0\,)\ \in\ \C^{\,2}$ such that there exists a local solution $u$ to \eqref{eq:psysS} which can be written in the following form:
\begin{equation*}
  u\,(\,t,\,x\,)\ \coloneq\ \partial_{\,q}^{\,-\,(\,m\,-\,1\,)}\comp\Bigl(\,\sum_{j\,=\,0}^{2}\,a_{\,j}\,(\,t,\,x\,)\,z^{\,j}\,\Bigr)\,,
\end{equation*}
where the function $z$ satisfies the following algebraic relation \eqref{eq:cube}:
\begin{equation*}
  z^{\,3}\ =\ p\,(\,t,\,x\,)\,z\ +\ q\,(\,t,\,x\,)\,.
\end{equation*}
The coefficients $p\,$, $q$ verify the `initial' conditions:
\begin{equation*}
  p\,(\,0,\,x\,)\ =\ \fO\,, \qquad
  q\,(\,0,\,x\,)\ =\ x\,.
\end{equation*}
\end{conj}
\begin{proof}
Open problem.
\end{proof}

\begin{remark}
The holomorphic functions $p\,(\,t,\,x\,)\,$, $q\,(\,t,\,x\,)$ and $\bigl\{\,a_{\,j}\,(\,t,\,x\,)\,\bigr\}_{j\,=\,0}^{\,2}$ whose existence is conjectured above depend also on $c_{\,1,\,2,\,3}$ and on the function $u_{\,0}\,(\,x\,)$ appearing in the formulation of the \nm{Cauchy} problem \eqref{eq:psysb}.
\end{remark}


\section{Numerical illustrations}
\label{sec:num}

In the absence of rigorous theoretical proof, in this Section, we would like to illustrate the practical convergence of the fixed point algorithm described in \cref{sec:it}. Thus, it will provide rational computational support towards various conjectures formulated in this study (in particular, the \cref{conj:2} and indirectly \cref{conj:m}).


\subsection{Test 1}
\label{sec:test1}

In order to come up with a practical algorithm, we have to truncate the expansion \eqref{eq:ser}:
\begin{equation*}
  u\,(\,t,\,x\,)\ \approx\ \sum_{k\,=\,1}^{N}\,\frac{1}{k}\;\Bigl(\,-\,p\,(\,t,\,x\,)\,b_{\,k-1}\,(\,t,\,x\,)\ +\ 3\,b_{\,k-3}\,(\,t,\,x\,)\,\Bigr)\,z^{\,k}\,,
\end{equation*}
Since we work with holomorphic functions, we shall expand all the functions in double truncated series in the vicinity of the origin $(\,0,\,0\,)\ \in\ \C^{\,2}$ as follows:
\begin{equation*}
  \psi\,(\,t,\,x\,)\ =\ \sum_{\substack{l,\,m\,=\,1 \\ l\,+\,m\ \leq\ M\,-\,1}}^{M\,-\,1} \psi_{\,l\,m}\;t^{\,l}\,x^{\,m}\ +\ \O\,\bigl(\,t^{\,M}\ +\ t^{\,M\,-\,1}\,x\ +\ \ldots\ t\,x^{\,M\,-\,1}\ +\ x^{\,M}\,\bigr)\,,
\end{equation*}
with $\psi\ \in\ \set*{p,\,q,\,b_{\,0},\,\ldots,\,b_{\,N\,-\,1}}\,$. So, in practice, we perform all our computations with such double truncated power series in independent variables $(\,t,\,x\,)\ \in\ \C^{\,2}\,$. We have to describe another technical issue towards the practical implementation of our algorithm here. Since we work with truncated data, our finitary mapping $\FFt$ by construction returns a slightly reduced set of data:
\begin{equation*}
  \Bigl(\,p\,(\,t,\,x\,),\,q\,(\,t,\,x\,),\,\set*{b_{\,k}\,(\,t,\,x\,)}_{\,k\,=\,0}^{\,N}\,\Bigr)\ \stackrel{\FFt}{\mapsto}\ \Bigl(\,\ph\,(\,t,\,x\,),\,\qh\,(\,t,\,x\,),\,\set{\bh_{\,k}\,(\,t,\,x\,)}_{\,k\,=\,0}^{\,N\,-\,3}\,\Bigr)\,.
\end{equation*}
Thus, at every iteration, we are practically losing three components of the data. Consequently, if we are willing to obtain $N$ components of data vector $\set*{b_{\,k}\,(\,t,\,x\,)}_{\,k\,=\,0}^{\,N}$ after $I$ iterations, we have to initiate the iterative process with $\set*{b_{\,k}\,(\,t,\,x\,)}_{\,k\,=\,0}^{\,N\,+\,3\, I}$ initial components. We took into account this observation in our computations. We are going to illustrate the work of the algorithm starting from the following initial guess:
\begin{subequations}\label{eq:data0}
\begin{align}\label{eq:data}
  p\,(\,t,\,x\,)\ &\coloneq\ \frac{t}{2}\,, \qquad
  q\,(\,t,\,x\,)\ \coloneq\ t\ +\ x\,, \\
  b_{\,0}\ &\coloneq\ \fW\,, \qquad b_{\,1}\ \coloneq\ \fW\,, \qquad b_{\,j}\ \coloneq\ \fO\,, \qquad \forall\,j\ \geq\ 2\,.
\end{align}
\end{subequations}
It corresponds to the \nm{Cauchy} datum \eqref{eq:init1} as it is explained at the beginning of \cref{sec:it}. Throughout all this \cref{sec:num} we shall assume $c_{\,1}\ =\ 1$ and $c_{\,2}\ =\ \frac{3}{4}\,$. However, three different possibilities are taken for the rest of coefficients\footnote{For the precise meaning of coefficients $\set*{c_{\,j}}_{\,j\,\geq\,0}$ see \cref{eq:ic3}.} $\set*{c_{\,j}}_{\,j\,\geq\,3}\ \subseteq\ \C\,$. For example, in this \cref{sec:test1} we take all $\set*{c_{\,j}}_{\,j\,\geq\,3}$ to be zero.

\begin{remark}
The choice of the initial data \eqref{eq:data} implies that $q_{\,t}\,(\,0,\,0\,)\ =\ 1\,$. This specifies one of two roots in $q_{\,t}\,(\,0,\,0\,)$ in the eikonal equation. Consequently, in all iterations we shall also have that $\qh_{\,t}\,(\,0,\,0\,)\ =\ 1\,$. Similarly, we shall also have $\ph_{\,t}\,(\,0,\,0\,)\ =\ \dfrac{1}{2}\,$.
\end{remark}

\begin{table}
  \centering
  \caption{\small\em Various numerical parameters used in our computations.}
  \medskip
  \label{tab:params}
  \begin{tabular}{lc}
  \toprule
  \textit{Parameter} & \textit{Value} \\
  \midrule
  Floating point arithmetics, significant digits & $30$ \\
  Truncation degree in $t\,$, $M$ & $25$ \\
  Truncation degree in $x\,$, $M$ & $25$ \\
  Total number of iterations, $I$ & $25$ \\
  Initial length of data vector, $N$ & $80\ =\ 3\times I\ +\ 5$ \\
  Observation point in $\C\,$, $x_{\,0}$ & $0.1\,\ui$ \\
  Temporal segment for the error control, $[\,a,\,b\,]$ & $[\,0,\,0.1\,]$ \\
  \bottomrule
  \end{tabular}
  \bigskip
\end{table}

In order to speed up the computations, we also turn to the floating-point arithmetics with, possibly, extended precision. The values of all other parameters are reported in \cref{tab:params}. The proposed algorithm was implemented in the computer algebra system \texttt{Maple}${}^{\text{\texttrademark}}$. The essential part of the employed code is reported in Appendix~\ref{app:b}. The complete program can be shared upon a simple request by email.

\begin{figure}
  \centering
  \subfigure[$t\ \mapsto\ p\,(\,t,\,x_{\,0}\,)$]{\includegraphics[width=0.48\textwidth]{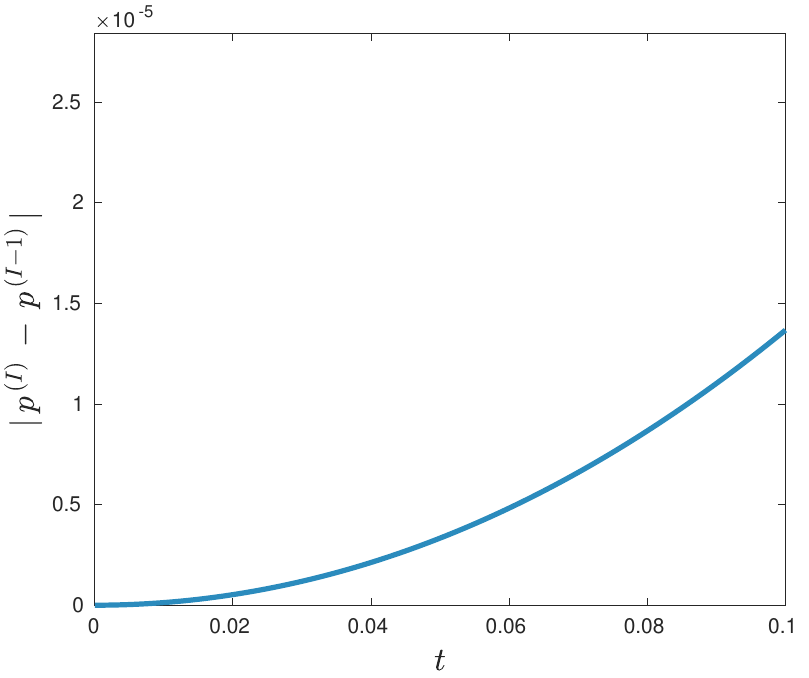}}
  \subfigure[$t\ \mapsto\ q\,(\,t,\,x_{\,0}\,)$]{\includegraphics[width=0.50\textwidth]{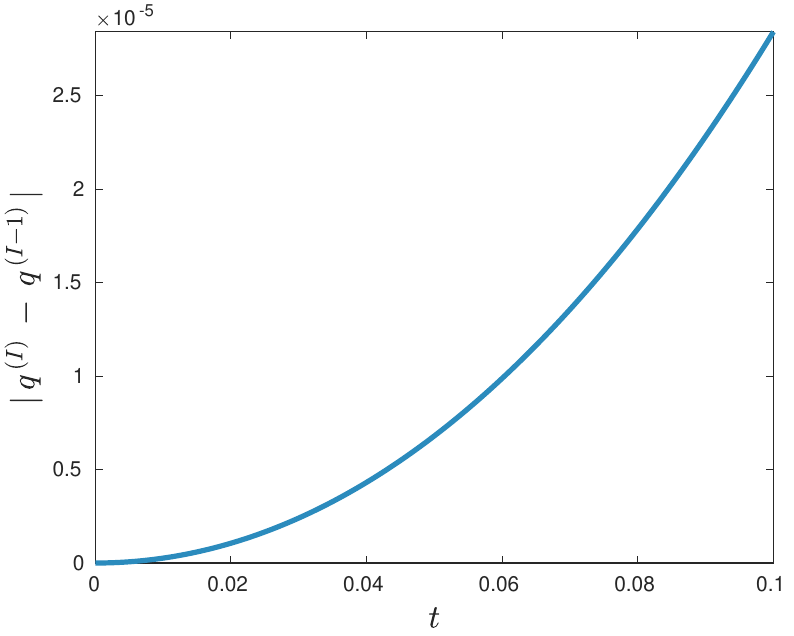}}
  \subfigure[$t\ \mapsto\ b_{\,0}\,(\,t,\,x_{\,0}\,)$]{\includegraphics[width=0.48\textwidth]{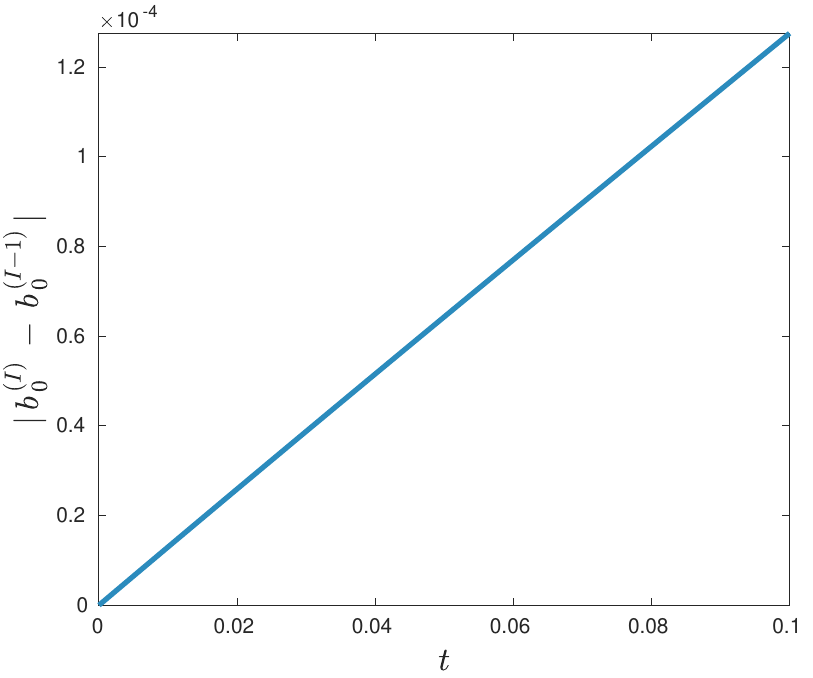}}
  \subfigure[$t\ \mapsto\ b_{\,1}\,(\,t,\,x_{\,0}\,)$]{\includegraphics[width=0.48\textwidth]{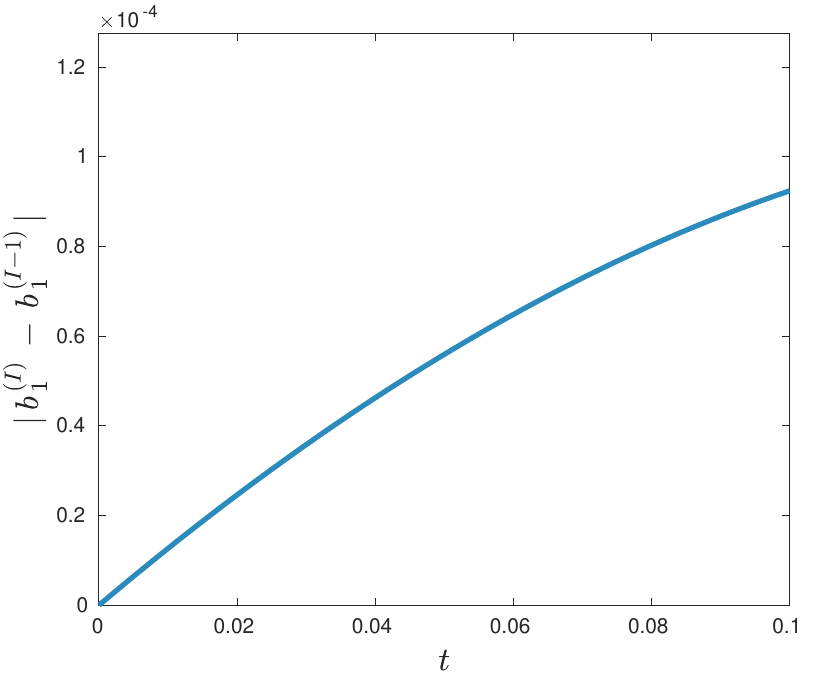}}
  \caption{\small\em The absolute value of the difference between two successive iterations $(I)$\up{th} and $(I\,-\,1)$\up{th} of solution data for $I\ =\ 10\,$. The values of several numerical parameters are reported in \cref{tab:params}.}
  \label{fig:figs}
\end{figure}

\begin{figure}
  \centering
  \subfigure[$\norm*{p^{\,(\,I\,)\,}\,(\,t,\,x_{\,0}\,)\ -\ p^{\,(\,I\,-\,1\,)}\,(\,t,\,x_{\,0}\,)}_{\,\infty}$]{\includegraphics[width=0.48\textwidth]{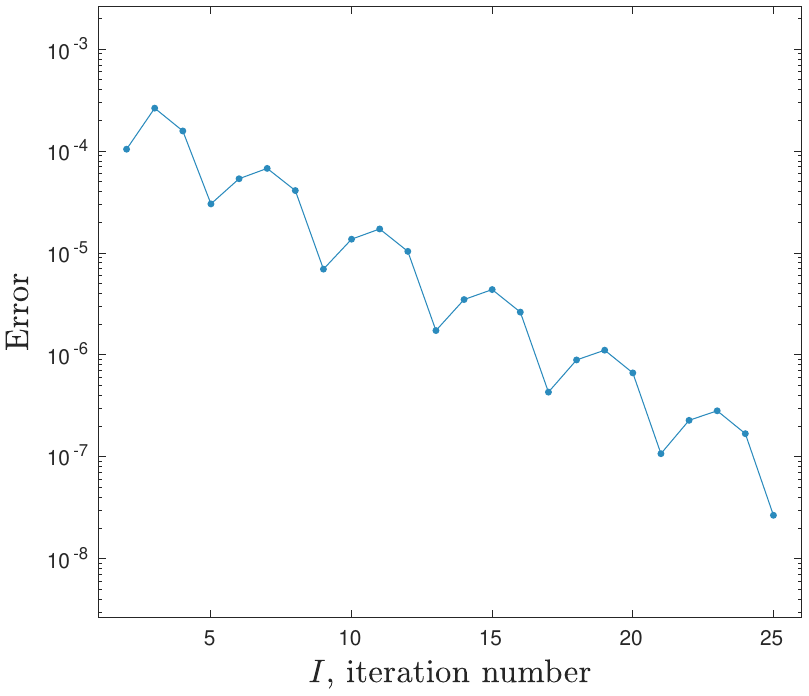}}
  \subfigure[$\norm*{q^{\,(\,I\,)\,}\,(\,t,\,x_{\,0}\,)\ -\ q^{\,(\,I\,-\,1\,)}\,(\,t,\,x_{\,0}\,)}_{\,\infty}$]{\includegraphics[width=0.50\textwidth]{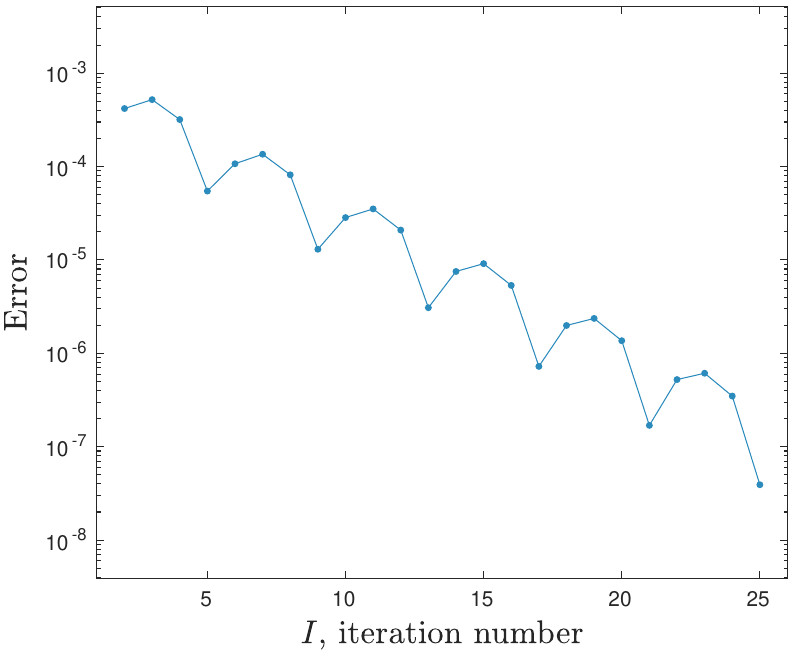}}
  \subfigure[$\norm*{b_{\,0}^{\,(\,I\,)\,}\,(\,t,\,x_{\,0}\,)\ -\ b_{\,0}^{\,(\,I\,-\,1\,)}\,(\,t,\,x_{\,0}\,)}_{\,\infty}$]{\includegraphics[width=0.48\textwidth]{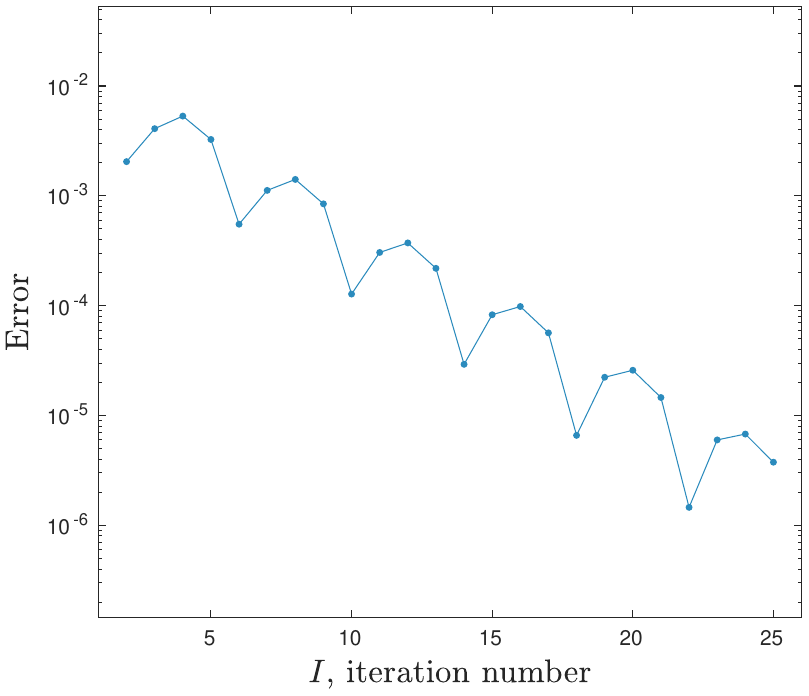}}
  \subfigure[$\norm*{b_{\,1}^{\,(\,I\,)\,}\,(\,t,\,x_{\,0}\,)\ -\ b_{\,1}^{\,(\,I\,-\,1\,)}\,(\,t,\,x_{\,0}\,)}_{\,\infty}$]{\includegraphics[width=0.48\textwidth]{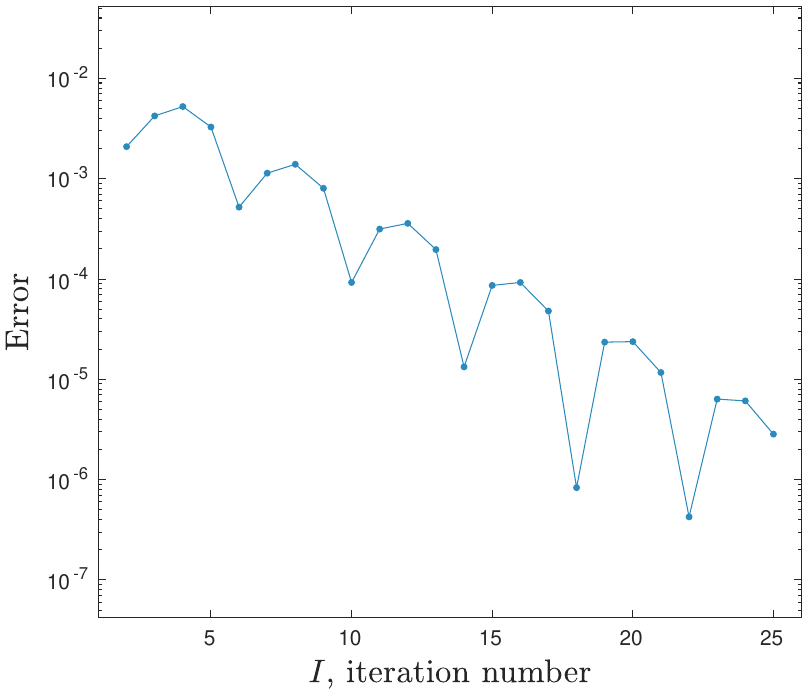}}
  \caption{\small\em The maximum norm $L_{\,\infty}\,\bigl(\,[\,a,\,b\,]\,\bigr)$ of the difference between two successive iterations $(I)$\up{th} and $(I\,-\,1)$\up{th} as a function of the iteration number in semi-logarithmic co\"ordinates. The initial data is given in \cref{eq:data0}. The values of several numerical parameters are reported in \cref{tab:params}.}
  \label{fig:err}
\end{figure}

To summarize, we can say that we perform numeric-symbolic computations in order to observe the convergence towards a fixed point in practice. If we denote by $\psi^{(\,j\,)}$ the value of the variable $\psi\ \in\ \set*{p,\,q,\,b_{\,0},\,\ldots,\,b_{\,N\,-\,3\,j}}$ after $j$ iterations, the closeness to the fixed point may be appreciated by looking at the norm of the difference between two successive iterations:
\begin{equation*}
  \norm{\psi^{\,(\,j\,)}\ -\ \psi^{\,(\,j\,-\,1\,)}}\,.
\end{equation*}
We illustrate the behaviour of the proposed algorithm on the initial data \eqref{eq:data0}. We perform $I$ iterations of our mapping. It turns out that even very moderate values of the parameter $I$ (\cf \cref{tab:params}) are enough to appreciate the convergence of the iterative fixed point process. The whole symbolic/numeric computation lasted about five minutes of the \acrshort{cpu} time (on our computers). The precise information is not very important because it may vary from one system to another. In four panels of \cref{fig:figs} we show the differences between two last iterations of functions $\set*{p,\,q,\,b_{\,0},\,b_{\,1}}\,$, constituting the problem data. It can be clearly seen that the absolute value of the difference is quite small already after $I\ =\ 10$ iterations. We checked that subsequent iterations reduced this difference further. The reduction of the $L_{\,\infty}$ norm of the difference between two successive iterations is shown in \cref{fig:err}. This Figure presents strong empirical evidence for the convergence towards a fixed point, which solves the underlying \acrshort{ivp} according to \cref{eq:prop}, even if the convergence seems to be far from being linear.

\begin{figure}
  \centering
  \subfigure[$\norm*{p^{\,(\,I\,)\,}\,(\,t,\,x_{\,0}\,)\ -\ p^{\,(\,I\,-\,1\,)}\,(\,t,\,x_{\,0}\,)}_{\,\infty}$]{\includegraphics[width=0.48\textwidth]{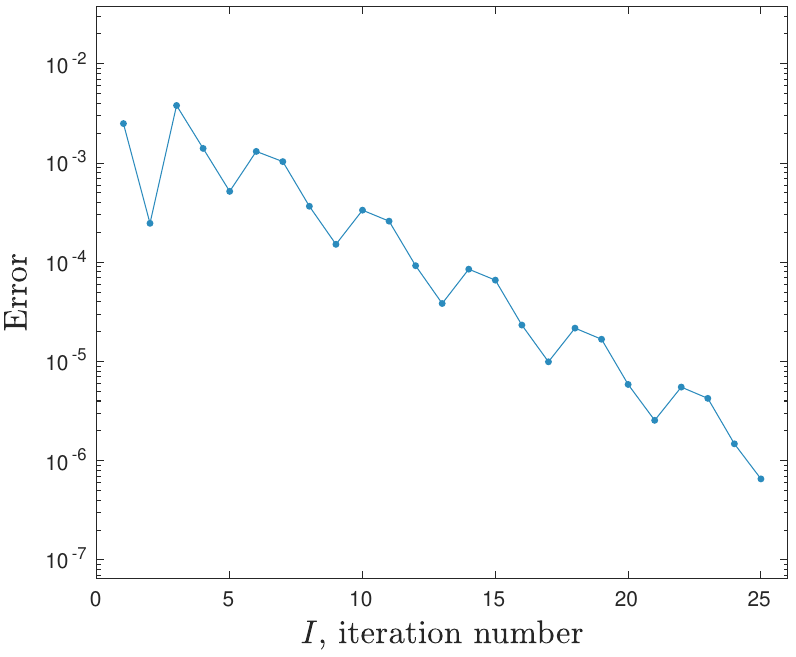}}
  \subfigure[$\norm*{q^{\,(\,I\,)\,}\,(\,t,\,x_{\,0}\,)\ -\ q^{\,(\,I\,-\,1\,)}\,(\,t,\,x_{\,0}\,)}_{\,\infty}$]{\includegraphics[width=0.46\textwidth]{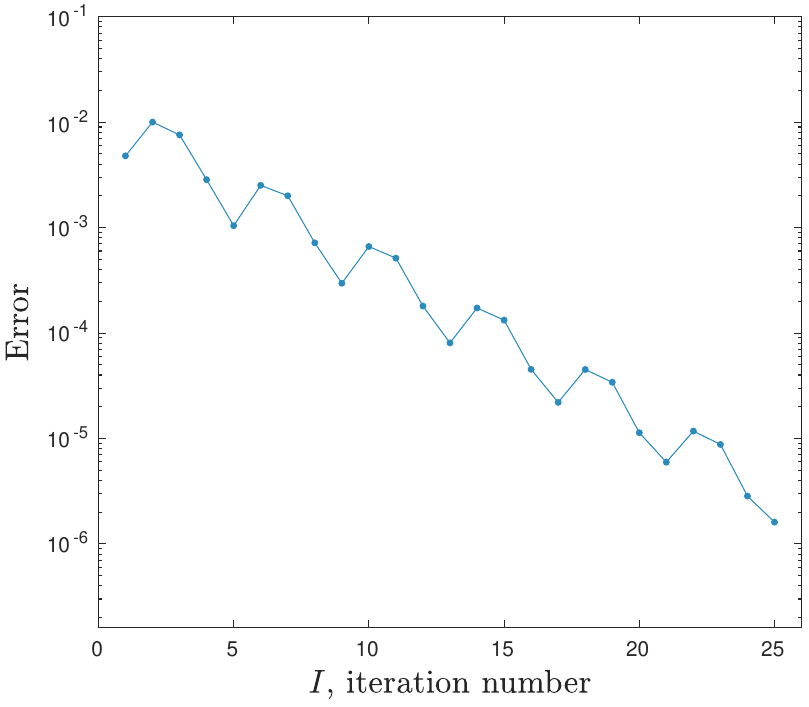}}
  \subfigure[$\norm*{b_{\,0}^{\,(\,I\,)\,}\,(\,t,\,x_{\,0}\,)\ -\ b_{\,0}^{\,(\,I\,-\,1\,)}\,(\,t,\,x_{\,0}\,)}_{\,\infty}$]{\includegraphics[width=0.48\textwidth]{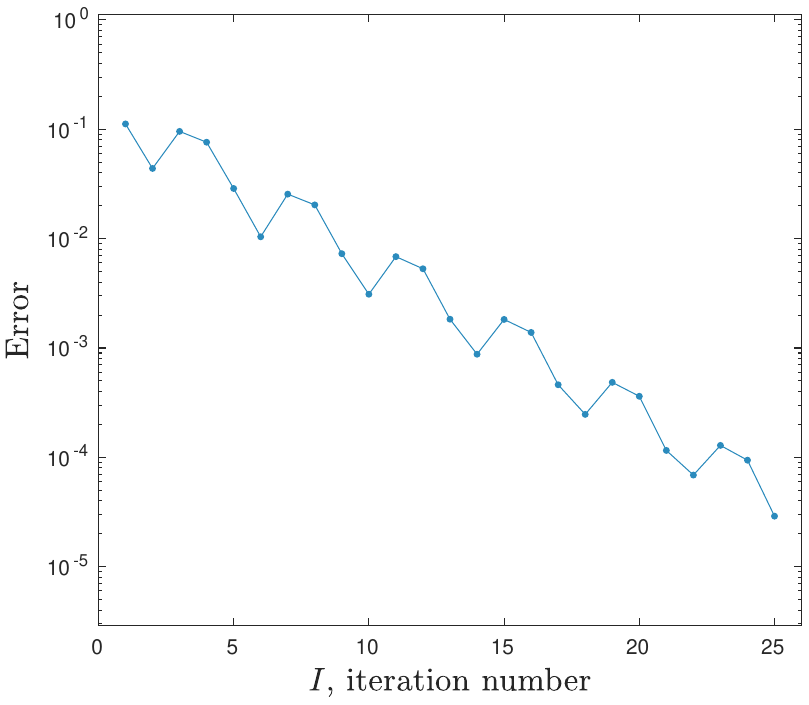}}
  \subfigure[$\norm*{b_{\,1}^{\,(\,I\,)\,}\,(\,t,\,x_{\,0}\,)\ -\ b_{\,1}^{\,(\,I\,-\,1\,)}\,(\,t,\,x_{\,0}\,)}_{\,\infty}$]{\includegraphics[width=0.50\textwidth]{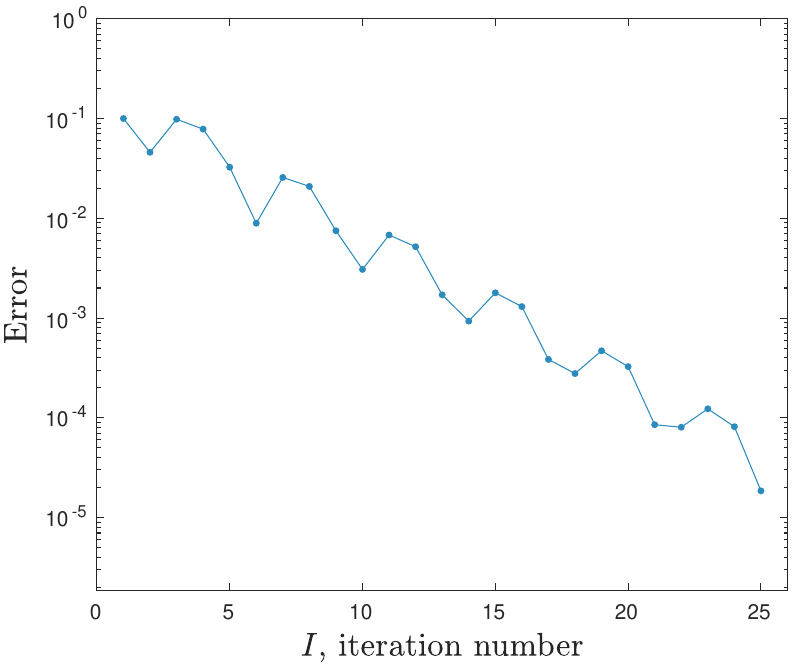}}
  \caption{\small\em The maximum norm $L_{\,\infty}\,\bigl(\,[\,a,\,b\,]\,\bigr)$ of the difference between two successive iterations $(I)$\up{th} and $(I\,-\,1)$\up{th} as a function of the iteration number in semi-logarithmic co\"ordinates. The initial data is given in \cref{eq:data1}. The values of several numerical parameters are reported in \cref{tab:params}.}
  \label{fig:err1}
\end{figure}

\begin{figure}
  \centering
  \subfigure[$t\ \mapsto\ \abs*{p_{\,1}^{\,(\,I\,)}\ -\ p_{\,2}^{\,(\,I\,)}}$]{\includegraphics[width=0.48\textwidth]{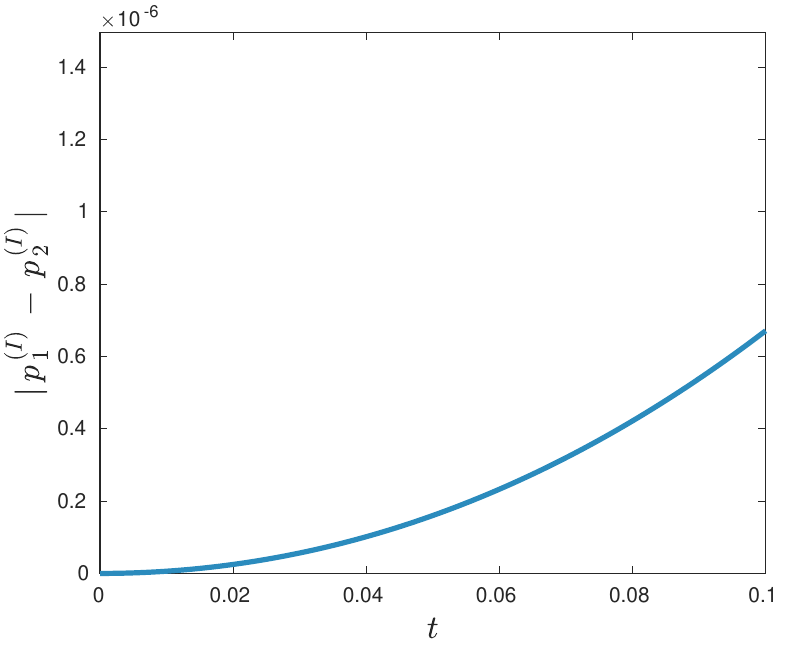}}
  \subfigure[$t\ \mapsto\ \abs*{q_{\,1}^{\,(\,I\,)}\ -\ q_{\,2}^{\,(\,I\,)}}$]{\includegraphics[width=0.50\textwidth]{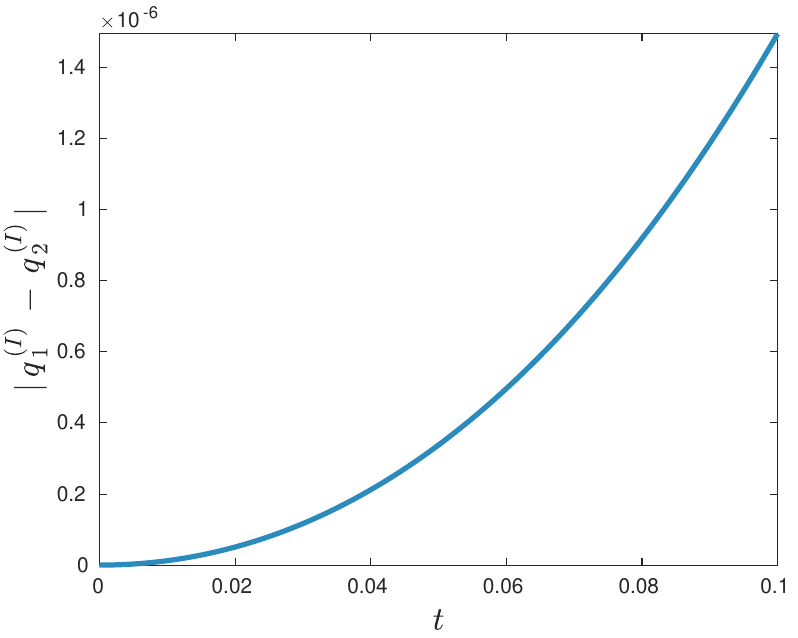}}
  \subfigure[$t\ \mapsto\ \abs*{b_{\,1,\,0}^{\,(\,I\,)}\ -\ b_{\,2,\,0}^{\,(\,I\,)}}$]{\includegraphics[width=0.48\textwidth]{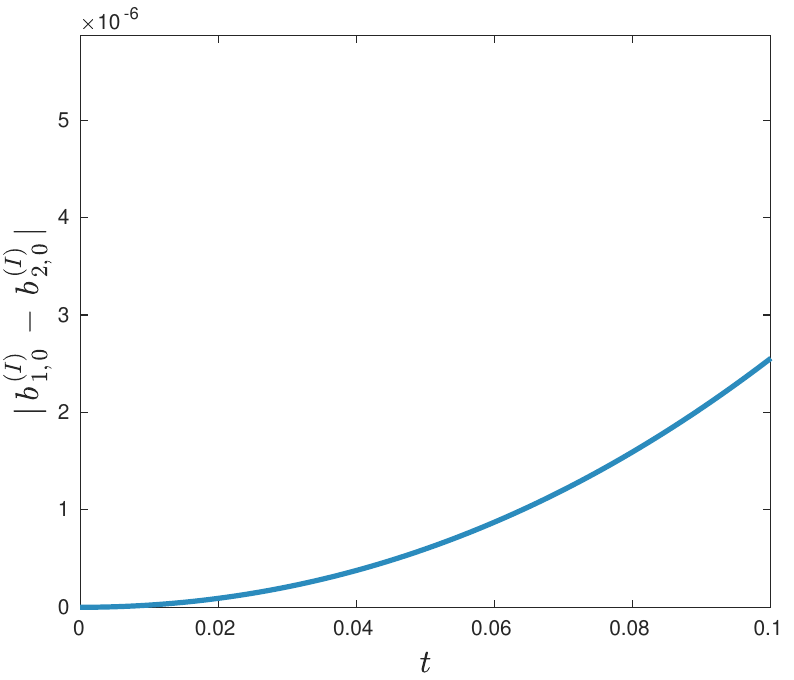}}
  \subfigure[$t\ \mapsto\ \abs*{b_{\,1,\,1}^{\,(\,I\,)}\ -\ b_{\,2,\,1}^{\,(\,I\,)}}$]{\includegraphics[width=0.48\textwidth]{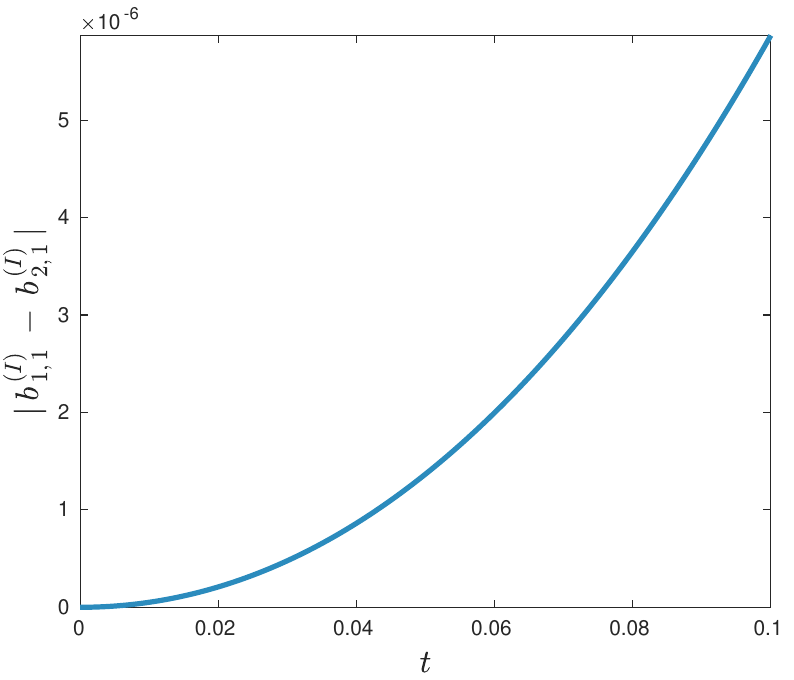}}
  \caption{\small\em Differences between corresponding components of the solution data at $I\ =\ 25$\up{th} iteration. Both iterative processes were initialized with data \eqref{eq:data0} and \eqref{eq:data1} respectively. The values of all numerical parameters are given in \cref{tab:params}.}
  \label{fig:diff}
\end{figure}

\begin{figure}
  \centering
  \includegraphics[width=0.59\textwidth]{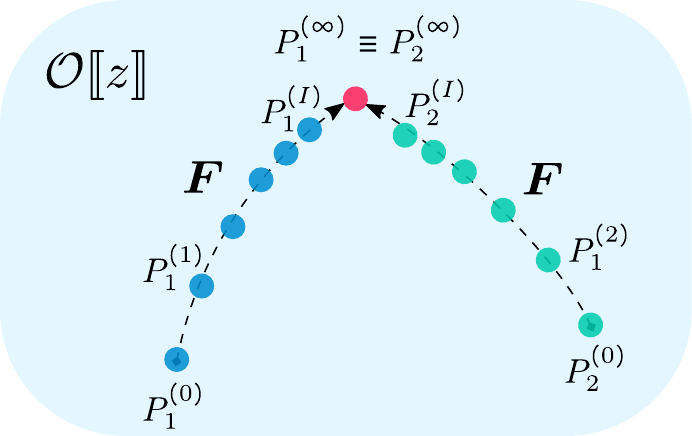}
  \caption{\small\em A schematic representation of the experimental findings of \cref{sec:num}: $P_{\,1}^{\,(\,0\,)}$ stands for the initial guess \eqref{eq:data0} while $P_{\,2}^{\,(\,0\,)}$ stands for \eqref{eq:data1}. The convergence should be understood in the sense of a numerically-supported conjecture.}
  \label{fig:conv}
\end{figure}

In order to present even more convincing numerical evidence, we include in our study another test case with the following initial guess:
\begin{subequations}\label{eq:data1}
\begin{align}
  p\,(\,t,\,x\,)\ &\coloneq\ \frac{t}{2}\,, \qquad
  q\,(\,t,\,x\,)\ \coloneq\ t\ +\ x\,, \\
  b_{\,0}\ &\coloneq\ \fW\ +\ t\,, \qquad b_{\,1}\ \coloneq\ \fW\ -\ t\,, \qquad b_{\,j}\ \coloneq\ \fO\,, \qquad \forall\,j\ \geq\ 2\,.
\end{align}
\end{subequations}
It is not difficult to see that these initial data \eqref{eq:data1} (which correspond to a different initialization with respect to \eqref{eq:data0}) verify the same initial conditions \eqref{eq:iconds}, \ie the \nm{Cauchy} datum \eqref{eq:init1}. Thus, if the solution we seek is unique\footnote{We may reasonably assume that at least within the solution ansatz we consider. Perhaps the best reason to believe in the uniqueness would be the \nm{Lax}--\nm{Glimm} theory \cite{Glimm1970}, but ``we sadly lack a local uniqueness theorem'' as the authors of \cite{Glimm1970} put it themselves (in the real case).}, we may expect the iterations to converge to the same fixed point. This explains the motivation behind this second numerical study. We are using precisely the same numerical parameters as reported above in \cref{tab:params} (and, thus, the same parameters as in the previous computations). The convergence of the fixed point iterations is reported in Figure~\ref{fig:err1}. We can observe the same (\cf Figure~\ref{fig:err}) roughly linear, but non-monotonic convergence to a fixed point. Now one may ask a legitimate question whether iterations starting at \eqref{eq:data0} and \eqref{eq:data1} converge to the same point in $\O\,\llbracket\,z\,\rrbracket\,$? We may answer this question by looking at the difference between obtained solutions data at the final iteration $I\,$. We can have a look at these differences $t\ \mapsto\ \abs{\psi_{\,1}^{\,(\,I\,)}\,(\,t,\,x_{\,0}\,)\ -\ \psi_{\,2}^{\,(\,I\,)}\,(\,t,\,x_{\,0}\,)}$ in \cref{fig:diff} where they are represented as functions of time\footnote{It turns out that in this particular case, these differences are only functions of time, \ie they do not depend on $x\,$.} for $t\ \in\ [\,a,\,b\,]\ \equiv\ \Bigl[\,0,\,\dfrac{1}{10}\,\Bigr]\,$. In particular, we would like to underline the fact that the magnitude of these differences is consistent with the results presented in Figures~\ref{fig:err} and \ref{fig:err1}. All these observations provide strong empirical evidence for the convergence of the iterative process described above. The experimental findings of this Section are schematically depicted and summarized in \cref{fig:conv}.

\subsection{Test 2}
\label{sec:test2}

\begin{figure}
  \centering
  \subfigure[$\norm*{p^{\,(\,I\,)\,}\,(\,t,\,x_{\,0}\,)\ -\ p^{\,(\,I\,-\,1\,)}\,(\,t,\,x_{\,0}\,)}_{\,\infty}$]{\includegraphics[width=0.48\textwidth]{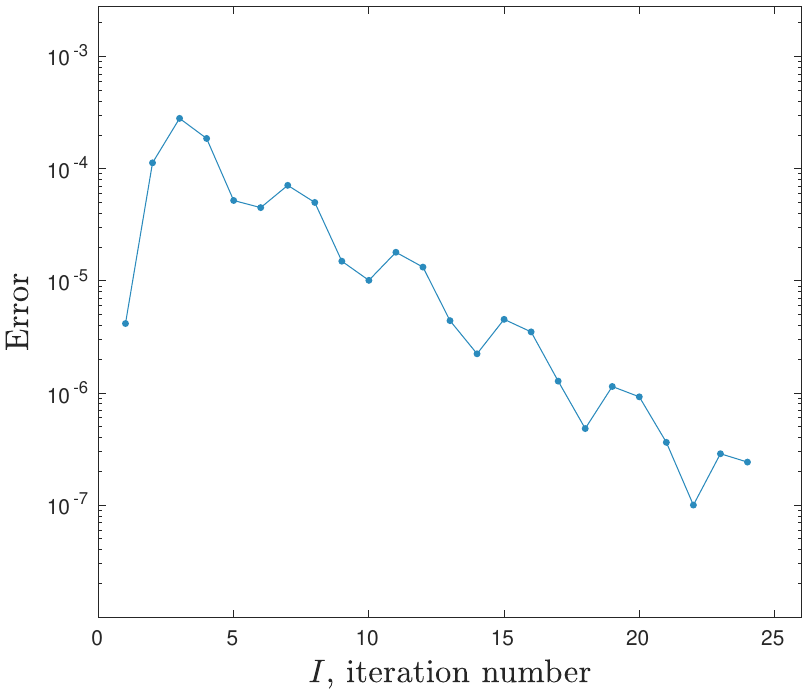}}
  \subfigure[$\norm*{q^{\,(\,I\,)\,}\,(\,t,\,x_{\,0}\,)\ -\ q^{\,(\,I\,-\,1\,)}\,(\,t,\,x_{\,0}\,)}_{\,\infty}$]{\includegraphics[width=0.48\textwidth]{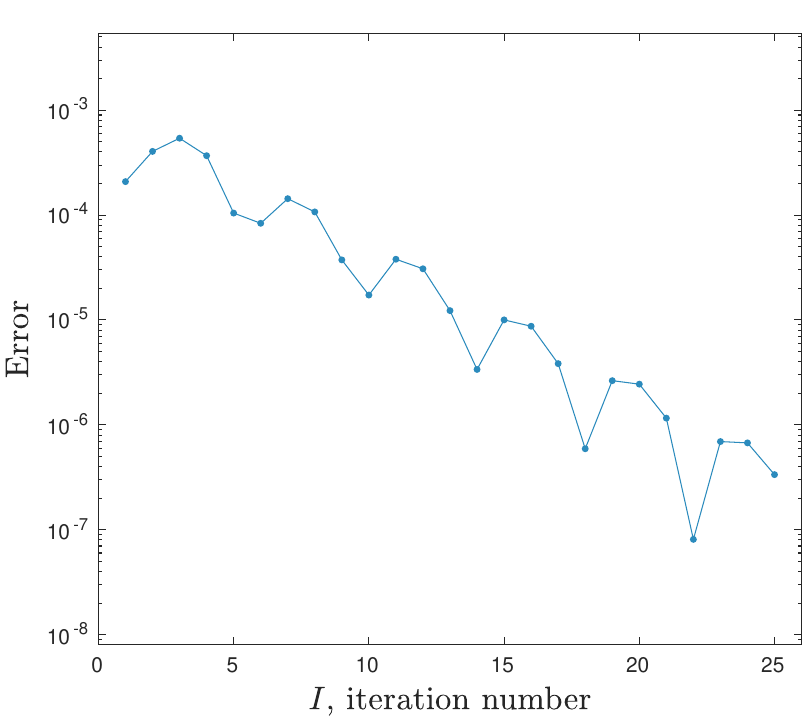}}
  \subfigure[$\norm*{b_{\,0}^{\,(\,I\,)\,}\,(\,t,\,x_{\,0}\,)\ -\ b_{\,0}^{\,(\,I\,-\,1\,)}\,(\,t,\,x_{\,0}\,)}_{\,\infty}$]{\includegraphics[width=0.48\textwidth]{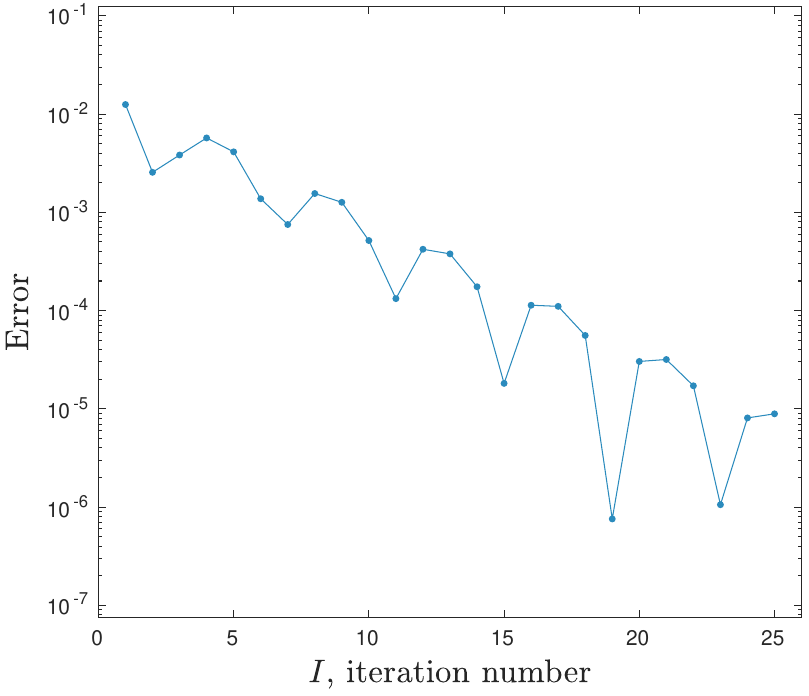}}
  \subfigure[$\norm*{b_{\,1}^{\,(\,I\,)\,}\,(\,t,\,x_{\,0}\,)\ -\ b_{\,1}^{\,(\,I\,-\,1\,)}\,(\,t,\,x_{\,0}\,)}_{\,\infty}$]{\includegraphics[width=0.48\textwidth]{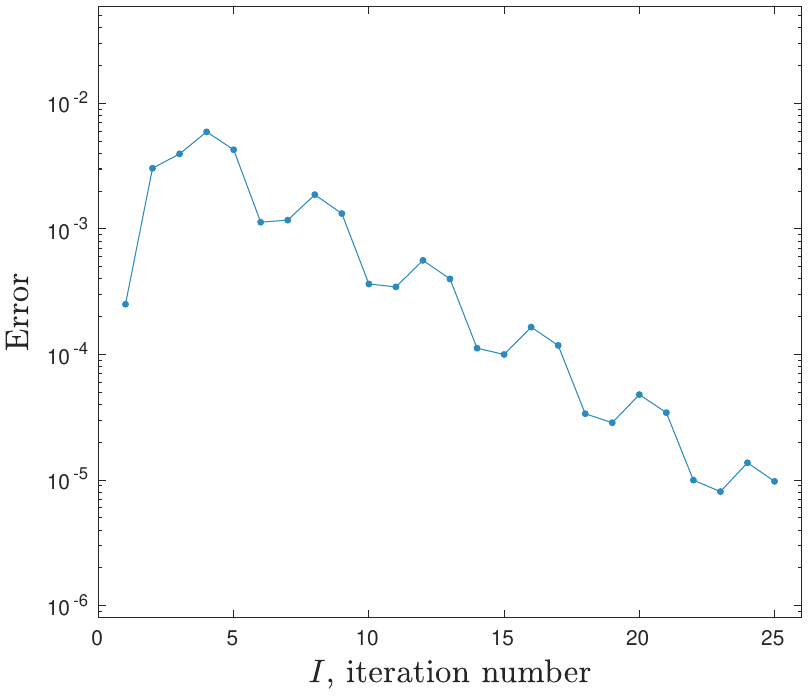}}
  \subfigure[$\norm*{b_{\,2}^{\,(\,I\,)\,}\,(\,t,\,x_{\,0}\,)\ -\ b_{\,2}^{\,(\,I\,-\,1\,)}\,(\,t,\,x_{\,0}\,)}_{\,\infty}$]{\includegraphics[width=0.48\textwidth]{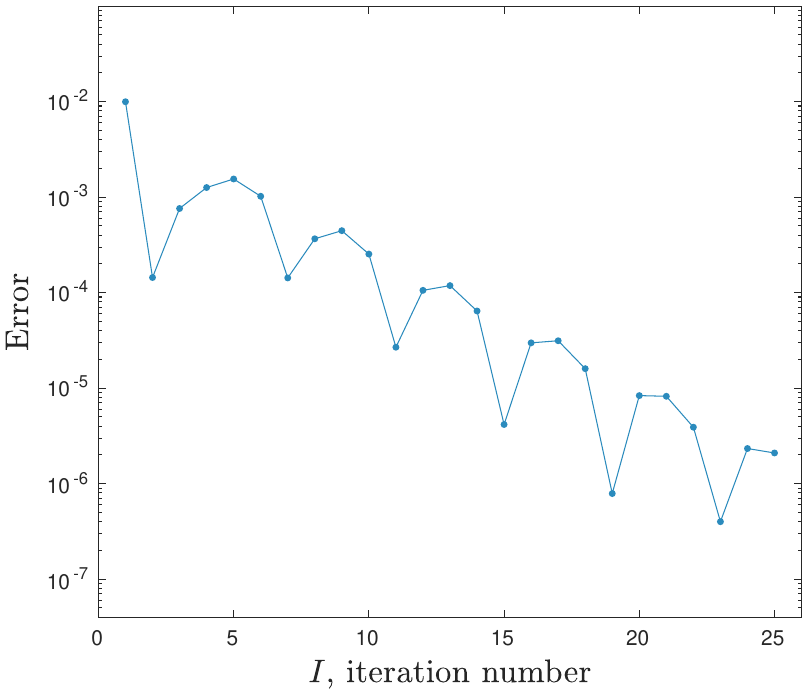}}
  \subfigure[$\norm*{b_{\,3}^{\,(\,I\,)\,}\,(\,t,\,x_{\,0}\,)\ -\ b_{\,3}^{\,(\,I\,-\,1\,)}\,(\,t,\,x_{\,0}\,)}_{\,\infty}$]{\includegraphics[width=0.48\textwidth]{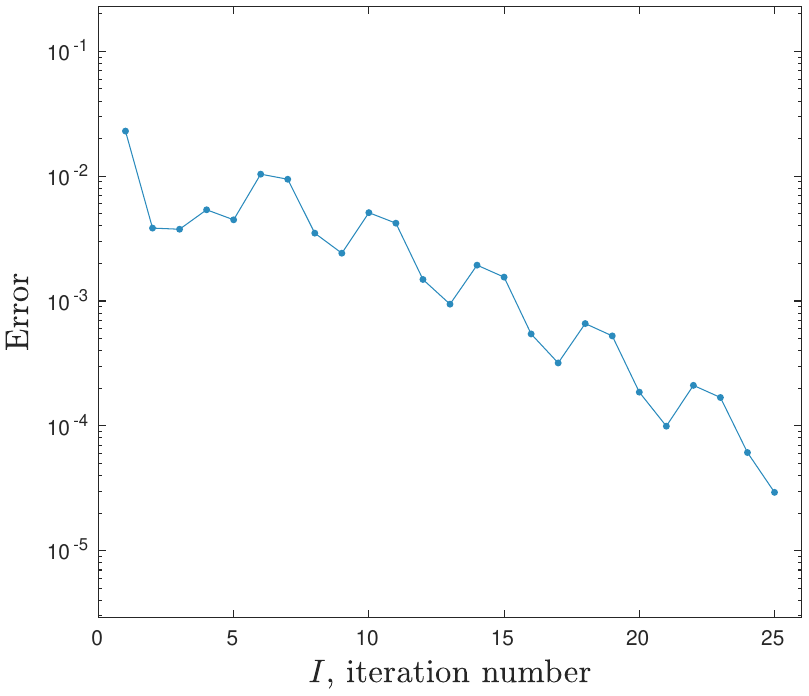}}
  \caption{\small\em The maximum norm $L_{\,\infty}\,\bigl(\,[\,a,\,b\,]\,\bigr)$ of the difference between two successive iterations $(I)$\up{th} and $(I\,-\,1)$\up{th}. The initial data is given in \cref{eq:data2}.}
  \label{fig:final}
\end{figure}

We shall consider a completely different initial condition from \eqref{eq:ic3}:
\begin{equation}\label{eq:init2}
  u\,(\,0,\,x\,)\ =\ x\ +\ \frac{3}{4}\;x^{\;\fourthirds}\ +\ \frac{3}{50}\;x^{\,\eightthirds}\ +\ \frac{1}{20}\;x^{\,3}\,.
\end{equation}
In terms of the solution data, the last initial datum translates in the following coefficients:
\begin{subequations}\label{eq:data2}
\begin{align}
  p\,(\,t,\,x\,)\ &\coloneq\ \frac{t}{2}\,, \qquad
  q\,(\,t,\,x\,)\ \coloneq\ t\ +\ x\,, \\
  b_{\,0}\ &\coloneq\ \fW\,, \qquad b_{\,1}\ \coloneq\ \fW\,, \\
  b_{\,2}\,(\,t,\,x\,)\ &\coloneq\ \frac{x}{10}\,, \qquad b_{\,3}\,(\,t,\,x\,)\ \coloneq\ \frac{x}{10}\,, \\
  \qquad b_{\,j}\ &\coloneq\ \fO\,, \qquad \forall\,j\ \geq\ 4\,.
\end{align}
\end{subequations}
The complexity of this test case is much higher than what we did before. In order to make the computational times reasonable, we had to reduce the parameter $M\ \coloneq\ 10\,$. All other numerical parameters were kept as in \cref{tab:params}. We studied the convergence of fixed point iterations under the map $\FFt$ for this \acrshort{ivp}. Namely, we monitored the differences between successive iterations for variables $p\,$, $q\,$, $b_{\,0}\,$, \dots, $b_{\,3}\,$. The results are reported in Figure~\ref{fig:final}. As we can see, the \emph{practical} linear non-monotonic convergence can be clearly witnessed even in this test case.

\begin{figure}
  \centering
  \subfigure[$\norm*{p^{\,(\,I\,)\,}\,(\,t,\,x_{\,0}\,)\ -\ p^{\,(\,I\,-\,1\,)}\,(\,t,\,x_{\,0}\,)}_{\,\infty}$]{\includegraphics[width=0.48\textwidth]{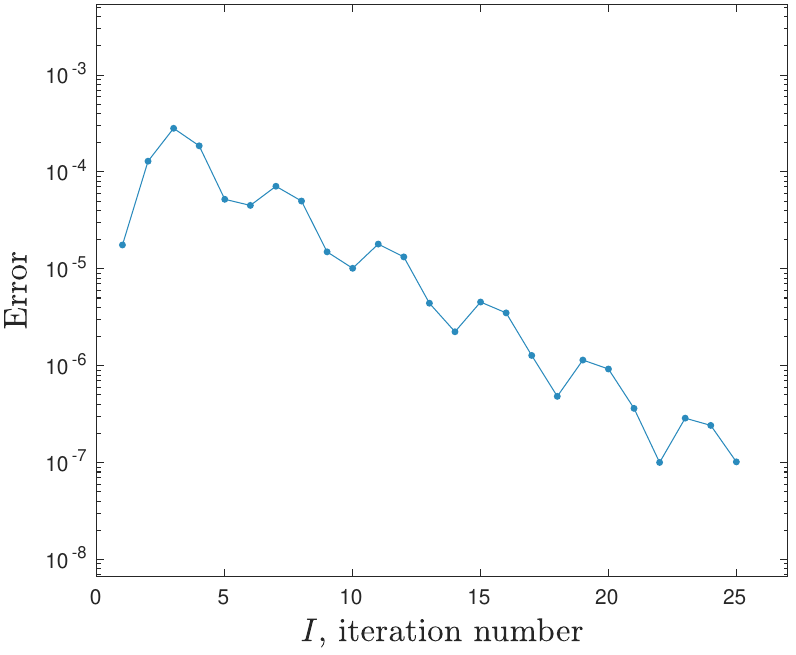}}
  \subfigure[$\norm*{q^{\,(\,I\,)\,}\,(\,t,\,x_{\,0}\,)\ -\ q^{\,(\,I\,-\,1\,)}\,(\,t,\,x_{\,0}\,)}_{\,\infty}$]{\includegraphics[width=0.48\textwidth]{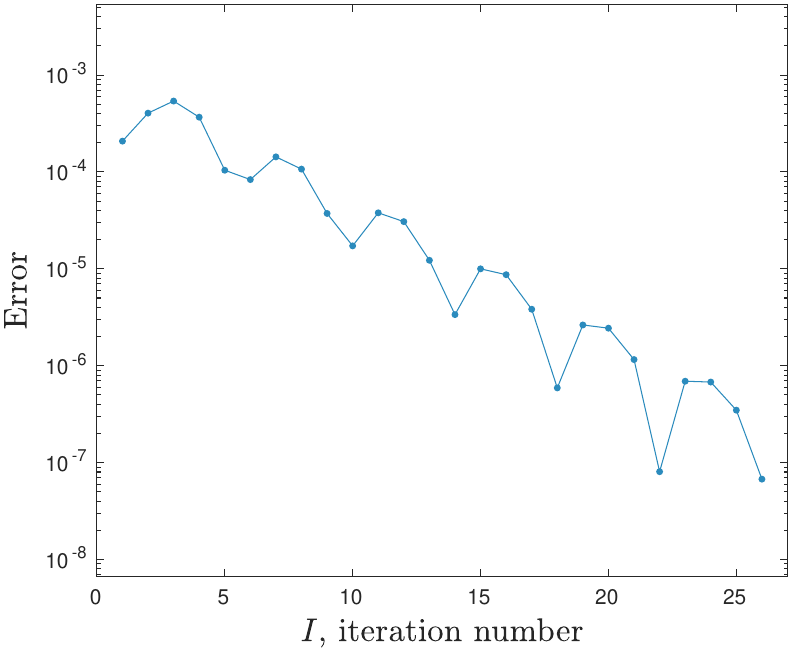}}
  \subfigure[$\norm*{b_{\,0}^{\,(\,I\,)\,}\,(\,t,\,x_{\,0}\,)\ -\ b_{\,0}^{\,(\,I\,-\,1\,)}\,(\,t,\,x_{\,0}\,)}_{\,\infty}$]{\includegraphics[width=0.48\textwidth]{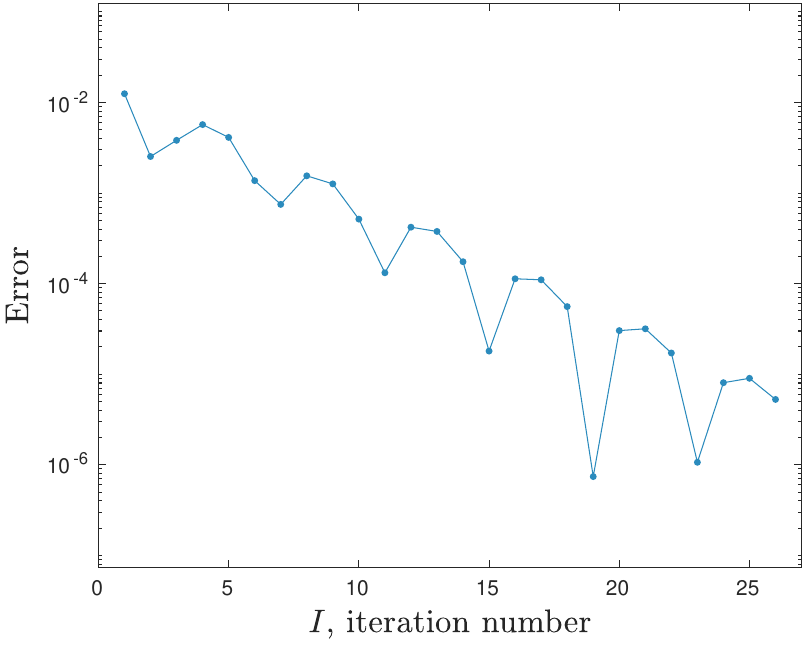}}
  \subfigure[$\norm*{b_{\,1}^{\,(\,I\,)\,}\,(\,t,\,x_{\,0}\,)\ -\ b_{\,1}^{\,(\,I\,-\,1\,)}\,(\,t,\,x_{\,0}\,)}_{\,\infty}$]{\includegraphics[width=0.48\textwidth]{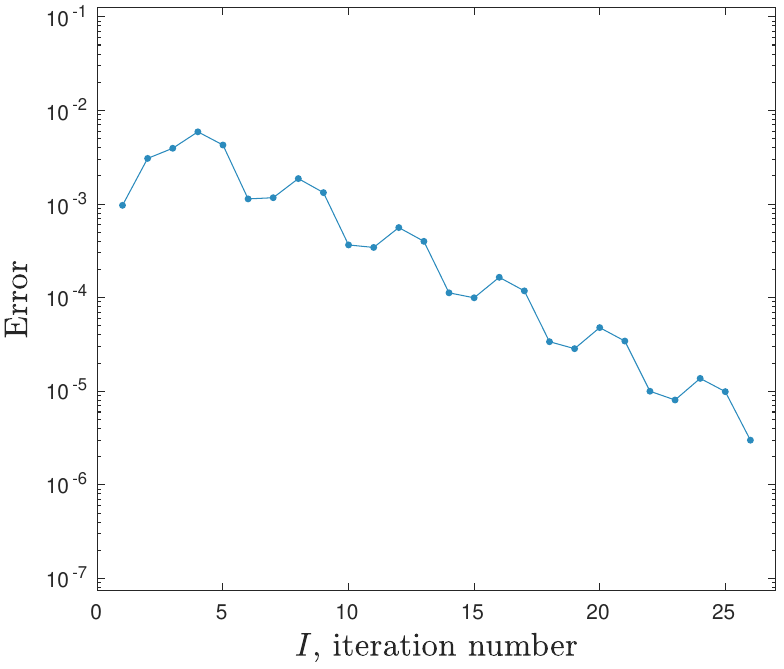}}
  \subfigure[$\norm*{b_{\,2}^{\,(\,I\,)\,}\,(\,t,\,x_{\,0}\,)\ -\ b_{\,2}^{\,(\,I\,-\,1\,)}\,(\,t,\,x_{\,0}\,)}_{\,\infty}$]{\includegraphics[width=0.48\textwidth]{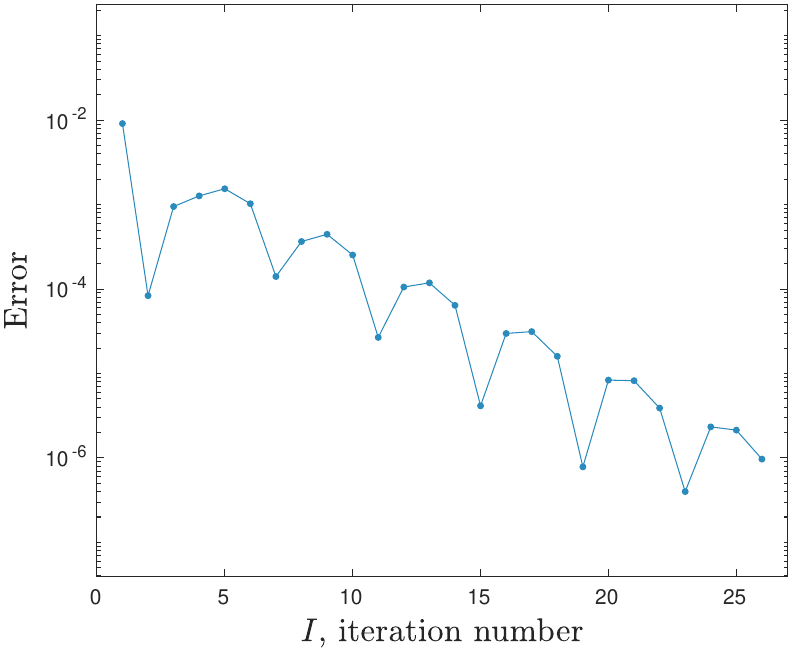}}
  \subfigure[$\norm*{b_{\,3}^{\,(\,I\,)\,}\,(\,t,\,x_{\,0}\,)\ -\ b_{\,3}^{\,(\,I\,-\,1\,)}\,(\,t,\,x_{\,0}\,)}_{\,\infty}$]{\includegraphics[width=0.48\textwidth]{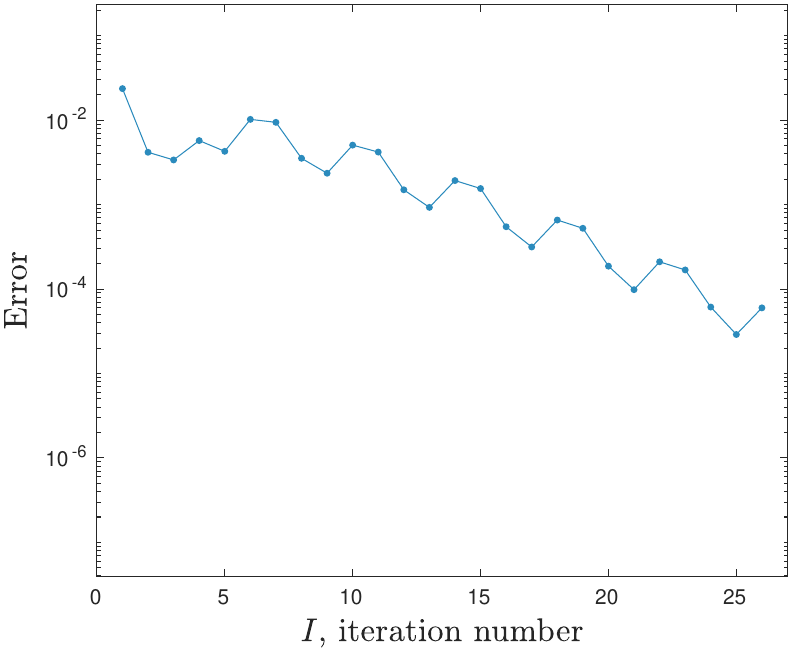}}
  \caption{\small\em The maximum norm $L_{\,\infty}\,\bigl(\,[\,a,\,b\,]\,\bigr)$ of the difference between two successive iterations $(I)$\up{th} and $(I\,-\,1)$\up{th}. The initial data is given in \cref{eq:data3}.}
  \label{fig:err2bis}
\end{figure}

\begin{figure}
  \centering
  \subfigure[$(\,t,\,x\,)\ \mapsto\ \abs*{p_{\,1}^{\,(\,I\,)}\ -\ p_{\,2}^{\,(\,I\,)}}$]{\includegraphics[width=0.48\textwidth]{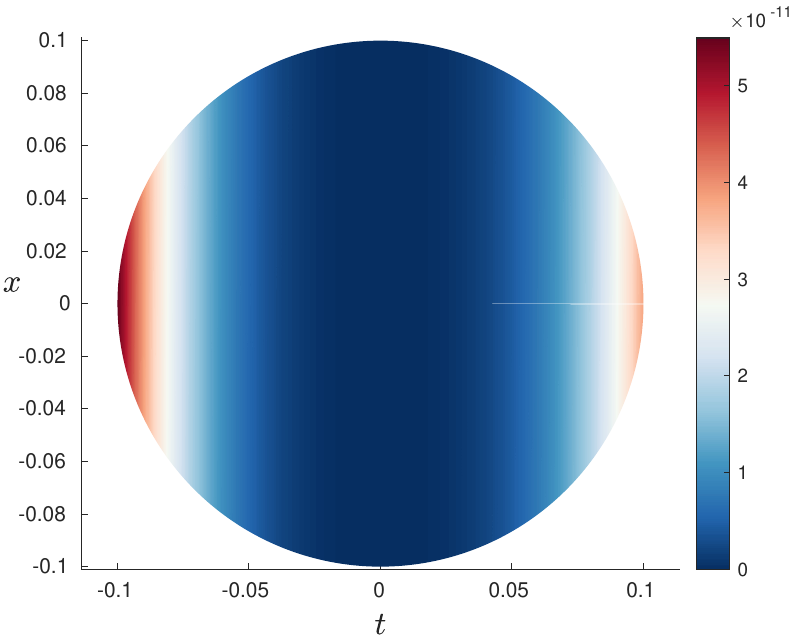}}
  \subfigure[$(\,t,\,x\,)\ \mapsto\ \abs*{q_{\,1}^{\,(\,I\,)}\ -\ q_{\,2}^{\,(\,I\,)}}$]{\includegraphics[width=0.48\textwidth]{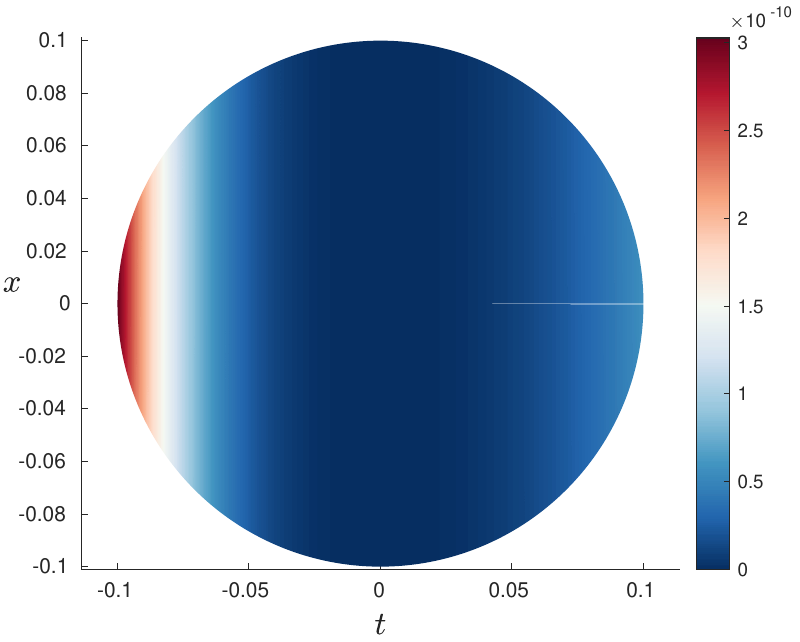}}
  \subfigure[$(\,t,\,x\,)\ \mapsto\ \abs*{b_{\,1,\,0}^{\,(\,I\,)}\ -\ b_{\,2,\,0}^{\,(\,I\,)}}$]{\includegraphics[width=0.48\textwidth]{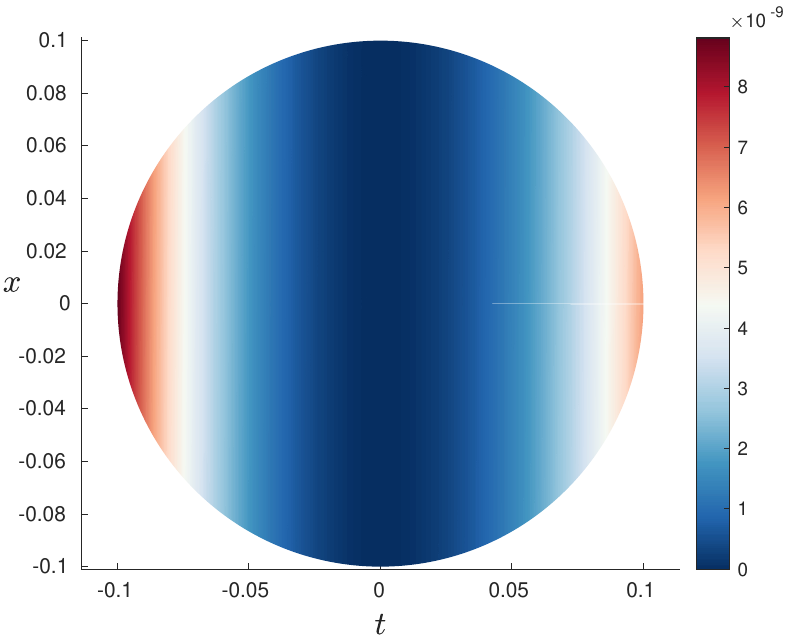}}
  \subfigure[$(\,t,\,x\,)\ \mapsto\ \abs*{b_{\,1,\,1}^{\,(\,I\,)}\ -\ b_{\,2,\,1}^{\,(\,I\,)}}$]{\includegraphics[width=0.48\textwidth]{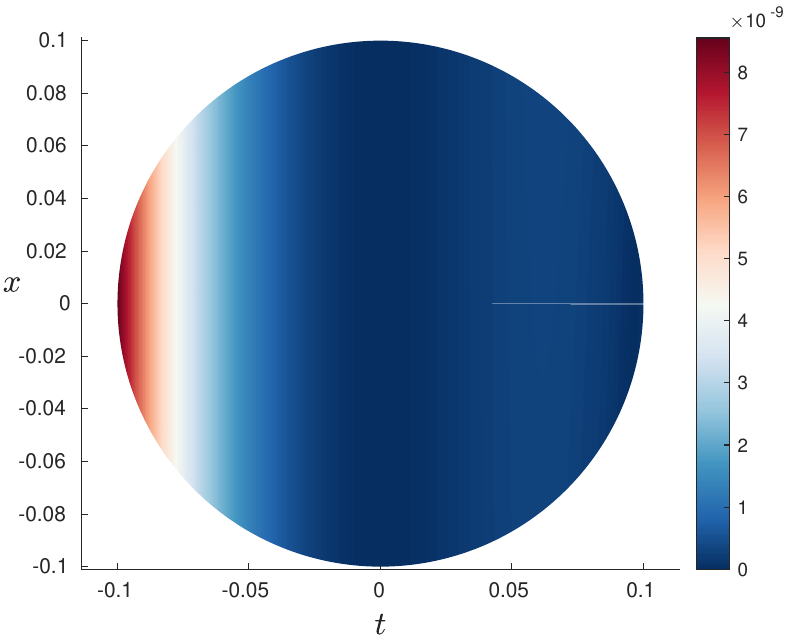}}
  \subfigure[$(\,t,\,x\,)\ \mapsto\ \abs*{b_{\,1,\,2}^{\,(\,I\,)}\ -\ b_{\,2,\,2}^{\,(\,I\,)}}$]{\includegraphics[width=0.48\textwidth]{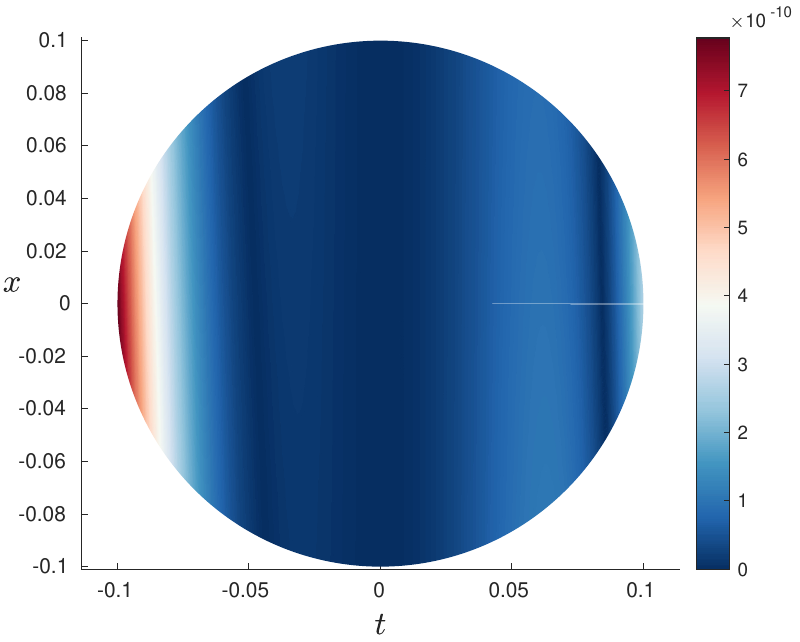}}
  \subfigure[$(\,t,\,x\,)\ \mapsto\ \abs*{b_{\,1,\,3}^{\,(\,I\,)}\ -\ b_{\,2,\,3}^{\,(\,I\,)}}$]{\includegraphics[width=0.48\textwidth]{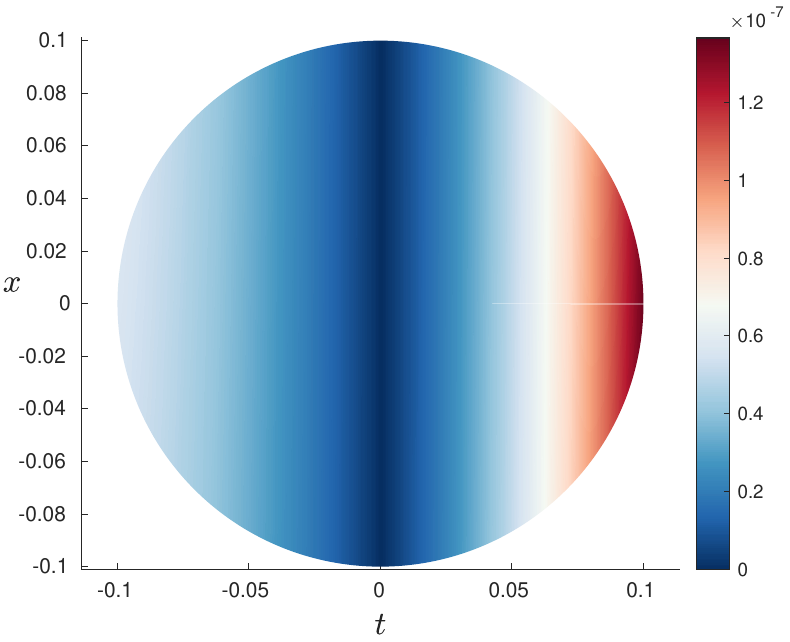}}
  \caption{\small\em The absolute value of the differences between corresponding components of the solution data at $I\ =\ 25$\up{th} iteration. Both iterative processes were initialized with data \eqref{eq:data2} and \eqref{eq:data3} respectively. The values of all numerical parameters are given in \cref{tab:params}.}
  \label{fig:close2}
\end{figure}

As we do it systematically, we also consider another initialization of the fixed point algorithm, which corresponds to the same \nm{Cauchy} datum \eqref{eq:init2} as before:
\begin{subequations}\label{eq:data3}
\begin{align*}
  p\,(\,t,\,x\,)\ &\coloneq\ \frac{t}{2}\,, \qquad
  q\,(\,t,\,x\,)\ \coloneq\ t\ +\ x\,, \\
  b_{\,0}\ &\coloneq\ \fW\,, \qquad b_{\,1}\ \coloneq\ \fW\ -\ \frac{t^{\,2}}{10}\,, \\
  b_{\,2}\,(\,t,\,x\,)\ &\coloneq\ \frac{x}{10}\ +\ \frac{t}{100}\,, \qquad b_{\,3}\,(\,t,\,x\,)\ \coloneq\ \frac{x}{10}\ -\ \frac{t}{100}\,, \\[10pt]
  \qquad b_{\,j}\ &\coloneq\ \fO\,, \qquad \forall\,j\ \geq\ 4\,.
\end{align*}
\end{subequations}
The algorithm was run with this initialization and all other parameters as above. The convergence of the iterations initialized with the initial date \eqref{eq:data3} is reported in \cref{fig:err2bis}. We can observe a very similar behaviour to \cref{fig:final} attesting one more time to the convergence of the proposed algorithm. One may ask the legitimate question of whether two different initializations \eqref{eq:data2} and \eqref{eq:data3} converge to the same point. This theoretical question is very complicated. However, what we can do in practice is to measure how close are the corresponding iterates from each other in the spirit of mental representation from \cref{fig:conv}. This is done in \cref{fig:close2} where the difference $(\,t,\,x\,)\ \mapsto\ \abs{\psi_{\,1}^{\,(\,I\,)}\,(\,t,\,x\,)\ -\ \psi_{\,2}^{\,(\,I\,)}\,(\,t,\,x\,)}$ with $\psi\ \in\ \set*{p,\,q,\,b_{\,0},\,\ldots,\,b_{\,3}}\,$. In order to make the error plot possible, we considered the real values of independent variables varying in some finite neighbourhood of the origin $\vO\ \in\ \R^{\,2}\,$:
\begin{equation*}
  \Set*{(\,t,\,x\,)\ \in\ \R^{\,2}}{t^{\,2}\ +\ x^{\,2}\ \leq\ \frac{1}{100}}\,.
\end{equation*}
In particular, \cref{fig:close2} shows that the convergence is the slowest for the component $b_{\,3}$ of the solution data but even for $b_{\,3}$ the results are perfectly acceptable.

The Section below contains another couple of computational investigations of the proposed fixed point algorithm.


\subsection{Test 3}
\label{sec:test3}

\begin{figure}
  \centering
  \subfigure[$\norm*{p^{\,(\,I\,)\,}\,(\,t,\,x_{\,0}\,)\ -\ p^{\,(\,I\,-\,1\,)}\,(\,t,\,x_{\,0}\,)}_{\,\infty}$]{\includegraphics[width=0.48\textwidth]{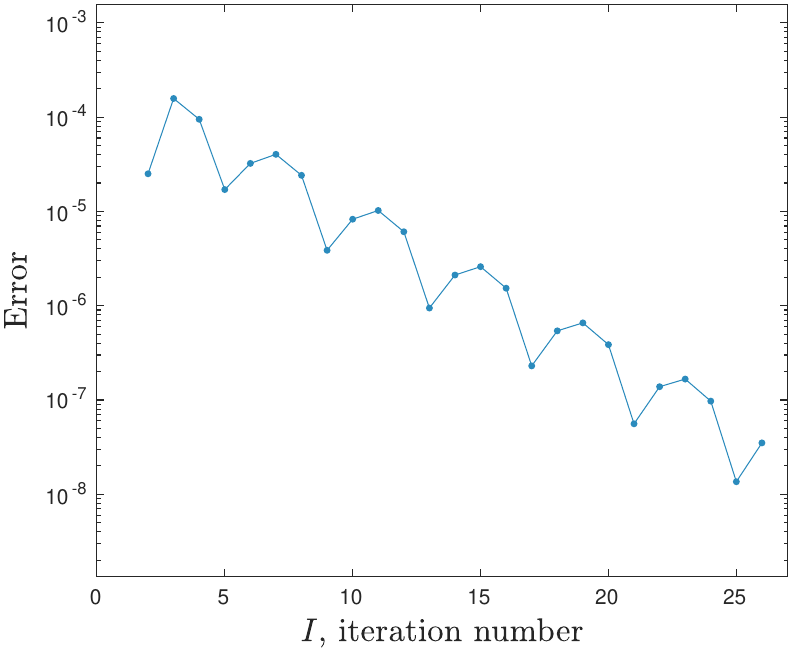}}
  \subfigure[$\norm*{q^{\,(\,I\,)\,}\,(\,t,\,x_{\,0}\,)\ -\ q^{\,(\,I\,-\,1\,)}\,(\,t,\,x_{\,0}\,)}_{\,\infty}$]{\includegraphics[width=0.48\textwidth]{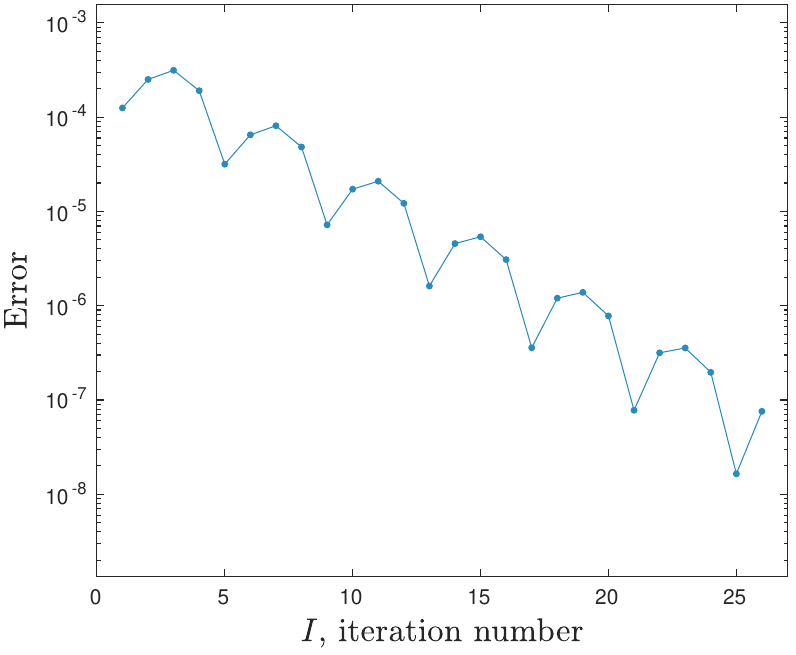}}
  \subfigure[$\norm*{b_{\,0}^{\,(\,I\,)\,}\,(\,t,\,x_{\,0}\,)\ -\ b_{\,0}^{\,(\,I\,-\,1\,)}\,(\,t,\,x_{\,0}\,)}_{\,\infty}$]{\includegraphics[width=0.48\textwidth]{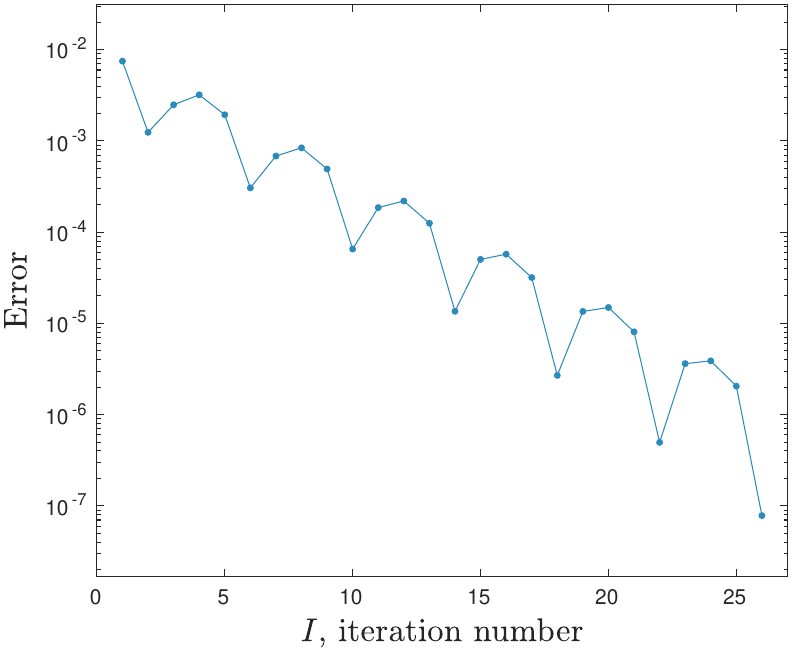}}
  \subfigure[$\norm*{b_{\,1}^{\,(\,I\,)\,}\,(\,t,\,x_{\,0}\,)\ -\ b_{\,1}^{\,(\,I\,-\,1\,)}\,(\,t,\,x_{\,0}\,)}_{\,\infty}$]{\includegraphics[width=0.48\textwidth]{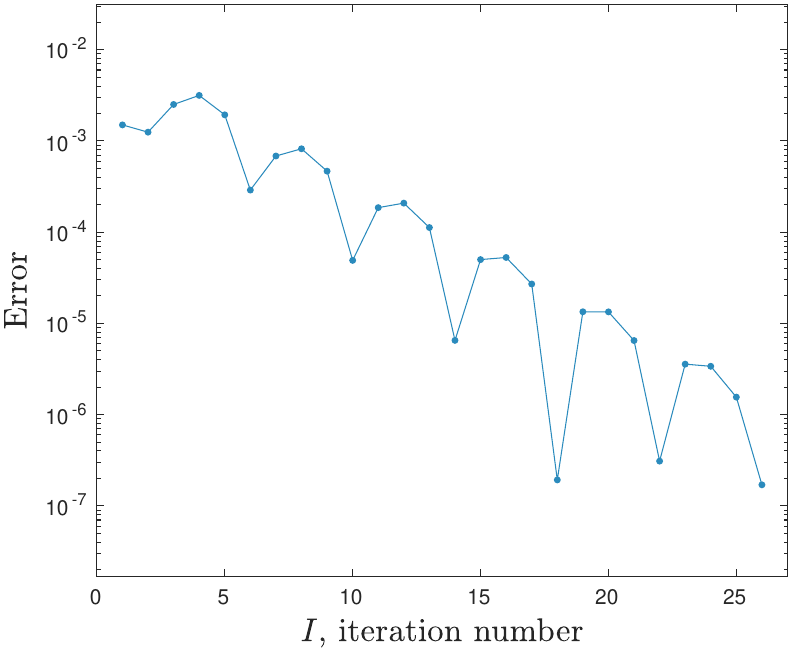}}
  \subfigure[$\norm*{b_{\,2}^{\,(\,I\,)\,}\,(\,t,\,x_{\,0}\,)\ -\ b_{\,2}^{\,(\,I\,-\,1\,)}\,(\,t,\,x_{\,0}\,)}_{\,\infty}$]{\includegraphics[width=0.48\textwidth]{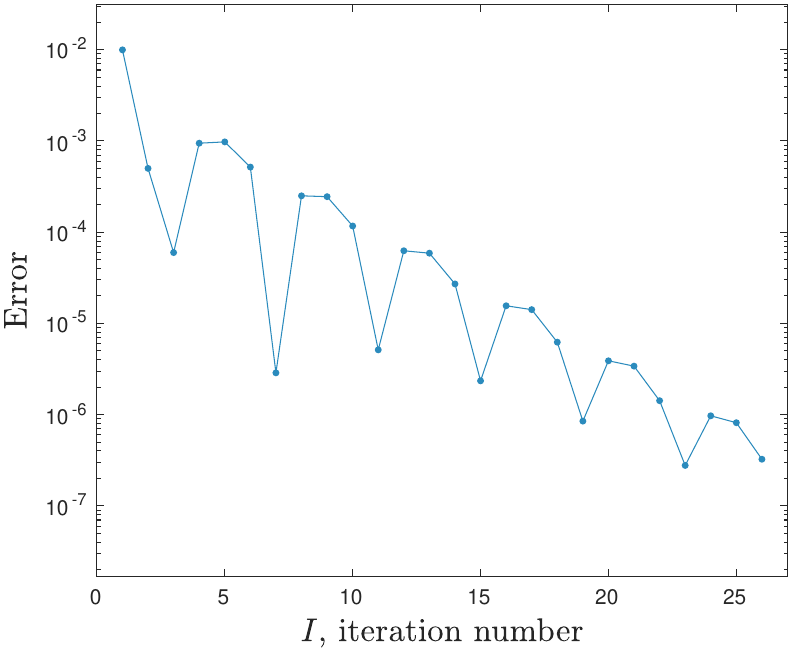}}
  \caption{\small\em The maximum norm $L_{\,\infty}\,\bigl(\,[\,a,\,b\,]\,\bigr)$ of the difference between two successive iterations $(I)$\up{th} and $(I\,-\,1)$\up{th}. The initial data is given in \cref{eq:data5}.}
  \label{fig:err3}
\end{figure}

\begin{figure}
  \centering
  \subfigure[$\norm*{p^{\,(\,I\,)\,}\,(\,t,\,x_{\,0}\,)\ -\ p^{\,(\,I\,-\,1\,)}\,(\,t,\,x_{\,0}\,)}_{\,\infty}$]{\includegraphics[width=0.48\textwidth]{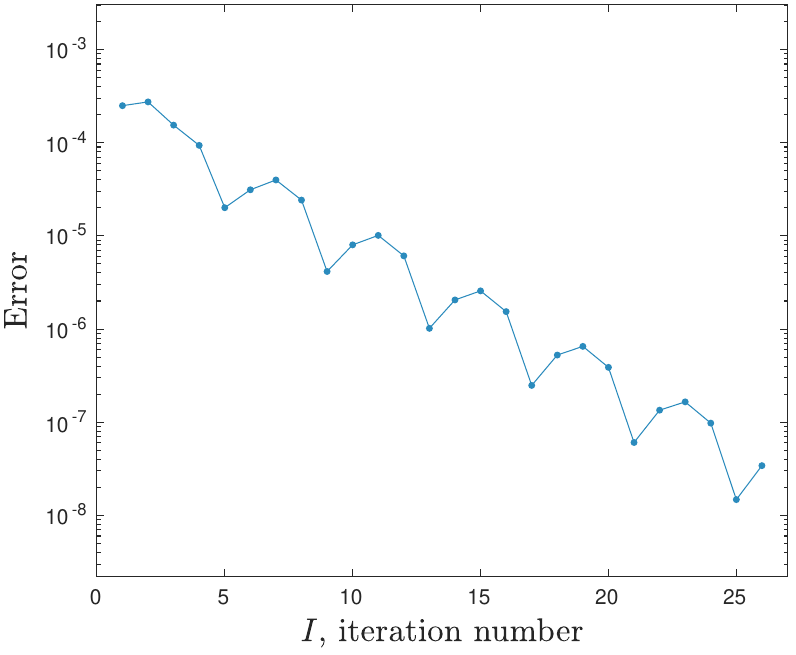}}
  \subfigure[$\norm*{q^{\,(\,I\,)\,}\,(\,t,\,x_{\,0}\,)\ -\ q^{\,(\,I\,-\,1\,)}\,(\,t,\,x_{\,0}\,)}_{\,\infty}$]{\includegraphics[width=0.48\textwidth]{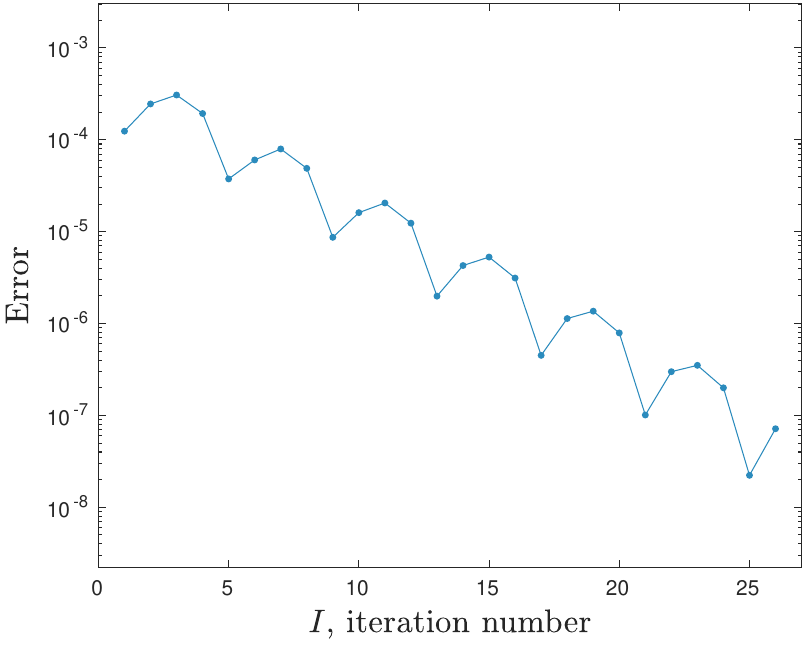}}
  \subfigure[$\norm*{b_{\,0}^{\,(\,I\,)\,}\,(\,t,\,x_{\,0}\,)\ -\ b_{\,0}^{\,(\,I\,-\,1\,)}\,(\,t,\,x_{\,0}\,)}_{\,\infty}$]{\includegraphics[width=0.48\textwidth]{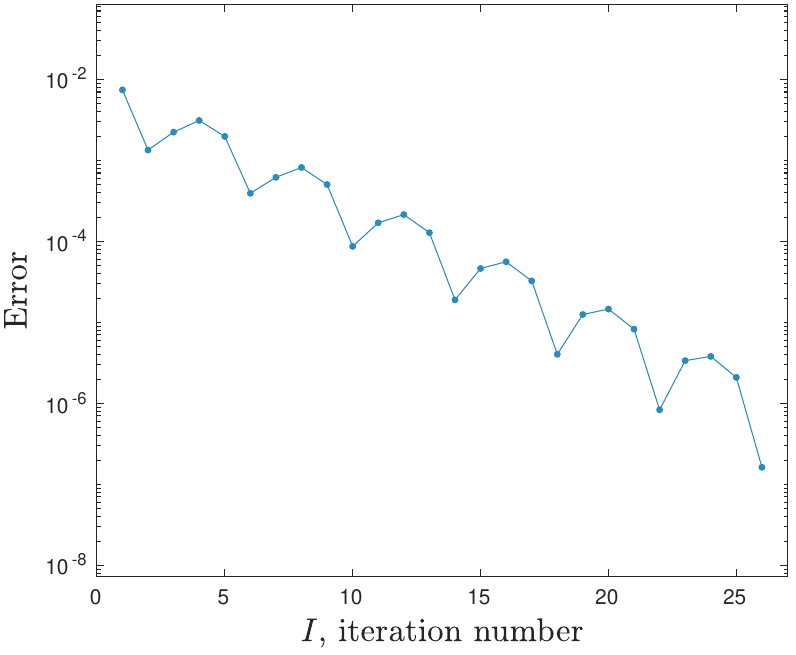}}
  \subfigure[$\norm*{b_{\,1}^{\,(\,I\,)\,}\,(\,t,\,x_{\,0}\,)\ -\ b_{\,1}^{\,(\,I\,-\,1\,)}\,(\,t,\,x_{\,0}\,)}_{\,\infty}$]{\includegraphics[width=0.48\textwidth]{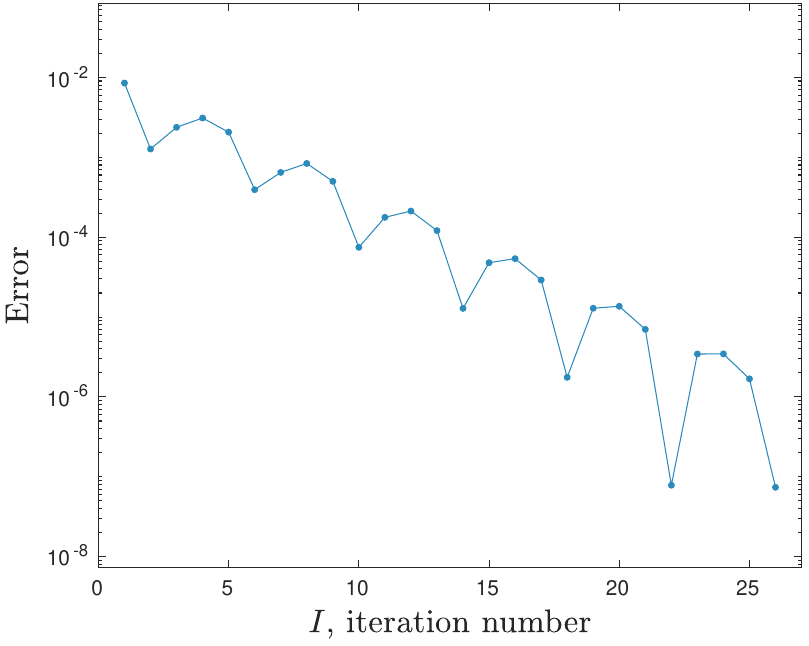}}
  \subfigure[$\norm*{b_{\,2}^{\,(\,I\,)\,}\,(\,t,\,x_{\,0}\,)\ -\ b_{\,2}^{\,(\,I\,-\,1\,)}\,(\,t,\,x_{\,0}\,)}_{\,\infty}$]{\includegraphics[width=0.48\textwidth]{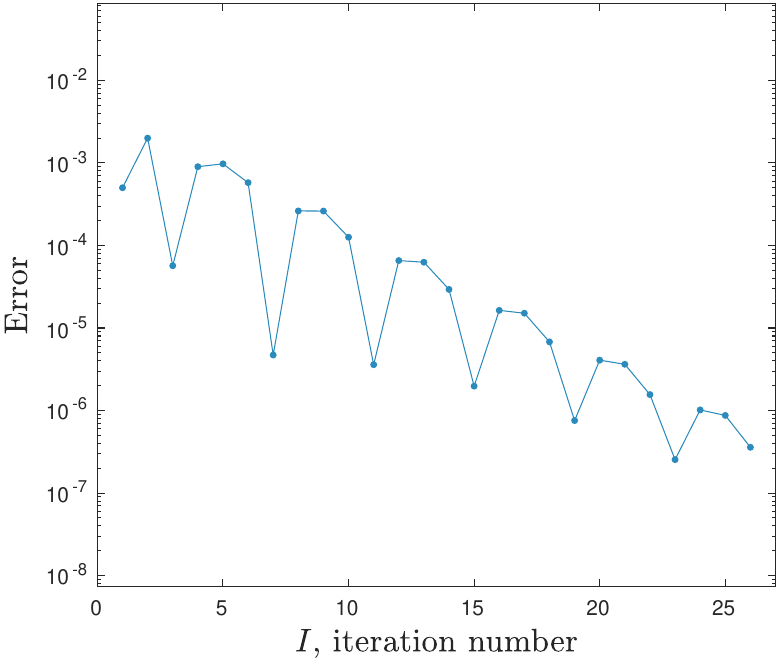}}
  \caption{\small\em The maximum norm $L_{\,\infty}\,\bigl(\,[\,a,\,b\,]\,\bigr)$ of the difference between two successive iterations $(I)$\up{th} and $(I\,-\,1)$\up{th}. The initial data is given in \cref{eq:data6}.}
  \label{fig:err3bis}
\end{figure}

\begin{figure}
  \centering
  \subfigure[$t\ \mapsto\ \abs*{p_{\,1}^{\,(\,I\,)}\ -\ p_{\,2}^{\,(\,I\,)}}$]{\includegraphics[width=0.48\textwidth]{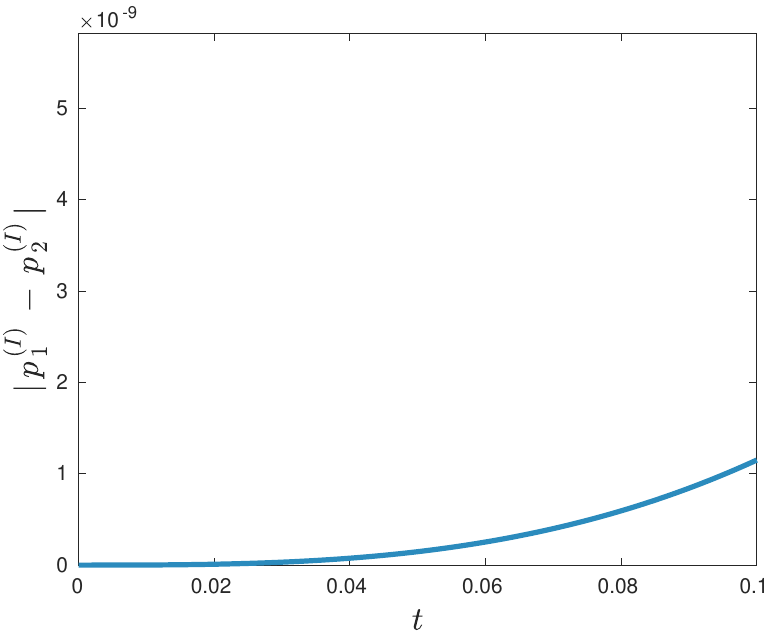}}
  \subfigure[$t\ \mapsto\ \abs*{q_{\,1}^{\,(\,I\,)}\ -\ q_{\,2}^{\,(\,I\,)}}$]{\includegraphics[width=0.50\textwidth]{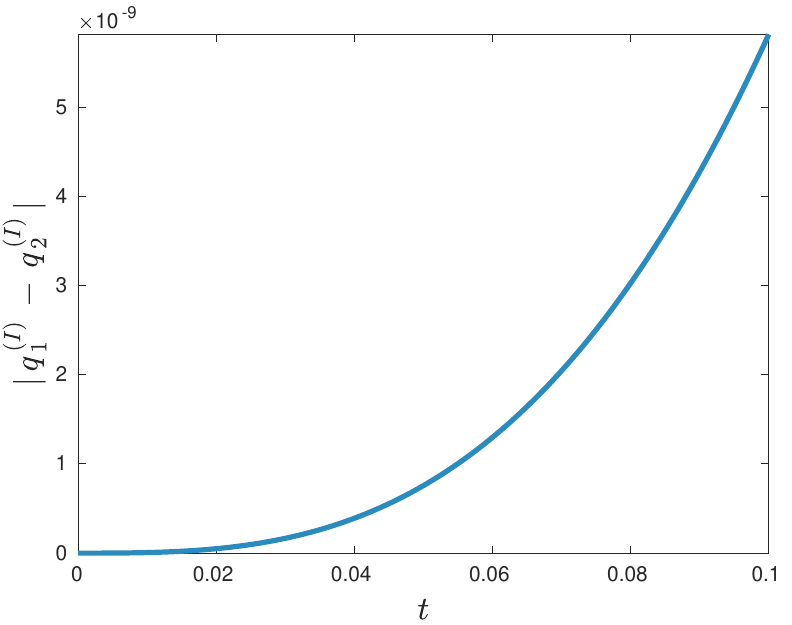}}
  \subfigure[$t\ \mapsto\ \abs*{b_{\,1,\,0}^{\,(\,I\,)}\ -\ b_{\,2,\,0}^{\,(\,I\,)}}$]{\includegraphics[width=0.48\textwidth]{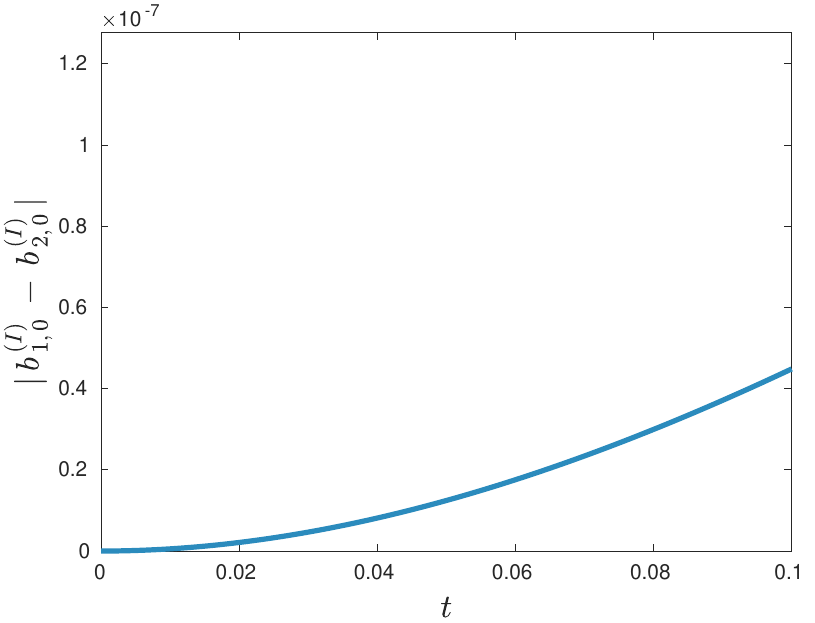}}
  \subfigure[$t\ \mapsto\ \abs*{b_{\,1,\,1}^{\,(\,I\,)}\ -\ b_{\,2,\,1}^{\,(\,I\,)}}$]{\includegraphics[width=0.48\textwidth]{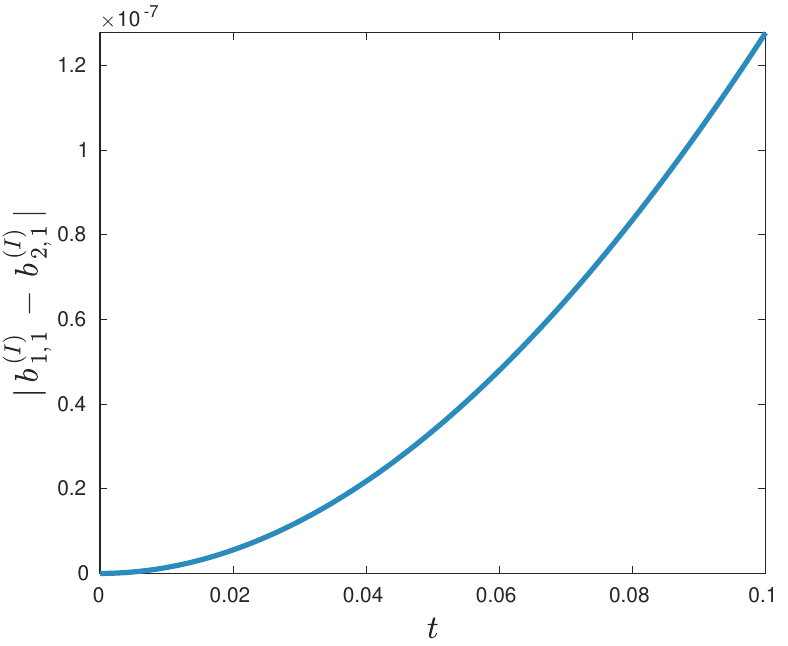}}
  \subfigure[$t\ \mapsto\ \abs*{b_{\,1,\,2}^{\,(\,I\,)}\ -\ b_{\,2,\,2}^{\,(\,I\,)}}$]{\includegraphics[width=0.48\textwidth]{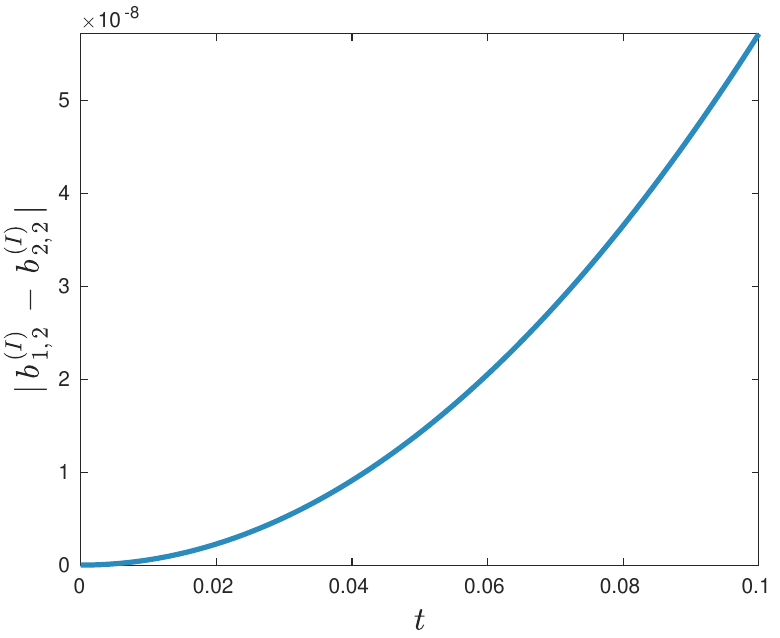}}
  \caption{\small\em Differences between corresponding components of the solution data at $I\ =\ 25$\up{th} iteration. Both iterative processes were initialized with data \eqref{eq:data5} and \eqref{eq:data6} respectively. The values of all numerical parameters are given in \cref{tab:params}.}
  \label{fig:diff3}
\end{figure}

As a final numerical test case, we consider another \acrshort{ivp} \eqref{eq:second}, represented by two different initial guesses, which contains explicitly the fractional power $x^{\,\fourthirds}$ along with $x^{\;\fivethirds}$ in perfect agreement with the generalized \cref{conj:m}. Namely, in the \nm{Cauchy} datum \eqref{eq:psysb} there is a term $c_{\,3}\,x^{\,m\,-\,1\,+\,\frac{2}{3}}$ (with $c_{\,3}\ \in\ \C$), which for $m\ =\ 2$ gives $c_{\,3}\,x^{\;\fivethirds}$ and it was absent in previous numerical Tests $1$ (\cref{sec:test1}) and $2$ (\cref{sec:test2}). As a result, we decided to include an additional numerical test which contains explicitly this term. Namely, we consider the following \nm{Cauchy} datum:
\begin{equation}\label{eq:init3}
  u\,(\,0,\,x\,)\ =\ x\ +\ \frac{3}{4}\;x^{\;\fourthirds}\ +\ \frac{3}{50}\;x^{\,\fivethirds}\,.
\end{equation}
Two slightly different initial data can represent this \nm{Cauchy} datum. The first initial guess reads:
\begin{subequations}\label{eq:data5}
\begin{align}
  p\,(\,t,\,x\,)\ &\coloneq\ \frac{t}{2}\,, \qquad
  q\,(\,t,\,x\,)\ \coloneq\ t\ +\ x\,, \\
  b_{\,0}\ &\coloneq\ \fW\,, \qquad b_{\,1}\ \coloneq\ \fW\,, \\
  b_{\,2}\,(\,t,\,x\,)\ &\coloneq\ \frac{1}{10}\,, \qquad b_{\,j}\ \coloneq\ \fO\,, \qquad \forall\,j\ \geq\ 3\,.
\end{align}
\end{subequations}
The second one is
\begin{subequations}\label{eq:data6}
\begin{align}
  p\,(\,t,\,x\,)\ &\coloneq\ \frac{t}{2}\,, \qquad
  q\,(\,t,\,x\,)\ \coloneq\ t\ +\ x\,, \\
  b_{\,0}\ &\coloneq\ \fW\,, \qquad b_{\,1}\ \coloneq\ \fW\ -\ \frac{t}{10}\,, \\
  b_{\,2}\,(\,t,\,x\,)\ &\coloneq\ \frac{1}{10}\ +\ \frac{t^{\,2}}{5}\,, \qquad b_{\,j}\ \coloneq\ \fO\,, \qquad \forall\,j\ \geq\ 3\,.
\end{align}
\end{subequations}
The convergence of both initializations has been studied symbolically and numerically using the code provided in Appendix~\ref{app:b}. The norms of differences between two successive iterations are reported in Figures~\ref{fig:err3} and \ref{fig:err3bis}. Moreover, \cref{fig:diff3} indicates that both iterations converge to the same element of $\O\,\llbracket\,z\,\rrbracket\,$.

To make an intermediate conclusion, we would like to mention that the Authors tested several other (and even gradually more complicated) configurations, and the proposed algorithm's practical convergence was invariably observed. This constatation gives us a good hope that Conjectures~\ref{conj:2}, \ref{conj:m} are true. More precisely, these numerical tests make quite plausible the following (informal) Conjecture, which, thanks to \nm{Picard} theorem, would imply the \cref{conj:2}:

\begin{minipage}{\textwidth}
\begin{conj}
There exists a \nm{Banach} space $E$ with current point $\Bigl(\,p\,(\,t,\,x\,),\allowbreak\,q\,(\,t,\,x\,),\allowbreak\,\set*{b_{\,k}\,(\,t,\,x\,)}_{\,k\,=\,0}^{\,\infty}\,\Bigr)$ and a closed subset $B\ \subseteq\ E$ satisfying the three following properties:
\begin{enumerate}
  \item The elements of $B$ satisfy for $t\ =\ 0$ the initial condition corresponding to the value $u\,(\,0,\,x\,)\,$.
  \item $\FF_{\,>}\,(\,\set{B}\,)\ \subseteq\ B\,$.
  \item The restriction $\FF\vert_{\,B}$ is a contraction.
\end{enumerate}
\end{conj}
\end{minipage}


\subsection{The singular locus $4\,p^{\,3}\ -\ 27\,q^{\,2}\ =\ 0$'s dependence on the \nm{Cauchy} datum}

Recall that the \nm{Leray} principle states that singularities of a ramified \textbf{linear} \nm{Cauchy} problem are determined by the singularity locus of the \nm{Cauchy} data. In our genuinely nonlinear case, this principle seems not to be valid anymore. 

Above, we solved (approximately) three different ramified \nm{Cauchy} problems \eqref{eq:init1}, \eqref{eq:init2} and \eqref{eq:init3}. The singularities of the solution are determined by functions $p$ and $q\,$. Let us compare them at the end of respective iterative processes:
\begin{align*}
  p^{\,(\,I\,)}_{\,1}\ &\approx\ \half\;t\ +\ 0.0250\,t^{\,2}\ +\ 0.0068\,t^{\,3}\ +\ \ldots \\
  p^{\,(\,I\,)}_{\,2}\ &\approx\ \half\;t\ +\ (\,0.0250\,+\,0.0004\,x\,+\,0.0025\,x^{\,2}\,)\,t^{\,2}\ +\\ 
  & \qquad (\,0.0175\,+\,0.0070\,x\,-\,0.0009\,x^{\,2}\,+\,0.0003\,x^{\,3}\,)\,t^{\,3}\ +\ \ldots \\
  p^{\,(\,I\,)}_{\,3}\ &\approx\ \half\;t\ +\ 0.0113\,t^{\,2}\ +\ 0.0019\,t^{\,3}\ +\ \ldots
\end{align*}
Let us have a look at the corresponding $q$ component approximations:
\begin{align*}
  q^{\,(\,I\,)}_{\,1}\ &\approx\ t\ +\ x\ -\ 0.0333\,t^{\,2}\ -\ 0.0050\,t^{\,3}\ +\ \ldots \\
  q^{\,(\,I\,)}_{\,2}\ &\approx\ t\ +\ x\ +\ (\,-0.0333\,+\,0.0133\,x\,)\,t^{\,2}\ +\\ 
  & \qquad (\,0.0082\,-\,0.0010\,x\,+\,0.0006\,x^{\,2}\,)\,t^{\,3}\ +\ \ldots \\
  q^{\,(\,I\,)}_{\,3}\ &\approx\ t\ +\ x\ -\ 0.0200\,t^{\,2}\ -\ 0.0017\,t^{\,3}\ +\ \ldots
\end{align*}
The lower indices show the test case number ($1\,$, $2$ or $3$ corresponding to Sections \ref{sec:test1}, \ref{sec:test2} and \ref{sec:test3} respectively) and we show only four significant digits for the sake of notation compactness. Just from the visual inspection of expressions provided above, we can see that solution data $p$ and $q$ \textbf{depend} on the \nm{Cauchy} datum. Thus, the \nm{Leray} principle seems to hold in our problem and it is another sign of the genuinely nonlinear problem.


\section{Further generalizations}
\label{sec:gen}

In this Section, we would like to announce the general conjecture that we would like to formulate as a result of our investigations on this topic. Let us denote by $(\,t,\,\x\,)$ a generic point of $\C^{\,n\,+\,1}\,$, where $\x\ =\ (\,x_{\,1},\,x_{\,2},\,\ldots,\,x_{\,n}\,)\ \in\ \C^{\,n}\,$. We fix also a polydisc $D\ \subset\ \C^{\,n\,+\,1}$ centered around the origin $\vO\ \in\ \C^{\,n\,+\,1}$ and we denote by $\O_{\,D}$ the ring of holomorphic functions $A\,(\,t,\,\x\,)$ on $D\,$.

\begin{definition}
Denote by $\Ss_{\,m\,-\,1}$ the vector space of polynomial functions in several variables on $\O_{\,D}$ of the following form:
\begin{equation*}
  F\,(\,t,\,\x,\,v_{\,1,\,\alpha_{\,1}},\,v_{\,2,\,\alpha_{\,2}},\,\ldots,\,,\,v_{\,n,\,\alpha_{\,n}}\,)\ \eqdef\ \sum_{a\ \in\ \ah}\,A_{\,a}\,(\,t,\,\x\,)\,\prod_{\ell\,=\,1}^{n}\,\prod_{k\,=\,1}^{m\,-\,1}\,v_{\,\ell,\,k}^{a\,(\,\ell,\,k)}\,,
\end{equation*}
where $A_{\,a}\ \in\ \O_{\,D}\,$, $\forall\,a\ \in\ \ah$ and the collection of maps $\ah$ is a finite set with elements of the type
\begin{equation*}
  a\,:\ n^{\,\sqsupset}\,\times\,(\,m\,-\,1\,)^{\,\sqsupset}\ \longrightarrow\ \N\,,
\end{equation*}
which do not vanish identically, so that we always have
\begin{equation*}
  F\,(\,t,\,\x,\,0,\,\ldots,\,0\,)\ \equiv\ \fO\,.
\end{equation*}
\end{definition}
In the Definition above, $v_{\,\ell,\,k}$ (respectively, $v_{\,\ell,\,k}^{a\,(\,\ell,\,k)}$) should correspond morally to $\partial_{\,x_{\,\ell}}^{\,k}$ (respectively, to the nonlinear operator $u\ \mapsto\ (\,\partial_{\,x_{\,\ell}}^{\,k}\,u\,)^{\,a\,(\,\ell,\,k)}$).

Let $u\,(\,t,\,\x\,)$ be a germ of a holomorphic function at a point of the polydisc $D\,$. For any $F\ \in\ \Ss_{\,m\,-\,1}$ we shall write:
\begin{equation*}
  F\,\bigl(\,t,\,\x,\,\partial_{\,x_{\,1}}^{\,\alpha_{\,1}}\,u,\,\partial_{\,x_{\,2}}^{\,\alpha_{\,2}}\,u,\,\ldots,\,,\,\partial_{\,x_{\,n}}^{\,\alpha_{\,n}}\,u\,\bigr)\ \eqdef\ \sum_{a\ \in\ \ah}\,A_{\,a}\,(\,t,\,\x\,)\,\prod_{\ell\,=\,1}^{n}\,\prod_{k\,=\,1}^{m\,-\,1}\,(\partial_{\,x_{\,\ell}}^{\,k}\,u)^{\,a\,(\,\ell,\,k)}\,.
\end{equation*}
Notice that in the previous expression, the involved partial derivatives of the germ $u$ are of order at most $m\ -\ 1\,$. Also, as above, we shall use extensively the multi-index notation. If $\beta\ \coloneq\ (\,\beta_{\,1},\,\beta_{\,2},\,\ldots,\,\beta_{\,n}\,)\ \in\ \N^{\,n}\,$, we shall denote $\abs{\beta}\ \bydef\ \beta_{\,1}\ +\ \beta_{\,2}\ +\ \ldots\ \beta_{\,n}$ and $\partial_{\,\x}^{\,\beta}\ \bydef\ \partial_{\,x_{\,1}}^{\,\beta_{\,1}}\comp\partial_{\,x_{\,2}}^{\,\beta_{\,2}}\comp\ldots\comp\partial_{\,x_{\,n}}^{\,\beta_{\,n}}\,$.

We also consider the following scalar quasi-linear differential operator of order $m\,$:
\begin{multline*}
  P\,(\,u\,)\ \eqdef\ \partial_{\,t}^{\,m}\,u\ +\ \sum_{j\,=\,0}^{m\,-\,1}\,\sum_{\substack{\beta\ \in\ \N^{\,n} \\ \abs{\beta}\ =\ m\,-\,j}}\,F_{\,j,\,\beta}\,\bigl(\,t,\,\x,\,\partial_{\,x_{\,1}}^{\,\alpha_{\,1}}\,u,\,\partial_{\,x_{\,2}}^{\,\alpha_{\,2}}\,u,\,\ldots,\,,\,\partial_{\,x_{\,n}}^{\,\alpha_{\,n}}\,u\,\bigr)\,\partial_{\,t}^{\,j}\comp\partial_{\,\x}^{\,\beta}\,u\ +\\ 
  G\,\bigl(\,t,\,\x,\,\partial_{\,x_{\,1}}^{\,\alpha_{\,1}}\,u,\,\partial_{\,x_{\,2}}^{\,\alpha_{\,2}}\,u,\,\ldots,\,,\,\partial_{\,x_{\,n}}^{\,\alpha_{\,n}}\,u\,\bigr)\,,
\end{multline*}
where $F_{\,j,\,\beta}$ and $G$ both belong to the vector space $\Ss_{\,m\,-\,1}\,$. We are finally led to formulate the following
\begin{conj}\label{conj:main}
Consider $p\ \in\ \N^{\,\times}$ and $p\ +\ 1$ holomorphic functions $\set*{\upsilon_{\,j}\,(\,\x\,)}_{j\,=\,0}^{\,p}$ on $D\ \cap\ \bigl(\,\set*{0}\,\times\,\C^{\,n}\,\bigr)$ such that both $\upsilon_{\,0}\,(\,\vO\,)$ and $\upsilon_{\,1}\,(\,\vO\,)$ are not zero. Let
\begin{equation}\label{eq:u0}
  u_{\,0}\,(\,\x\,)\ =\ \sum_{j\,=\,0}^{p}\,\upsilon_{\,j}\,(\,\x\,)\,x_{\,1}^{\,m\,-\,1\,+\,\frac{j}{p\,+\,1}}
\end{equation}
be a germ of a holomorphic function ramified around $x_{\,1}\ =\ 0\,$. We make the following two assumptions:
\begin{enumerate}
  \item The following polynomial equation in variable $\xi$ admits $m$ distinct roots\footnote{This is precisely the place where we see that the initial data enters explicitly into the definition of a suitable differential operator in the context of our study.}:
  \begin{multline}\label{eq:11}
    \xi^{\,m}\ +\\ 
    \sum_{j\,=\,0}^{m\,-\,1}\,\xi^{\,j}\,F_{\,j,\,(\,m-j,\,0,\,\ldots,\,0\,)}\,\bigl(\,0,\,\vO,\,\partial_{\,x_{\,1}}^{\,\alpha_{\,1}}\,u_{\,0}\,(\,\vO\,),\,\partial_{\,x_{\,2}}^{\,\alpha_{\,2}}\,u_{\,0}\,(\,\vO\,),\,\ldots,\,\partial_{\,x_{\,n}}^{\,\alpha_{\,n}}\,u_{\,0}\,(\,\vO\,)\,\bigr)\ =\ 0\,.
  \end{multline}
  \item $\partial_{\,v_{\,1,\,m\,-\,1}}\,F_{\,0,\,(\,m,\,0,\,\ldots,\,0\,)}\,\bigl(\,0,\,\vO,\,\partial_{\,x_{\,1}}^{\,\alpha_{\,1}}\,u_{\,0}\,(\,\vO\,),\,\partial_{\,x_{\,2}}^{\,\alpha_{\,2}}\,u_{\,0}\,(\,\vO\,),\,\ldots,\,\partial_{\,x_{\,n}}^{\,\alpha_{\,n}}\,u_{\,0}\,(\,\vO\,)\,\bigr)\ \neq\ 0\,.$
\end{enumerate}
Then, the following statements are true. For each choice of a root to the algebraic \cref{eq:11}, there exist holomorphic functions $\set*{q_{\,\ell}\,(\,t,\,\x\,)}_{\,\ell\,=\,0}^{\,p\,-\,1}$ defined in a neighbourhood of the origin $\vO\ \in\ \C^{\,n\,+\,1}$ such that
\begin{equation*}
  q_{\,0}\,(\,0,\,\x\,)\ =\ x_{\,1}\,, \qquad
  q_{\,\ell}\,(\,0,\,\x\,)\ =\ \fO\,, \qquad
  \ell\ \in\ (\,p\,-\,1\,)^{\,\sqsupset}\,.
\end{equation*}
We denote by $z\,(\,t,\,\x\,)$ a holomorphic germ satisfying the following algebraic equation outside the swallowtail singularity:
\begin{equation*}
  z^{\,p\,+\,1}\,(\,t,\,\x\,)\ -\ q_{\,p\,-\,1}\,(\,t,\,\x\,)\,z^{\,p\,-\,1}\,(\,t,\,\x\,)\ -\ \ldots\ q_{\,0}\,(\,t,\,\x\,)\ =\ \fO\,.
\end{equation*}
There exist also the holomorphic functions $\set*{a_{\,\ell}\,(\,t,\,\x\,)}_{\,\ell\,=\,0}^{\,p}$ defined on a neighbourhood of the origin $\vO\ \in\ \C^{\,n\,+\,1}$ such that the function
\begin{equation*}
  u\,(\,t,\,\x\,)\ =\ \partial_{\,q_{\,0}}^{\,-\,m\,+\,1}\bigl(\,a_{\,p}\,(\,t,\,\x\,)\,z^{\,p}\,(\,t,\,\x\,)\ +\ \ldots\ +\ a_{\,1}\,(\,t,\,\x\,)\,z\,(\,t,\,\x\,)\ +\ a_{\,0}\,(\,t,\,\x\,)\,\bigr)
\end{equation*}
satisfies:
\begin{align*}
  P\,(\,u\,)\,(\,t,\,\x\,)\ &\equiv\ 0\,, \\
  u\,(\,0,\,\x\,)\ &\equiv\ u_{\,0}\,(\,\x\,)\,,
\end{align*}
where $u_{\,0}\,(\,\x\,)$ was defined in \eqref{eq:u0}.
\end{conj}

Let us now discuss the motivation behind this general \cref{conj:main}. The $p-$tuple $(\,z,\,z^{\,2},\,\ldots,\,z^{\,p}\,)$ is the solution of a holonomic $\D-$module whose characteristic variety $\Vv$ is included in the union of the zero section of $\T^{\,\ast}\,\C^{\,p}$ and of the co-normal to the swallow tail singularity associated to the algebraic equation:
\begin{equation*}
  z^{\,p\,+\,1}\ -\ q_{\,p\,-\,1}\,z^{\,p\,-\,1}\ -\ \ldots\ q_{\,1}\,z\ -\ q_{\,0}\ =\ \fO\,.
\end{equation*}
The fact that this singularity is \emph{stable} plays the crucial r\^ole in our constructions. Moreover, $\T_{\,0}^{\,\ast}\,\C^{\,p}\ \bigcap\ \Vv$ coincides with the line co-normal to $q_{\,0}\ =\ 0$ and $q_{\,0}\,(\,0,\,\x\,)\ =\ x_{\,1}\,$. These important facts explain why we consider only the contributions of the monomials $\partial_{\,t}^{\,j}\comp\partial_{\,x_{\,1}}^{\,m\,-\,j}$ in \cref{eq:11}.

The assumption (1) means that all the characteristics associated with this geometry are simple. The proof of \cref{conj:main} should allow to construct $m$ solutions of this type, each corresponding to a choice of a root in \cref{eq:11}. It is quite clear that holomorphic functions $\set*{q_{\,\ell}\,(\,t,\,\x\,)}_{\,\ell\,=\,0}^{\,p\,-\,1}$ highly depend on the initial datum $u_{\,0}\,(\,\x\,)$ and, even more precisely, on $\set*{\upsilon_{\,j}\,(\,\x\,)}_{\,j\,=\,0}^{\,p}\,$. The assumption (2) means that the \acrshort{pde} $P\,(\,u\,)\ =\ \fO$ is genuinely nonlinear with respect to this geometry. The simplest example of such an operator $P\,(\cdot)$ from the \cref{conj:main} statement is given by the following expression:
\begin{equation*}
  P\,(\,u\,)\ \coloneq\ \partial_{\,t}^{\,m}\,u\ +\ \bigl(\,\partial_{\,x_{\,1}}^{\,m\,-\,1}\,u\,\bigr)^{\,n}\,\partial_{\,x_{\,1}}^{\,m}\,u\,, \qquad n\ \in\ \N_{\,\geq\,1}\,.
\end{equation*}
To make things clearer, let us emphasize that, for instance, the two following quasi-linear operators do not satisfy the assumption (2) of \cref{conj:main}:
\begin{itemize}
  \item $P_{\,1}\,(\,u\,)\ =\ \partial_{\,t}^{\,m}\,u\ +\ \bigl(\,\partial_{\,x_{\,1}}^{\,m\,-\,2}\,u\,\bigr)^{\,n}\,\partial_{\,x_{\,1}}^{\,m}\,u\,$,
  \item $P_{\,2}\,(\,u\,)\ =\ \partial_{\,t}^{\,m}\,u\ +\ \bigl(\,\partial_{\,x_{\,1}}^{\,m}\,u\,\bigr)^{\,n}\,\partial_{\,x_{\,2}}^{\,m}\,u\,$.
\end{itemize}

We summarize below the proposed strategy to attack this problem:
\begin{enumerate}
  \item We seek a solution of $P\,(\,u\,)\ =\ \fO$ under the form
  \begin{equation*}
    u\,(\,t,\,\x\,)\ =\ \partial_{\,q_{\,0}}^{\,-\,m\,+\,1}\,\Bigl(\,\sum_{j\,=\,0}^{+\,\infty}\,b_{\,j}\,(\,t,\,\x\,)\,z^{\,j}\,\Bigr)\,.
  \end{equation*}
  Since it does not seem to be possible to construct a solution directly, we proceed as in \cref{sec:sec}, \ie we try to find an algorithm with an underlying map $\FF$ such that a fixed point of $\FF$ gives a solution. The fact that the swallow tail singularity should be characteristic for the \acrshort{pde} will lead to the eikonal equation. So, in some sense, we use an ansatz and the method of the geometric optics\footnote{This method is also known as the \nm{Hamilton}--\nm{Jacobi} theory.}. A heuristic (but crucial) underlying idea is that the \nm{Hamiltonian} flow of the principal symbol of the linearized operator should propagate the singularities along the co-normal of the swallow tail.
  \item Find suitable semi-norms allowing to construct a \nm{Banach} space so that one might apply the fixed point theorem to the mapping $\FF\,$.
\end{enumerate}
We hope to motivate and stimulate the research in this direction. Of course, this programme is extremely difficult. The best rational reasons to believe that it should work are the validity of \cref{thm:2} (the case of algebraic equations of the second degree) and the numerical convergence results for the algorithm of \cref{sec:sec}.


\section{Conclusions and perspectives}
\label{sec:concl}

Above, we presented the main results of the present manuscript. The main conclusions and perspectives of this study are outlined below.


\subsection{Conclusions}

In this article, we considered several ramified \acrfull{cp}s for the first and second-order genuinely nonlinear \acrshort{pde}s in the full complex setting. The term \acrshort{cp} should be understood in the context because we specify only one initial condition even for \acrshort{pde}s \eqref{eq:fam} with $m\ >\ 1\,$. We succeeded in understanding completely the case of the \acrshort{ibe} \eqref{eq:burg} using the methods of \nm{Cauchy}--\nm{Kovalevskaya} and those of the contact geometry. Both theories give the same singularities for the \acrshort{ibe} solutions to the \acrshort{ivp} \eqref{eq:cauchy2}. Then, we switched to the second-order \acrshort{ivp} \eqref{eq:second} for the genuinely nonlinear \acrshort{pde} \eqref{eq:sec}. However, this problem cannot be addressed with the same methods. This observation led us to consider some holonomic $\D-$modules and their geometry, which gave us the solution ansatz. The closed-form expression is not possible even within the proposed ansatz. Consequently, we devised an iterative method whose fixed point gives the desired solution. Unfortunately, we were not able to prove rigorously its convergence. Nevertheless, we performed the practical computations with the celebrated computer algebra software \texttt{Maple}${}^{\text{\texttrademark}}$ in order to check in practice the behaviour of the proposed algorithm. Indeed, our symbolic computations clearly indicate the convergence of the proposed method. The main part of the employed code is provided in Appendix~\ref{app:b}. On the basis of this empirical evidence, we were able to formulate several new conjectures which remain open for the moment, the main one being the \cref{conj:main}.


\subsection{Perspectives}

We already mentioned in \cref{rem:wf} the need to study the micro-local singularity of distributions $z$ and $z^{\,2}\,$. The Authors are currently working to establish a rigorous convergence proof for the iterative scheme proposed in Section~\ref{sec:sec}. The underlying theoretical setting was briefly mentioned in \cref{sec:hope}. This result will allow us to establish the existence and uniqueness result for the considered \acrshort{ivp} to Equation~\eqref{eq:sec}. For the moment, the appropriate \nm{Banach} spaces seem to have been identified (see \cref{sec:hope}). Of course, the application of these methods to other fully nonlinear and even higher order \acrshort{pde}s (such as the family of \acrshort{pde}s \eqref{eq:fam} with $m\ >\ 2$) is to be expected. Finally, we would like to see the connections between our theory and the more classical theory of shock waves for $p-$systems as presented in \cite[Chapter~12]{Serre2000}. We shall investigate the interplay and connections between both theories. As the first glance at possible connections, we invite the reader to consult \cref{app:a}. We would like to explore other types of singularities as a separate research direction. According to the classification of catastrophes \cite{Arnold1985}, the cusp and swallowtail belong to the class $A_{\,n}$ for some $n\ \in\ \N\,$. We believe that it would be interesting to explore other types of singularities $D_{\,n}$ and $E_{\,n}$ and the solutions they generate in the current genuinely nonlinear \acrshort{pde} framework.


\subsection*{Acknowledgments}
\addcontentsline{toc}{subsection}{Acknowledgments}

E.~\nm{Leichtnam} would like to acknowledge the hospitality of LAMA UMR \#5127 (University \nm{Savoie Mont Blanc}), where this work took a physical place. Both authors would like to thank F.~\nm{Golse} (CMLS, \'Ecole Polytechnique) and D.~\nm{Serre} (UMPA, ENS de \nm{Lyon}) for very helpful discussions. The second Author (EL) thanks G.~\nm{Lebeau} for fruitful discussions about this topic some time ago (during the mid-eighties).


\subsection*{Funding}
\addcontentsline{toc}{subsection}{Funding}

This publication is based upon work supported by the Khalifa University of Science and Technology under Award No. FSU-2023-014.


\subsection*{Author contributions}
\addcontentsline{toc}{subsection}{Author contributions}

All the Authors contributed equally to this work.


\subsection*{Conflict of interest statement}
\addcontentsline{toc}{subsection}{Conflict of interest statement}

The Authors certify that they have \textbf{no} affiliations with or involvement in any organization or entity with any financial interest (such as honoraria; educational grants; participation in speakers' bureaus; membership, employment, consultancies, stock ownership, or other equity interest; and expert testimony or patent-licensing arrangements), or non-financial interest (such as personal or professional relationships, affiliations, knowledge or beliefs) in the subject matter or materials discussed in this manuscript.


\subsection*{Disclaimer statement}
\addcontentsline{toc}{subsection}{Disclaimer statement}

The opinions expressed in this manuscript are those of the Authors and do not necessarily reflect the views of their employers or any other affiliated organizations.



\bigskip\bigskip
\invisiblesection{References}
\bibliographystyle{acm}
\bibliography{mybiblio}
\bigskip\bigskip


\appendix
\section{Nomenclature}

\begin{description}
  \item[$f\,:\ A\ \longrightarrow\ B$\ ] A map from the domain $\dom\,f\ =\ A$ to the co-domain (or target) $\cod\,f\ =\ B\,$.
  \item[$f_{\,>}$\ ] Image of a map. Consider a map $f\,:\ A\ \longrightarrow\ B\,$. Then, the image $f_{\,>}\,:\ \powerset\,(\,A\,)\ \longrightarrow\ \powerset\,(\,B\,)$ is defined in the usual way, where $\powerset\,(\cdot)$ denotes the power set.
  \item[$f^{\,<}$\ ] Pre-image of a map. Consider a map $f\,:\ A\ \longrightarrow\ B\,$. Then, the pre-image $f^{\,<}\,:\ \powerset\,(\,B\,)\ \longrightarrow\ \powerset\,(\,A\,)$ is defined in the usual way, where $\powerset\,(\cdot)$ denotes the power set.
  \item[$\comp$\ ] Composition operation for functions of differential operators.
  \item[$f^{m\,\comp}\,(\,-\,)$\ ] The $m$\up{th} iterate of a mapping $f\,:\ A\ \longrightarrow\ A$ for some $m\ \in\ \N^{\,\times}\,$, \ie $f^{m\,\comp}\,(\,-\,)\ \eqdef\ \underbrace{(\,f\comp f\comp\ldots\comp f\,)}_{m \text{times}}\,(\,-\,)$
  \item[$\ud\,f$\ ] The differential of a smooth function $f\,$.
  \item[$\jj^{\,1}\,(\,u\,)$\ ] $1-$jet of a smooth function $u\,$.
  \item[$\N$\ ] The set of natural numbers starting from $0\,$.
  \item[$\N^{\,\times}$\ ] The set of strictly positive natural numbers.
  \item[$\R$\ ] The set of real numbers.
  \item[$\R^{\,\times}$\ ] The multiplicative group of real numbers.
  \item[$\R_{\,\geq\,0}$\ ] The set of non-negative real numbers.
  \item[$\C$\ ] The set of complex numbers.
  \item[$\C^{\,\times}$\ ] The multiplicative group of complex numbers.
  \item[$\abs{\cdot}$\ ] The absolute value of a real of complex number.
  \item[$\norm{\cdot}$\ ] The norm function on a \nm{Banach} space.
  \item[$\O_{\,D}$\ ] The ring of holomorphic functions defined on the polydisc $D\,$.
  \item[$\O\,{[\,z\,]}$\ ] The ring of polynomials in the formal variable $z$ with coefficients being holomorphic germs at the origin.
  \item[$\O\,\llbracket\,z\,\rrbracket$\ ] The ring of formal power series in the formal variable $z$ with coefficients being holomorphic germs at the origin.
  \item[$n^{\,\sqsupset}$\ ] Finite set $\set{1,\,2,\,\ldots,\,n}\,$, $\forall\,n\ \in\ \N_{\,\geq\,1}\,$.
  \item[$n^{\,\sqsubset}$\ ] Finite set $\set{0,\,2,\,\ldots,\,n\,-\,1}\,$, $\forall\,n\ \in\ \N_{\,\geq\,1}\,$.
  \item[$\underline{r}$\ ] A function taking the constant value $r\ \in\ \R\,$. The domain and co-domain of this function should be clear from the context (usually, from the left-hand side of the equation, where it appears).
  \item[$\underline{\bm{r}}$\ ] A function taking the constant vector value $(\,r,\,r,\,\ldots,\,r\,)\ \in\ \R^{\,m}$ with $r\ \in\ \R$ and $m\ \in\ \N_{\,\geq\,2}\,$. The domain and co-domain of this function should be clear from the context (usually, from the left-hand side of the equation, where it appears).
  \item[$(\,-\,)^{\,\top}$\ ] The transposition operator acting on a vector or on a linear operator. Strictly speaking, the result of this operation belongs to the dual vector space.
  \item[$\x\ $] Element of the vector space $\x\ \eqdef\ (\,x_{1},\,x_{2},\,\ldots,\,x_{n}\,)\ \in\ \C^{\,n}\,$.
  \item[$\partial_{\,x_{\,j}}^{\,m}$\ ] Partial derivative operator of the order $m$ with respect to the independent variable $x_{\,j}\,$. Sometimes, we also use the index notation $(\cdot)_{\,x_{\,i}}$ to denote the first-order partial derivative.
  \item[$E\ $] A \nm{Banach} space.
  \item[$\FF\ $] Fixed point mapping.
  \item[$\FFt\ $] Finitary truncation of the fixed point mapping.
  \item[$\psi^{\,(\,I\,)}\ $] Function $\psi$ on the $(I)$\up{th} iteration.
  \item[$\coloneq$\ ] Assignment of the right-hand side to the left-hand side.
  \item[$\eqdef$\ ] The left-hand side is defined by the right-hand side.
  \item[$\defeq$\ ] The right-hand side is defined by the left-hand side.
  \item[$\bydef$\ ] By definition.
  \item[$\equiv$\ ] Equal identically.
  \item[$\cong$\ ] Isomorphic.
  \item[$\kk^{\,n}$\ ] The standard $n-$dimensional vector space with the base field $\kk\,$.
  \item[$\U_{\,\z}$\ ] Some open neighbourhood of a point $\z\ \in\ \C^{\,m}\,$.
  \item[$\xib$\ ] Bold version of the Greek letter $\xi\,$.
  \item[$Z\,(\,G\,)$\ ] The center of the group $G\,$.
  \item[$U_{\,n}$\ ] The multiplicative group of the roots of unity of degree $n\,$.
  \item[$S_{\,n}$\ ] The symmetric group with $n$ symbols.
  \item[$\GL_{n}\,(\,\kk\,)$\ ] The general linear group of the square $n\times n$ matrices with entries belonging to the field $\kk\,$.
  \item[$\SL_{n}\,(\,\kk\,)$\ ] The special linear group of the square $n\times n$ matrices with entries belonging to the field $\kk\,$.
  \item[$\PSL_{n}\,(\,\kk\,)$\ ] The projective special linear group defined as $\SL_{n}\,(\,\kk\,)\,/\,Z\,\bigl(\,\SL_{n}\,(\,\kk\,)\,\bigr)\,$. For $n\ =\ 2$ and $\kk\ =\ \Z$ it is called the modular group.
  \item[$\Br_{n}$\ ] The braid group with $n$ strands.
  \item[$\M$\ ] Holomorphic, smooth or topological manifold, depending on the context.
  \item[$\pi_{\,1}\,(\,\M\,)$\ ] The fundamental (or the first homotopy) group of a topological manifold $\M\,$.
  \item[$\T_{\,x}\,\M$\ ] Tangent space to a smooth manifold $\M$ at point $x\ \in\ \M\,$.
  \item[$\omega$\ ] The standard symplectic $2-$form on $\T_{\,x}\,\M\,$.
  \item[$\pb*{-}{-}$\ ] The standard \nm{Poisson} bracket on $C^{\,\infty}\,(\,\M\,)\,$.
  \item[$\T\,\M$\ ] Tangent bundle of a smooth manifold $\M\,$.
  \item[$\T_{\,x}^{\,\ast}\,\M$\ ] Co-tangent space to a smooth manifold $\M$ at point $x\ \in\ \M\,$.
  \item[$\T^{\,\ast}\,\M$\ ] Co-tangent bundle of a smooth manifold $\M\,$.
  \item[$\Nn\,(\,\M\,)$\ ] Co-normal to a smooth sub-manifold $\M\,$.
  \item[$\singsupp\,(\,\cdot\,)$\ ] The singular support of a distribution.
  \item[$\WF\,(\,\cdot\,)$\ ] The wave front set of a distribution.
  \item[$\Xx\ $] The complex manifold.
  \item[$\set*{\Xx_{\,\lambda}}_{\,\lambda\,\in\,\Lambda}\ $] The \nm{Whitney} stratification of the complex manifold $\Xx\,$.
  \item[$\D_{\,\Xx}\ $] The sheaf of differential operators of \emph{finite order} on the complex manifold $\Xx\,$.
  \item[$\mathfrak{M}\ $] The coherent $\D_{\,\Xx}-$module.
  \item[$\O_{\,\Xx}\ $] The sheaf of holomorphic functions on $\Xx\,$, which is a left coherent $\D_{\,\Xx}-$module.
  \item[$\mathsf{Ext}\ $] The right derived functor of the hom-functor.
  \item[$\Gal\,(\,-\,)\ $] \nm{Galois} group of a field extension.
\end{description}


\section{Details of computations}
\label{app:0}

In this Appendix, we explain more precisely how we construct the iterations from \cref{sec:it}. Namely, we specify how the expression $(\,\widehat{\maltese}\,)$ is constructed from $(\,\maltese\,)\,$.

Recall, that the solution ansatz is given by
\begin{equation*}
  u\,(\,t,\,x\,)\ =\ \sum_{k\,=\,1}^{\infty}\,\frac{1}{k}\;\Bigl(\,-\,p\,(\,t,\,x\,)\,b_{\,k-1}\,(\,t,\,x\,)\ +\ 3\,b_{\,k-3}\,(\,t,\,x\,)\,\Bigr)\,z^{\,k}\,.
\end{equation*}
One checks by direct computation that
\begin{multline}\label{eq:one}
  u_{\,t\,t}\ =\ \sum_{k\,=\,0}^{+\,\infty}\,\biggl[\,\Bigl(\,2\ -\ \frac{2}{k}\,\Bigr)\,b_{\,k\,-\,1,\,t}\,p_{\,t}\ +\ \Bigl(\,1\ -\ \frac{1}{k}\,\Bigr)\,b_{\,k\,-\,1}\,p_{\,t\,t}\ +\ b_{\,k}\,q_{\,t\,t}\ +\ \underline{2\,b_{\,k,\,t}\,q_{\,t}}\ +\\ 
  \frac{1}{k}\;\bigl(\-\,p\,b_{\,k\,-\,1,\,t\,t}\ +\ 3\,b_{\,k\,-\,3,\,t\,t}\,\bigr)\,\biggr]\,z^{\,k}\ +\\ 
  \frac{1}{3\,z^{\,2}\,-\,p}\;\sum_{k\,=\,0}^{+\,\infty}\,\biggl[\,(\,k\,-\,1\,)\,b_{\,k\,-\,1}\,p_{t}^{\,2}\ +\ 2\,k\,b_{\,k}\,p_{\,t}\,q_{\,t}\ +\ (\,k\,+\,1\,)\,b_{\,k\,+\,1}\,q_{t}^{\,2}\,\biggr]\,z^{\,k}\,,
\end{multline}
Moreover, one checks that the coefficient in front of $\dfrac{1}{3\,z^{\,2}\ -\ p}$ in $(\,\maltese\,)\ \bydef\ u_{\,t\,t}\ -\ u_{\,x}\,u_{\,x\,x}$ is equal to
\begin{multline}\label{eq:two}
  \sum_{k\,=\,0}^{+\,\infty}\,\biggl[\,(\,k\,-\,1\,)\,p_{t}^{2}\,b_{\,k\,-\,1}\ +\ \underline{2\,k\,b_{\,k}\,p_{t}\,q_{t}}\ +\ \underline{(\,k\,+\,1\,)\,b_{\,k\,+\,1}\,q_{t}^{2}}\\ 
  -\ \sum_{j\,=\,0}^{\,k}\biggl(\,\Bigl(\,1\,-\,\frac{1}{j}\,\Bigr)\,p_{x}\,b_{\,j\,-\,1}\ +\ q_{x}\,b_{\,j}\ +\ \frac{1}{j}\;\bigl(\,-\,p\,b_{\,j\,-\,1,\,x}\ +\ 3\,b_{\,j\,-\,3,\,x}\,\bigr)\,\biggr)\,\cdot\\
  \Bigl(\,(\,k\,-\,j\,-\,1\,)\,p_{x}^{2}\,b_{\,k\,-\,j\,-\,1}\ +\ 2\,(\,k\,-\,j\,)p_{x}\,q_{x}\,b_{\,k\,-\,j}\\ 
  +\ (\,k\,-\,j\,+\,1\,)\,q_{x}^{2}\,b_{\,k\,-\,j\,+\,1}\,\Bigr)\,\biggr]\,z^{\,k}\,.
\end{multline}
Now, we construct $(\,\widehat{\maltese}\,)$ from $(\,\maltese\,)\ \equiv\ u_{\,t\,t}\ -\ u_{\,x}\,u_{\,x\,x}$ in the following way. First, we replace in \eqref{eq:one} the term $2\,b_{\,k,\,t}\,q_{t}$ by $2\,\bh_{\,k,\,t}\,q_{t}\,$. Then, we replace in \eqref{eq:two} the term $2\,k\,b_{\,k} p_{t}\,q_{t}$ by $2\,k\,b_{\,k} \ph_{t}\,q_{t}$ and $(\,k\,+\,1\,)\,b_{\,k\,+\,1}\,q_{t}^{2}$ by $(\,k\,+\,1\,)\,b_{\,k\,+\,1}\,q_{t}\,\qh_{t}\,$. Then, the eikonal \cref{eq:eiko} is obtained by writing that $(\,\widehat{\maltese}\,)\ \in\ \O\,\llbracket\,z\,\rrbracket\,$, namely, the coefficient \eqref{eq:two} of $\dfrac{1}{3\,z^{\,2}\ -\ p}$ in $(\,\widehat{\maltese}\,)$ is divisible by $3\,z^{\,2}\ -\ p\,$. Finally, \cref{eq:b} is obtained by fulfilling the requirement that $(\,\widehat{\maltese}\,)\ \equiv\ \fO\,$.

\begin{remark}
In \cref{eq:two} we could have replaced the term $2\,k\,b_{\,k}\,p_{t}\,q_{t}$ (respectively, $(\,k\,+\,1\,)\,b_{\,k\,+\,1}\,q_{t}^{2}$) by $2\,k\,b_{\,k}\,\ph_{t}\,\qh_{t}$ (respectively $(\,k\,+\,1\,)\,b_{\,k\,+\,1}\,\qh_{t}^{2}$) and thus get a nonlinear eikonal equation. What we did is roughly a simple linearisation of this, which is more convenient for implementing in practice. In other words, no matter the eikonal equation's form, one has to perform the infinity of fixed point iterations. We took special care to make each such iterative step as easy as possible.
\end{remark}


\section{A simple explicit example of shock wave formation in the \acrshort{ibe}}
\label{app:a}

In this Appendix, we would like to provide one simple explicit example and discuss some relations between our theory of ramified solutions and the classical theory of shock waves. This Appendix should be perceived only as an appetizer. The deeper connections will be revealed in the subsequent works.

Let us consider the following \acrshort{ivp} for the \acrshort{ibe}:
\begin{align}\label{eq:quad1}
  u_{\,t}\ -\ u\,u_{\,x}\ &=\ \fO\,, \\
  u\,(\,0,\,x\,)\ &=\ x^{\,\half}\,, \label{eq:quad2}
\end{align}
where we took the square root singularity to simplify algebraic computations which follow. In the real case, we would consider the \acrshort{ivp} on a half-line $(\,t,\,x\,)\ \in\ \R_{\,\geq\,0}\times\R_{\,\geq\,0}\,$, but there are no such restrictions in the complex case. So, in agreement with our line of thinking, we assume that $(\,t,\,x\,)\ \in\ \C^{\,2}\,$. Using the methods of contact geometry described in Section~\ref{sec:cg}, we introduce the auxiliary variable
\begin{equation}\label{eq:char}
  y\ \eqdef\ x\ -\ x^{\,\half}\,t\,.
\end{equation}
The solution $u$ verifies the following second-order algebraic equation:
\begin{equation}\label{eq:quad}
  u^{\,2}\ -\ t\,u\ -\ y\ =\ \fO\,.
\end{equation}
The last algebraic \cref{eq:quad} can be also rewritten in an equivalent form:
\begin{equation*}
  \Bigl(\,u\ -\ \frac{t}{2}\,\Bigr)^{\,2}\ -\ \Bigl(\,y\ +\ \frac{t^{\,2}}{4}\,\Bigr)\ =\ \fO\,.
\end{equation*}
Now it is obvious that the discriminant of \cref{eq:quad} is
\begin{equation*}
  \delta\ \eqdef\ y\ +\ \frac{t^{\,2}}{4}\,.
\end{equation*}
Our theory states that the complex solution $u$ will be ramified around the locus $\delta\ =\ \fO\,$. Taking into account \eqref{eq:char}, two solutions to \cref{eq:quad} read
\begin{equation*}
  r_{\,\pm}\ \eqdef\ \frac{t}{2}\ \pm\ \sqrt{\Bigl(\,\frac{t}{2}\ -\ \sqrt{x}\,\Bigr)^{\,2}}\,.
\end{equation*}

\begin{figure}[H]
  \centering
  \includegraphics[width=0.51\textwidth]{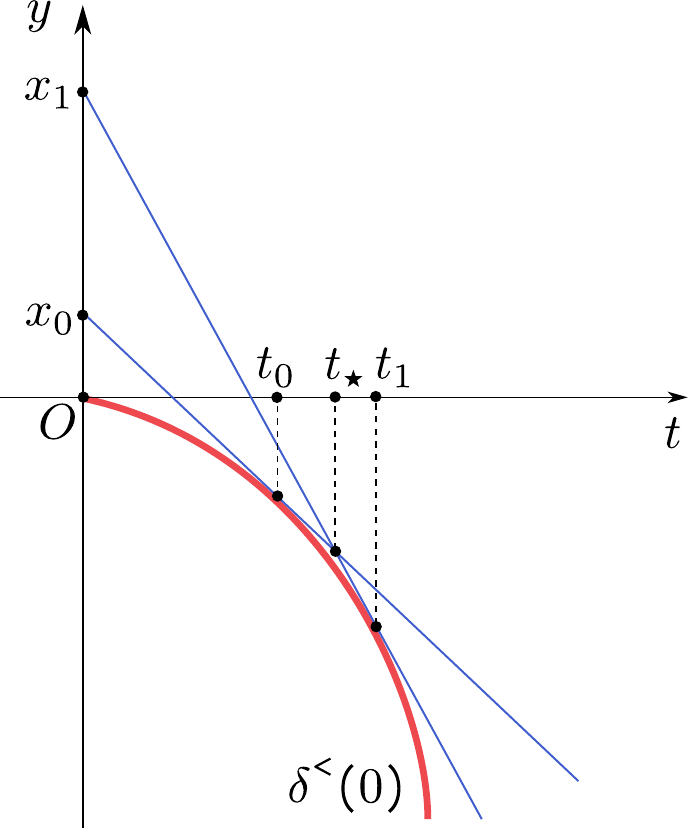}
  \caption{\small\em A schematic graphical illustration of the shock wave formation in the \acrshort{ibe} for an \acrshort{ivp} \eqref{eq:quad1}, \eqref{eq:quad2} with the square root algebraic singularity.}
  \label{fig:contact}
\end{figure}

We take two points $x_{\,0}\ \neq\ x_{\,1}$ and we consider two characteristics \eqref{eq:char} passing through these points at $t\ =\ 0\,$. If $x_{\,0,\,1}$ are real, we may assume that $0\ <\ x_{\,0}\ <\ x_{\,1}$ to fix the ideas. Along the characteristic $y_{\,j}\ \eqdef\ x_{\,j}\ -\ x_{\,j}^{\,\half}\,t$ the solution $u$ takes the value $x_{\,j}^{\,\half}\,$, $j\ \in\ 2^{\,\sqsubset}\,$. A shock wave forms when two characteristics cross each other since in the crossing point, the solution value is obviously contradictory:
\begin{equation*}
  x_{\,0}\ \neq\ x_{\,1}\ \qquad\ \Longrightarrow\ \qquad\ x_{\,0}^{\,\half}\ \neq\ x_{\,1}^{\,\half}\,.
\end{equation*}
The crossing will take place when
\begin{equation*}
  x_{\,0}\ -\ x_{\,0}^{\,\half}\,t\ =\ x_{\,1}\ -\ x_{\,1}^{\,\half}\,t\,.
\end{equation*}
It can be easily checked that the characteristic $y\ =\ x_{\,0}\ -\ x_{\,0}^{\,\half}\,t$ meets the singular locus $\delta\ =\ \fO$ in the real plane at time instance $t_{\,0}\ \eqdef\ 2\,x_{\,0}^{\,\half}\,$. Indeed,
\begin{equation*}
  y_{\,0}\ =\ x_{\,0}\ -\ x_{\,0}^{\,\half}\,t_{\,0}\ \equiv\ x_{\,0}\ -\ 2\,x_{\,0}\ =\ -\,x_{\,0}\,.
\end{equation*}
Henceforth,
\begin{equation*}
  \delta\,(\,t_{\,0},\,x_{\,0}\,)\ \bydef\ y_{\,0}\ +\ \frac{t_{\,0}^{2}}{4}\ =\ -\,x_{\,0}\ +\ \frac{\bigl(\,2\,x_{\,0}^{\,\half}\,\bigr)^{\,2}}{4}\ \equiv\ 0
\end{equation*}
and we deduce that $(\,t_{\,0},\,x_{\,0}\,)\ \in\ \delta^{\,<}\,(\,\set{0}\,)\,$. Moreover, it is not difficult to compute exactly the crossing time $t_{\,\star}\,$:
\begin{equation*}
  t_{\,\star}\ =\ x_{\,0}^{\,\half}\ +\ x_{\,1}^{\,\half}\ >\ 2\,x_{\,0}^{\,\half}\ \bydef\ t_{\,0}\,.
\end{equation*}
We would like to explain what happens precisely at time $t\ =\ t_{\,0}\,$. Until the crossing time $t\ =\ t_{\,\star}\,$, the solution value along the characteristic $y_{\,0}$ is $x_{\,0}^{\,\half}\,$. The point is that at time $t\ =\ t_{\,0}$, we change the solution branch. Indeed, for $0\ <\ t\ <\ t_{\,0}\,$:
\begin{equation*}
  r_{\,+}\ \bydef\ \frac{t}{2}\ +\ x_{\,0}^{\,\half}\ -\ \frac{t}{2}\ \equiv\ x_{\,0}^{\,\half}\,,
\end{equation*}
while for $t_{\,0}\ <\ t\ <\ t^{\,\star}$ we have
\begin{equation*}
  r_{\,-}\ \bydef\ \frac{t}{2}\ -\ \Bigl(\,\frac{t}{2}\ -\ x_{\,0}^{\,\half}\,\Bigr)\ \equiv\ x_{\,0}^{\,\half}\,.
\end{equation*}
In other words, at the discriminant contact point $t\ =\ t_{\,0}\,$, we have no other choice except to jump from one solution branch to another in order to keep the constant solution value along the given characteristic. The changing of the solution branch precedes and explains the apparition of the shock. This situation is illustrated in Figure~\ref{fig:contact}.

\section{Maple code}
\label{app:b}

In this Appendix we provide the \texttt{Maple}${}^{\text{\texttrademark}}$ code which implements one iteration of the mapping $\FFt$ used in \cref{sec:num}:
\bigskip
\begin{lstlisting}[language=Maple]
restart:
with(LinearAlgebra):
Digits := 30;
DEG := 25;

[...]

OneIter := proc(v::Array)
  description "One iteration of the fixed point algorithm";
  local N::integer,K::integer,l::integer,p::symbol,q::symbol;
  local px::symbol,pt::symbol,qx::symbol,qt::symbol; 
  local qtt::symbol,qxx::symbol,b::Array,A::Array,B::Array;
  local BB::Array,C::Array,M::Matrix,Minv::Matrix,B0::symbol;
  local B1::symbol,ptt::symbol,B0rhs::symbol,B1rhs::symbol;
  local k::integer,j::integer,P::symbol,Q::symbol,Pt::symbol;
  local Qt::symbol,RHS::symbol,rhs1::symbol,rhs2::symbol;
  local M11::symbol,M12::symbol,M21::symbol,M22::symbol;
  p := convert(v[1], float);
  q := convert(v[2], float);
  pt := convert(diff(p, t), float);
  px := convert(diff(p, x), float);
  qt := convert(diff(q, t), float);
  qx := convert(diff(q, x), float);
  ptt := convert(diff(p, t$2), float);
  qtt := convert(diff(q, t$2), float);
  qxx := convert(diff(q, x$2), float);
  N := ArrayNumElems(v) - 2;
  b := Array(-3 .. N - 1);
  b[-3 .. -1] := Array([0., 0., 0.]);
  b[0 .. N - 1] := convert(v[3 .. N + 2], float);
  K := ceil((N - 2)/2);
  B := Array(0 .. K);
  for k from 0 to K do
    B[k] := convert(b[0]*diff(q, x)*((k - 1.0)*px^2*b[k - 1] + 2.0*k*px*qx*b[k] + (k + 1)*qx^2*b[k + 1]), float);
    for j from 1 to k do
      B[k] := B[k] + convert(((1.0 - 1.0/j)*px*b[j - 1] + qx*b[j] + (-p*diff(b[j - 1], x) + 3.0*diff(b[j - 3], x))/j)*((k - j - 1.0)*px^2*b[k - j - 1] + 2.0*(k - j)*px*qx*b[k - j] + (k - j + 1.0)*qx^2*b[k - j + 1]), float);
    end do;
    B[k] := convert(B[k] - (k - 1.0)*pt^2*b[k - 1], float);
  end do;
  rhs1 := convert(add((p/3.0)^k*B[2*k], k = 0 .. floor(K/2)), float);
  rhs2 := convert(add((p/3.0)^k*B[2*k + 1], k = 0 .. floor((K - 1)/2)), float);
  M11 := convert(4.0*qt*add((p/3.0)^k*k*b[2*k], k = 0 .. floor((N - 1)/2)), float);
  M12 := convert(qt*add((p/3.0)^k*(2.0*k + 1.0)*b[2*k + 1], k = 0 .. floor((N - 2)/2)), float);
  M21 := convert(2.0*qt*add((p/3.0)^k*(2.0*k + 1.0)*b[2*k + 1], k = 0 .. floor((N - 2)/2)), float);
  M22 := convert(2.0*qt*add((p/3.0)^k*(k + 1.0)*b[2*k + 2], k = 0 .. floor((N - 3)/2)), float);
  M := convert(Matrix([[M11, M12], [M21, M22]]), float);
  Minv := MatrixInverse(M);
  Pt := mtaylor(rhs1*Minv[1, 1] + rhs2*Minv[1, 2], [t, x], DEG);
  Qt := mtaylor(rhs1*Minv[2, 1] + rhs2*Minv[2, 2], [t, x], DEG);
  P := dsolve({P(0.) = 0., diff(P(t), t) = Pt});
  Q := dsolve({Q(0.) = 1.0*x, diff(Q(t), t) = Qt});
  A := Array(0 .. N - 2);
  for k from 0 to N - 2 do
    A[k] := -b[0]*diff(q, x)*((k - 1.0)*px^2*b[k - 1] + 2.0*k*px*qx*b[k] + (k + 1.0)*qx^2*b[k + 1]);
    for j from 1 to k do
      A[k] := A[k] - ((1.0 - 1.0/j)*px*b[j - 1] + qx*b[j] + (-p*diff(b[j - 1], x) + 3.0*diff(b[j - 3], x))/j)*((k - j - 1.0)*px^2*b[k - j - 1] + 2.0*(k - j)*px*qx*b[k - j] + (k - j + 1.0)*qx^2*b[k - j + 1]);
    end do;
    A[k] := convert(A[k] + (k - 1.0)*pt^2*b[k - 1] + 2.0*k*Pt*qt*b[k] + (k + 1.0)*Qt*qt*b[k + 1], float);
  end do;
  C := Array(0 .. N - 4);
  for k from 0 to N - 4 do
    C[k] := convert(add((p/3.0)^l*A[k + 2 + 2*l], l = 0 .. floor((N - k - 4)/2)), float);
  end do;
  B0rhs := mtaylor(1/(2.0*qt)*(-b[0]*qtt - C[0] + b[0]*qx*(b[0]*qxx + 2.0*qx*diff(b[0], x))), [t, x], DEG);
  B0 := dsolve({B0(0.0) = 1.0, diff(B0(t), t) = B0rhs});
  B1rhs := mtaylor(1/(2.0*qt)*(-b[1]*qtt + p*diff(b[0], t $ 2) - C[1] + b[0]*qx*(b[1]*qxx + 2.0*qx*diff(b[1], x) - p*diff(b[0], x $ 2)) + (b[1]*qx - p*diff(b[0], x))*(b[0]*qxx + 2*qx*diff(b[0], x))), [t, x], DEG);
  B1 := dsolve({B1(0.0) = 1.0, diff(B1(t), t) = B1rhs});
  BB := Array(2 .. N - 4);
  for k from 2 to N - 4 do
    RHS := -(2.0 - 2.0/k)*pt*diff(b[k - 1], t) - (1.0 - 1.0/k)*ptt*b[k - 1] - qtt*b[k] - 1.0/k*(3.0*diff(b[k - 3], t $ 2) - p*diff(b[k - 1], t $ 2)) - C[k] + qx*b[0]*((2.0 - 2.0/k)*px*diff(b[k - 1], x) + (1.0 - 1.0/k)*ptt*b[k - 1] + qxx*b[k] + 2.0*qx*diff(b[k], x) + 1.0/k*(-p*diff(b[k - 1], x $ 2) + 3.0*diff(b[k - 3], x $ 2)));
    for j from 1 to k-1 do
      RHS := RHS + ((1.0 - 1.0/j)*px*b[j - 1] + qx*b[j] + 1.0/j*(-p*diff(b[j - 1], x) + 3.0*diff(b[j - 3], x)))*((2.0 - 2.0/(k - j))*px*diff(b[k - j - 1], x) + (1.0 - 1.0/(k - j))*ptt*b[k - j - 1] + qxx*b[k - j] + 2.0*qx*diff(b[k - j], x) + 1.0/(k - j)*(-p*diff(b[k - j - 1], x $ 2) + 3.0*diff(b[k - j - 3], x $ 2)));
    end do;
    RHS := RHS + ((1.0 - 1.0/k)*px*b[k - 1] + qx*b[k] + 1.0/k*(-p*diff(b[k - 1], x) + 3.0*diff(b[k - 3], x)))*(qxx*b[0] + 2.0*qx*diff(b[0], x));
    RHS := mtaylor(RHS/(2.0*qt), [t, x], DEG);
    BB[k] := dsolve({beta(0.0) = 0.0, diff(beta(t), t) = RHS});
  end do;
  return convert(Array([convert(series(rhs(P), t = 0, DEG), polynom), convert(series(rhs(Q), t = 0, DEG), polynom), convert(series(rhs(B0), t = 0, DEG), polynom), convert(series(rhs(B1), t = 0, DEG), polynom), seq(convert(series(rhs(BB[j]), t = 0, DEG), polynom), j = 2 .. N - 4)]), float);
end proc;
\end{lstlisting}
The code provided above was used to study the convergence of initial data \eqref{eq:data0}, \eqref{eq:data1} and all the others from \cref{sec:num}. A few modifications are needed to run the \acrshort{ivp} \eqref{eq:data2}. The complete programs can be obtained under a simple request by email.


\section{An example of a genuinely nonlinear problem}

Sometimes, an example is worth one thousand words. In our terminology, the following \nm{Cauchy} problem is nonlinear (quasi-linear) but not genuinely nonlinear:
\begin{equation*}
  u_{\,t\,t}\ -\ u\,u_{\,x\,x}\ =\ \fO\,, \qquad u\,(\,0,\,x\,)\ =\ c_{\,1}\,x\ +\ c_{\,2}\,x^{\,1\,+\,\third}\,,
\end{equation*}
for some real non-zero constants $c_{\,1,\,2}\ \in\ \R^{\,\times}\,$. In contrast, the following \acrshort{pde} is genuinely nonlinear:
\begin{equation*}
  u_{\,t\,t}\ -\ u_{\,x}\,u_{\,x\,x}\ =\ \fO\,, \qquad u\,(\,0,\,x\,)\ =\ c_{\,1}\,x\ +\ c_{\,2}\,x^{\,1\,+\,\third}\,.
\end{equation*}
In reality, this notion depends on the operator \emph{and} on the initial condition as well. That is why we speak above about the \nm{Cauchy} problem instead of just a \acrshort{pde}. The notion of a genuinely nonlinear problem was made more precise in \cref{sec:gen}.


\printglossary[type=\acronymtype]


\TheEnd

\end{document}